\documentclass{cspmA2}

\usepackage[figuresright]{rotating}
  \usepackage{floatpag}
  \rotfloatpagestyle{empty}

\usepackage{amsthm}
\usepackage{graphicx}

\usepackage{newtxtext}
\usepackage{newtxmath}

\usepackage{tikz-cd}

\newcommand{\inputtikz}[1]{\includegraphics{figures/#1}}

\usepackage[scaled=0.9]{couriers}

\usepackage{hyperref}

\usepackage{makeidx}
\makeindex

\newcounter{mparcnt}

\newcommand{\alttext}[3]{\paragraph{#1, page~{#2}} \quad #3\\[6pt]}
 \newcommand{\commentAlt}[2]{\immediate\write\tempfile{\noexpand\alttext{#1}{\thepage}{#2}}}

\newcommand{\commentLongAlt}[2]{\immediate\write\tempfile{\noexpand\alttext{#1}{\thepage\ (Long-Alt)}{#2}}}
\newwrite\tempfile
\immediate\openout\tempfile=alt-text.tex

\theoremstyle{plain}
\newtheorem{theorem}{Theorem}[chapter]
\newtheorem{lemma}[theorem]{Lemma}
\newtheorem{proposition}[theorem]{Proposition}
\newtheorem{corollary}[theorem]{Corollary}
\newtheorem{Exercise}[theorem]{Exercise}

  \def\makeRRlabeldot#1{\hss\llap{#1}}
  \renewcommand\theenumi{{\rm (\roman{enumi})}}
  \renewcommand\theenumii{{\rm (\alph{enumii})}}
  \renewcommand\theenumiii{{\rm (\arabic{enumiii})}}
  \renewcommand\theenumiv{{\rm (\Alph{enumiv})}}

\newcounter{newfoot}
\newenvironment{newfn}[1][]{\refstepcounter{newfoot}#1}

\begin{document}
\newcommand{\dis}{\displaystyle}
\newcommand{\txt}{\textstyle}

\newcommand{\noi}{\noindent}
\newcommand{\blok}{\hspace{-1pt}\tikz[scale=1,baseline=-4pt]{\draw[line width=5pt] (0,0) -- ++ (0.4,0);}\hspace{1pt}}
\newcommand{\med}{\medskip}
\newcommand{\quand}{\quad\mbox{and}\quad}


\newcommand{\halmos}{\rule{1ex}{1.4ex}}
\newcommand{\QED}{\nopagebreak{\hspace*{\fill}$\halmos$\medskip}}
\newenvironment{Proof}[1][]{\noi\textbf{Proof #1}}{\QED}
\newcommand{\bpro}{\begin{Proof}}
\newcommand{\epro}{\end{Proof}}


\newcommand{\al}{\alpha}
\newcommand{\bet}{\beta}
\newcommand{\ga}{\gamma}
\newcommand{\Ga}{\Gamma}
\newcommand{\de}{\delta}
\newcommand{\De}{\Delta}
\newcommand{\eps}{\varepsilon}
\newcommand{\la}{\lambda}
\newcommand{\La}{\Lambda}
\newcommand{\sig}{\sigma}
\newcommand{\tet}{\theta}
\newcommand{\om}{\omega}
\newcommand{\Om}{\Omega}

\newcommand{\si}{\ensuremath{\sigma}}

\newcommand{\Ai}{{\cal A}}
\newcommand{\Bi}{{\cal B}}
\newcommand{\Ci}{{\cal C}}
\newcommand{\Di}{{\cal D}}
\newcommand{\Ei}{{\cal E}}
\newcommand{\Fi}{{\cal F}}
\newcommand{\Gi}{{\cal G}}
\newcommand{\Hi}{{\cal H}}
\newcommand{\Ii}{{\cal I}}
\newcommand{\Ki}{{\cal K}}
\newcommand{\Li}{{\cal L}}
\newcommand{\Mi}{{\cal M}}
\newcommand{\Ni}{{\cal N}}
\newcommand{\Pc}{{\cal P}}
\newcommand{\Ri}{{\cal R}}
\newcommand{\Si}{{\cal S}}
\newcommand{\Ti}{{\cal T}}
\newcommand{\Vi}{{\cal V}}
\newcommand{\Xc}{{\cal X}}

\newcommand{\R}{{\mathbb R}}
\newcommand{\N}{{\mathbb N}}
\newcommand{\Z}{{\mathbb Z}}
\newcommand{\C}{{\mathbb C}}
\newcommand{\E}{{\mathbb E}}
\renewcommand{\P}{{\mathbb P}}

\newcommand{\desd}{\ensuremath{\Leftrightarrow}}
\newcommand{\volgt}{\ensuremath{\Rightarrow}}
\newcommand{\up}{\uparrow}
\newcommand{\down}{\downarrow}
\newcommand{\sub}{\subset}
\newcommand{\beh}{\backslash}
\newcommand{\asto}[1]{\underset{{#1}\to\infty}{\longrightarrow}}
\newcommand{\Asto}[1]{\underset{{#1}\to\infty}{\Longrightarrow}}

\newcommand{\ti}{\tilde}
\newcommand{\dgg}{\dagger}
\newcommand{\ov}{\overline}
\newcommand{\un}{\underline}
\newcommand{\subb}[2]{_{\begin{array}{c}\ \\[-.65cm]\scriptstyle{#1}\\[-.15cm]\scriptstyle{#2}\end{array}}}
\newcommand{\subbb}[3]{_{\begin{array}{c}\scriptstyle{#1}\\[-.15cm]\scriptstyle{#2}\\[-.15cm]\scriptstyle{#3}\end{array}}}

\newcommand{\ffrac}[2]{{\textstyle\frac{{#1}}{{#2}}}}
\newcommand{\dif}[1]{\ffrac{\partial}{\partial{#1}}}
\newcommand{\diff}[1]{\ffrac{\partial^2}{{\partial{#1}}^2}}

\newcommand{\cn}{\colon}
\newcommand{\di}{\mathrm{d}}
\newcommand{\half}{{[0,\infty)}}
\newcommand{\expo}{\mbox{\large\it e}}
\newcommand{\ex}[1]{\expo^{\,\textstyle{#1}}}
\newcommand{\var}{{\rm Var}}
\newcommand{\cov}{{\rm Cov}}
\newcommand{\ha}{\ffrac{1}{2}}
\newcommand{\haa}{\frac{1}{2}}

\newcommand{\Fb}{{\mathbf F}}
\newcommand{\Xb}{{\mathbf X}}
\newcommand{\Yb}{{\mathbf Y}}
\newcommand{\Zb}{{\mathbf Z}}

\newcommand{\mk}{\mathfrak{m}}

\newcommand{\cube}[3]{
\begin{scope}[cm={0.70711,-0.40825,0.70711,0.40825,(#1)},transform shape]
\clip (-#2/2,0) rectangle (0,#2/2);
\node[anchor=center,inner sep=0pt] at (0,0) {\includegraphics[width=#2cm]{figdata/#3.pdf}};
\draw[fill=white,opacity=0.3] (-#2/2,0) rectangle (0,#2/2);
\end{scope}
\begin{scope}[cm={0.70711,-0.40825,0,0.81650,(#1)},transform shape]
\clip (-#2/2,-#2/2) rectangle (0,0);
\node[anchor=center,inner sep=0pt] at (0,0) {\includegraphics[width=#2cm]{figdata/#3.pdf}};
\end{scope}
\begin{scope}[cm={0,-0.81650,0.70711,0.40825,(#1)},transform shape]
\clip (0,0) rectangle (#2/2,#2/2);
\node[anchor=center,inner sep=0pt] at (0,0) {\includegraphics[width=#2cm]{figdata/#3.pdf}};
\draw[fill=black,opacity=0.3] (0,0) rectangle (#2/2,#2/2);
\end{scope}
}

\newenvironment{dedication}
  {\clearpage           
   \thispagestyle{empty}
   \vspace*{\stretch{1}}
   \itshape             
   \raggedleft          
  }
  {\par 
   \vspace{\stretch{3}} 
   \clearpage           
  }

\newcommand{\rvot}[1]{
\draw[very thick,<-] (#1)++(0.22,0) -- ++ (0.78,0);
\draw[line width=5pt] (#1)++(-0.2,0) -- ++ (0.4,0);
}
\newcommand{\lvot}[1]{
\draw[very thick,<-] (#1)++(-0.22,0) -- ++ (-0.78,0);
\draw[line width=5pt] (#1)++(0.2,0) -- ++ (-0.4,0);
}
\newcommand{\rrw}[1]{
\draw[very thick,->] (#1)++(0.22,0) -- ++ (0.78,0);
\draw[line width=5pt] (#1)++(-0.2,0) -- ++ (0.4,0);
}
\newcommand{\lrw}[1]{
\draw[very thick,->] (#1)++(-0.22,0) -- ++ (-0.78,0);
\draw[line width=5pt] (#1)++(0.2,0) -- ++ (-0.4,0);
}

\title{A Course in Interacting Particle Systems}
\author{Jan M.\ Swart}

\frontmatter
\maketitle

\begin{dedication}
Everyone has the right to life, liberty and the security of person.\\
(Universal Declaration of Human Rights, Article 3.)
\vspace{14cm}

\parbox{13cm}{\rm This material will be published by Cambridge University Press as \emph{A Course in Interacting Particle Systems} by Jan M.\ Swart. This pre-publication version is free to view and download for personal use only. Not for re-distribution, re-sale or use in derivative works. \textcopyright\ Institute of Information Theory and Automation of the Czech Academy of Sciences, Pod Vod\'arenskou v\v{e}\v{z}\'i 1143/4, 182 00 Prague 8, Czech Republic.}
\end{dedication}

\tableofcontents
\listofcontributors

\subsection*{Author}

\begin{center}
Jan M.\ Swart,\\
The Czech Academy of Sciences,\\
Institute of Information Theory and Automation
\end{center}

\noi
Jan Swart is a research fellow at the Institute of Information Theory
and Automation (UTIA) of the Czech Academy of Sciences, Prague. He
regularly teaches at Charles University and has coauthored over forty
papers in probability theory with an emphasis on interacting particle
systems. He was born (in 1970) in the Netherlands. After obtaining his
PhD (in 1999) under the supervision of Frank den Hollander he first
moved to Germany where as a young postdoc he worked with Klaus
Fleichmann and others and then moved to the Czech republic where he
has been employed at UTIA since 2005.

\chapter*{Preface}

Interacting particle systems, in the sense we will be using the word in this book, are countable systems of locally interacting Markov processes. Each interacting particle system is defined on a lattice: a countable set with (usually) some concept of distance defined on it; the canonical choice is the $d$-dimensional integer lattice $\Z^d$. Situated on each point in this lattice, there is a continuous-time Markov process with a finite state space (often even of cardinality two) whose jump rates depend on the states of the Markov processes on near-by sites. Interacting particle systems are often used as extremely simplified ``toy models'' for stochastic phenomena that involve a spatial structure.

An attractive property of interacting particle systems is that they are easy to simulate on a computer.\footnote{To get started doing this yourself, look at my simulation library that is available from {\tt http://staff.utia.cas.cz/swart/simulate.html}.} Although the definition of an interacting particle system often looks very simple, and problems of existence and uniqueness have long been settled, it is often surprisingly difficult to prove anything nontrivial about its behavior. With a few exceptions, explicit calculations tend not to be feasible, so one has to be satisfied with qualitative statements and some explicit bounds. Despite intensive research over more than fifty years, some easy-to-formulate problems still remain open while the solutions of others have required the development of nontrivial and complicated techniques.

Luckily, as a reward for all this, it turns out that despite their simple rules, interacting particle systems are often remarkably subtle models that capture the sort of phenomena one is interested in much better than might initially be expected. Thus, while it may seem outrageous to assume that ``Plants of a certain type occupy points in the square lattice $\Z^2$, live for an exponential time with mean one, and place seeds on unoccupied neighboring sites with rate $\lambda$'' it turns out that making the model more realistic often does not change much in its overall behavior. Indeed, there is a general philosophy in the field, that is still insufficiently understood, that says that interacting particle systems come in ``universality classes'' with the property that all models in one class have roughly the same behavior.

As a mathematical discipline, the subject of interacting particle systems is still relatively young. It started around 1970 with the work of F.~Spitzer\index{Spitzer} \cite{Spi70} and R.L.~Dobrushin\index{Dobrushin} \cite{Dob71}, with many other authors joining in during the next few years. By 1975, general existence and uniqueness questions had been settled, four classical models had been introduced (the exclusion process, the stochastic Ising model, the voter model and the contact process), and elementary (and less elementary) properties of these models had been proved. In 1985, when Liggett\index{Liggett} published his famous book \cite{Lig85}, the subject had established itself as a mature field of study. Since then, it has continued to grow rapidly, to the point where it is impossible to accurately capture the state of the art in a single book. Liggett's second book \cite{Lig99} focuses on three of the four classical models only. Such is the sophistication of modern methods that by now it would be possible to write a book on each of the four classical models alone.

While interacting particle systems, in the narrow sense we defined them above, have apparently not been the subject of mathematical study before 1970, the subject has close links to some problems that are considerably older. In particular, the Ising model (without time evolution) has been studied since 1925 while both the Ising model and the contact process have close connections to percolation, which has been studied since the late 1950-ies. In recent years, more links between interacting particle systems and other, older subjects of mathematical research have been established, and the field continues to receive new impulses not only from the applied, but also from the more theoretical side. Until 1990 most of the work concentrated on the $d$-dimensional integer lattice. Since then other lattices such as trees have gained popularity. Interacting particle systems on random graphs are a hot topic. Scaling limits, both deterministic (hydrodynamic limits) and random (SPDE's, super Brownian motion, the Brownian web) are an old subject that continues to see exciting developments.

Apart from Liggett's books \cite{Lig85,Lig99}, there exist a number of other books that treat interacting particle systems in one way or another. Durrett's lecture notes from 1988 \cite{Dur88} and his St.\ Flour lecture notes \cite{Dur95} still make great reading. Both contain many pictures of numerical simulations, discuss the mean field limit, and use percolation theory as a unifying idea behind the proofs. While his 1988 book focuses mainly on the classical models, his St.\ Flour lecture notes show how block arguments together with assumptions such as sufficiently long range interaction or rapid stirring can be used to treat a plethora of models. Both books focus on developing the great ideas and getting one's hands dirty doing actual calculations, as opposed to carefully developing the basic technical framework. The St.\ Flour lecture notes are quite high level, so Master students and even beginning Ph.D.\ students may find it hard to see that all the technical details can really be filled in.

Two books that don't have interacting particle systems as their primary topic but nevertheless say something useful about them are Liggett's \emph{Continuous time Markov processes} \cite{Lig10} and Grimmett's \emph{Probability on Graphs} \cite{Gri18}. Apart from giving a short introduction to particle systems Liggett's book is also a useful reference for the classical theory of continuous-time Markov chains. Grimmett's book contains a lot material on percolation theory, which includes the contact process (as a form of oriented percolation) and the random cluster model with its application to the Ising and Potts models. Recently, Lanchier \cite{Lan24} has made an impressive attempt to give a complete overview of the existing literature on interacting particle systems in the life and social sciences, which he managed at the cost of being only able to sketch the main proof ideas.

The present book grew out of lecture notes I wrote for courses I have been giving at Charles University in Prague at regular intervals starting in 2009. Since the idea was that it should be possible to cover most of the material in a one semester course, a lot of basic material that can be found in Liggett's classic book \cite{Lig85} is omitted here. What is new in the book has less to do with new results but more to do with how the material is presented. The introductory chapter shows a wide variety of models that reflect the present state of the subject. There is an informal discussion of phase transitions and of critical behavior, which is still poorly understood from a mathematical side but nevertheless important to get the full picture. And there is a whole chapter devoted to the mean-field limit, which from the mathematical techniques involved is a bit off topic but nevertheless essential to get a good complete understanding.

One of the most fundamental novelties is how graphical representations are given a central role in the construction of all kinds of interacting particle systems, instead of being viewed as a tool for the study of certain specific models only. This reflects the modern state of the art that indeed uses graphical representations all of the time. It also allows the basic existence and uniqueness results (presented in Chapter~\ref{C:construct}) to be proved in a way that prepares for the discussion of duality in Chapter~\ref{C:dual}. The basic existence and uniqueness result (Theorem~\ref{T:Poispart}) is a pathwise result, that despite being based on well-known methods has not appeared in print in this form before. Most of the duality in this book is pathwise duality (which is a modern word for an old concept). Stochastic flows, both forward and backward in time, are given a central place.

The material is meant to be presentable (with minor omissions here and there) during a one semester course. As a preparation for the book, the students need a basic course in measure theory and probability. It is also preferable if they have at least some prior experience with continuous-time Markov chains, so that some of the more standard sections of Chapter~\ref{C:Markov} can be skipped over quickly in favor of sections containing material that is less widely known. I have tried to make the book reasonably self-contained, but not at all costs, so material about differential equations or about Feller semigroups is cited without proof. Chapter~\ref{C:Markov} and Sections \ref{S:conintro}--\ref{S:Gencon} contain the core technical results that one needs in order to understand the rest of the book. Chapters \ref{C:meanfield}, \ref{C:monot}, \ref{C:dual}, and \ref{C:percol} have been written in such a way that they do not depend too much on each other mutually.

People that I am indebted to for their comments and suggestions include Tibor Mach, Aernout van Enter, Sam Olesker-Taylor, Jan Niklas Latz, Jim Fill, Cristina Toninelli, and I am sure more people whom I am now forgetting. Work on this book was sponsored by GA\v{C}R grant 25-16267S.




\mainmatter

\chapter{Introduction}\label{C:intro}

\section{General set-up}\label{S:setup}

Let $S$ be a finite set, called the 
\emph{local state space},\index{local state space} and let $\La$ be a
countable set, called the \emph{lattice}.\index{lattice} We let
$S^\La$ denote the Cartesian product space of $\La$ copies of $S$,
that is, elements $x$ of $S^\La$ are of the form
\[
x=\big(x(i)\big)_{i\in\La}\quad\mbox{with}\quad x(i)\in S\ \forall\ i\in\La.
\]
Equivalently, $S^\La$ is nothing else than the set of all functions
$x\cn\La\to S$.

\emph{Interacting particle systems}\index{interacting particle system}
are continuous-time Markov processes $X=(X_t)_{t\geq 0}$ with a state
space \index{state space} of the form $S^\La$. Thus, $(X_t)_{t\geq 0}$
is a Markov process such that at each time $t\geq 0$, the state of $X$
is of the form
\[
X_t=\big(X_t(i)\big)_{i\in\La}\quad\mbox{with}\quad X_t(i)\in S\ \forall\ i\in\La.
\]
We call $X_t(i)$ the \emph{local state}\index{local state} of $X$ at time $t$
and at the \emph{position} $i$. Positions $i\in\La$ are also often called
\emph{sites}.\index{site}

The time evolution of continuous-time Markov processes is usually
characterized by their \emph{generator}\index{generator} $G$, which is an
operator acting on functions $f\cn\Si\to\R$, where $\Si$ is the state space. For
example, in the case of Brownian motion, the state space is $\Si=\R$ and the
generator is the differential operator $G=\ha\diff{x}$. In the case of an
interacting particle system, the state space is of the form $\Si=S^\La$ and
the generator can usually be written in the form
\begin{equation}\label{Gdef}
Gf(x)=\sum_{m\in\Gi}r_m\big\{f\big(m(x)\big)-f\big(x\big)\big\}
\qquad(x\in S^\La).
\end{equation}
Here $\Gi$ is a set whose elements are \emph{local maps} $m\cn S^\La\to S^\La$
and $(r_m)_{m\in\Gi}$ is a collection of nonnegative constants called
\emph{rates}, that say with which Poisson intensity the local map $m$ should
be applied to the configuration $X_t$. The precise definitions will be given
in later chapters, but at the moment it suffices to say that if we approximate
$(X_t)_{t\geq 0}$ by a discrete-time Markov chain where time is increased in
steps of size $\di t$, then
\[\begin{array}{rl}
\dis r_m\,\di t\ &\dis\mbox{is the probability that the map $m$}\\
&\dis\mbox{is applied during the time interval $(t,t+\di t]$.}
\end{array}\]

Often, the lattice $\La$ has the structure of an (undirected)
graph. In this case, we let $E$ denote the corresponding \emph{edge
set}.\index{edge set} This is a set of unordered pairs $\{i,j\}$ with
$i,j\in\La$ and $i\neq j$, that are called \emph{edges}. In drawings
of the graph, the fact that $\{i,j\}\in E$ is indicated by connecting
the points representing $i$ and $j$ by a line segment. We let
\[\index{0E@$\Ei$}
\Ei:=\big\{(i,j):\{i,j\}\in E\big\}
\]
denote the corresponding set of all \emph{ordered} pairs $(i,j)$ that
correspond to an edge. We call
\begin{equation}\label{Ni}\index{0Ni@$\Ni_i$}
\Ni_i:=\big\{j\in\La:\{i,j\}\in E\big\}
\end{equation}
the \emph{neighborhood} of the site $i$.

Many well-known and well-studied interacting particle systems are
defined on the \emph{$d$-dimensional integer lattice}\index{integer
  lattice} $\Z^d$. We denote the origin by
$0=(0,\ldots,0)\in\Z^d$. For any $i=(i_1,\ldots,i_d)\in\Z^d$, we let
\[\index{00norm1@$\Vert\,\cdot\,\Vert_1$}\index{00normi@$\Vert\,\cdot\,\Vert_\infty$}
\|i\|_1:=\sum_{k=1}^d|i_k|
\quand
\|i\|_\infty:=\max_{k=1,\ldots,d}|i_k|\qquad(i\in\Z^d)
\]
denote the $\ell_1$-norm and supremum-norm, respectively. For $R\geq 1$, we set
\begin{equation}\label{Edef}
E^d:=\big\{\{i,j\}:\|i-j\|_1=1\big\}\quand
E^d_R:=\big\{\{i,j\}:0<\|i-j\|_\infty\leq R\big\}.
\end{equation}
Then $(\Z^d,E^d)$ is the integer lattice equipped with the
\emph{nearest neighbor}\index{nearest neighbor} graph structure and
$(\Z^d,E^d_R)$ is the graph obtained by connecting all points within
$\|\,\cdot\,\|_\infty$-distance $R$ with an edge. We let $\Ei^d$ and
$\Ei^d_R$ denote the corresponding sets of ordered pairs $(i,j)$.

The graphs we have just introduced have the property that they ``look
the same'' seen from any vertex.
An \emph{automorphism}\index{automorphism!of graphs} of a graph $(\La,E)$
is a bijection $\psi\cn\La\to\La$ that ``preserves the edges'' in the
sense that $\{\psi(i),\psi(j)\}\in E$ if and only if $\{i,j\}\in E$. A
graph $(\La,E)$ is called \emph{vertex
transitive}\index{vertex transitive!graph}\index{transitivity!of graphs}
if for each pair of vertices $i,j\in\La$, there exists an automorphism
$\psi$ such that $\psi(i)=j$. The graphs $(\Z^d,E^d)$ and
$(\Z^d,E^d_R)$ we have just introduced are clearly vertex transitive
(just take $\psi(k):=k+j-i$). Another example of vertex transitive
graphs are $d$-regular trees, that is, infinite trees in which each
vertex has precisely $d$ neighbors.

Before we turn to rigorous mathematical theory, it is good to see a
number of examples. It is easy to simulate interacting particle
systems on a computer. In simulations, the infinite graphs
$(\Z^d,E^d)$ or $(\Z^d,E^d_R)$ are replaced by a finite piece of
$\Z^d$, with some choice of the boundary conditions (for example
periodic boundary conditions).

\section{The voter model}\label{S:voter}

For each $i,j\in\La$, the \emph{voter model map}\index{voter model map}
${\rm vot}_{ij}\cn S^\La\to S^\La$ is defined as
\begin{equation}\label{votmap}\index{0vot@${\tt vot}_{ij}$}
{\tt vot}_{ij}(x)(k):=\left\{\begin{array}{ll}
x(i)\quad&\mbox{if }k=j,\\[5pt]
x(k)\quad&\mbox{otherwise.}
\end{array}\right.
\end{equation}
Applying ${\tt vot}_{ij}$ to a configuration $x$ has the effect that
local state of the site $i$ is copied onto the site $j$. The
\emph{nearest neighbor voter model}\index{nearest neighbor} on $\Z^d$
is the interacting particle system with generator
\begin{equation}\label{Gvot}
G_{\rm vot}f(x)
:=\frac{1}{|\Ni_0|}\sum_{(i,j)\in\Ei^d}
\big\{f\big({\tt vot}_{ij}(x)\big)-f\big(x\big)\big\}
\qquad(x\in S^{\Z^d}).
\end{equation}
Here $\Ni_0$ is the neighborhood of the origin and $|\Ni_0|=2d$
denotes its cardinality.\index{00normj@$\vert A\vert$} Similarly,
replacing the set of directed edges $\Ei^d$ by $\Ei^d_R$ and replacing
$\Ni_0$ by the appropriate set of neighbors in this new graph, we
obtain the \emph{range $R$} \index{range of a particle system} voter
model.

\begin{figure}[htb]
\begin{center}
\includegraphics[width=4.5cm]{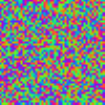} \quad
\includegraphics[width=4.5cm]{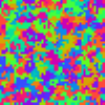}

\vspace{10pt}

\includegraphics[width=4.5cm]{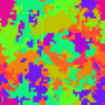}  \quad
\includegraphics[width=4.5cm]{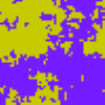}
\caption{Four snapshots of a two-dimensional voter model with periodic
  boundary conditions. Initially, the types of sites are i.i.d. Time evolved
  in these pictures is 0, 1, 32, and 500.}
\label{fig:vot2d}
\commentAlt{Figure~\ref{fig:vot2d}}{Four pictures showing a
  checkerboard of colors. Progressively, the number of colors is
  reduced and the areas with the same color are becoming larger.}
\end{center}
\end{figure}

In the context of the voter model, the local state $x(i)$ at a site
$i$ is often called the \emph{type}\index{type of a site} at $i$. The
voter model is often used to model biological populations, where
organisms with different genetic types occupy sites in space. Note
that since each site $j$ has $|\Ni_j|=|\Ni_0|$ neighbors, the total
rate of all maps ${\tt vot}_{ij}$ with $i\in\Ni_j$ is one. In view of
this, an alternative way to describe the dynamics in (\ref{Gvot}) is
to say that with rate 1, the organism living at a given site dies, and
is replaced by a descendant chosen with equal probability from its
neighbors.

An alternative interpretation, that has given the voter model its name, is
that sites represent people and types represent political opinions. With rate
one, an individual becomes unsure what political party to vote for, asks a
randomly chosen neighbor, and copies his/her opinion.

In Figure~\ref{fig:vot2d}, we see four snapshots of the time evolution of
a two-dimensional nearest-neighbor voter model. The initial state is
constructed by assigning i.i.d.\ types to the sites. Due to the copying
dynamics, we see patches appear where every site in a local neighborhood has
the same type. As time proceeds, these patches, usually called
\emph{clusters},\index{cluster} grow in size, so that eventually, for any
$N\geq 1$, the probability that all sites within distance $N$ of the origin
are of the same type tends to one.\footnote{In spite of this, for the model on
  the infinite lattice, it is still true that the origin changes its type
  infinitely often.}

\begin{figure}[htb]
\begin{center}
\inputtikz{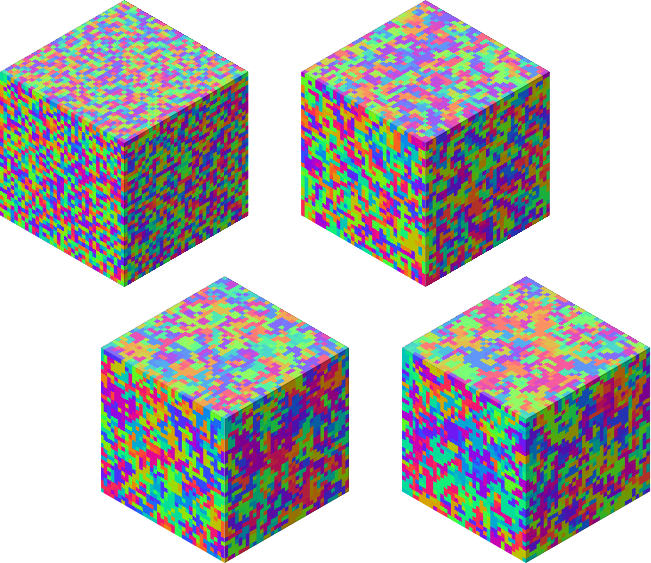}
\caption{Four snapshots of a three-dimensional voter model
  with periodic boundary conditions. Initially, the types of sites are
  i.i.d. Time evolved in these pictures is 0, 4, 32, and 250.}
\label{fig:vot3d}
\commentAlt{Figure~\ref{fig:vot3d}}{Four pictures showing a cube made
  up of small colored cubes. Progressively, the areas with the same
  color are becoming larger, but this process comes to a halt and the
  last two cubes look similar.}
\end{center}
\end{figure}

It turns out that this sort of behavior, called
\emph{clustering},\index{clustering} is dimension dependent. The voter model
clusters in dimensions 1 and 2, but not in dimensions 3 and more.
In Figure~\ref{fig:vot3d}, we see four snapshots of the time evolution of
a three-dimensional voter model. The model is simulated on a cube with
periodic boundary conditions. In this case, we see that even after a
long time, there are still many different types near the
origin.\footnote{On a finite lattice, such as we use in our
  simulations, one would eventually see one type take over, but the
  time one has to wait for this is very long compared to dimensions 1
  and 2. On the infinite lattice, the probability that the origin has
  a different type from its right neighbor tends to a positive limit
  as time tends to infinity.}

\section{The contact process}\index{contact process}

The contact process is another interacting particle system with a biological
interpretation. For this process, we choose the local state space $S=\{0,1\}$.
We interpret a site such that $X_t(i)=1$ as \emph{occupied} by an organism,
and a site such that $X_t(i)=0$ as \emph{empty}. Alternatively, the contact
process can be seen as a model for the spread of an infection. In this case,
sites with $X_t(i)=1$ are called \emph{infected} and sites with $X_t(i)=0$ are
called \emph{healthy}.

For each $i,j\in\La$, we define a \emph{branching map}\index{branching map}
${\tt bra}_{ij}\cn\{0,1\}^\La\to\{0,1\}^\La$ as
\begin{equation}\label{bramap}\index{0bra@${\tt bra}_{ij}$}
{\tt bra}_{ij}(x)(k):=\left\{\begin{array}{ll}
x(i)\vee x(j)\quad&\mbox{if }k=j,\\[5pt]
x(k)\quad&\mbox{otherwise.}
\end{array}\right.
\end{equation}
Note that this says that if prior to the application of ${\tt bra}_{ij}$, the
site $i$ is occupied, then after the application of ${\tt bra}_{ij}$, the site
$j$ will also be occupied, regardless of its previous state. If initially $i$
is empty, then nothing happens. We interpret this as the organism at $i$
giving \emph{birth} to a new organism at $j$, or the infected site $i$
\emph{infecting} the site $j$. If $j$ is already occupied/infected, then
nothing happens.

For each $i\in\La$, we also define a \emph{death map}\index{death map}
${\tt death}_i\cn\{0,1\}^\La\to\{0,1\}^\La$ as
\begin{equation}\label{deathmap}\index{0death@${\tt death}_i$}
{\tt death}_i(x)(k):=\left\{\begin{array}{ll}
0\quad&\mbox{if }k=i,\\[5pt]
x(k)\quad&\mbox{otherwise.}
\end{array}\right.
\end{equation}
If the map ${\tt death}_i$ is applied, then an organism at $i$, if there is
any, dies, respectively, the site $i$, if it is infected, \emph{recovers} from
the infection.

\begin{figure}[htb]
\begin{center}
\includegraphics[width=4.5cm]{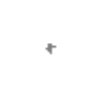} \quad
\includegraphics[width=4.5cm]{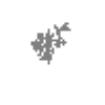}

\vspace{10pt}

\includegraphics[width=4.5cm]{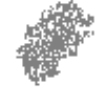}  \quad
\includegraphics[width=4.5cm]{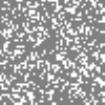}
\caption{Four snapshots of a two-dimensional contact process. Initially, only
  a single site is infected. The infection rate is 2, the death rate is 1, and
  time evolved in these pictures is 1, 5, 10, and 20.}
\label{fig:contact}
\commentAlt{Figure~\ref{fig:contact}}{Four pictures showing a subset
  of a checkerboard, increasing in size, inside of which most, but not
  all small squares are dark grey. Finally, this area fills out the
  whole picture.}
\end{center}
\end{figure}

Recalling (\ref{Edef}), the (nearest neighbor) contact process with
\emph{infection rate}\index{infection rate} $\la\geq 0$ and
\emph{death rate}\index{death rate} $\de\geq 0$ is the interacting
particle system with generator
\begin{equation}\begin{array}{r@{\,}c@{\,}l}\label{Gcontact}
\dis G_{\rm cont}f(x)
&:=&\dis\la\sum_{(i,j)\in\Ei^d}
\big\{f\big({\tt bra}_{ij}(x)\big)-f\big(x\big)\big\}\\[5pt]
&&\dis+\de\sum_{i\in\Z^d}\big\{f\big({\tt death}_i(x)\big)-f\big(x\big)\big\}
\qquad(x\in\{0,1\}^{\Z^d}).
\end{array}\end{equation}
This says that infected sites infect each healthy neighbor with rate $\la$,
and infected sites recover with rate $\de$.

In Figure~\ref{fig:contact}, we see four snapshots of the time evolution
of a two-dimensional contact process. Occupied sites are black and empty sites
are white. Initially, only the origin is occupied. The infection rate is 2 and
the death rate is 1. In this example, the infection spreads through the whole
population, eventually reaching a steady state\footnote{In fact, on the finite
  square used in our simulations, one can prove that the infection dies out
  a.s. However, the time one has to wait for this is exponentially large in
  the system size. For the size of system shown in Figure~\ref{fig:contact},
  this time is already too long to be numerically observable.} where a positive
fraction of the population is infected. Of course, starting from a single
infected site, there is always a positive probability that the infection dies
out in the initial stages of the epidemic.

\begin{figure}
\begin{center}
\inputtikz{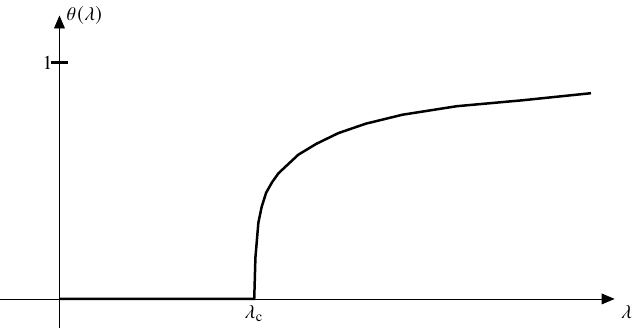}
\end{center}
\caption{Survival probability of the one-dimensional contact process.}
\label{fig:surprob}
\commentAlt{Figure~\ref{fig:surprob}}{Graph of a continuous function
  that is constantly zero left of a critical point called
  $\lambda_c$ and positive on the right of this point. The
  derivative at the critical point is infinite.}
\end{figure}

Unlike the voter model, the behavior of the contact process is roughly
similar in different dimensions. Instead of the dimension, this time,
the proportion $\la/\de$ of the infection rate to the death rate
determines the long-time behavior. By changing the speed of time, we
can without loss of generality choose one of the constants $\la$ and
$\de$ to be one, and it is customary to set $\de:=1$.  Let
$e_i\in\{0,1\}^\La$ be defined by $e_i(j):=1$ if $i=j$ and $:=0$
otherwise.\index{0ei@$e_i$} In Figure~\ref{fig:surprob}, we have
plotted the \emph{survival probability}
\begin{equation}\label{surprob}\index{0zzhthetalambda@$\tet(\la)$}
\tet(\la):=\P^{e_0}[X_t\neq 0\ \forall t\geq 0]
\end{equation}
of the one-dimensional contact process, started in $X_0=e_0$, that is,
with a single infected site at the origin, as a function of the
infection rate $\la$. Note that since $(\Z^d,E^d)$ is vertex
transitive, there is nothing special about the origin here: we could
have picked any other site instead. For reasons that we cannot explain
here, $\tet(\la)$ is in fact the same as the probability that the
origin is infected in equilibrium; this will be proved in
Lemma~\ref{L:theta}.

It turns out that for the nearest-neighbor contact process on $\Z^d$,
there exists a \emph{critical
 value}\index{critical!infection rate}\index{critical!point}
$\la_{\rm c}=\la_{\rm c}(d)$ with $0<\la_{\rm c}<\infty$ such that
$\tet(\la)=0$ for $\la\leq\la_{\rm c}$ and $\tet(\la)>0$ for $\la>\la_{\rm c}$.
The function $\tet$ is continuous, strictly increasing and concave on
$[\la_{\rm c},\infty)$ and satisfies $\lim_{\la\to\infty}\tet(\la)=1$. One has
\cite[Table~3.2]{HHL08}
\begin{equation}\label{lacrit}
\la_{\rm c}(1)=1.648924\pm0.00011.
\end{equation}
Proving these statements is not easy, however. For example, continuity of the
function $\tet$ in the point $\la_{\rm c}$ was proved only in 1990
\cite{BG90}, seventeen years after the introduction of the model in
\cite{CS73,Har74}. The best\footnote{There exists a sequence of rigorous upper
  bounds on the constant from (\ref{lacrit}) that is known to converge to the
  real value, but these bounds are so difficult to calculate that the best
  bound that has really been achieved by this method is much worse than the one
  in \cite{Lig95}.} rigorous upper bound on the constant from (\ref{lacrit})
is $\la_{\rm c}(1)\leq 1.942$ which is proved in \cite{Lig95}.

Krone \cite{Kro99} introduced a \emph{two-stage contact
process}\index{two-stage contact process}. In this model, the local
state space is $\{0,1,2\}$ where $0$ represents an empty site, $1$ a
young organism, and $2$ an adult organism. In a branching event, an
adult organism produces a young organism on an empty neighboring
site. In addition, young organisms can grow up. Both young and adults
can die, the young possibly at a higher rate than the adults. The
behavior of this model is similar to that of the contact process.

\section{Ising and Potts models}\label{S:IsingPotts}
\index{Ising model}\index{Potts model}

In a \emph{stochastic Ising model}, sites in the lattice $\Z^d$ are
interpreted as atoms in a crystal, that can have two possible local
states, usually denoted by $-1$ and $+1$. In the traditional
interpretation, these states describe the direction of the magnetic
field of the atom, and because of this, the local state $x(i)$ of a
site $i$ is usually called the \emph{spin} at $i$. More generally, one
can consider \emph{stochastic Potts models} where each ``spin'' can
have $q\geq 2$ possible values. In this case, the local state space is
traditionally denoted as $S=\{1,\ldots,q\}$, the special case $q=2$
corresponding to the Ising model (except for a small difference in
notation between $S=\{-1,+1\}$ and $S=\{1,2\}$).

Given a state $x$ and site $i$, we let
\begin{equation}\index{0Nxi@$N_{x,i}$}\label{Mxi}
N_{x,i}(\sig):=\sum_{j\in\Ni_i}1_{\txt\{x(j)=\sig\}}\qquad(\sig\in S)
\end{equation}
denote the number of neighbors of the site $i$ that have the spin
value $\sig\in S$. In the Ising and Potts models, sites like or
dislike to have the same spin value as their neighbors, depending on a
parameter $\bet\in\R$ called the
\emph{inverse temperature}.\index{inverse temperature}
In the physical interpretation of the model, $1/\bet$ corresponds (up
to a multiplicative constant) to the temperature (in degrees Kelvin
above the absolute zero). Adding a so-called \emph{Glauber
dynamics}\index{Glauber dynamics} \cite{Gla63} to the
model,\footnote{The terms \emph{Ising model} and \emph{Potts model}
refer only to certain Gibbs measures. A \emph{stochastic} Ising model
or Potts model is any interacting particle system that has these Gibbs
measures as its invariant laws (usually reversible). There exist
several different ways to invent a dynamics with this property. This
will be explained in a bit more detail in Section~\ref{S:oth}. In this
section, we stick to Glauber dynamics.} sites update their spin values
with rate one, and at such an event choose a new spin value with
probabilities that depend on the values of their neighbors. More
precisely, the \emph{stochastic Potts model} with \emph{Glauber
dynamics} is the interacting particle system that evolves in such a
way that
\begin{equation}\label{Glauber}
\mbox{site $i$ flips to the value $\sig$ with rate}\quad
r^\sig_i(x):=\frac{e^{\bet N_{x,i}(\sig)}}{\sum_{\tau\in S}e^{\bet N_{x,i}(\tau)}}.
\end{equation}
More formally, we can write the generator as
\begin{equation}\label{GPotts}\begin{array}{r@{\,}c@{\,}l}
\dis G_{\rm Potts}f(x)&:=&\dis\sum_{i\in\Z^d}\sum_{\sig\in S}
r^\sig_i(x)\big\{f\big(m^\sig_i(x)\big)-f\big(x\big)\big\},
\end{array}\end{equation}
where $m^\sig_i\cn S^\La\to S^\La$ are maps defined by
\begin{equation}\label{msig}
m^\sig_i(x)(j):=\left\{\begin{array}{ll}
\sig\quad&\mbox{if }j=i,\\[5pt]
x(j)\quad&\mbox{otherwise.}\end{array}\right.
\end{equation}
The attentive reader may notice that the way we have written the
generator in (\ref{GPotts}) is different from the way we have written
our generators so far, since unlike the rates $r_m$ in (\ref{Gdef}),
the rates $r^\sig_i(x)$ depend on the state $x$. This will be
explained in more detail in Chapter~\ref{C:construct}. In particular,
in Section~\ref{S:Isap}, we will see that it is possible to rewrite
the generator in (\ref{GPotts}) in a way that fits the general form
(\ref{Gdef}) (with rates that do not depend on the state $x$) but for
the Potts model, unlike the models we have seen so far, this way of
writing the generator is less natural and more complicated.

Returning to our informal description in (\ref{Glauber}), we notice
that for $\bet>0$, sites prefer to have spin values that agree with as
many neighbors as possible, that is, the model is
\emph{ferromagnetic}. \index{ferromagnetic} For $\bet<0$, the model is
\emph{antiferromagnetic}. \index{antiferromagnetic} These terms
reflect the situation that in some materials, neighboring spins like
to line up, which can lead to long-range order that has the effect
that the material can be magnetized. Antiferromagnetic materials, on
the other hand, lack this effect.

Alternatively, Potts models can also be interpreted as social or economic
models, where sites represent people or firms and spin values represent
opinions or the state (financially healthy or not) of a firm \cite{BD01}.

\begin{figure}[htb]
\begin{center}
\includegraphics[width=4.5cm]{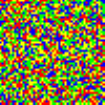} \quad
\includegraphics[width=4.5cm]{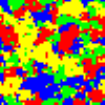}

\vspace{10pt}

\includegraphics[width=4.5cm]{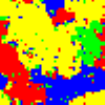}  \quad
\includegraphics[width=4.5cm]{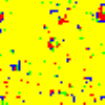}
\caption{Four snapshots of a $q=4$, $\bet=1.2$ Potts model with Glauber
  dynamics and periodic boundary conditions. Initially, the types of sites are
  i.i.d. Time evolved in these pictures is 0, 4, 32, 500.}
\label{fig:Potts}
\commentAlt{Figure~\ref{fig:Potts}}{Four pictures showing a
  checkerboard with four colors. Progressively, the areas with the
  same color are becoming larger, but inside large areas of one color
  there are small specks of different colors.}
\end{center}
\end{figure}

In Figure~\ref{fig:Potts} we see four snapshots of a two-dimensional
nearest-neighbor Potts model with four possible spin values. We have used
periodic boundary conditions, and the value of the parameter $\bet$ is
$1.2$. Superficially, the behavior is similar to that of a voter model, in the
sense that the system forms clusters of growing size that in the end take over
any finite neighborhood of the origin. Contrary to the voter model, however,
even in the middle of a large cluster that is predominantly of one color, sites
can still flip to other values as is clear from (\ref{Glauber}), so in the
simulations we see many small islands of different colors inside large
clusters where one color dominates. Another difference is that
clustering\index{clustering} happens only when the value of the
parameter $\bet$ is large enough. For small values of $\bet$, the behavior is
roughly similar to the voter model in dimensions $d\geq 3$. There is a
critical value $0<\bet_{\rm c}<\infty$ \index{critical!point} where the model
changes from one type of behavior to the other type of behavior. In this
respect, the model is similar to the contact process.

\begin{figure}
\begin{center}
\inputtikz{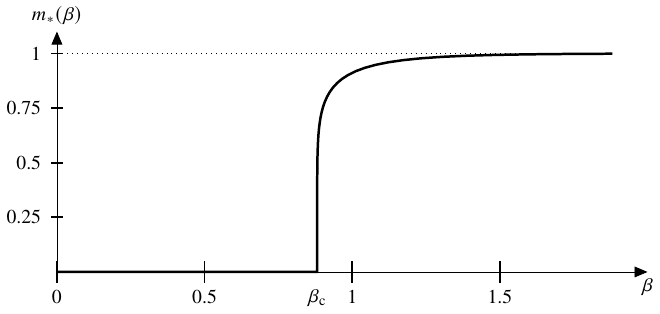}
\caption{The spontaneous magnetization of the two-dimensional Ising model.}
\label{fig:Ismag}
\commentAlt{Figure~\ref{fig:Ismag}}{Graph of a continuous function
  that is constantly zero left of a critical point called $\beta_c$
  and positive on the right of this point. The derivative at the
  critical point is infinite. The graph looks steeper than in
  Figure~\ref{fig:surprob}.}
\end{center}
\end{figure}

To make this critical value visible, imagine that instead of periodic boundary
conditions, we would use frozen boundary conditions where the sites at the
boundary are kept fixed at one chosen color, say color~1. Then the system has
a unique invariant law (equilibrium), in which for sufficiently large values
of $\bet$ the color 1 is (much) more frequent than the other colors, but for
low values of $\bet$ all colors occur with (almost) the same frequency. In
particular, for the Ising model, where the set of possible spin values is
$\{-1,+1\}$, we let
\begin{equation}\begin{array}{rl}\index{0mbet@$m_\ast(\bet)$}\label{mast}
\dis m_\ast(\bet):=
&\mbox{the expectation of $x(0)$ with $+1$ boundary}\\
&\mbox{conditions, in the limit of large system size.}
\end{array}\end{equation}
This function is called the \emph{spontaneous magnetization}.
\index{spontaneous magnetization}\index{magnetization}
For the Ising model in two dimensions, the spontaneous magnetization can be
explicitly calculated, as was first done by Onsager \cite{Ons44}.
The formula is
\begin{equation}\label{Onsager}
m_\ast(\bet)=\left\{\begin{array}{ll}
\big(1-{\rm sinh}(\bet)^{-4}\big)^{1/8}
\quad&\mbox{for }\bet\geq\bet_{\rm c}:=\log(1+\sqrt{2}),\\[5pt]
0\quad&\mbox{for }\bet\leq\bet_{\rm c}.
\end{array}\right.
\end{equation}
This function is plotted in Figure~\ref{fig:Ismag}. In this case, the critical
point $\bet_{\rm c}$ is known explicitly.

For Ising models in dimensions $d\geq 3$, the graph of $m_\ast(\bet)$ looks
roughly similar to Figure~\ref{fig:Ismag}, with $\bet_{\rm c}\approx 0.442$
in dimension $3$ \cite{GPA01}, but no explicit formulas are known.

\begin{figure}[htb]
\begin{center}
\inputtikz{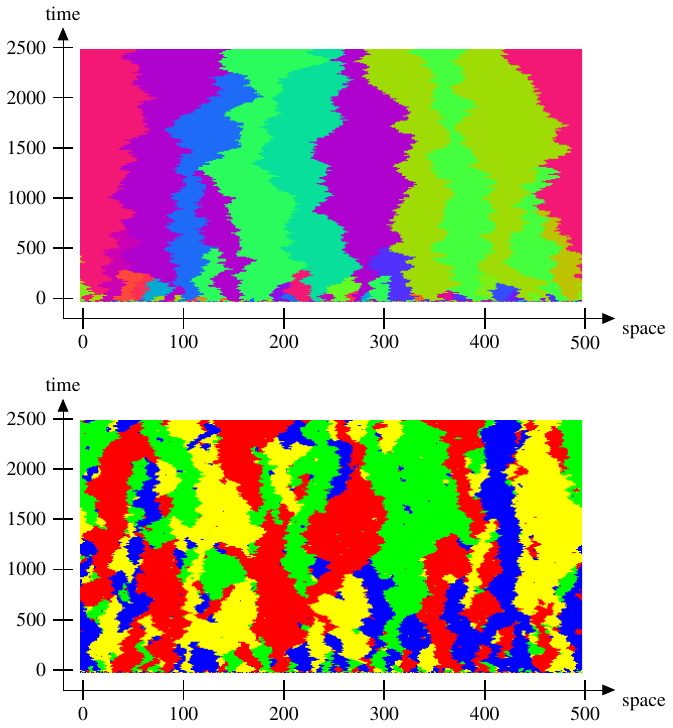}
\caption{Time evolution of a one-dimensional voter model (above) and a
  one-dimensional Potts model with a high value of $\bet$ (below).}
\label{fig:Potvot1D}
\commentAlt{Figure~\ref{fig:Potvot1D}}{Pictures of a voter and Potts
  model in which space is plotted horizontally, time vertically, and
  the color of a space-time point indicates its type. See long
  description.}
\commentLongAlt{Figure~\ref{fig:Potvot1D}}{In the voter model, the
  boundaries between regions of the same color are coalescing random
  walk paths, that are erratic and on a large scale look like Brownian
  motions. After two paths coalesce, the color that was previously
  between them has died out. The picture for the Potts model is
  similar but there are only four colors and colors that were locally
  extinct can sometimes reappear. Because of this, intervals with the
  same color stop growing in size after some time. An attentive
  observer may notice that new colors appear more frequently on the
  boundary between two existing colors than inside areas of one
  color.}
\end{center}
\end{figure}

In dimension one, one has $m^\ast(\bet)=0$ for \emph{all} $\bet\geq
0$. More generally, one-dimensional Potts models do not show long
range order, even if $\bet$ is very large.\footnote{This was first
noticed by Ising \cite{Isi25}, who introduced the model but noticed
that it was uninteresting, since an incorrect heuristic reasoning led
him to believe that what he had proved in dimension~1 would probably
hold in any dimension. Peierls \cite{Pei36} realized that dimension
matters and proved that the Ising model in higher dimensions does show
long range order.} By this we mean that in equilibrium, the
correlation between the spin values at $0$ and a point $i\in\Z$ tends
to zero as $i\to\infty$ for any value of $\bet$ (even though the decay
is slow if $\bet$ is large). In Figure~\ref{fig:Potvot1D}, we compare
the time evolution of a one-dimensional Potts model (with a large
value of $\bet$) with the time evolution of a one-dimensional voter
model. In the voter model, the cluster size keeps growing, but in the
Potts model, the typical cluster size converges to a finite limit.

\section{Phase transitions}\index{phase transition}\label{S:phase}

Figures~\ref{fig:surprob} and \ref{fig:Ismag} are examples of a phenomenon
that is often observed in interacting particle systems. As a parameter
governing the dynamics crosses a particular value, the system goes through
an abrupt change in behavior. This is called a \emph{phase transition} and the
value of the parameter is called the \emph{point of the phase transition} or,
in the mathematical literature, \emph{critical point}.\index{critical!point}
As we will see in a moment, in the physics literature, the term critical point
has a more restricted meaning. The term ``phase transition'' of course also
describes the behavior that certain materials change from a gas, fluid, or
solid phase into another phase at a particular value of the temperature,
pressure etc., and from the theoretical physicist's point of view, this is
indeed the same phenomenon.

In both Figure~\ref{fig:surprob} and \ref{fig:Ismag}, the point of the
phase transition in fact separates two regimes, one where the
interacting particle systems (on the infinite lattice) has a unique
invariant law (below $\la_{\rm c}$ and $\bet_{\rm c}$) and another
regime where there are more invariant laws (above $\la_{\rm c}$ and
$\bet_{\rm c}$). Indeed, for the contact process, the delta measure on
the empty configuration is always an invariant law, but above
$\la_{\rm c}$, a second, nontrivial invariant law also appears. Potts
models have $q$ invariant laws (one corresponding to each color) above
the critical point.\footnote{More precisely, they have $q$ invariant
laws that have the additional property that they are also translation
invariant in space. Depending on the dimension, there may exist
additional invariant laws that are not translation invariant.}
Multiple invariant laws are a general phenomenon associated with phase
transitions.

Phase transitions are classified into \emph{first order} and
\emph{second order} phase transitions.\footnote{This terminology was
introduced by Paul Ehrenfest. The idea is that in first order phase
transitions, the first derivative of the free energy has a
discontinuity, while in a second order phase transitions, the first
derivative of the free energy is continuous and only the second
derivative makes a jump.} Second order phase transitions are also
called \emph{continuous} phase transitions. The phase transitions in
Figures~\ref{fig:surprob} and \ref{fig:Ismag} are both second order,
since the functions $\tet$ and $m_\ast$ are continuous at the critical
points $\la_{\rm c}$ and $\bet_{\rm c}$, respectively. Also, second
order phase transitions are characterized by the fact that at the
critical point, there is only one invariant law. By contrast, if we
would draw the function $m_\ast(\bet)$ of a Potts model for
sufficiently large values of $q$ (in dimension two, for $q>4$), then
the plot of $m_\ast$ would make a jump at $\bet_{\rm c}$ and the
system would have multiple invariant laws at this point, which means
that this phase transition is first order.

It can be difficult to prove that a given phase transition is first or
second order. While for the two-dimensional Ising model, continuity of the
magnetization follows from Onsager's solution \cite{Ons44}, the analogous
statement for the three-dimensional Ising model was only proved recently
\cite{ADS15} (70 years after Onsager!).

For the Ising model, it is known (but only partially proved) that
\[
m_\ast(\bet)\propto(\bet-\bet_{\rm c})^c\quad\mbox{as }\bet\down\bet_{\rm c},
\]
where $c$ is a \emph{critical exponent}, which is given
by\footnote{This exponent is traditionally called $\bet$. The values
given here are taken from the Wikipedia page \emph{Ising critical
exponents} (retrieved 9.2.2025).}
\[
c=1/8\mbox{ in dim 2},\quad
c\approx 0.326\mbox{ in dim 3},\quand
c=1/2\mbox{ in dim $\geq 4$.}
\]
For the contact process, it has numerically been observed that
\[
\tet(\la)\propto(\la-\la_{\rm c})^c\qquad\mbox{as }\la\down\la_{\rm c},
\]
with a critical exponent \cite[Table~4.3]{HHL08}
\[\begin{array}{l}
\dis c\approx 0.276\mbox{ in dim 1},\quad
c\approx 0.583\mbox{ in dim 2},\\[5pt]
c\approx 0.813\mbox{ in dim 3},\quand
c=1\mbox{ in dim $\geq 4$.}
\end{array}\]
In theoretical physics, (nonrigorous) \emph{renormalization group theory}
\index{renormalization group theory} is used to explain these critical
exponents and calculate them. According to this theory, critical
exponents are \emph{universal}. For example, the nearest-neighbor
model and the range $R$ models with different values of $R$ all have
different values of the critical point, but the critical exponent $c$
has the same value for all these models.\footnote{Universality in the
range $R$ does not always hold.  It has been proved that the $q=3$
ferromagnetic Potts model in dimension two has a first order phase
transition for large $R$ \cite{GB07}, while the model with $R=1$ is
known to have a second order phase transition \cite{DST17}.}  Also,
changing from the square lattice to, for example, the triangular
lattice has no effect on $c$.

Critical exponents are associated only with second order phase transitions. At
the critical point of a second order phase transition, one observes
\emph{critical behavior},\index{critical behavior} which involves, for example,
power-law decay of correlations. For this reason, physicists use the term
``critical point'' only for second order phase transitions.

So far, there is no mathematical theory that can explain critical
behavior, except in high dimensions (where one uses a technique called
the \emph{lace expansion}) and in a few two-dimensional models (that
have a conformally invariant scaling limit that can be described using the
Schramm-Loewner equation).

\section{Variations on the voter model}

Apart from the models discussed so far, lots of other interacting particle
systems have been introduced and studied in the literature to model a
plethora of phenomena. Some of these behave very similarly to the models we
have already seen (and even appear to have the same critical exponents), while
others are completely different. In this and the next sections, we take a brief
look at some of these models to get an impression of the possibilities.

\begin{figure}[htb]
\begin{center}
\inputtikz{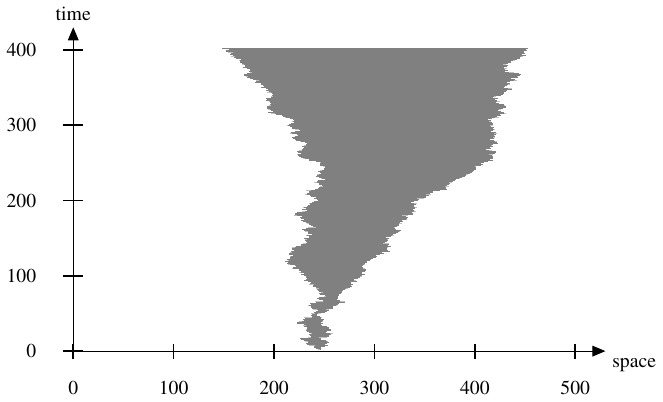}
\caption{Time evolution of a one-dimensional biased voter model with
 bias $s=0.2$.}
\label{fig:biasvoter}
\commentAlt{Figure~\ref{fig:biasvoter}}{Space is plotted horizontally
  and time vertically. An interval of dark grey color grows in size
  and is bounded by two random walk paths that look like drifted
  Brownian motions.}
\end{center}
\end{figure}

The \emph{biased voter model}\index{biased voter model} with \emph{bias}
$s\geq 0$ is the interacting particle system with state space $\{0,1\}^{\Z^d}$
and generator (compare (\ref{Gvot}))
\begin{equation}\begin{array}{r@{\,}c@{\,}l}\label{Gbiasvot}
G_{\rm bias}f(x)
&:=&\dis\frac{1}{|\Ni_0|}\sum_{(i,j)\in\Ei^d}
\big\{f\big({\tt vot}_{ij}(x)\big)-f\big(x\big)\big\}\\[5pt]
&&\dis+\frac{s}{|\Ni_0|}\sum_{(i,j)\in\Ei^d}
\big\{f\big({\tt bra}_{ij}(x)\big)-f\big(x\big)\big\},
\end{array}\end{equation}
where ${\tt vot}_{ij}$ and ${\tt bra}_{ij}$ are the voter and branching maps
defined in (\ref{votmap}) and (\ref{bramap}). The biased voter model describes
a situation where one genetic type of an organism (in this case, type~1) is
more fit than the other type, and hence reproduces at a larger
rate. Alternatively, this type may represent a new idea or opinion that is more
attractive than the current opinion. Contrary to the normal voter model, even
if we start with just a single individual of type~1, there is a positive
probability that type~1 never dies out and indeed takes over the whole
population, as can be seen in Figure~\ref{fig:biasvoter}.

Fix $i\in\Z^d$ and for any $x\in\{0,1\}^{\Z^d}$, let
\[
f_\tau(x):=\frac{1}{|\Ni_i|}\sum_{j\in\Ni_i}1_{\txt\{x(j)=\tau\}}
\qquad(\tau=0,1)
\]
be the frequency of type $\tau$ in the neighborhood $\Ni_i$. In the standard
voter model, if the present state is $x$, then the site $i$ changes its type
with the following rates:
\[\begin{array}{ll}
\dis 0\mapsto 1\quad&\dis\mbox{with rate }f_1(x),\\[5pt]
\dis 1\mapsto 0\quad&\dis\mbox{with rate }f_0(x).
\end{array}\]
In the biased voter model, this is changed to
\[\begin{array}{ll}
\dis 0\mapsto 1\quad&\dis\mbox{with rate }(1+s)f_1(x),\\[5pt]
\dis 1\mapsto 0\quad&\dis\mbox{with rate }f_0(x).
\end{array}\]
Another generalization of the voter model, introduced by Neuhauser and
Pacala\index{Neuhauser--Pacala model} in \cite{NP99}, is defined by
the rates
\begin{equation}\begin{array}{ll}\label{NP99}
\dis 0\mapsto 1\quad&\dis\mbox{with rate }
f_1(x)\big(f_0(x)+\al f_1(x)\big),\\[5pt]
\dis 1\mapsto 0\quad&\dis\mbox{with rate }
f_0(x)\big(f_1(x)+\al f_0(x)\big),
\end{array}\end{equation}
where $0\leq\al\leq 1$ is a model parameter. Another way of expressing this is
to say that if the individual at $i$ is of type $\tau$, then this individual
dies with rate
\begin{equation}\label{rebde}
f_\tau(x)+\al f_{1-\tau}(x),
\end{equation}
and once an individual has died, just as in the normal contact process, it is
replaced by a descendant of a uniformly chosen neighbor.

If $\al=1$, then the rate of dying in (\ref{rebde}) is one and we are back at
the standard voter model, but for $\al<1$, individuals die less often if they
are surrounded by a lot of individuals of the other type. In biology, this
models \emph{balancing selection}. This is the effect that individuals that
differ from their neighbors experience less competition, which results in a
selective drive for high biodiversity.

In the social interpretation of the voter model, we may interpret
(\ref{rebde}) as saying that persons change their mind \emph{less} often if
they disagree with a lot of neighbors, that is, the model in (\ref{NP99}) has
``rebellious'' behavior.

\begin{figure}[htb]
\begin{center}
\inputtikz{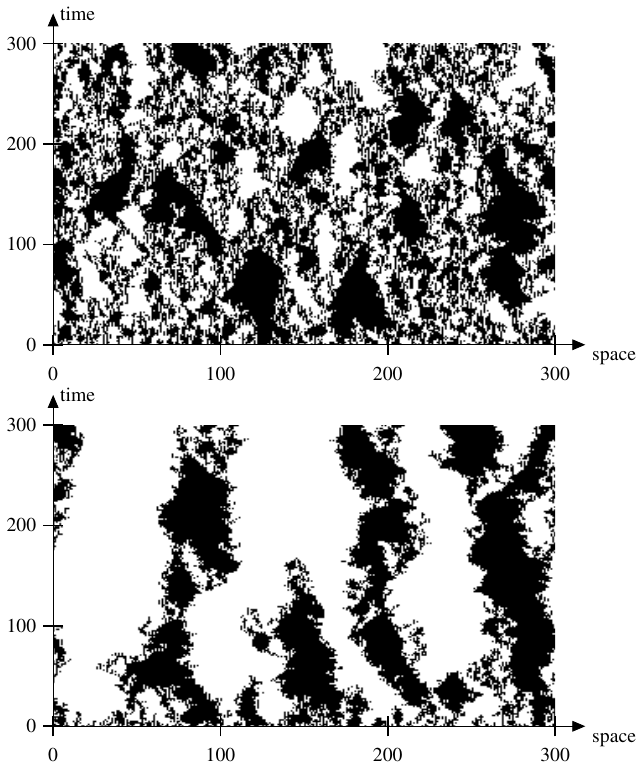}
\caption{Evolution of the Neuhauser--Pacala model with $R=2$ and $\al=0.2$
(top picture) and $\al=0.5$ (bottom picture).}
\label{fig:rebel}
\commentAlt{Figure~\ref{fig:rebel}}{The bottom picture looks like the
  voter model in Figure~\ref{fig:Potvot1D} with only two colors, black
  and white, and with more noisy boundaries. The top picture shows an
  apparent stationary state of smaller and larger black and white
  patches.}
\end{center}
\end{figure}


Numerical simulations, shown in Figure~\ref{fig:rebel}, suggest that
in one dimension, the model in (\ref{NP99}) with range $R\geq 2$
exhibits a phase transition in $\al$. For $\al$ sufficiently close to
1, the model behaves essentially as a voter model, with clusters
growing in time, but for small values of $\al$ (which represent strong
rebellious behavior), the cluster size tends to a finite limit. The
latter statement has been rigorously proved, but proving voter-like
behavior for $\al$ sufficiently close to one is an open problem.

\section{The exclusion process}\label{S:excl}

The exclusion process is a model for traffic or other forms of
transport. The local state space is $S=\{0,1\}$. Sites $i\in\La$ with
$x(i)=1$ are interpreted as being \emph{occupied} by a particle. Sites
with $x(i)=0$ are \emph{empty}. For each $i,j\in\La$ with $i\neq j$,
we define an \emph{asymmetric exclusion map} ${\tt asep}_{ij}\cn
S^\La\to S^\La$ by
\begin{equation}\label{asep}\index{0asep@${\tt asep}_{ij}$}
{\tt asep}_{ij}(x)(k):=\left\{\begin{array}{ll}
0\quad&\mbox{if $k=i$ and $x(j)=0$,}\\[5pt]
x(i)\quad&\mbox{if $k=j$ and $x(j)=0$,}\\[5pt]
x(k)\quad&\mbox{otherwise.}
\end{array}\right.
\end{equation}
Applying ${\tt asep}_{ij}$ to a configuration $x$ has the effect that
if there is a particle at $i$ and the site $j$ is empty, then the
particle at $i$ jumps to $j$. If there is no particle at $i$ or the
site $j$ is already occupied, then nothing happens. Note that these
dynamics preserve the number of particles. The one-dimensional lattice
$\La=\Z$ is of particular interest. The \emph{asymmetric simple
exclusion process} (ASEP)\index{exclusion process}\index{ASEP}
on $\Z$ with parameter $p\in[0,1]$ is the interacting particle system
with generator
\begin{equation}\begin{array}{r@{\,}c@{\,}l}\label{Gasep}
\dis G_{\rm asep}f(x)&:=&\dis(1-p)\sum_{i\in\Z}
\big\{f\big({\tt asep}_{i,i-1}\big)-f\big(x\big)\big\}\\[5pt]
&&\dis+p\sum_{i\in\Z}\big\{f\big({\tt asep}_{i,i+1}\big)-f\big(x\big)\big\}.
\end{array}
\end{equation} 
The process with $p=1$ is called the \emph{totally asymmetric simple
exclusion process} (TASEP)\index{TASEP} and the process with $p=\ha$ is
called the \emph{symmetric exclusion process} (SEP), or simply the
\emph{exclusion process}. The latter can alternatively also be defined
in a different way. For each $i,j\in\La$, we define an \emph{exclusion
map} \index{exclusion map} ${\tt excl}_{ij}\cn S^\La\to S^\La$ by
\begin{equation}\label{exclmap}\index{0excl@${\tt excl}_{ij}$}
{\tt excl}_{ij}(x)(k):=\left\{\begin{array}{ll}
x(j)\quad&\mbox{if }k=i,\\[5pt]
x(i)\quad&\mbox{if }k=j,\\[5pt]
x(k)\quad&\mbox{otherwise.}
\end{array}\right.
\end{equation}
Applying ${\tt excl}_{ij}$ to a configuration $x$ has the effect of
interchanging the types of $i$ and $j$. The interacting particle
system with state space $\{0,1\}^{\Z^d}$ and generator
\begin{equation}\label{Gexcl}
G_{\rm excl}f(x)
=\frac{1}{|\Ni_0|}\sum_{\{i,j\}\in E^d}
\big\{f\big({\tt excl}_{ij}(x)\big)-f\big(x\big)\big\}
\qquad(x\in\{0,1\}^{\Z^d})
\end{equation}
is called the (symmetric) \emph{exclusion process} on
$\Z^d$.\index{exclusion process}\index{SEP} One can check that in the
one-dimensional case, this is the same process as the one with the
generator in (\ref{Gasep}) for $p=\ha$. This follows from the fact
that in both processes, the same transitions happen at the same
rates. Indeed, from the point of view of how many particles there are
on each site, if a particle tries to jump to an already occupied site,
then it does not matter if the jump does not take place or the two
particles interchange their positions. Mathematically, the equality
(in law) of both processes follows from the fact that setting $p=\ha$
in (\ref{Gasep}) and $d=1$ in (\ref{Gexcl}), one has that
$G_{\rm asep}f=G_{\rm excl}f$ for all functions $f\cn\{0,1\}^\Z\to\R$
that depend on finitely many coordinates. This will be proved
rigorously in Sections \ref{S:Poiscon} and \ref{S:Gencon} of
Chapter~\ref{C:construct}. In the symmetric exclusion process,
individual particles move according to random walks, that are
independent as long as the particles are sufficiently far
apart. Particles never meet, and the total number of particles is
preserved.

If the totally asymmetric simple exclusion process is started in a
deterministic initial state, then its distribution at any later time
is a determinantal point process. This means that TASEP is one of the
rare examples of an interacting particle system that is, in some
sense, explicitly solvable. There are close connections between TASEP,
the so-called KPZ universality class (after the Kardar–Parisi–Zhang
equation), and random matrix theory. For this reason, TASEP is one of
the most studied interacting particle systems, see \cite{Fer13,MQ17}.

\section{Branching and coalescing particles}

For each $i,j\in\Z^d$, we define a \emph{coalescing random walk map}
\index{coalescing random walk map}
${\rm rw}_{ij}\cn\{0,1\}^{\Z^d}\to\{0,1\}^{\Z^d}$ by
\begin{equation}\label{rwmap}\index{0rw@${\tt rw}_{ij}$}
{\tt rw}_{ij}(x)(k):=\left\{\begin{array}{ll}
0\quad&\mbox{if }k=i,\\[5pt]
x(i)\vee x(j)\quad&\mbox{if }k=j,\\[5pt]
x(k)\quad&\mbox{otherwise.}
\end{array}\right.
\end{equation}
Applying ${\tt rw}_{ij}$ to a configuration $x$ has the effect that if the
site $i$ is occupied by a particle, then this particle jumps to the site $j$.
If there is already a particle at $j$, then the two particles coalesce.

\begin{figure}[htb]
\begin{center}
\inputtikz{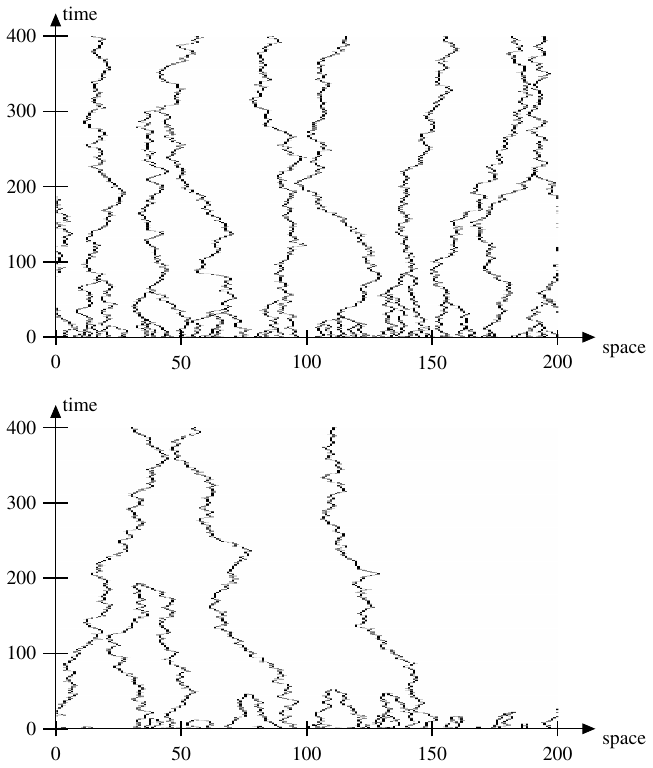}
\caption{Systems of coalescing random walks (above) and annihilating random
  walks (below).}
\label{fig:coalannih}
\commentAlt{Figure~\ref{fig:coalannih}}{Random walk paths with time
  running upwards. When two paths meet, in the top picture they
  continue as one while in the bottom picture both paths disappear. As
  a result, the density is reduced, initially quickly, then more
  slowly.}
\end{center}
\end{figure}

The interacting particle system with generator
\begin{equation}\label{Grw}
G_{\rm rw}f(x)
=\frac{1}{|\Ni_0|}\sum_{(i,j)\in\Ei^d}
\big\{f\big({\tt rw}_{ij}(x)\big)-f\big(x\big)\big\}
\qquad(x\in\{0,1\}^{\Z^d})
\end{equation}
describes a system of coalescing random walks, where each particle jumps with
rate 1 to a uniformly chosen neighboring site, and two particles on the same
site coalesce; see Figure~\ref{fig:coalannih}. Likewise, replacing the 
coalescing random walk map by the \emph{annihilating random walk map}
\index{annihilating random walk map} defined as
\begin{equation}\label{annmap}\index{0ann@${\tt arw}_{ij}$}
{\tt arw}_{ij}(x)(k):=\left\{\begin{array}{ll}
0\quad&\mbox{if }k=i,\\[5pt]
x(i)+x(j)\quad{\rm mod}(2)\quad&\mbox{if }k=j,\\[5pt]
x(k)\quad&\mbox{otherwise,}
\end{array}\right.
\end{equation}
yields a system of annihilating random walks, that kill each other as
soon as two particles land on the same site; see
Figure~\ref{fig:coalannih}. If a system of one-dimensional
nearest-neighbor coalescing or annihilating random walks is started in
a deterministic initial state, then its law at any positive time is a
Pfaffian point process \cite{GP+18}. Thus, coalescing or annihilating
random walks on $\Z$ are in some sense solvable, similar to TASEP.

\begin{figure}[htb]
\begin{center}
\inputtikz{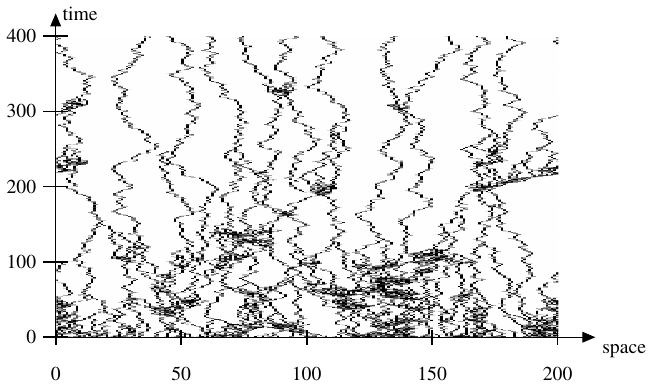}
\caption{A one-dimensional interacting particle system
with cooperative branching and coalescing random walk dynamics.}
\label{fig:coop}
\commentAlt{Figure~\ref{fig:coop}}{A rather messy picture showing
  coalescing random walks as well as occasional clusters containing
  many paths that through coalescence are then quicky reduced to just
  a few.}
\end{center}
\end{figure}

The previous two maps (the coalescing and annihilating random walk
maps) as well as the exclusion map can be combined with, for example,
the branching map and death map from (\ref{bramap}) and
(\ref{deathmap}). In particular, adding coalescing random walk or
exclusion dynamics to a contact process models displacement
(migration) of organisms. Since in many organisms, you actually need
two parents to produce offspring, several authors
\cite{Nob92,Dur92,Neu94,SS15} have studied particle systems where the
branching map is replaced by the \emph{cooperative branching map}
\index{cooperative branching map}
\begin{equation}\label{coopmap}\index{0coop@${\tt coop}_{ii'j}$}
{\tt coop}_{ii'j}(x)(k):=\left\{\begin{array}{ll}
1\quad&\mbox{if }k=j,\ x(i)=1,\ x(i')=1,\\[5pt]
x(k)\quad&\mbox{otherwise.}
\end{array}\right.
\end{equation}
See Figure~\ref{fig:coop} for a one-dimensional interacting particle system
involving cooperative branching and coalescing random walks.

\begin{figure}[htb]
\begin{center}
\inputtikz{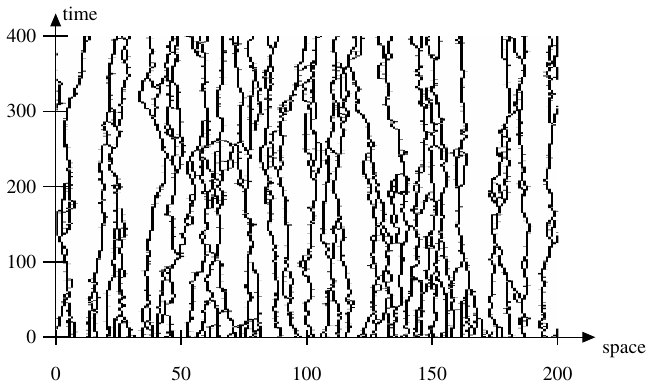}
\caption{A system with branching and killing.}
\label{fig:kill}
\commentAlt{Figure~\ref{fig:kill}}{From a distance this looks like
  branching and coalescing random walks. Upon closer inspection, there
  are small spikes on the paths which are caused by particles creating
  a new particle that is almost immediately killed.}
\end{center}
\end{figure}

We define a \emph{killing map}\index{killing map} by
\begin{equation}\label{killmap}\index{0kill@${\tt kill}_{ij}$}
{\tt kill}_{ij}(x)(k):=\left\{\begin{array}{ll}
0\quad&\mbox{if }k=j,\ x(i)=1,\ x(j)=1,\\[5pt]
x(k)\quad&\mbox{otherwise.}
\end{array}\right.
\end{equation}
In words, this says that if there are particles at $i$ and $j$, then
the particle at $i$ kills the particle at $j$. Sudbury
\cite{Sud97,Sud99} has studied a ``biased annihilating branching
process'' with generator of the form\index{biased annihilating branching process}
\begin{equation}\begin{array}{r@{\,}c@{\,}l}\label{BABP}
\dis G_{\rm babp}f(x) &:=&\dis\la\sum_{(i,j)\in\Ei^1}
\big\{f\big({\tt bra}_{ij}(x)\big)-f\big(x\big)\big\}\\[5pt]
&&\dis+\sum_{(i,j)\in\Ei^1} \big\{f\big({\tt kill}_{ij}(x)\big)-f\big(x\big)\big\}
\qquad(x\in\{0,1\}^\Z).
\end{array}\end{equation}
In the physics literature, this model is known as the
Fredrickson--Andersen one spin facilitated
model,\index{Fredrickson--Andersen model} see formula (26) in
\cite{RS03} (with $f=1$). It is part of the class of \emph{kinetically
constrained models}\index{kinetically constrained models} \cite{HT25}.
In the mathematical literature on this subject, a slight variant of
the model has been studied \cite{BDT19}. Figure~\ref{fig:kill} shows a
simulation of such a system when $\la=0.2$. When $\la$ is small, in
the simulations, the process seems to behave similarly to systems of
branching and coalescing random walks.

\section{Periodic behavior}\label{S:period}

The previous sections served to give a short introduction to some of
the main lines of research in interacting particle systems and to
showcase how much is known. By contrast, the present section is about
a subject that is certainly not, at present, a main line of research
and mainly serves to demonstrate how much there still is that we know
very little about.

An invariant law of an interacting particle system is a probability
distribution on the space $S^\La$ of all possible configurations with
the property that if the system at time zero is distributed according
to this law, then at all later times it is also distributed according
to this law. Invariant laws need not be unique. For example, Potts
models above the critical point have $q$ different invariant laws,
that are characterized by the color that occupies the majority of the
sites.

For all the interacting particle systems and initial states we have
considered so far, the system has the property that as time tends to
infinity, the distribution of the system converges to an invariant
law. This need not always be the case. Perhaps the simplest way in
which this can fail is if the system has a \emph{periodic law}, that
is, a law that has the property that if the system at time zero is
distributed according to this law, then it returns to this law after a
finite time $T>0$ (the \emph{period}), but the system has a different
distribution at all intermediate times $0<t<T$.

Very little is known rigorously about interacting particle systems
with periodic laws. Jahnel and K\"ulske \cite{JK14} have constructed a
three dimensional interacting particle system that has a periodic
law. A general result due to Mountford \cite{Mou95} implies that one
dimensional systems with finite range interactions cannot have
periodic laws. Beyond this, almost nothing is known rigorously for
spatial models, although there are some studies of periodic behavior
in the mean-field limit (see Chapter~\ref{C:meanfield}). In
particular, it is not known whether periodic laws are possible in two
dimensions. The construction in \cite{JK14} is rather abstract since
they do not write down the dynamics of their system explicitly but
only prove that such a system exists. Their system also does not have
finite range interactions, although the strength of the interaction
decays exponentially in the distance, which is almost as good.

Numerical simulations suggest that periodic behavior is not a rare
phenomenon. Several interacting particle systems with explicit
dynamics are known to exhibit periodic behavior in simulations. All
known examples seem to work only in dimensions three and higher,
however, which suggests that, perhaps, periodic behavior is not
possible in two dimensions. This would be in line with (though not
rigorously follow from) the Mermin--Wagner theorem from statistical
physics that states, informally, that continuous symmetry breaking is
not possible in two dimensions \cite{MW66}.

Some mechanisms that can lead to periodic behavior are investigated in
\cite{DFR13,CFT16,Tov19,CDFT20}. A particularly simple model that
numerically seems to exhibit periodic behavior is the \emph{cycle
conform model}\index{cycle conform model} that we describe now. Its
local state space is $S=\{0,1,2\}$ and it is based on two maps, a
\emph{cycle map} and a \emph{conform map}. For each $i\in\La$, the
cycle map ${\tt cyc}_i\cn S^\La\to S^\La$ is defined by
\begin{equation}\index{0coop@${\tt cyc}_i$}
{\tt cyc}_i(x)(k):=\left\{\begin{array}{ll}
x(i)+1\quad\mbox{mod}(3)\quad&\mbox{if }k=i,\\[5pt]
x(k)\quad&\mbox{otherwise,}
\end{array}\right.
\end{equation}
and for each $i,i',j\in\La$, the conform map ${\tt conf}_{ii'j}\cn
S^\La\to S^\La$ is defined by
\begin{equation}\index{0coop@${\tt conf}_{ii'j}$}
{\tt conf}_{ii'j}(x)(k):=\left\{\begin{array}{ll}
x(i)\quad&\mbox{if }k=j\mbox{ and }x(i)=x(i'),\\[5pt]
x(k)\quad&\mbox{otherwise.}
\end{array}\right.
\end{equation}
We assume that the lattice $\La$ has the structure of a graph with set
of edges $E$ and that each site has at least two neighbors. For each
$j\in\La$, we set (compare (\ref{Ni}))
\begin{equation}\label{Ni2}\index{0Ni2@$\Ni^2_i$}
\Ni^2_j:=\big\{(i,i'):i,i'\in\Ni_j,\ i\neq i'\big\}.
\end{equation}
The \emph{cycle conform model} with parameter $\al\in[0,1]$ is the
interacting particle system with generator
\begin{equation}\begin{array}{r@{\,}c@{\,}l}
\dis G_{\rm cc}f(x)&:=&\dis(1-\al)\sum_{j\in\La}
\big\{f\big({\tt cyc}_j(x)\big)-f\big(x\big)\big\}\\[5pt]
&&\dis+\al\sum_{j\in\La}\frac{1}{|\Ni^2_j|}\sum_{(i,i')\in\Ni^2_j}
\big\{f\big({\tt conf}_{ii'j}(x)\big)-f\big(x\big)\big\}.
\end{array}\end{equation}
In words, the dynamics can be described as follows. Each site
$j\in\La$ becomes active at rate one. With probability $1-\al$ the
site $j$ cycles, and with probability $\al$ it conforms. If the site
cycles, then it just changes its type to the next type modulo 3. If
the site $j$ conforms, then it samples two neighboring sites $i$ and
$i'$ at random, and if these happen to have the same type, then the
site $j$ copies their type. The idea of this sort of dynamics is to
give sites a tendency to conform to the type that is locally in the
majority. If we would sample just one neighbor $i$, then we would
obtain voter model dynamics which gives each type the same chance to
spread. By sampling two sites, we introduce a nonlinearity that favors
conformation to the local majority.

\begin{figure}
\begin{center}
\inputtikz{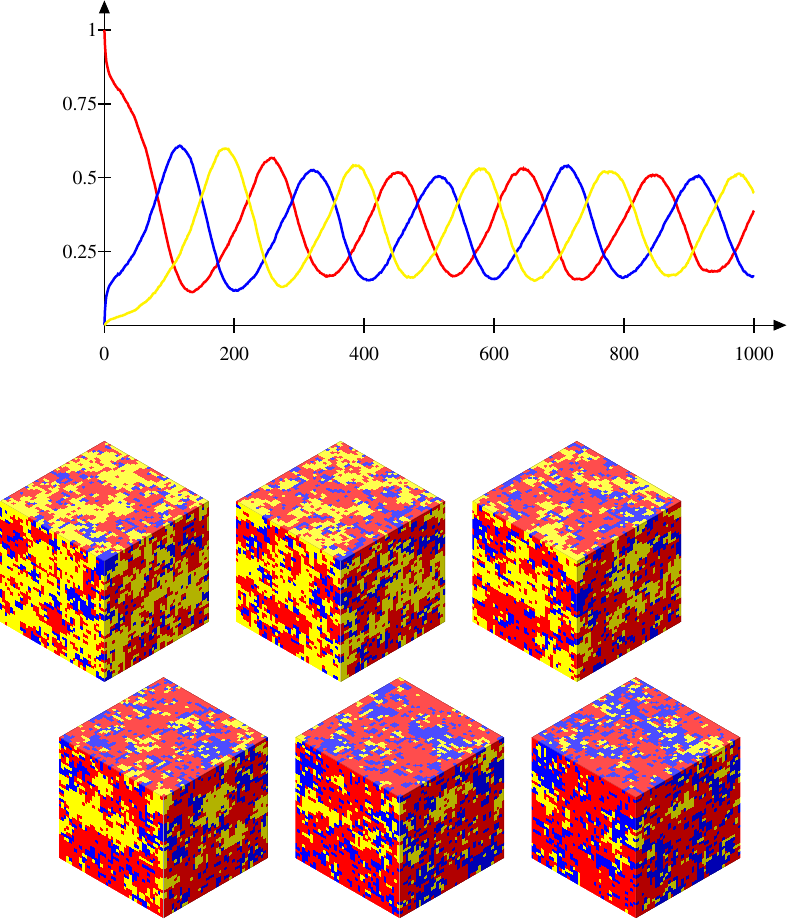}
\end{center}
\caption{Periodic behavior of the cycle conform model with $\al=0.915$
  on $\Z^3$ with nearest-neighbor edges. Simulation on a cube with
  sidelength 40 and periodic boundary conditions. Colors indicate the
  three states $0,1,2$. Shown are the frequency of each type as a
  function of time and the state of the cube at times 601, 613, 626,
  638, 651, and 663. Numerically, one sees periodic behavior roughly
  in the regime $0.903<\al<0.933$.}
\label{fig:cyconf}
\commentAlt{Figure~\ref{fig:cyconf}}{Three graphs of sinus-like
  functions in the colors red, blue, and yellow. Below that six cubes
  in which the frequency of colors reflects the graphs. See long
  description.}
\commentLongAlt{Figure~\ref{fig:cyconf}}{Three graphs in the colors
  red, blue, and yellow plotted together. The red graph starts at one
  and quickly drops. The yellow and blue graphs start at zero. Blue
  quickly increases, yellow more slowly. Soon, the graphs look like
  sinus functions that oscillate between 0.15 and 0.55. First the blue
  graph has a maximum, then yellow, then red, then blue again and so
  on. Below six cubes consisting of many small cubes in the colors
  red, blue, and yellow. In the first cube yellow dominates, red comes
  second, and blue is rare. Progressively yellow becomes rarer, red
  becomes more common, and and the end also blue starts to
  increase. In the final picture red dominates, blue comes second, and
  yellow is rare. The colors are not randomly (independently)
  distributed but come in small clusters of the same color.}
\end{figure}

Numerical simulations on $\Z^d$ in dimensions $d=1,2,3$ suggest the
following picture. In dimensions 1 and 2 there is a unique invariant
law for each $0<\al<1$. In dimension 3 there are two
critical values $0<\al_{\rm c}<\al'_{\rm c}<1$. For $\al<\al_{\rm c}$,
there is a unique invariant law. For $\al>\al'_{\rm c}$, there are three
invariant laws, that are moreover invariant under translations, in
which one of the three local states has a majority. In the
intermediate regime $\al_{\rm c}<\al<\al'_{\rm c}$ the system exhibits
periodic behavior, see Figure~\ref{fig:cyconf}.

\chapter{Continuous-time Markov chains}\label{C:Markov}

\section{Finite state space}\label{S:finMark}

Let $S$ be any finite set. A (real) matrix indexed by $S$ is a
collection of real constants $A=(A(x,y))_{x,y\in S}$. We calculate
with such matrices in the same way as with normal square
matrices. Thus, the product $AB$ of two matrices is defined as
\[
(AB)(x,z):=\sum_{y\in S}A(x,y)B(y,z)\qquad(x,z\in S).
\]
We let $1$ denote the identity matrix which has $1(x,y)=1$ if $x=y$
and $=0$ if $x\neq y$. We will sometimes denote this as
$1(x,y)=1_{\{x=y\}}$, where quite generally we let
$1_{\{\ldots\}}$\index{00yyz@$1_{\{\ldots\}}$} denote the indicator of
the event $\{\ldots\}$, that is, this is one if the conditions inside
$\{\ldots\}$ are satisfied and zero otherwise. We define the $n$-th
power $A^n$ of a matrix indexed by $S$ in the obvious way, with
$A^0:=1$. If $f\cn S\to\R$ is a function, then we also define
\begin{equation}\label{fAf}
Af(x):=\sum_{y\in S}A(x,y)f(y)
\quand
fA(y):=\sum_{x\in S}f(x)A(x,y).
\end{equation}
A \emph{matrix semigroup}\index{matrix semigroup} is a collection of
matrices $(A_t)_{t\geq 0}$ such that
\[
\lim_{t\down 0}A_t=A_0=1\quand A_sA_t=A_{s+t}\qquad(s,t\geq 0).
\]
If $G$ is a matrix indexed by $S$, then setting
\[
A_t=\ex{tG}:=\sum_{n=0}^\infty\frac{1}{n!}(tG)^n\qquad(t\geq 0).
\]
defines a matrix semigroup $(A_t)_{t\geq 0}$, and conversely every
matrix semigroup $(A_t)_{t\geq 0}$ is of this form. The matrix $G$ is
called the \emph{generator}\index{generator!of a semigroup} of
$(A_t)_{t\geq 0}$. The function $t\mapsto A_t$ is continuously
differentiable and one has
\begin{equation}\label{difAt}
\dif{t}A_t=GA_t=A_tG\qquad(t\geq 0)
\end{equation}

A \emph{probability kernel}\index{probability kernel} on $S$ is a
matrix $K=(K(x,y))_{x,y\in S}$ such that $K(x,y)\geq 0$ $(x,y\in S)$
and $\sum_{y\in S}K(x,y)=1$ $(x\in S)$. Clearly, the composition of
two probability kernels yields a third probability kernel. A
probability kernel is \emph{deterministic}
\index{probability kernel!deterministic} if it is of the form
\[
K_m(x,y):=\left\{\begin{array}{ll}
1\quad&\mbox{if }m(x)=y,\\[5pt]
0\quad&\mbox{otherwise,}
\end{array}\right.
\]
for some function $m\cn S\to S$. It is easy to see that the space of
all probability kernels on a finite set $S$ is convex, and the
deterministic probability kernels are exactly the extremal elements of
this set. It follows that each probability kernel can be written as a
convex combination of deterministic probability kernels. Another way
to say this is that for each probability kernel $K$ on $S$, it is
possible to find a random map $M\cn S\to S$ such that\footnote{Indeed,
this formula says nothing else than $K=\sum_m\P[M=m]K_m$, where the
sum runs over all maps $m\cn S\to S$, and $K_m$ is the deterministic
kernel defined by the map $m$.}
\begin{equation}\label{Krep}
K(x,y)=\P\big[M(x)=y\big]\qquad(x,y\in S).
\end{equation}
A formula of this form is called a \emph{random mapping
representation}\index{random mapping representation!of a probability kernel}
of the probability kernel $K$.

A \emph{Markov semigroup}\index{Markov semigroup} is a matrix
semigroup $(P_t)_{t\geq 0}$ consisting of probability kernels.

\begin{Exercise}[Markov generators]
Let $G$ be a matrix indexed by a finite set $S$. Show that $G$
generates a Markov semigroup if and only if
\begin{equation}\label{Gform}
G(x,y)\geq 0\quad(x\neq y)\quand\sum_yG(x,y)=0.
\end{equation}
\end{Exercise}

By definition, we say that a function $f$ that is defined on an
interval $I\sub\R$ is \emph{piecewise constant}\index{piecewise constant}
if each compact subinterval of $I$ can be divided into
finitely many subintervals, so that $f$ is constant on each
subinterval. By definition, a \emph{Markov process} with semigroup
$(P_t)_{t\geq 0}$ is a stochastic process $X=(X_t)_{t\geq 0}$ with
values in $S$ and piecewise constant, right-continuous sample paths,
such that
\begin{equation}\label{Mark}
\P\big[X_u\in\,\cdot\,\big|\,(X_s)_{0\leq s\leq t}\big]
=P_{u-t}(X_t,\,\cdot\,)\quad{\rm a.s.}
\qquad(0\leq t\leq u).
\end{equation}
Here, on the left-hand side, we condition on the \si-field generated
by the random variables $(X_s)_{0\leq s\leq t}$. One can prove that
formula (\ref{Mark}) is equivalent to the statement that
\begin{equation}\label{Markfdd}
\P\big[X_0=x_0,\ldots,X_{t_n}=x_n\big]=\P[X_0=x_0]P_{t_1-t_0}(x_0,x_1)\cdots P_{t_n-t_{n-1}}(x_{n-1},x_n)
\end{equation}
for all $0=t_0<t_1<\cdots<t_n$. From this last formula, we see that
for each initial law $\P[X_0=\,\cdot\,]=\mu$, there is a unique Markov
process with semigroup $(P_t)_{t\geq 0}$ and this initial law. We say
that $(P_t)_{t\geq 0}$ are the \emph{transition
kernels}\index{transition kernel} of the Markov process. It is custom
to let $\P^\mu$ denote the law of the Markov process with initial law
$\mu$, and to let $\P^x$\index{0Px@$\P^x$} denote the law of the
Markov process with deterministic initial state $X_0=x$ a.s. We let
$\E^\mu$ and $\E^x$\index{0Ex@$\E^x$} denote expectation with respect
to $\P^\mu$ and $\P^x$, respectively. Recalling our notation
(\ref{fAf}), we see that for any probability law $\mu$ on $S$ and
function $f\cn S\to\R$,
\[\begin{array}{r@{\,}c@{\,}l}
\dis\mu P_t(x)&=&\dis\P^\mu[X_t=x],\\[5pt]
\dis P_tf(x)&=&\dis\E^x[f(X_t)]
\end{array}\qquad(t\geq 0,\ x\in S).
\]
In particular, $\mu P_t$ is the law of the process at time $t$. We note that
\[
\P^x[X_t=y]=P_t(x,y)=1_{\{x=y\}}+tG(x,y)+O(t^2)\quad\mbox{as }t\down 0.
\]
For $x\neq y$, we call $G(x,y)$ the \emph{rate}\index{rate of jumps}
of jumps from $x$ to $y$. Intuitively, if the process is in $x$, then
in the next infinitesimal time interval of length $\di t$ it has a
probability $G(x,y)\di t$ to jump to $y$, independently for all $y\neq
x$.

Let $S$ be a finite set. Let $\Ki$ be a finite set whose elements are
probability kernels on $S$ and let $(r_K)_{K\in\Ki}$ be nonnegative
constants. Then it is straightforward to check that setting
\begin{equation}\label{GK}
Gf:=\sum_{K\in\Ki}r_K\,\big\{Kf-f\big\}
\end{equation}
defines a Markov generator. The following exercise says that
conversely, each Markov generator can be written in this form, where
we can even choose the set $\Ki$ so that it has only one element.

\begin{Exercise}
Let $S$ be a finite set. Show that each Markov generator $G$ on $S$
can be written in the form $Gf=r\{Kf-f\}$, where $r\geq 0$ is a
constant and $K$ is a probability kernel on $S$. \emph{Hint:} first
add a multiple of the identity matrix to $G$ to make all diagonal
entries nonnegative and then normalize.
\end{Exercise}

If all kernels in the set $\Ki$ are deterministic, then our expression
(\ref{GK}) for $G$ takes the form
\begin{equation}\label{Grep0}
Gf(x)=\sum_{m\in\Gi}r_m\big\{f\big(m(x)\big)-f\big(x\big)\big\},
\end{equation}
where $\Gi$ is a finite set whose elements are maps $m\cn S\to S$ and
$(r_m)_{m\in\Gi}$ are nonnegative constants. This way of writing a
generator will come back in formula (\ref{Grep}) of
Section~\ref{S:Poichain}.

If $(X_t)_{t\geq 0}$ is a Markov process with finite state space $S$,
semigroup $(P_t)_{t\geq 0}$, generator $G$, and initial law $\mu$, and
$f\cn S\to\R$ is a function, then $\mu P_tf=\E^\mu[f(X_t)]$ is the
mean of $f$ at time $t$. But what about the variance of $f$? It turns
out that there is a nice formula for this quantity, that is sometimes
useful.\footnote{We will actually only need Proposition~\ref{P:cov}
and its infinite dimensional analogue Proposition~\ref{P:cov2} in
Section~\ref{S:poscor} below, which is not used anywhere else, so the
material here can be skipped at a first reading.} For any probability
measure $\mu$ on $S$ and functions $f,g\cn S\to\R$ we adopt the
notation
\begin{equation}\label{cov}\index{0Ctov@$\cov_\mu(f,g)$}\index{0Var@$\var_\mu(f)$}
\cov_\mu(f,g):=\mu(fg)-(\mu f)(\mu g)\quand\var_\mu(f):=\cov_\mu(f,f).
\end{equation}
In words, $\cov_\mu(f,g)$ is the covariance of $f$ and $g$ under $\mu$
and $\var_\mu(f)$ is the variance of $f$. We define a function
$\Ga_G(f,g)\cn S\to\R$ by\footnote{This is called the carr\'e du champ
operator. Often, a factor $\ha$ is added to its definition because of
its relation to Dirichlet forms.}
\begin{equation}\label{Gadef}\index{0zzC@$\Ga_G(f,g)$}
\Ga_G(f,g)(x):=\sum_{y\in S}G(x,y)\big\{f(y)-f(x)\big\}\big\{g(y)-g(x)\big\}.
\end{equation}
An alternative formula for $\Ga_G(f,g)$ is
\[
\Ga_G(f,g)=G(fg)-(Gf)g-f(Gg).
\]
To see that both formulas are equivalent, we calculate
\[\begin{array}{c}
\dis G(fg)(x)=\sum_yG(x,y)f(y)g(y),\\[5pt]
\dis\big((Gf)g\big)(x)=\sum_yG(x,y)f(y)g(x),\quad\big(f(Gg)\big)(x)=\sum_yG(x,y)f(x)g(y).
\end{array}\]
Using the fact that $\sum_yG(x,y)=0$, this gives
\[\begin{array}{l}
\big[G(fg)-(Gf)g-f(Gg)\big](x)\\[5pt]
\dis\quad=\sum_yG(x,y)\big[f(y)g(y)-f(y)g(x)-f(x)g(y)+f(x)g(x)\big],
\end{array}\]
which agrees with our first formula for $\Ga_G(f,g)$. The following formula is well-known (see, for example, the proofs of \cite[Thm~2.1]{Led00} and \cite[Lemma~3.3]{Jou07}).

\begin{proposition}[Covariance formula]
Let\label{P:cov} $(P_t)_{t\geq 0}$ be the semigroup of a Markov
process with finite state space $S$ and generator $G$. Then for each
probability measure $\mu$ on $S$ and functions $f,g\cn S\to\R$, one
has
\[
\cov_{\mu P_t}(f,g)=\cov_\mu(P_tf,P_tg)+\int_0^t\!\di s\,\mu P_{t-s}\Ga_G(P_sf,P_sg).
\]
\end{proposition}

\begin{Proof}
For each $s\in[0,t]$, define $U_s\cn S\to\R$ by
\[
U_s:=P_{t-s}\big((P_sf)(P_sg)\big).
\]
Then
\[\begin{array}{r@{\,}c@{\,}l}
\dis\mu(U_t-U_0)&=&\dis\mu\big((P_tf)(P_tg)\big)-\mu P_t(fg)\\[5pt]
&=&\dis\mu\big((P_tf)(P_tg)\big)-(\mu P_tf)(\mu P_tg)+(\mu P_tf)(\mu P_tg)-\mu P_t(fg)\\[5pt]
&=&\dis\cov_\mu(P_tf,P_tg)-\cov_{\mu P_t}(f,g),
\end{array}\]
so to complete the proof it suffices to show that
\[
U_t-U_0=\int_0^t\!\di s\,\dif{s}U_s=-\int_0^t\!\di s\,P_{t-s}\Ga_G(P_sf,P_sg).
\]
Using (\ref{difAt}), we see that
\[\begin{array}{r@{\,}c@{\,}l}
\dis\dif{t_1}P_{t_1}\big((P_{t_2}f)(P_{t_3}g)\big)&=&\dis P_{t_1}G\big((P_{t_2}f)(P_{t_3}g)\big),\\[5pt]
\dis\dif{t_2}P_{t_1}\big((P_{t_2}f)(P_{t_3}g)\big)&=&\dis P_{t_1}\big((GP_{t_2}f)(P_{t_3}g)\big),\\[5pt]
\dis\dif{t_3}P_{t_1}\big((P_{t_2}f)(P_{t_3}g)\big)&=&\dis P_{t_1}\big((P_{t_2}f)(GP_{t_3}g)\big).
\end{array}\]
It follows that
\[
\dif{s}U_s=P_{t-s}\big\{\big((GP_sf)(P_sg)\big)+\big((P_sf)(GP_sg)\big)-G\big((P_sf)(P_sg)\big)\big\}
\]
which equals $-P_{t-s}\Ga_G(P_sf,P_sg)$, as required.
\end{Proof}

Let $X=(X_t)_{t\geq 0}$ be a continuous-time Markov process with
finite state space $S$, generator $G$, and semigroup $(P_t)_{t\geq
  0}$. By definition, an \emph{invariant law}\index{invariant law} of
$X$ is a probability measure $\nu$ on $S$ such that
\[
\nu P_t=\nu\qquad(t\geq 0).
\]
This says that if we start the process in the initial law
$\P[X_0\in\,\cdot\,]=\nu$, then $\P[X_t\in\,\cdot\,]=\nu$ for all
$t\geq 0$. By definition, the Markov process is \emph{irreducible}
\index{irreducible Markov chain} if for each $x,y\in S$, there exist
$x_0,\ldots,x_n\in S$ with $x=x_0$ and $y=x_n$, such that
$G(x_{k-1},x_k)>0$ for all $1\leq k\leq n$. The basic result about
invariant laws for continuous-time Markov process with finite state
space is the following theorem. A proof can be found in many places,
such as, for example, \cite[Thm~2.66]{Lig10}.

\begin{theorem}[Convergence to equilibrium]
Let\label{T:finergo} $X=(X_t)_{t\geq 0}$ be a continuous-time Markov
process with finite state space $S$, generator $G$, and semigroup
$(P_t)_{t\geq 0}$. If $X$ is irreducible, then it has a unique
invariant law $\nu$. Moreover, one has
\[
\mu P_t\asto{t}\nu\quad\mbox{for each probability law $\mu$ on $S$.}
\]
\end{theorem}

If $(P_t)_{t\geq 0}$ is Markov semigroup on a finite set $S$ and $\nu$
is an invariant law for $(P_t)_{t\geq 0}$, then it is possible to
construct a process $(X_t)_{t\in\R}$ whose finite-dimensional
distributions are given by
\begin{equation}\label{nustat}
\P\big[X_0=x_0,\ldots,X_{t_n}=x_n\big]=\nu(x_0)P_{t_1-t_0}(x_0,x_1)\cdots P_{t_n-t_{n-1}}(x_{n-1},x_n)
\end{equation}
for all $t_0<\cdots<t_n$. Such a process is
\emph{stationary}.\footnote{Recall that a process $(X_t)_{t\in\R}$ is
stationary if for each $s\in\R$, it is equal in distribution to
$(X'_t)_{t\in\R}$ defined as $X'_t:=X_{t-s}$
$(t\in\R)$.\index{stationary!process}} For this reason, invariant laws
are sometimes called \emph{stationary laws}. By definition, an
invariant law $\nu$ is \emph{reversible}\index{reversible law} if the
stationary process $(X_t)_{t\in\R}$ is equal in law to the
time-reversed process $(X_{-t})_{t\in\R}$. It is well-known
\cite[Exercise~2.44]{Lig10} that this is equivalent to $\nu$
satisfying the \emph{detailed balance}\index{detailed balance}
equations
\[
\nu(x)G(x,y)=\nu(y)G(y,x)\qquad(x,y\in S,\ x\neq y).
\]
The left-hand side of this equation can be interpreted as the
frequency with which the stationary process jumps from $x$ to
$y$. Detailed balance then says that jumps from $x$ to $y$ happen at
the same frequency as jumps from $y$ to $x$.

An irreducible continuous-time Markov process with finite state space
is called \emph{reversible} if its unique invariant law is
reversible. Examples are the Ising and Potts models with Glauber
dynamics on finite lattices, which have Gibbs measures as their
reversible laws. Other examples of processes that have a reversible
law are the Fredrickson--Andersen one spin facilitated model for which
a product measure with a suitably chosen intensity is a reversible
law, and the symmetric exclusion process for which product measures
with arbitrary intensities are reversible laws. Reversibility is a
useful property that allows for the use of techniques that are not
available for irreversible models, such as Dirichlet form techniques
and Poincar\'e or log-Sobolev inequalities. We refer to \cite{Sal97}
as a general introduction to this material and more specifically to
\cite{Mar99} for spin systems with Glauber dynamics. Most interacting
particle systems we will consider in this book are not reversible.

\begin{Exercise}
Let $X=(X_t)_{t\geq 0}$ be a continuous-time Markov process with
finite state space $S$, generator $G$, and semigroup $(P_t)_{t\geq 0}$.
Let $\nu$ be an invariant law and let $(X_t)_{t\in\R}$ be the
stationary process from (\ref{nustat}). Assume that $\nu(x)>0$ for all
$x\in S$. Show that the time-reversed process $(X'_t)_{t\in\R}$
defined as $X'_t:=X_{-t}$ is a stationary Markov process and calculate
its generator $G'$.
\end{Exercise}

\section{The embedded Markov chain}\label{S:embed}

Continuous-time Markov processes with countable state space (also
known as continuous-time Markov chains) can in many ways be treated in
the same way as those with a finite state space, but there are some
complications. The first complication one has to deal with is that
such processes may explode. The second complication is that their
long-time behavior is more complicated than in the finite case: they
may fail to have invariant laws, or have invariant (in particular
reversible) measures that are infinite. For this reason, they need to
be distinguished into positive recurrent, null recurrent, and
transient processes. A good general reference for this material is
\cite[Chapter~2]{Lig10}.

We calculate with matrices indexed by a countably infinite set $S$ in
the same way as for finite $S$, provided the infinite sums are
well-defined (that is, not of the form $\infty-\infty$). Generalizing
our earlier definition, we say that $K$ is a \emph{subprobability
kernel}\index{subprobability kernel} if $\sum_yK(x,y)\leq 1$ for all
$x\in S$. Also when $S$ is infinite, we define generators as in
(\ref{Gform}). Note that
\[
\sum_{y\in S}G(x,y)=G(x,x)+\sum_{y:\,y\neq x}G(x,y).
\]
Since $G(x,y)\geq 0$ for $x\neq y$, the infinite sum on the right-hand
side is always well-defined, though a priori it may be infinite. The
condition $\sum_yG(x,y)=0$ says that it must be finite and equal to
$-G(x,x)$, however. In the special context of continuous-time Markov
chains, a generator is traditionally called a
\emph{Q-matrix}\index{Q-matrix} (and denoted as $Q$) but we will stick
to the term generator.

It is well-known \cite[Section~2.5.2]{Lig10} that one can construct a
continuous-time Markov chain with generator $G$ from its associated
embedded discrete Markov chain and independent, exponentially
distributed holding times. We now recall this construction. Let $G$ be
a generator, let $c(x):=-G(x,x)$ $(x\in S)$, and let $K$ be the
probability kernel on $S$ defined by
\[
K(x,y):=\left\{\begin{array}{ll}
\dis c(x)^{-1}G(x,y)\quad&\mbox{if }c(x)>0,\ x\neq y,\\[5pt]
\dis 1\quad&\mbox{if }c(x)=0,\ x=y,\\[5pt]
\dis 0\quad&\mbox{otherwise.}
\end{array}\right.
\]
For each $x\in S$, let $(Y^x_k)_{k\geq 0}$ be the discrete-time Markov
chain with initial state $Y^x_0=x$ and transition kernel $K$. Set
$N:=\inf\{n\geq 0:c(Y^x_n)=0\}$, which may be infinite, and let
$(\sig_k)_{k\geq 0}$ be i.i.d.\ exponentially distributed random
variables with mean one, independent of $(Y^x_k)_{k\geq 0}$. We define
$(\tau_k)_{0\leq k\leq N+1}$ by
\[
\tau_0:=0,\quad\tau_n:=\sum_{k=0}^{n-1}\sig_k/c(Y^x_k)\quad(1\leq n\leq N+1),
\]
where we use the conventions that $\sig_k/c(Y^x_k):=\infty$ if
$c(Y^x_k)=0$ and $N+1:=\infty$ if $N=\infty$. We set
$\tau:=\tau_{N+1}$. Note that $\tau=\infty$ on the event that
$N<\infty$, but $\tau$ may be finite on the event that $N=\infty$. We
define a stochastic process $(X^x_t)_{t\geq 0}$ with values in
$S_\infty:=S\cup\{\infty\}$ by
\[
X^x_t:=\left\{\begin{array}{ll}
\dis Y^x_k\quad&\mbox{if }t\in[\tau_k,\tau_{k+1}),\ 0\leq k<N+1,\\[5pt]
\dis\infty\quad&\mbox{if }t\geq\tau.
\end{array}\right.
\]
We call $\tau$ the \emph{explosion time}.\index{explosion time} We set
$X^\infty_t:=\infty$ $(t\geq 0)$ and define probability kernels $(\ov
P_t)_{t\geq 0}$ on $S_\infty$ by
\[
\ov P_t(x,y):=\P\big[X^x_t=y\big]\qquad(t\geq 0,\ x,y\in S_\infty).
\]
We call $(X^x_t)_{t\geq 0}$ the \emph{continuous-time Markov
chain}\index{continuous-time Markov chain} with generator $G$ and we
call $(Y^x_k)_{k\geq 0}$ its associated \emph{embedded discrete-time
Markov chain}.\index{embedded Markov chain} It is well-known
\cite[Section~2.5.2]{Lig10} that $(X^x_t)_{t\geq 0}$ is a Markov
process (in the sense of (\ref{Mark}) and (\ref{Markfdd})) with state
space $S_\infty$ and transition kernels $(\ov P_t)_{t\geq 0}$. In
Exercise~\ref{E:Markprop} below, you will be be asked to prove
this. An alternative proof will be suggested in
Exercise~\ref{E:Markprop2}. The random times
\[
\eta_k:=\sig_k/c(Y^x_k)\quad(0\leq k<N)
\]
are called the \emph{holding times}.\index{holding times} Note that
conditional on the embedded chain $(Y^x_k)_{k\geq 0}$, the holding
times $(\eta_k)_{0\leq k<N}$ are independent exponentially distributed
such that $\eta_k$ has mean $1/c(Y^x_k)$. If $\tau=\infty$ a.s.\ for
each initial state $x\in S$, then we say that the continuous-time
Markov chain with generator $G$ is \emph{nonexplosive}. In the
opposite case, it is \emph{explosive}.\index{explosive Markov process}
We let
\begin{equation}\label{subP}
P_t(x,y):=\ov P_t(x,y)\qquad(t\geq 0,\ x,y\in S)
\end{equation}
denote the restrictions of the transition kernels $(\ov P_t)_{t\geq 0}$
to $S$. If $G$ is explosive, then these are only subprobability
kernels.

\begin{Exercise}
Show\label{E:Markprop} that $(X^x_t)_{t\geq 0}$ is a Markov process in
the sense of (\ref{Mark}) with state space $S_\infty$ and transition
kernels $(\ov P_t)_{t\geq 0}$. \emph{Hint:} let $M$ be the number of
jumps of the process $(X^x_s)_{0\leq s\leq t}$, which may be
infinite. Then after conditioning on $(X^x_s)_{0\leq s\leq t}$, you
know $M$ as well as $(Y^x_k)_{0\leq k<M+1}$ and $(\sig_k)_{0\leq k<M}$,
plus in the case that $M<\infty$ you have the information
that $\sum_{k=0}^M\sig_k/c(Y^x_k)>t$. Given all this information, what
do you know about the process $(X^x_u)_{u\geq t}$?
\end{Exercise}

\section{Generator construction}\label{S:genchain}

Let $G$ be the generator of a continuous-time Markov chain with
countable state space $S$. Generalizing our earlier definition to
countable state spaces, we call the collection of subprobability
kernels $(P_t)_{t\geq 0}$ defined in (\ref{subP}) the \emph{Markov
semigroup}\index{Markov semigroup} with generator $G$. In this section
we make a more direct link between $(P_t)_{t\geq 0}$ and $G$.

We calculate with infinite matrices as in the finite case. We observe
that if $f\cn S\to\R$ is nonnegative, then all terms in the infinite
sum $Gf(x):=\sum_{y\in S}G(x,y)f(y)$ except one are nonnegative so
\begin{equation}\label{Gfdef}
Gf\cn S\to(-\infty,\infty]\mbox{ is well-defined for all }f\cn S\to\half.
\end{equation}
Let $u\cn S\times\half\to\half$ be a function. We say that $u$ solves
the \emph{Kolmogorov backward equation}\index{Kolmogorov backward equation}
\index{backward!Kolmogorov equation}
\begin{equation}\label{backeq}
\dif{t}u_t(x)=\sum_yG(x,y)u_t(y)\qquad(t\geq 0,\ x\in S),
\end{equation}
if the function $t\mapsto u_t(x)$ is continuously differentiable for
each $x\in S$ and (\ref{backeq}) holds. By (\ref{Gfdef}), the
right-hand side of (\ref{backeq}) is well-defined, and the equality in
(\ref{backeq}) implies that it must be finite. We say that $u$ is a
\emph{minimal}\index{minimal solution to backward equation} solution
to (\ref{backeq}) if any other solution $u'$ with the same initial
condition $u'_0=u_0$ satisfies $u_t(x)\leq u'_t(x)$ for all $t\geq 0$
and $x\in S$. Note that for a given initial condition, there can be at
most one minimal solution to (\ref{backeq}). The following theorem is
the main result of this section.

\begin{theorem}[Generator construction]
Let\label{T:Gchain} $G$ be the generator of a continuous-time Markov
chain with countable state space $S$ and let $(P_t)_{t\geq 0}$ be the
Markov semigroup with generator $G$. Then for each bounded function
$f\cn S\to\half$, the function
\[
u_t(x):=P_tf(x)\qquad(t\geq 0,\ x\in S)
\]
is the minimal solution to the Kolmogorov backward equation
(\ref{backeq}) with initial condition~$f$.
\end{theorem}

To prepare for the proof of Theorem~\ref{T:Gchain}, as a first step,
we set $c(x):=-G(x,x)$ $(x\in S)$, and we consider the equation
\begin{equation}\label{intKol}
u_t(x)=u_0(x)\ex{-c(x)t}+\int_0^t\!\di s\,\ex{-c(x)s}\sum_{y:\,y\neq x}G(x,y)u_{t-s}(y).
\end{equation}
By definition, a solution to (\ref{intKol}) is a function $u\cn
S\times\half\to\half$ such that $t\mapsto u_t(x)$ is measurable for
all $x\in S$ and (\ref{intKol}) holds for all $t\geq 0$ and $x\in
S$. Note that since $G(x,y)\geq 0$ for $x\neq y$, the sum over $y$ and
consequently also the integral over $s$ are well-defined, even though
a priori the outcome may be $\infty$ (a posteriori, of course,
(\ref{intKol}) implies that the outcome must be finite).

\begin{lemma}[First jump decomposition]
Under\label{L:fijump} the assumptions of Theorem~\ref{T:Gchain}, the
function $u$ solves (\ref{intKol}) with $u_0=f$.
\end{lemma}

\begin{Proof}
Since $f$ is bounded $u_t(x):=P_tf(x)<\infty$ for all $t\geq 0$ and
$x\in S$. (This is the only place in the proof where we use the
boundedness of $f$.) Let $(X^x_t)_{t\geq 0}$ be the continuous-time
Markov chain with generator $G$ and initial state $x\in S$,
constructed from the embedded discrete-time Markov chain
$(Y^x_k)_{k\geq 0}$ and i.i.d.\ standard exponential random variables
$(\sig_k)_{k\geq 0}$ as in the previous section. We extend $f$ to
$S_\infty$ by setting $f(\infty):=0$ so that
$P_tf(x)=\E[f(X^x_t)]$. Let $\tau_k$ denote the time when
$(X^x_t)_{t\geq 0}$ makes its $k$-th jump. If $c(x)=0$, then
$\tau_1=\infty$ a.s.\ and $X_t=x$ a.s.\ for all $t\geq 0$, which
implies $u_t(x)=P_tf(x)=\E[f(X_t)]=f(x)$ $(t\geq 0)$, so
(\ref{intKol}) is trivially satisfied with $u_0(x)=f(x)$. We assume,
therefore, from now on that $c(x)>0$. In this case, conditional on
$\tau_1=s$ and $X^x_{\tau_1}=y$, the process $(X^x_{\tau_1+t})_{t\geq 0}$
is equally distributed with $(X^y_t)_{t\geq 0}$, which allows us
to write
\[\begin{array}{l}
\dis\E[f(X^x_t)]
=f(x)\P[\tau_1>t]+\int_0^t\P[\tau_1\in\di t]\sum_{y:\,y\neq x}\P[X^x_{\tau_1}=y]\P\big[f(X^y_{t-s})\big]\\[5pt]
\dis\quad=f(x)\ex{-c(x)t}+\int_0^tc(x)\ex{-c(x)s}\di s\,c(x)^{-1}\sum_{y:\,y\neq x}G(x,y)P_{t-s}f(y),
\end{array}\]
which shows that $u_t:=P_tf$ solves (\ref{intKol}) with $u_0=f$.
\end{Proof}

\begin{lemma}[The backward equation]
If\label{L:difint} a function $u\cn S\times\half\to\half$ solves
(\ref{backeq}), then it solves (\ref{intKol}). Conversely, each
bounded solution to (\ref{intKol}) also solves (\ref{backeq}).
\end{lemma}

\begin{Proof}
If $u$ solves (\ref{backeq}), then
\[
\dif{t}u_t(x)+c(x)u_t(x)=\sum_{y:\,y\neq x}G(x,y)u_t(y),
\]
which implies
\[
\dif{t}\Big(\ex{c(x)t}u_t(x)\Big)=\ex{c(x)t}\sum_{y:\,y\neq x}G(x,y)u_t(y)\qquad(t\geq 0,\ x\in S).
\]
Integrating and then multiplying both sides of the equation by
$e^{-c(x)t}$, we obtain
\[
u_t(x)=u_0(x)\ex{-c(x)t}+\int_0^t\!\di s\,\ex{-c(x)(t-s)}\sum_{y:\,y\neq x}G(x,y)u_s(y),
\]
which after the substitution $s\mapsto t-s$ yields (\ref{intKol}).

Conversely, if (\ref{intKol}) holds and $u$ is bounded, then the
right-hand side is continuous in $t$ for each $x$ and hence so is the
left-hand side. But then the right-hand side must actually be
continuously differentiable as a function of $t$ and the same must be
true for the left-hand side. We can then reverse the argument above
(differentiating instead of integrating) to obtain (\ref{backeq}).
\end{Proof}

\begin{lemma}[Comparison principle]
Let\label{L:compar} $G$ be the generator of a con\-tinu\-ous-time
Markov chain with countable state space $S$ and let $(P_t)_{t\geq 0}$
be the Markov semigroup with generator $G$. Assume that $u\cn
S\times\half\to\half$ satisfies
\[
\dif{t}u_t(x)\geq\sum_yG(x,y)u_t(y)\qquad(t\geq 0,\ x\in S),
\]
where $t\mapsto u_t(x)$ is continuously differentiable for each $x\in S$. Then
\[
P_tu_0(x)\leq u_t(x)\qquad(t\geq 0,\ x\in S).
\]
\end{lemma}

\begin{Proof}
Let $(X^x_t)_{t\geq 0}$ be the continuous-time Markov chain with
generator $G$ and initial state $x$ and as in Section~\ref{S:embed},
let $\tau_k$ denote the time of its $k$-th jump, for $1\leq k\leq N$,
where $N$ denotes the total number of jumps, which may be finite or
infinite, and let $\tau$ denote the explosion time. Then
\[
P_tu_0(x)=\E\big[u_0(X^x_t)1_{\{t<\tau\}}\big]\qquad(t\geq 0).
\]
On the event that $N<\infty$ we set $\tau_k:=\infty$ for $k>N$ and we define
\[
u^{(n)}_0(x)=\E\big[u_0(X^x_t)1_{\{t<\tau_n\}}\big]\qquad(t\geq 0).
\]
Then
\[
u^{(n)}_t(x)\asto{n}P_tu_0(x)\qquad(t\geq 0).
\]
Using the same argument as in the proof of Lemma~\ref{L:fijump}, we see that
\[
u^{(n+1)}_t(x)=u_0(x)\ex{-c(x)t}+\int_0^t\!\di s\,\ex{-c(x)s}\sum_{y:\,y\neq x}G(x,y)u^{(n)}_{t-s}(y).
\]
By the same argument as in the proof of Lemma~\ref{L:difint}, with all
equalities replaced by inequalities,
\[
u_t(x)\geq u_0(x)\ex{-c(x)t}+\int_0^t\!\di s\,\ex{-c(x)s}\sum_{y:\,y\neq x}G(x,y)u_{t-s}(y).
\]
We claim that $u^{(n)}_t(x)\leq u_t(x)$ $(n\geq 0,\ t\geq 0)$. The
proof is by induction. Clearly $u^{(0)}_t(x)=0\leq u_t(x)$ $(t\geq
0)$. Assuming that the statement holds for $n$, we have
\[\begin{array}{l}
\dis u_t(x)\geq u_0(x)\ex{-c(x)t}+\int_0^t\!\di s\,\ex{-c(x)s}\sum_{y:\,y\neq x}G(x,y)u_{t-s}(y)\\[5pt]
\dis\ \geq u_0(x)\ex{-c(x)t}+\int_0^t\!\di s\,\ex{-c(x)s}\sum_{y:\,y\neq x}G(x,y)u^{(n)}_{t-s}(y)=u^{(n+1)}_t(x).
\end{array}\]
Letting $n\to\infty$ we obtain $u_t(x)\geq P_t(x)$.
\end{Proof}

\begin{Proof}[of Theorem~\ref{T:Gchain}]
By Lemmas \ref{L:fijump} and \ref{L:difint}, $u$ solves the Kolmogorov
backward equation (\ref{backeq}) with initial condition $f$. If $u'$
is another solution, then Lemma~\ref{L:compar} implies that $u\leq
u'$, showing that $u$ is minimal.
\end{Proof}

We conclude this section with the following observation.

\begin{proposition}[Uniqueness of solutions]
The Kolmogorov backward equation (\ref{backeq}) has a unique bounded
solution $u$ with initial condition $u_0=f$ for each bounded function
$f\cn S\to\half$ if and only if the continuous-time Markov chain with
generator $G$ is nonexplosive.
\end{proposition}

\begin{Proof}
For each $r\in\R$, let $\un r\cn S\to\R$ denote the function that is
constantly equal to $r$. If $G$ is explosive, then $u_t:=P_t\un 1$ and
$u'_t:=\un 1$ $(t\geq 0)$ are two different bounded solutions of the
Kolmogorov backward equation (\ref{backeq}) with initial condition
$u_0=u'_0=\un 1$, proving that solutions are not unique.

On the other hand, assume that $G$ is nonexplosive and that $u$ is a
bounded solution with initial condition $u_0=f$. Since $u$ is bounded,
there exist an $r\geq 0$ such that $u_t\leq\un r$ $(t\geq
0)$. Lemma~\ref{L:compar} tells us that $P_tf\leq u_t$ $(t\geq
0)$. Also, since $\un r-u$ solves (\ref{backeq}) with initial
condition $\un r-f$, Lemma~\ref{L:compar} tells us that $\un r-u_t\geq
P_t(\un r-f)=\un r-P_tf$ $(t\geq 0)$ where in the last step we have
used that $G$ is nonexplosive. Combining these inequalities, we see
that $u_t=P_tf$ $(t\geq 0)$.
\end{Proof}

\section{Lyapunov functions}

It is tempting to think of explosive continuous-time Markov chains as
pathological, but there exist very natural chains that are
explosive. In fact, each transient chain can with a suitable random
time transformation be transformed into an explosive chain, so from
this point of view the distinction between transient and recurrent
chains would appear to be more fundamental than the distinction
between explosive and nonexplosive chains. Nevertheless, it is useful
to have at our disposal a technique for proving that a given chain is
nonexplosive. In the present section, we will show how
nonexplosiveness can be proved with the help of Lyapunov
functions. Below is the main result of this section. The term
``Lyapunov function'' originates in the stability theory of ordinary
differential equations but is sometimes also used for certain
functions occurring in Foster's theorem, that gives necessary and
sufficient conditions for positive recurrence of a Markov chain. The
role of the function $L$ in the following theorem is similar, so using
the term in a general sense, we may call it a Lyapunov function
too. Note that $GL$ in condition~(ii) below is well-defined by
(\ref{Gfdef}).

\begin{theorem}[Sufficient conditions for nonexplosiveness]
Let\label{T:Lyap} $G$ be the generator of a continuous-time Markov
chain with countable state space $S$ and let $c(x):=-G(x,x)$ $(x\in
S)$. Assume that there exists a function $L\cn S\to\half$ and constant
$\la\in\R$ such that:
\begin{enumerate}[(ii)]
\item $\sup\{c(x):x\in S,\ L(x)<C\}$ is finite for all $C<\infty$,
\item $GL\leq\la L$.
\end{enumerate}
Then the continuous-time Markov chain $(X_t)_{t\geq 0}$ with generator
$G$ is nonexplosive and
\[
\E^x\big[L(X_t)\big]\leq\ex{\la t}L(x)\qquad(t\geq 0,\ x\in S).
\]
\end{theorem}

The proof of Theorem~\ref{T:Lyap} depends on two lemmas.

\begin{lemma}[Exponential bound]
Let\label{L:Lex} $(X_t)_{t\geq 0}$ be a continuous-time Mar\-kov chain
with generator $G$, started in $X_0=x$. Assume that $L\cn S\to\half$
satisfies $GL\leq\la L$ for some $\la\in\R$. Then
\begin{equation}\label{Lex}
\E^x\big[L(X_t)1_{\{t<\tau\}}\big]\leq\ex{\la t}L(x)\qquad(t\geq 0),
\end{equation}
where $\tau$ denotes the explosion time of $(X_t)_{t\geq 0}$.
\end{lemma}

\begin{Proof}
The function $u_t(x):=L(x)e^{\la t}$ satisfies $\dif{t}u_t\geq Gu_t$
$(t\geq 0)$, so Lemma~\ref{L:compar} tells us that $P_tu_0\leq u_t$
$(t\geq 0)$, which is the same as (\ref{Lex}).
\end{Proof}

\begin{lemma}[Bounded jump rates]
Let\label{L:bdrate} $G$ be the generator of a con\-tinu\-ous-time
Markov chain with countable state space $S$ and let $c(x):=-G(x,x)$
$(x\in S)$. Assume that $\sup_{x\in S}c(x)<\infty$. Then $G$ is
nonexplosive.
\end{lemma}

\begin{Proof}
Let $(X^x_t)_{t\geq 0}$ be the continuous-time Markov chain with
generator $G$ and initial state $x$, constructed from its embedded
discrete-time Markov chain $(Y^x_k)_{k\geq 0}$ and i.i.d.\ standard
exponential random variables $(\sig_k)_{k\geq 0}$ as in
Section~\ref{S:embed}. Let $N$ be the total number of jumps, which may
be finite or infinite, and let $\tau$ denote the explosion time. On
the event that $N<\infty$ we have $\tau=\infty$ while on the event
that $N=\infty$ we have
\[
\tau=\sum_{k=0}^\infty\sig_k/c(Y^x_k).
\]
By our assumption that $C:=\sup_{x\in S}c(x)<\infty$ we can estimate
this from below by
\[
\tau\geq C^{-1}\sum_{k=1}^\infty\sig_k
\]
which is $\infty$ a.s.\ by the strong law of large numbers.
\end{Proof}

\begin{Proof}[of Theorem~\ref{T:Lyap}]
It suffices to prove that $G$ is nonexplosive, since the statement
about the expectation of $L(X_t)$ then follows from
Lemma~\ref{L:Lex}. We set
\[
S_C:=\big\{x\in S:L(x)<C\big\}
\]
and define a generator $G_C$ by
\[
G_C(x,y):=\left\{\begin{array}{ll}
\dis G(x,y)\quad&\mbox{if }x\in S_C,\\[5pt]
\dis 0\quad&\mbox{if }x\not\in S_C
\end{array}\right.
\]
We let $(X^x_t)_{t\geq 0}$ and $(X^{x,C}_t)_{t\geq 0}$ denote the
continuous-time Markov chains with generators $G$ and $G_C$
respectively. It follows from the construction of these processes in
terms of their embedded Markov chains that we can naturally couple
these processes such that
\[
X^x_t=X^{x,C}_t\qquad\forall t\leq\tau_C:=\inf\big\{t\geq 0:X^x_t\not\in S_C\big\}.
\]
In fact, we then have
\[
X^{x,C}_t=X^x_{t\wedge\tau_C}\qquad(t\geq 0),
\]
that is, $(X^{x,C}_t)_{t\geq 0}$ corresponds to the process
$(X^x_t)_{t\geq 0}$ stopped as soon as it leaves $S_C$.

Let $\tau$ denote the explosion time of $(X^x_t)_{t\geq 0}$.
Lemma~\ref{L:bdrate} and assumption~(i) of the theorem imply
that $(X^{x,C}_t)_{t\geq 0}$ is nonexplosive. Since the processes
$(X^x_t)_{t\geq 0}$ and $(X^{x,C}_t)_{t\geq 0}$ are equal up to time
$\tau_C$ it follows that $\tau_C\leq\tau$. Making $\la$ larger if
necessary, we can without loss of generality assume that $\la\geq 0$.
Assumption~(ii) of the theorem then implies that also $G_CL\leq\la L$.
Indeed, $G_CL(x)=GL(x)$ if $L(x)<C$ and $G_CL(x)=0$ otherwise. We
can therefore use Lemma~\ref{L:Lex} and the fact that $G_C$ is
nonexplosive to conclude that
\[
C\P^x[\tau_C\leq t]\leq\E\big[L(X^{x,C}_t)\big]\leq\ex{\la t}L(x)\qquad(t\geq 0).
\]
Since $\tau_C\leq\tau$, it follows that
\[
\P^x[\tau\leq t]\leq\ex{\la t}L(x)/C\qquad(t\geq 0),
\]
so letting $C\to\infty$ we see that $\tau=\infty$ a.s.
\end{Proof}

\section{Poisson point sets}\label{S:Pois}

The construction of a continuous-time Markov chain from its embedded
discrete-time Markov chain is useful for theoretical purposes, but for
the purpose of studying interacting particle systems a different
construction, that is based on Poisson point sets, will turn out to be
much more useful. To prepare for this, in the present section, we
recall the definition of Poisson point sets and some of their basic
properties.

Let $S$ be a \si-compact\footnote{This means that there exists a
countable collection of compact sets $S_i\sub S$ such that
$\bigcup_iS_i=S$.} metrizable space. We will mainly be interested in
the case that $S=\Gi\times\R$ where $\Gi$ is a countable set. We let
$\Si$ denote the Borel-\si-field on $S$. A \emph{locally finite
measure}\index{locally finite measure} on $(S,\Si)$ is a measure $\mu$
such that $\mu(C)<\infty$ for all compact $C\sub S$.

Let $(\Om,\Fi,\P)$ be our underlying probability space. A random measure on
$S$ is a function $\xi\cn\Om\times\Si\to[0,\infty]$ such that for fixed
$\om\in\Om$, the function $\xi(\om,\,\cdot\,)$ is a locally finite measure
on $(S,\Si)$, and for fixed $A\in\Si$, the function $\xi(\,\cdot\,,A)$ is
measurable. By \cite[Lemma~1.37]{Kal97}, we can think of $\xi$ as a random
variable with values in the space of locally finite measures on $(S,\Si)$,
equipped with the \si-field generated by the maps $\mu\mapsto\mu(A)$ with
$A\in\Si$. Then the integral $\int f\di\xi$ defines a $[0,\infty]$-valued
random variable for all measurable $f\cn S\to[0,\infty]$. There exists a unique
measure, denoted by $\E[\xi]$, such that
\[
\int\!\! f\,\di\E[\xi]=\E\bigl[\int\!\! f\,\di\xi\bigr]
\]
for all measurable $f\cn S\to[0,\infty]$. The measure $\E[\xi]$ is called the
\emph{intensity}\index{intensity!of a random measure} of $\xi$.

The following result follows from \cite[Lemma~10.1 and Prop.~10.4]{Kal97}.
\footnote{In fact, \cite[Prop.~10.4]{Kal97} shows that it is possible
to construct Poisson point measures on arbitrary measurable spaces,
assuming only that the intensity measure is \si-finite, but we will
not need this generality.} Below, $\Si_{\rm loc}:=\{A\in\Si:\ov A
\mbox{ is compact}\}$\index{0SEloc@${\cal S}_{\rm loc}$} denotes the
set of measurable subsets of $S$ whose closure is compact.

\begin{proposition}[Poisson point measures]
Let\label{P:Pois} $\mu$ be a locally finite measure on $(S,\Si)$. Then
there exists a random measure $\xi$, unique in distribution, such that
for any disjoint $A_1,\ldots,A_n\in\Si_{\rm loc}$, the random
variables $\xi(A_1),\ldots,\xi(A_n)$ are independent and $\xi(A_i)$ is
Poisson distributed with mean $\mu(A_i)$.
\end{proposition}

We call a random measure $\xi$ satisfying the conditions of
Proposition~\ref{P:Pois} a \emph{Poisson point measure}
\index{Poisson!point measure} with \emph{intensity}
\index{intensity!of a Poisson point measure} $\mu$. Indeed, one can
check that $\E[\xi]=\mu$. We note that $\xi(A)\in\N$ for all
$A\in\Si_{\rm loc}$. Such measures are called (locally finite)
\emph{counting measures}\index{counting measure}. Each locally finite
counting measure $\nu$ on $S$ is of the form
\[
\nu=\sum_{x\in{\rm supp}(\nu)}n_x\de_x,
\]
where ${\rm supp}(\nu)$, the support of $\nu$, is a locally finite subset
of $S$, the $n_x$ are positive integers, and $\de_x$ denotes the delta-measure
at $x$. We say that $\nu$ is \emph{simple}\index{simple counting measure} if
$n_x=1$ for all $x\in{\rm supp}(\nu)$. Recall that a measure $\mu$ has an
\emph{atom}\index{atom of a measure} at $x$ if $\mu(\{x\})>0$. A measure $\mu$
is called \emph{atomless}\index{atomless measure} if it has no atoms, that is,
$\mu(\{x\})=0$ for all $x\in S$. The already mentioned
\cite[Prop.~10.4]{Kal97} tells us the following.

\begin{lemma}[Simple Poisson point measures]
Let\label{L:simple} $\xi$ be a Poisson point measure with locally
finite intensity $\mu$. Then $\xi$ is a.s.\ simple if and only if
$\mu$ is atomless.
\end{lemma}

If $\mu$ is atomless, then a Poisson point measure $\xi$ with intensity $\mu$
is characterized by its support $\om:={\rm supp}(\xi)$. We call $\om$ a
\emph{Poisson point set}\index{Poisson!point set} with intensity
$\mu$. Intuitively, $\om$ is a set such that
$\P[\om\cap\di x\neq\emptyset]=\mu(\di x)$, independently for each
infinitesimal subset $\di x\sub S$.

For any counting measure $\nu$ on $S$ and measurable function $f\cn
S\to[0,1]$ we introduce the notation
\[
f^{\txt\nu}:=\prod_{i\in I}f(x_i)\qquad\mbox{with}\quad\nu=\sum_{i\in I}\de_{x_i},
\]
where the index set $I$ is either finite or countably infinite and, by
definition, a product of zero factors is one. Thus $f^0:=1$, where
$0$ denotes the counting measure that is identically
zero. Alternatively, our definition says that
\[
f^{\txt\nu}=\ex{\int(\log f)\di\nu},
\]
where $\log 0:=-\infty$ and $e^{-\infty}:=0$. It is easy to see that
$f^\nu f^{\nu'}=f^{\nu+\nu'}$.

\begin{lemma}[Laplace functionals]\label{L:PoisLap}
Let $\mu$ be a locally finite measure on $(S,\Si)$ and let $\xi$ be a Poisson
point measure with intensity $\mu$. Then
\begin{equation}\label{Pois}
\E\big[(1-f)^{\txt\xi}\big]=\ex{-\int f\di\mu}
\end{equation}
for each measurable $f\cn S\to[0,1]$. Conversely, if $\xi$ is a random counting
measure and (\ref{Pois}) holds for all continuous, compactly supported $f$,
then $\xi$ is a Poisson point measure with intensity $\mu$.
\end{lemma}

\begin{Proof}
The fact that Poisson point measures satisfy (\ref{Pois}) is proved in
\cite[Lemma~10.2]{Kal97}, which is written in terms of $-\log f$,
rather than $f$. The fact that knowing (\ref{Pois}) for all continuous,
compactly supported $f$ determines the law of $\xi$ uniquely follows
from \cite[Lemma~10.1]{Kal97}.
\end{Proof}

Formula (\ref{Pois}) can be interpreted in terms of thinning. Consider
a counting measure $\nu=\sum_i\de_{x_i}$, let $f\cn S\to[0,1]$ be
measurable, and let $\chi_i$ be independent \emph{Bernoulli random
variables}\index{Bernoulli random variable} (that is, random variables
with values in $\{0,1\}$) with $\P[\chi_i=1]=f(x_i)$. Then the random
counting measure
\[
\nu':=\sum_i\chi_i\de_{x_i}
\]
is called an \emph{$f$-thinning}\index{thinning} of the counting
measure $\nu$. Note that
\[
\P[\nu'=0]=\prod_i\P[\chi_i=0]=(1-f)^\nu.
\]
In view of this, the left-hand side of (\ref{Pois}) can be interpreted
as the probability that after thinning the random counting measure
$\xi$ with $f$, no points remain.

We cite the following fact from \cite[Lemma~10.17]{Kal97}.

\begin{lemma}[Poisson points on the half-line]
Let\label{L:Poisexp} $0<c<\infty$ and let $\ell$\index{0l@$\ell$}
denote the Lebesgue measure on $\half$. Let $(\tau_k)_{k\geq 0}$ be
real random variables such that $\tau_0=0$ and
$\sig_k:=\tau_k-\tau_{k-1}>0$ $(k\geq 1)$. Then $\om:=\{\tau_k:k\geq
1\}$ is a Poisson point set on $\half$ with intensity $c\ell$ if and
only if the random variables $(\sig_k)_{k\geq 1}$ are
i.i.d.\ exponentially distributed with mean $c^{-1}$.
\end{lemma}

We will need the following property of Poisson sets.

\begin{proposition}[Markov property of Poisson sets]
Let\label{P:PoisMark} $S$ be a countable set, let $\mu$ be a locally
finite measure on $S$, and let $\om$ be a Poisson point set on
$S\times\half$ with intensity measure
\[
\rho\big(\{x\}\times[s,t]\big):=\mu(\{x\})(t-s)\qquad(x\in S,\ 0\leq s\leq t).
\]
Let $S'\sub S$ and assume that $0<\mu(S')<\infty$. Set
\[
\tau:=\inf\{\tau_x:x\in S'\}\quad\mbox{with}\quad
\tau_x:=\inf\big\{t\geq 0:(x,t)\in\om\big\}\qquad(x\in S').
\]
Then a.s., there exists a unique $X\in S'$ such that $\tau_X=\tau$. Setting
\[
\om':=\big\{(x,t-\tau):(x,t)\in\om,\ t>\tau\big\},
\]
one has that the random variables $\tau$, $X$, and $\om'$ are independent,
\begin{equation}\label{tauX}
\P[\tau\geq t]=\ex{-\mu(S')t}\quad(t\geq 0),
\qquad
\P[X=x]=\frac{\mu(\{x\})}{\mu(S')}\quad(x\in S'),
\end{equation}
and $\om'$ is equally distributed with $\om$.
\end{proposition}

\begin{Proof}
We first prove the statement if $S'=S$. Let $\la:=\mu(S)$ and
$\pi:=\mu/\la$. Let $(\sig_k)_{k\geq 1}$ be i.i.d.\ exponentially
distributed random variables with parameter $\la$ and let
$\tau_k:=\sum_{i=1}^k\sig_i$ $(k\geq 1)$. Let $(X_k)_{k\geq 1}$ be
i.i.d.\ with law $\pi$ and independent of $(\sig_k)_{k\geq 1}$. Set
\[
\xi:=\sum_{k=1}^\infty\de_{\tau_k}
\quand
\eta:=\sum_{k=1}^\infty\de_{(X_k,\tau_k)}.
\]
By Lemma~\ref{L:Poisexp}, $\xi$ is a Poisson point measure on $\half$
with intensity measure $\la\ell$. We claim that $\eta$ is a Poisson
point measure on $S\times\half$ with intensity measure $\rho$ as in
the proposition. To see this, we apply Lemma~\ref{L:PoisLap}. Let
$f\cn S\times\R\to[0,1]$ be continuous and compactly supported. Define
$\ov f\cn\half\to[0,1]$ by
\begin{equation}
\ov f(t):=\sum_{x\in S}\pi(x)f(x,t)\qquad(t\geq 0).
\end{equation}
Then, using the fact that $\xi$ is a Poisson point measure on $\half$
with intensity measure $\la\ell$, we see that
\[\begin{array}{l}
\dis\E\big[(1-f)^{\txt\eta}\big]=\E\big[\prod_{k=1}^\infty\big(1-f(X_k,\tau_k)\big)\big]
=\prod_{k=1}^\infty\E\big[1-f(X_k,\tau_k)\big]\\[5pt]
\dis\quad=\prod_{k=1}^\infty\E\big[1-\ov f(\tau_k)\big]=\E\big[(1-\ov f)^{\txt\xi}\big]=\ex{-\la\int_0^\infty\ov f(t)\,\di t}=\ex{-\int f\,\di\rho}.
\end{array}\]
It follows that
\[
\om:=\big\{(X_k,\tau_k):k\geq 1\big\}
\]
is a Poisson point set on $S\times\half$ with intensity measure
$\rho$. Now clearly $\tau=\tau_1$ and $X=X_1$ are distributed as in
(\ref{tauX}) while
\[
\om'=\big\{(X_k,\tau_k-\tau_1):k\geq 2\big\}
\]
is independent of $(X_1,\tau_1)$ and equally distributed with
$\om$. This completes the proof in the special case that $S'=S$. The
general case follows immediately by applying what we have already
proved to the restriction of $\om$ to $S'\times\half$ and then using
that this is independent of to the restriction of $\om$ to $(S\beh
S')\times\half$.
\end{Proof}

\begin{Exercise}\label{E:sumPois}
Let $\xi_1,\xi_2$ be independent Poisson point measures with intensities $\mu_1,\mu_2$. Show that $\xi_1+\xi_2$ is a Poisson point measures with intensity $\mu_1+\mu_2$. \emph{Hint:} Lemma~\ref{L:PoisLap}.
\end{Exercise}

\section{Poisson construction of Markov processes}\label{S:Poichain}

In the present section we will show how a continuous-time Markov chain
with countable state space $S$ can be constructed by applying certain
maps $m\cn S\to S$ at the times of a Poisson point process. We start
with the following observation.

\begin{lemma}[Random mapping representation]
Let\label{L:randmap} $S$ be a countable set, let $\Gi$ be a countable
collection of maps $m\cn S\to S$, and let $(r_m)_{m\in\Gi}$ be
nonnegative real numbers such that
\begin{equation}\label{locfin}
c(x):=\sum_{m:\,m(x)\neq x}r_m<\infty\quad\mbox{for all }x\in S.
\end{equation}
Then there exists a generator $G$ of a continuous-time Markov chain
with state space $S$ such that
\begin{equation}\label{Grep}
Gf(x)=\sum_{m\in\Gi}r_m\big\{f\big(m(x)\big)-f\big(x\big)\big\}
\end{equation}
for all bounded $f\cn S\to\half$. Conversely, each generator $G$ of a
con\-tinu\-ous-time Markov chain with state space $S$ can be written
in the form (\ref{Grep}) for a suitable choice of a collection $\Gi$
of maps $m\cn S\to S$ and nonnegative real numbers $(r_m)_{m\in\Gi}$
satisfying (\ref{locfin}).
\end{lemma}

\begin{Proof}
Let $\Gi$ be a countable collection of maps $m\cn S\to S$, let
$(r_m)_{m\in\Gi}$ be nonnegative real numbers satisfying
(\ref{locfin}), and let $c(x)$ be as defined in (\ref{locfin}). Then
it is straightforward to check that
\[
G(x,y):=\left\{\begin{array}{ll}
\dis\sum\subb{m\in\Gi}{m(x)=y}r_m\quad&\mbox{if }x\neq y,\\[5pt]
\dis-c(x)\quad&\mbox{if }x=y
\end{array}\right.\qquad(x,y\in S)
\]
defines a generator such that (\ref{Grep}) holds. To see that each
generator $G$ can be written in this form, we define for each $x,y\in
S$ with $x\neq y$ a map $m_{x,y}\cn S\to S$ by
\[
m_{x,y}(z):=\left\{\begin{array}{ll}
y\quad&\mbox{if }z=x,\\[5pt]
z\quad&\mbox{otherwise,}
\end{array}\right.
\]
we set $\Gi:=\{m_{x,y}:x,y\in S,\ x\neq y\}$ and
$r_{m_{x,y}}:=G(x,y)$. It is then straightforward to check that the
rates $(r_m)_{m\in\Gi}$ satisfy (\ref{locfin}) and (\ref{Grep}) holds.
\end{Proof}

We call the way of writing $G$ as in (\ref{Grep}) a \emph{random
mapping representation} of $G$.
\index{random mapping representation!of a generator} Recall that random mapping
representations of probability kernels have been defined in
Section~\ref{S:finMark}. There is a close connection between the
two. Indeed, we will see that if $(P_t)_{t\geq 0}$ is the Markov
semigroup with generator $G$, then using a random mapping
representation of $G$ we can for each $t\geq 0$ construct a random
mapping representation of $P_t$. We will do this by giving a Poisson
construction of the continuous-time Markov chain with generator
$G$. In the context of interacting particle systems, such Poisson
constructions are known as \emph{graphical representations} and they
have many applications.

From now on, we assume that $G$ is the generator of a continuous-time
Markov chain with countable state space $S$ and we fix a random
mapping representation of $G$ of the form (\ref{Grep}) in terms of
rates $(r_m)_{m\in\Gi}$ satisfying (\ref{locfin}). We equip the space
$\Gi\times\R$ with the measure
\begin{equation}\label{rho}
\rho\big(\{m\}\times[s,t]\big):=r_m(t-s)\qquad\big(m\in\Gi,\ s\leq t\big).
\end{equation}
Let $\om$ be a Poisson point set with intensity $\rho$. We call $\om$
the \emph{graphical representation} associated with the random mapping
representation (\ref{Grep}). We claim that for each $t\in\R$, there
exists at most one $m\in\Gi$ such that $(m,t)\in\om$. To see this, we
note that for each $m\in\Gi$, the set
\[
\xi_m:=\sum_{t:\,(m,t)\in\om}\de_t
\]
is a Poisson point measure on $\R$ with intensity $r_m\ell$, where
$\ell$\index{0l@$\ell$} denotes the Lebesgue measure. Since the sets
$\R\times\{m\}$ $(m\in\Gi)$ are disjoint, the random measures $\xi_m$
$(m\in\Gi)$ are independent, and hence by Exercise~\ref{E:sumPois},
for each $m\neq m'$, the measure $\xi_m+\xi_{m'}$ is a Poisson point
measure on $\R$ with intensity $(r_m+r_{m'})\ell$. Since the Lebesgue
measure is atomless, by Lemma~\ref{L:simple}, this Poisson point
measure is simple, so there are no times $t\in\R$ for which both
$(m,t)\in\om$ and $(m',t)\in\om$. In view of this, we can
unambiguously define a random function $\R\ni t\mapsto\mk^\om_t\in\Gi$ by
setting
\[\index{0mt@$\mk^\om_t$}
\mk^\om_t:=\left\{\begin{array}{ll}
m\quad&\mbox{if }(m,t)\in\om,\\[5pt]
1\quad&\mbox{otherwise,}
\end{array}\right.\]
where we write $1$ to denote the identity map.

As before, we write $S_\infty:=S\cup\{\infty\}$. We equip $S_\infty$
with a topology such that a set $A\sub S_\infty$ is closed if and only
if $A$ is either finite or $A$ is infinite and $\infty\in A$. One can
check that $S_\infty$ is compact in this topology and a sequence
$x_n\in S$ converges to $\infty$ if and only if it leaves every finite
subset of $S$, that is, for each finite $S'\sub S$, there exists an
$m$ such that $x_n\not\in S'$ for all $n\geq m$. The topological space
$S_\infty$ is known as the \emph{one-point
compactification}\index{one-point compactification} of $S$.
We extend the maps $m\in\Gi$ to $S_\infty$ by putting $m(\infty):=\infty$
$(m\in\Gi)$.

Fix $s\in\R$. By definition, we say that a random function
$X\cn[s,\infty)\to S_\infty$ solves the evolution equation
\begin{equation}\label{evolve}
X_t=\mk^\om_t(X_{t-})\qquad(t>s)
\end{equation}
if
\begin{enumerate}[(ii)]
\item $\lim_{r\down t}X_r=X_t$ $(t\geq s)$ and $\lim_{r\up t}X_r=:X_{t-}$
 exists $(t>s)$,
\item $X_t=\infty$ for all $t\geq\tau:=\inf\{r\geq s:X_r=\infty\}$,
\end{enumerate}
and (\ref{evolve}) holds. Below is the main result of this section.

\begin{theorem}[Stochastic flow]
Let\label{T:stochflow} $G$ be the generator of a con\-tinu\-ous-time
Markov chain with countable state space $S$ and let $\om$ be the
graphical representation associated with a random mapping
representation of $G$. Then almost surely, for each $s\in\R$ and $x\in
S_\infty$, there exists a unique solution $(X^{s,x}_t)_{t\geq s}$ to
the evolution equation (\ref{evolve}) with initial state
$X^{s,x}_s=x$. Setting
\begin{equation}\label{flowdef}\index{0Xb@$\Xb_{s,t}$}
\Xb_{s,t}(x):=X^{s,x}_t\qquad(s\leq t,\ x\in S_\infty)
\end{equation}
defines a collection of maps $(\Xb_{s,t})_{s\leq t}$ from $S_\infty$
into itself such that
\begin{equation}\label{stochflow}
\Xb_{s,s}=1
\quand\Xb_{t,u}\circ\Xb_{s,t}=\Xb_{s,u}\qquad(s\leq t\leq u).
\end{equation}
If $G$ is nonexplosive, then almost surely, $\Xb_{s,t}$ maps $S$ into
itself for all $s\leq t$. If $s\in\R$ and $X_0$ is an $S$-valued
random variable with law $\mu$, independent of $\om$, then the process
$(X_t)_{t\geq 0}$ defined as
\begin{equation}\label{Xflow}
X_t:=\Xb_{s,s+t}(X_0)\qquad(t\geq 0)
\end{equation}
is distributed as the continuous-time Markov chain with generator $G$
and initial law~$\mu$.
\end{theorem}

Formula (\ref{stochflow}) says that the random maps
$(\Xb_{s,t})_{s\leq t}$ form a \emph{stochastic
flow}.\index{stochastic flow} This stochastic flow is
\emph{stationary}\index{stationary!stochastic flow} in the sense that
\[
\Xb_{0,t}\mbox{ is equally distributed with }\Xb_{s,s+t}\quad(s\in\R).
\]
We note that since $\Xb_{s,t}$ is constructed using only Poisson
points of the form $(m,r)$ with $s<r\leq t$, and restrictions of a
Poisson point set to disjoint parts of the space are independent, it
follows that the stochastic flow $(\Xb_{s,t})_{s\leq t}$ has
\emph{independent increments}\index{independent increments of a stochastic flow}
in the sense that
\[
\Xb_{t_0,t_1},\ldots,\Xb_{t_{n-1},t_n}\quad\mbox{are independent}\quad
\forall\ t_0<\cdots<t_n.
\]
We note that (\ref{Xflow}) implies that
\[
P_t(x,y)=\P\big[\Xb_{s,s+t}(x)=y\big]\qquad(s,t\in\R,\ x,y\in S),
\]
so indeed, as announced, we have found a random mapping representation
of the subprobability kernels $(P_t)_{t\geq 0}$.\med

\begin{Proof}[of Theorem~\ref{T:stochflow}]
We start by proving that almost surely, for each $s\in\R$ and $x\in
S_\infty$, there exists a unique solution $(X^{s,x}_t)_{t\geq s}$ to
the evolution equation (\ref{evolve}) with initial state
$X^{s,x}_s=x$. If $x=\infty$, then clearly $X^\infty_t=\infty$ $(t\geq
s)$ is the unique solution of (\ref{evolve}) so without loss of
generality we assume from now on that $x\in S$. For each $x\in S$, the
set
\[
\big\{t\in\R:\mk^\om_t(x)\neq x\big\}
\]
is a Poisson point set on $\R$ with intensity $c(x)$ defined in
(\ref{locfin}), which is finite by assumption. This allows us to
inductively define times $(\tau_n)_{n\geq 0}$ and a discrete chain
$(Y^x_n)_{n\geq 0}$ by setting $\tau_0:=0$, $Y^x_0:=x$,
\[
\tau_{n+1}:=\left\{\begin{array}{ll}
\dis\inf\big\{t>\tau_n:\mk^\om_t(x_n)\neq x_n\big\}\quad&\mbox{if }\tau_n<\infty,\\[5pt]
\dis\infty\quad&\mbox{if }\tau_n=\infty,
\end{array}\right.
\]
and
\[
Y^x_{n+1}:=\left\{\begin{array}{ll}
\dis\mk^\om_{\tau_{n+1}}(Y^x_n)\quad&\mbox{if }\tau_{n+1}<\infty,\\[5pt]
\dis Y^x_n\quad&\mbox{if }\tau_{n+1}=\infty.
\end{array}\right.
\]
We set $\tau:=\lim_{n\to\infty}\tau_n$ and
\[
N:=\inf\big\{n\geq 0:\tau_{n+1}=\infty\big\}=\inf\big\{n\geq 0:c(Y^x_n)=0\big\}.
\]
We claim that $(X^{s,x}_t)_{t\geq s}$ defined as
\[
X^x_t:=\left\{\begin{array}{ll}
\dis Y^x_k\quad&\mbox{if }t\in[\tau_k,\tau_{k+1}),\ 0\leq k<N+1,\\[5pt]
\dis\infty\quad&\mbox{if }t\geq\tau
\end{array}\right.
\]
solves the evolution equation (\ref{evolve}). It is easy to see that
$t\mapsto X_t$ is right-continuous with left limits, that
$\tau=\inf\{r\geq s:X_r=\infty\}$ and $X_t=\infty$ for all
$t\geq\tau$, and that (\ref{evolve}) holds for all $t\neq\tau$. To see
that (\ref{evolve}) holds at the time $\tau$ should it be finite, we
need to show that $X_{\tau-}=\infty$ on the event that $\tau<\infty$. To see
this, imagine on the contrary that $\tau<\infty$ while $X_t$ does not
converge to $\infty$ as $t\to\tau$. By the definition of the one-point
compactification, this implies that there exist a finite set $S'\sub
S$ and times $s_n\to\tau$ such that $X_{s_n}\in S'$. This, in turn,
implies that during the finite time interval $[0,\tau)$, the function
$X_t$ makes infinitely many jumps that start at some point in $S'$
and end in some other point in $S$. But this is impossible, since by
(\ref{locfin}),
\[
\bigcup_{y\in S'}\big\{t\in\R:\mk^\om_t(y)\neq y\big\}
\]
is a locally finite subset of $\R$. This completes the proof that
$(X^{s,x}_t)_{t\geq s}$ solves (\ref{evolve}). If $(X'_t)_{t\geq s}$
is another solution, then we see by induction that $X'_t=X^{s,x}_t$
for all $0\leq t\leq\tau_n$ and for all $n\geq 0$. This implies that
$X'_{\tau-}=\infty$. By (\ref{evolve}) and the way we have defined $m(\infty)$
for $m\in\Gi$, we must have $X'_\tau=X'_{\tau-}$ if $\tau<\infty$ so
$X'_t=\infty$ for all $t\geq\tau$ by property~(ii), which shows that
solutions to (\ref{evolve}) are unique.

Note that our previous argument holds almost surely for all $x\in S$
and $s\in\R$ simultaneously, that is, this includes random times $s$
that may be chosen in dependence on the Poisson set $\om$. For
deterministic $x$ and $s$, we claim that $(X^{s,x}_{s+t})_{t\geq 0}$
is distributed as the continuous-time Markov chain with generator $G$
and initial state $x$. To see this, let $(\sig'_k)_{k\geq 0}$ be
i.i.d.\ exponentially distributed random variables with mean one,
independent of everything else. Define $(\sig_k)_{k\geq 0}$ by
$\sig_k:=c(Y^x_k)\tau_k$ if $c(Y^x_k)>0$ and $\sig_k:=\sig'_k$
otherwise. Using Proposition~\ref{P:PoisMark} we see by induction that
$(Y^x_k)_{k\geq 0}$ is the embedded Markov chain, $(\sig_k)_{k\geq 0}$
are i.i.d.\ exponentially distributed random variables with mean one,
independent of $(Y^x_k)_{k\geq 0}$, and $(X^{s,x}_{s+t})_{t\geq 0}$ is
the continuous-time Markov chain constructed in terms of its embedded
Markov chain and exponential holding times as in
Section~\ref{S:embed}.

Let $(\Xb_{s,t})_{s\leq t}$ be defined in (\ref{flowdef}). Then it is
straightforward to check that $(\Xb_{s,t})_{s\leq t}$ is a stochastic
flow in the sense of (\ref{stochflow}). If $X_0=x$ is deterministic,
then we have just seen that the process in (\ref{Xflow}) is
distributed as the continuous-time Markov chain with generator $G$ and
initial law $\mu$. The general case follows by conditioning on $X_0$,
which is independent of everything else.

To complete the proof, we must show that if $G$ is nonexplosive, then
almost surely, $\Xb_{s,t}$ maps $S$ into itself for all $s\leq t$. In
other words, we must show that
\[
\Xb_{s,t}(x)\in S\quad\forall s,t\in\R\mbox{ with }s\leq t\mbox{ and }x\in S\quad{\rm a.s.}
\]
If $s$ is deterministic, then by what we have proved
$X_t:=\Xb_{s,s+t}(x)$ $(t\geq 0)$ is the continuous-time Markov chain
with generator $G$ and initial state $x$, so if this process is
nonexplosive, then
\[
\Xb_{s,t}(x)\in S\quad\forall t\in[s,\infty)\quad{\rm a.s.}\quad(s\in\R,\ x\in S).
\]
Since $s\mapsto\Xb_{s,t}(x)$ is constant between the times of the
Poisson point process
\[
\big\{s\in\R:\mk^\om_s(x)\neq x\big\}
\]
which has finite intensity by (\ref{locfin}), we can improve our
previous statement to
\[
\Xb_{s,t}(x)\in S\quad\forall s,t\in\R\mbox{ with }s\leq t\quad{\rm a.s.}\quad(x\in S),
\]
and since $S$ is countable, we see that the statement holds for all
$s,t\in\R$ with $s\leq t$ and $x\in S$ simultaneously.
\end{Proof}

\begin{Exercise}
Use\label{E:Markprop2} formula (\ref{Xflow}) as well as the fact that
the stochastic flow $(\Xb_{s,t})_{s\leq t}$ is stationary with
independent increments to give an alternative proof of the fact that
the construction in Section~\ref{S:embed} of $(X^x_t)_{t\geq 0}$ via
the embedded Markov chain yields a Markov process in the sense of
(\ref{Mark}) with state space $S_\infty$ and transition kernels
$(\ov P_t)_{t\geq 0}$.
\end{Exercise}

For later use, we conclude this section with a theorem that is very
similar to Theorem~\ref{T:stochflow}. The starting point is again a
graphical representation $\om$ that is associated with a random
mapping representation of the generator $G$ of a continuous-time
Markov chain. The only difference is that this time, we will apply the
maps associated with elements of $\om$ in the reverse order. By
definition, we say that a random function $X\cn(\infty,u]\to S_\infty$
solves the evolution equation
\begin{equation}\label{backevolve}
X_{t-}=\mk^\om_t(X_t)\qquad(t\leq u)
\end{equation}
if
\begin{enumerate}[(ii)]
\item $\lim_{r\down t}X_r=X_t$ $(t<u)$ and $\lim_{r\up t}X_r=:X_{t-}$ exists
 $(t\leq u)$,
\item $X_{t-}=\infty$ for all $t\leq\tau:=\sup\{t\leq u:X_t=\infty\}$,
\end{enumerate}
and (\ref{backevolve}) holds. Note that we allow for the case that
$(m,u)\in\om$ for some $m\in\Gi$ and in this case it is possible that
$X_{u-}\neq X_u$. This is a difference with solutions of
(\ref{evolve}) which never make a jump at the initial time
$s$. Because we work ``backwards in time'', in (\ref{Xflow}) below we
obtain a Markov process with left-continuous sample paths. In spite of
these small differences, the proof of Theorem~\ref{T:backflow} is
practically identical to the proof of Theorem~\ref{T:stochflow}.

\begin{theorem}[Backward stochastic flow]
Let\label{T:backflow} $G$ be the generator of a con\-tinu\-ous-time
Markov chain with countable state space $S$ and let $\om$ be the
graphical representation associated with a random mapping
representation of $G$. Then almost surely, for each $u\in\R$ and $x\in
S_\infty$, there exists a unique solution $(X^{u,x}_t)_{t\leq u}$ to
the evolution equation (\ref{backevolve}) with final state
$X^{u,x}_u=x$. Setting
\begin{equation}\label{backdef}
\Xb_{u,t}(x):=X^{u,x}_t\qquad(u\geq t,\ x\in S_\infty)
\end{equation}
defines a collection of maps $(\Xb_{u,t})_{u\geq t}$ from $S_\infty$
into itself such that
\begin{equation}\label{backflow}
\Xb_{u,u}=1
\quand\Xb_{t,s}\circ\Xb_{u,t}=\Xb_{u,s}\qquad(u\geq t\geq s).
\end{equation}
If $G$ is nonexplosive, then almost surely, $\Xb_{u,t}$ maps $S$ into
itself for all $u\geq t$. If $u\in\R$ and $X_0$ is an $S$-valued
random variable with law $\mu$, independent of $\om$, then the process
$(X_t)_{t\geq 0}$ defined as
\begin{equation}\label{backXflow}
X_t:=\Xb_{u,u-t}(X_0)\qquad(t\geq 0)
\end{equation}
is distributed as the left-continuous modification of the
continuous-time Mar\-kov chain with generator $G$ and initial
law~$\mu$.
\end{theorem}

We call a collection of maps $(\Xb_{u,t})_{u\geq t}$ as in
(\ref{backflow}) a \emph{backward stochastic flow}.
\index{backward!stochastic flow} Stationarity and independent
increments are defined as in the forward case.

\begin{Exercise}
Prove Theorem~\ref{T:backflow}.
\end{Exercise}

\section{An example: ASEP}

In this section we look at an asymmetric simple exclusion process
(ASEP) on the natural numbers with a finite, fixed number of
particles. This demonstrates the theory developed so far and at the
same time serves as a warm-up for the final two sections of this
chapter, in which we show more generally how to construct interacting
particle systems in which the lattice may be infinite but the number
of particles is finite. The construction of interacting particle
systems with infinitely many particles will have to wait till
Chapter~\ref{C:construct}.

We fix an integer $n\geq 1$ and write
\[
\Si_n:=\big\{x\in\{0,1\}^\N:\sum_{i=0}^\infty x(i)=n\big\}.
\]
For each $i,j\in\N$ with $i\neq j$, we let ${\tt asep}_{ij}$ denote
the asymmetric exclusion map defined in Section~\ref{S:excl}. We
observe that ${\tt asep}_{ij}$ preserves the number of particles, that
is, it maps the space $\Si_n$ into itself. We will be interested in
the continuous-time Markov chain with countable state space $\Si_n$
and generator
\begin{equation}\begin{array}{r@{\,}c@{\,}l}\label{Gasepn}
\dis Gf(x)
&:=&\dis\sum_{i=1}^\infty r^-_i\big\{f\big({\tt asep}_{i,i-1}(x)\big)-f\big(x\big)\big\}\\[5pt]
&&\dis+\sum_{i=1}^\infty r^+_i\big\{f\big({\tt asep}_{i-1,i}(x)\big)-f\big(x\big)\big\},
\end{array}\end{equation}
where $r^\pm_i$ are nonnegative rates. To see that this is indeed the
generator of a continuous-time Markov chain with state space $\Si_n$,
we must check condition (\ref{locfin}) of Lemma~\ref{L:randmap}. Since
${\tt asep}_{ij}(x)\neq x$ if and only if $x(i)=1$ and $x(j)=0$, the
constant $c(x)$ from (\ref{locfin}) is given by
\begin{equation}\label{cSn}
c(x)=\sum_{i=1}^\infty 1_{\{x(i-1)=0,\ x(i)=1\}}r^-_i
+\sum_{i=1}^\infty 1_{\{x(i-1)=1,\ x(i)=0\}}r^+_i.
\end{equation}
For $x\in\Si_n$, the sums in (\ref{cSn}) have at most $2n$ nonzero
terms so clearly $c(x)<\infty$ for all $x\in\Si_n$. It follows that
$G$, defined in (\ref{Gasepn}) is the generator of a (possibly
explosive) continuous-time Markov chain with state space $\Si_n$.

\begin{lemma}[Nonexplosiveness]
Assume\label{L:asepexpl} that there exists a constant $K<\infty$ such
that $r^+_i\leq Ki$ $(i\geq 1)$. Then the continuous-time Markov chain
$(X_t)_{t\geq 0}$ with generator $G$ and state space $\Si_n$ is
nonexplosive.
\end{lemma}

\begin{Proof}
For $x\in\Si_n$, let $R(x):=\sup\{i\in\N:x(i)=1\}$ denote the position
of the right-most particle. We will apply Theorem~\ref{T:Lyap} to the
Lyapunov function
\[
L(x):=\big(R(x)+1\big)^2\qquad(x\in\Si_n).
\]
Since $L$ can only increase due to the right-most particle making a
jump to the right,
\[\begin{array}{r@{\,}c@{\,}l}
\dis GL(x)&\leq&\dis r^+_{R(x)+1}\big[\big(R(x)+2\big)^2-\big(R(x)+1\big)^2\big]\\[5pt]
&\leq&\dis K\big(R(x)+1\big)\big(2R(x)+3\big)\leq 3K\big(R(x)+1\big)^2,
\end{array}\]
so condition (ii) of Theorem~\ref{T:Lyap} is satisfied with
$\la=3K$. In view of (\ref{cSn}), condition (i) is also satisfied, so
we conclude that $(X_t)_{t\geq 0}$ is nonexplosive.
\end{Proof}

\begin{Exercise}
Fix\label{E:explos} $\ha<p\leq 1$ and $\al>1$ and assume that
$r^-_i=(1-p)i^\al$ and $r^+_i=p(i-1)^\al$ $(i\geq 1)$. Assume that
$n=1$ (there is only one particle) and write $X_t=e_{\xi_t}$ where
$e_i\in\{0,1\}^\La$ is defined as $e_i(j):=1$ if $i=j$ and $:=0$
otherwise. Then $(\xi_t)_{t\geq 0}$ is a continuous-time Markov chain
with state space $\N$ that jumps from $i$ to $i-1$ with rate
$(1-p)i^\al$ $(i\geq 1)$ and from $i$ to $i+1$ with rate $pi^\al$
$(i\geq 0)$. Let $(Y_k)_{k\geq 0}$ be the embedded Markov chain of
$(\xi_t)_{t\geq 0}$. Show that
\[
\sum_{k=0}^\infty Y_k^{-\al}<\infty\quad{\rm a.s.}
\]
Use this to conclude that the continuous-time Markov chain
$(X_t)_{t\geq 0}$ is explosive.
\end{Exercise}

Formula (\ref{Gasepn}) is a random mapping representation of the
generator $G$. Combining Lemma~\ref{L:asepexpl} with
Theorem~\ref{T:stochflow}, we can use this random mapping
representation to define a stochastic flow $(\Xb_{s,t})_{s\leq t}$ on
$\Si_n$, and then construct the continuous-time Markov chain
$(X_t)_{t\geq 0}$ in terms of this stochastic flow $(\Xb_{s,t})_{s\leq t}$
as in (\ref{Xflow}).

Random mapping representations are in general not unique, and
different random mapping representations lead to different stochastic
flows for the same continuous-time Markov chain. To demonstrate this,
we assume from now on that
\begin{equation}\label{rsym}
r_i:=r^-_i=r^+_i\qquad(i\geq 1).
\end{equation}
We recall the definition of the (symmetric) exclusion map ${\tt excl}_{ij}$
from Section~\ref{S:excl}. We claim that under the
symmetry assumption (\ref{rsym}), we can rewrite our definition of the
generator $G$ from (\ref{Gasepn}) as
\begin{equation}\label{Gsym}
Gf(x)=\sum_{i=1}^\infty r_i\big\{f\big({\tt excl}_{i-1,i}(x)\big)-f\big(x\big)\big\}.
\end{equation}
To check this, it suffices to check that for each $x,y\in\Si_n$ with
$x\neq y$, the rate $G(x,y)$ of jumps from $x$ to $y$ is the same
regardless of whether we define $G$ by (\ref{Gasepn}) or by
(\ref{Gsym}). This rate is zero unless we are in one of the following
two cases:
\renewcommand{\labelenumi}{{\rm\Roman{enumi}.}}
\begin{enumerate}[II.]
\item There exists an $i\geq 1$ such that $x(i-1)=0$, $x(i)=1$, $y(i-1)=1$, $y(i)=0$, and $x(j)=y(j)$ for all $j\not\in\{i-1,i\}$,
\item There exists an $i\geq 1$ such that $x(i-1)=1$, $x(i)=0$, $y(i-1)=0$, $y(i)=1$, and $x(j)=y(j)$ for all $j\not\in\{i-1,i\}$.
\end{enumerate}
\renewcommand{\labelenumi}{{\rm (\roman{enumi})}}
Defining $G$ as in (\ref{Gasepn}), we see that $G(x,y)=r^-_i$ in
case~I and $G(x,y)=r^+_i$ in case~II. On the other hand, defining $G$
as in (\ref{Gsym}), we have $G(x,y)=r_i$ in both cases. In particular,
under the symmetry assumption (\ref{rsym}), both definitions are
equivalent.

\begin{figure}
\begin{center}
\inputtikz{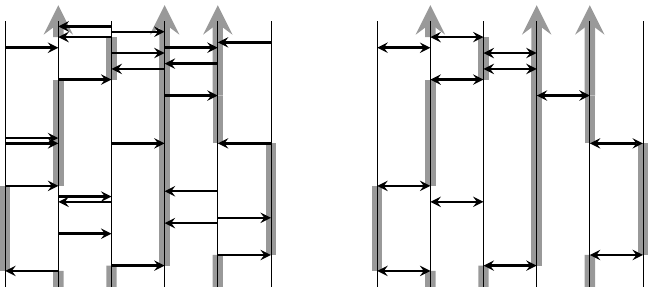}
\end{center}
\caption{Two different Poisson constructions of the same
  continuous-time Markov chain. Each picture shows a different
  graphical representation of the same symmetric exclusion
  process. Time is plotted upwards. The arrows in the picture on the
  left symbolize the application of the map ${\tt asep}_{ij}$, which
  has the effect that if there is a particle at $i$ and the site $j$
  is vacant, then the particle at $i$ jumps to $j$. The arrows in the
  picture on the right symbolize the application of the map ${\tt excl}_{ij}$,
  which has the effect that particles at $i$ and $j$
  exchange their positions. The Poisson density of arrows in the left
  picture is twice as high as in the right picture (or the same, if we
  count double arrows double).}
\label{fig:sep}
\commentAlt{Figure~\ref{fig:sep}}{Three gray paths indicate the
  movement of particles. The black arrows in the left picture have one
  arrowhead and can transport particles in one direction, those on the
  right have two arrowheads and can transport particles in both
  directions.}
\end{figure}

Even though under the assumption (\ref{rsym}) formulas (\ref{Gasepn})
and (\ref{Gsym}) define the same generator, they are different random
mapping representations that lead to different stochastic flows. This
is illustrated\footnote{In general, two different graphical
representations of the same Markov process only yield two Markov
processes that are equal in law. In Figure~\ref{fig:sep}, for didactic
purposes, these processes together with their graphical
representations have been coupled so that the Markov processes are
a.s.\ equal. In this example, can you think of a coupling that
achieves this?} in Figure~\ref{fig:sep}. To further illuminate this,
let us define
\[
\Si_{\rm fin}:=\bigcup_{n=0}^\infty\Si_n.
\]
Note that $\Si_{\rm fin}$ is countable. Our previous arguments show
that $G$ is the generator of a continuous-time Markov chain with state
space $\Si_{\rm fin}$. By Theorem~\ref{T:stochflow}, we can use the
random mapping representations (\ref{Gasepn}) and (\ref{Gsym}) to
construct two different stochastic flows $(\Xb_{s,t})_{s\leq t}$ and
$(\Xb'_{s,t})_{s\leq t}$. The next exercise demonstrates that these
stochastic flows have different properties.

\begin{Exercise}
Show that the random maps $\Xb'_{s,t}\cn\Si_{\rm fin}\to\Si_{\rm fin}$
are \emph{additive} in the sense that
\[
\Xb'_{s,t}(x\vee y)=\Xb'_{s,t}(x)\vee\Xb'_{s,t}(y)\qquad(s\leq t,\ x,y\in\Si_{\rm fin}),
\]
where $(x\vee y)(i)=x(i)\vee y(i)$ denotes the pointwise maximum of
two configurations $x,y\in\Si_{\rm fin}$. Show that the maps
$\Xb_{s,t}\cn\Si_{\rm fin}\to\Si_{\rm fin}$ are not additive.
\end{Exercise}

\begin{Exercise}
Show that if (\ref{rsym}) is weakened to $r^-_i\leq r^+_i$ $(i\geq 1)$,
then the generator $G$ in (\ref{Gasepn}) can be rewritten as
\[\begin{array}{r@{\,}c@{\,}l}
\dis Gf(x)
&=&\dis\sum_{i=1}^\infty r^-_i\big\{f\big({\tt excl}_{i-1,i}(x)\big)-f\big(x\big)\big\}\\[5pt]
&&\dis+\sum_{i=1}^\infty(r^+_i-r^-_i)\big\{f\big({\tt asep}_{i-1,i}(x)\big)-f\big(x\big)\big\}.
\end{array}\]
\end{Exercise}

\section{Local maps}\label{S:local}

In this section we return to the general set-up of
Section~\ref{S:setup}. Thus, $S$ is a finite set, called the
\emph{local state space}, $\La$ is a countable set, called the
\emph{lattice}, and $S^\La$ denotes the Cartesian product of $\La$
copies of $S$, that is, this is the space of all functions $x\cn\La\to
S$. The interacting particle systems we are interested in are Markov
processes with state space $S^\La$ and generator $G$ of the form
(\ref{Gdef}), that is,
\[
Gf(x)=\sum_{m\in\Gi}r_m\big\{f\big(m(x)\big)-f\big(x\big)\big\},
\]
were $\Gi$ is a set whose elements are maps $m\cn S^\La\to S^\La$ and
$(r_m)_{m\in\Gi}$ are nonnegative rates.

If the lattice $\La$ is finite, then so is $S^\La$ and we can
immediately apply Theorem~\ref{T:stochflow} to construct our
interacting particle system from a graphical representation $\om$. If
$\La$ is infinite, then $S^\La$ is uncountable (as long as $S$ has at
least two elements), so Theorem~\ref{T:stochflow} is not
applicable. In Chapter~\ref{C:construct} we will develop the methods
needed to construct interacting particle systems on infinite lattices
from a graphical representation. There is one important special case
where Theorem~\ref{T:stochflow} is applicable, however, even if the
lattice is finite. Many interacting particle systems with a local
state space of the form $S=\{0,1\}$ have the property that if we start
the system in an initial configuration with finitely many ones, then
the system stays in such states for all times. In the previous
section, we have shown this for the exclusion process. Other examples
are the voter model, the contact process, and various systems of
branching and coalescing particles. A counterexample are stochastic
Ising models.

In Section~\ref{S:finsys} we will show how Theorem~\ref{T:stochflow}
can be applied to construct interacting particle systems on infinite
lattices, but started in an initial configuration with finitely many
ones. To prepare for this, in the present section, we take a closer
look at the sort of maps that are typically used to construct an
interacting particle system, such as the voter map in (\ref{votmap})
and the branching and death maps in (\ref{bramap}) and
(\ref{deathmap}).

We will always equip the state space $S^\La$ of an interacting
particle system with the \emph{product topology},\index{product topology}
which says that a sequence $x_n\in S^\La$ converges to a
limit $x$ if and only if
\[
x_n(i)\asto{n}x(i)\qquad\forall i\in\La.
\]
Note that since $S$ is finite, this simply says that for each
$i\in\La$, there is an $N$ (which may depend on $i$) such that
$x_n(i)=x(i)$ for all $n\geq N$. Since $S$ is finite, it is in
particular compact, so by Tychonoff's theorem, the space $S^\La$ is
compact in the product topology.

Let $S$ and $T$ be finite sets, let $\La$ be a countable set, and let
$f\cn S^\La\to T$ be a function. Then we say that a point $j\in\La$ is
\emph{$f$-relevant}\index{f-relevant@$f$-relevant}
\index{relevant coordinate for function} if
\[
\exists x,y\in S^\La\mbox{ s.t.\ }f(x)\neq f(y)
\mbox{ and }x(k)=y(k)\ \forall k\neq j,
\]
that is, changing the value of $x$ in $j$ may change the value of
$f(x)$. We write
\[\index{0Rf@$\Ri(f)$}
\Ri(f):=\big\{j\in\La:j\mbox{ is $f$-relevant}\big\}.
\]
The following lemma (which we have taken from \cite[Lemma~24]{SS18})
says that a function $f\cn S^\La\to T$ is continuous with respect to
the product topology if and only if it depends on finitely many
coordinates.

\begin{lemma}[Continuous maps]
Let\label{L:contprod} $S$ and $T$ be finite sets and let $\La$ be a
countable set. Then a function $f\cn S^\La\to T$ is continuous with
respect to the product topology if and only if the following two
conditions are satisfied:
\begin{enumerate}[(ii)]
\item $\Ri(f)$ is finite,
\item If $x,y\in S^\La$ satisfy $x(j)=y(j)$ for all $j\in\Ri(f)$, then $f(x)=f(y)$.
\end{enumerate}
\end{lemma}

Before we give the proof of Lemma~\ref{L:contprod}, we first make some
observations. The following exercise shows how continuity can fail if
condition~(i) of Lemma~\ref{L:contprod} does not hold.

\begin{Exercise}[A discontinuous map]
Let $2\N:=\{2n:n\in\N\}$ and $2\N+1:=\{2n+1:n\in\N\}$. Define $f\cn\{0,1\}^\N\to\{0,1\}$ by
\begin{equation}
f(x):=\left\{\begin{array}{ll}
1\quad\mbox{if }\inf\{i\in\N:x(i)=1\}\in 2\N\cup\{\infty\},\\[5pt]
0\quad\mbox{if }\inf\{i\in\N:x(i)=1\}\in 2\N+1.
\end{array}\right.
\end{equation}
Show that $f$ satisfies condition~(ii) of Lemma~\ref{L:contprod} but
not condition~(i). Show that $f$ is not continuous.
\end{Exercise}

The following exercise shows that contrary to what one might initially
have guessed, condition~(ii) of Lemma~\ref{L:contprod} is not
automatically satisfied, even when condition~(i) holds.

\begin{Exercise}[Another discontinuous map]
Define $f\cn\{0,1\}^\N\to\{0,1\}$ by
\begin{equation}
f(x):=\left\{\begin{array}{ll}
1\quad\mbox{if }\{i\in\N:x(i)=1\}\mbox{ is finite},\\[5pt]
0\quad\mbox{if }\{i\in\N:x(i)=1\}\mbox{ is infinite}.
\end{array}\right.
\end{equation}
Show that $f$ satisfies condition~(i) of Lemma~\ref{L:contprod} but
not condition~(ii). Show that $f$ is not continuous.
\end{Exercise}

\begin{Proof}[of Lemma~\ref{L:contprod}]
Let $(\al_j)_{j\in\La}$ be strictly positive constants such that
$\sum_{j\in\La}\al_j<\infty$. Then the metric
\begin{equation}
d(x,y):=\sum_{j\in\La}\al_j1_{\txt\{x(j)\neq y(j)\}}
\qquad(x,y\in S^\La)
\end{equation}
generates the product topology on $S^\La$. By Tychonoff's theorem,
$S^\La$ is compact, so the function $f$ is uniformly continuous. Since
the target space $T$ is finite, this means that there exists an
$\eps>0$ such that $d(x,y)<\eps$ implies $f(x)=f(y)$. Since
$\sum_{j\in\La}\al_j<\infty$, there exists some finite $\La'\sub\La$
such that $\sum_{j\in\La\beh\La'}\al_j<\eps$. It follows that
\begin{description}
\item[{\rm(ii)'}] If $x,y\in S^\La$ satisfy $x(j)=y(j)$ for all $j\in\La'$,
  then $f(x)=f(y)$.
\end{description}
We conclude from this that $\Ri(f)\sub\La'$, proving (i). If this is a
strict inclusion, then we can inductively remove those points from
$\La'$ that are not elements of $\Ri(f)$ while preserving the property
(ii)', until in a finite number of steps we see that (ii) holds.


Conversely, if (i) and (ii) hold and $x_k\to x$ pointwise, then by (i)
there exists some $n$ such that $x_k(j)=x(j)$ for all $j\in\Ri(f)$ and
hence by (ii) $f(x_k)=f(x)$ for all $k\geq n$, proving that $f$ is
continuous.
\end{Proof}

For any map $m\cn S^\La\to S^\La$ and $i\in\La$, we define $m[i]\cn
S^\La\to S$ by $m[i](x):=m(x)(i)$\index{0mi@$m[i]$} $(x\in
S^\La,\ i\in\La)$. It follows immediately from the definition of the
product topology that $m$ is continuous if and only if $m[i]$ is
continuous for all $i\in\La$. We let
\[\index{0Dm@$\Di(m)$}
\Di(m):=\big\{i\in\La:\exists x\in S^\La\mbox{ s.t.\ }m(x)(i)\neq x(i)\big\}
\]
denote the set of lattice points $i$ for which $m[i]$ is not the trivial map
$m[i](x)=x(i)$ $(x\in S^\La)$. Note that $\Ri(m[i])=\{i\}$ if $i\not\in\Di(m)$.

By definition, a \emph{local map}\index{local map} is a function $m\cn
S^\La\to S^\La$ such that:
\begin{enumerate}[(ii)]
\item $m$ is continuous,
\item $\Di(m)$ is finite.
\end{enumerate}
In view of Lemma~\ref{L:contprod}, this says that $m$ is local if $m$
changes the values of at most finitely many lattice points using
information from finitely many lattice points only. The following
exercise describes yet another way to look at local maps.

\begin{Exercise}[Local maps]\label{E:locmap}
Show that a map $m\cn S^\La\to S^\La$ is local if and only if there
exists a finite set $\De\sub\La$ and a map $m'\cn S^\De\to S^\De$ such
that
\[
m(x)(k)=\left\{\begin{array}{ll}
m'\big((x(i))_{i\in\De}\big)(k)\quad&\mbox{if }k\in\De,\\[5pt]
x(k)\quad&\mbox{otherwise.}
\end{array}\right.
\]
\end{Exercise}

Before we continue, it is good to see a number of examples.
\begin{itemize}
\item The voter map ${\tt vot}_{ij}$ defined in (\ref{votmap}) satisfies
\[
\Di({\tt vot}_{ij})=\{j\}\quand\Ri({\tt vot}_{ij}[j])=\{i\},
\]
since only the type at $j$ changes, and it suffices to know the type
at $i$ to predict the new type at $j$.
\item The branching map ${\tt bra}_{ij}$ defined in (\ref{bramap}) satisfies
\[
\Di({\tt bra}_{ij})=\{j\}\quand\Ri({\tt bra}_{ij}[j])=\{i,j\},
\]
since only the type at $j$ changes, but we need to know both the type
at $i$ and at $j$ to predict the new type at $j$ since
${\tt bra}_{ij}(x)(j)=x(i)\vee x(j)$.
\item The death map ${\tt death}_i$ defined in (\ref{deathmap}) satisfies
\[
\Di({\tt death}_i)=\{i\}\quand\Ri({\tt death}_i[i])=\emptyset
\]
since only the type at $i$ changes, and the new type at $i$ is 0
regardless of $x$.
\item For each $i\in\La$, we can similarly define a \emph{birth
map}\index{birth map} ${\tt birth}_i\cn\{0,1\}^\La\to\{0,1\}^\La$ as
\begin{equation}\label{birthmap}\index{0birth@${\tt birth}_i$}
{\tt birth}_i(x)(k):=\left\{\begin{array}{ll}
1\quad&\mbox{if }k=i,\\[5pt]
x(k)\quad&\mbox{otherwise.}
\end{array}\right.
\end{equation}
Then
\[
\Di({\tt birth}_i)=\{i\}\quand\Ri({\tt birth}_i[i])=\emptyset.
\]
\item The coalescing random walk map ${\tt rw}_{ij}$ defined in
  (\ref{rwmap}) satisfies
\[
\Di({\tt rw}_{ij})=\{i,j\},\quad\Ri({\tt rw}_{ij}[i])=\emptyset,
\quand\Ri({\tt rw}_{ij}[j])=\{i,j\},
\]
since the types at both $i$ and $j$ can change, the new type at $i$ is 0
regardless of the previous state, but to calculate ${\tt rw}_{ij}(x)(j)$
we need to know both $x(i)$ and $x(j)$.
\end{itemize}

\begin{Exercise}[Exclusion and cooperative branching maps]
Recall the asymmetric and symmetric exclusion maps ${\tt asep}_{ij}$
and ${\tt excl}_{ij}$ defined in (\ref{asep}) and (\ref{exclmap}), and
the cooperative branching map ${\tt coop}_{ii'j}$ defined in
(\ref{coopmap}). When $m$ is any of these maps, determine $\Di(m)$,
and determine $\Ri(m[i])$ for all $i\in\Di(m)$.
\end{Exercise}

\section{Systems of finitely many particles}\label{S:finsys}

Throughout this section we assume that $S$ is a finite set containing
a special element that we denote by $0$. For $x\in S^\La$ we introduce
the notation
\[\index{00normx@$\vert x\vert$}
|x|:=\big|\{i\in\La:x(i)\neq 0\}\big|\qquad(x\in S^\La)
\]
and we write
\[\index{0SiLa@$\Si(\La),\Si_{\rm fin}(\La)$}
\Si(\La):=S^\La\quand\Si_{\rm fin}(\La):=\big\{x\in\Si(\La):|x|<\infty\big\}.
\]
It is easy to see that $\Si_{\rm fin}(\La)$ is countable. We let $\un
0\in\Si(\La)$ \index{00yyy@$\un 0,\un 1$} denote the configuration
that is identically zero, that is, this is the constant function
defined as $\un 0(i):=0$ $(i\in\La)$. We will be interested in local maps
$m\cn\Si(\La)\to\Si(\La)$ that satisfy
\begin{equation}\label{m0}
m(\un 0)=\un 0.
\end{equation}
Almost all the examples of local maps mentioned in the previous
section satisfy (\ref{m0}). Indeed, this holds for the local maps
\[
{\tt vot}_{ij},\quad{\tt bra}_{ij},\quad{\tt death}_i,\quad{\tt rw}_{ij},
\quad{\tt asep}_{ij},\quad{\tt excl}_{ij},\quand{\tt coop}_{ijk},
\]
while ${\tt birth}_i$ is the only local map mentioned in the previous
section that does not map $\un 0$ into itself.

Let $\Gi$ is a countable collection of local maps
$m\cn\Si(\La)\to\Si(\La)$ that all satisfy (\ref{m0}) and let
$(r_m)_{m\in\Gi}$ be rates. Then under suitable assumptions on the
rates, we may expect that
\begin{equation}\label{GIPS}
Gf(x):=\sum_{m\in\Gi}r_m\big\{f\big(m(x)\big)-f\big(x\big)\big\}
\end{equation}
is the generator of a nonexplosive continuous-time Markov chain with
countable state space $\Si_{\rm fin}(\La)$. To formulate sufficient
conditions for this to be true, for any local map $m\cn S^\La\to
S^\La$ we introduce the following notation:
\[\begin{array}{c}\index{0Rm@$\Ri(m)$}\index{0Rmup@$\Ri^\up_i(m)$}
\index{0Rmdw@$\Ri^\down_i(m)$}
\dis\Ri(m):=\big\{(i,j)\in\La^2:i\mbox{ is $m[j]$-relevant}\big\},\\[5pt]
\dis\Ri^\up_i(m):=\big\{j\in\La:(i,j)\in\Ri(m)\big\},\quad
\Ri^\down_j(m):=\big\{i\in\La:(i,j)\in\Ri(m)\big\}.
\end{array}\]
Here is the main result of this section. The form of condition
(\ref{upsum}) is inspired by \cite{Lat24}. Below, we let
$1_A$\index{00yyza@$1_A$} denote the indictor function of a set
$A\sub\La$, that is, $1_A(i):=1$ if $i\in A$ and $:=0$ if $i\in\La\beh
A$.

\begin{theorem}[Finite particle configurations]
Let\label{T:finIPS} $S$ and $\La$ be a finite and countable set,
respectively, and assume that $S$ contains a special element denoted
as $0$. Let $\Gi$ be a countable collection of local maps $m\cn
S^\La\to S^\La$ such that $m(\un 0)=\un 0$ for all $m\in\Gi$ and let
$(r_m)_{m\in\Gi}$ be nonnegative rates. Assume that
\begin{equation}\label{upsum}
{\rm(i)}\ \sup_{i\in\La}\sum_{m\in\Gi}r_m1_{\Di(m)}(i)<\infty,\quad
{\rm(ii)}\ \sup_{i\in\La}\sum_{m\in\Gi}r_m\big|\Ri^\up_i(m)\beh\{i\}\big|<\infty.
\end{equation}
Then $G$ defined in (\ref{GIPS}) is the generator of a nonexplosive
continuous-time Markov chain with state space $\Si_{\rm fin}(\La)$.
Moreover, the process started in $X_0=x\in\Si_{\rm fin}(\La)$ satisfies
\begin{equation}\label{expsiz}
\E^x\big[|X_t|\big]\leq|x|\ex{K_\up t}\quad(t\geq 0)\quad\mbox{with}\quad
K_\up:=\sup_{i\in\La}\sum_{m\in\Gi}r_m\big(|\Ri^\up_i(m)|-1\big).
\end{equation}
\end{theorem}

\begin{Proof}
We start by checking condition (\ref{locfin}) which is necessary and
sufficient for (\ref{GIPS}) to define the generator of a (possibly
explosive) continuous-time Markov chain. In our present setting,
(\ref{locfin}) reads
\[
\sum_{m:\,m(x)\neq x}r_m<\infty\quad\mbox{for all }x\in\Si_{\rm fin}(\La).
\]
Let $x\in\Si_{\rm fin}(\La)$ and let $A:=\{i\in\La:x(i)\neq 0\}$ which
is finite by the definition of $\Si_{\rm fin}(\La)$. If $m(x)\neq x$,
then $m(x)(j)\neq x(j)$ for some some $j\in\Di(m)$. If $j\not\in A$,
then by the fact that $m(\un 0)=\un 0$, there must exist an $i\in A$
such that $i\in\Ri(m[j])$. This allows us to estimate
\[\begin{array}{r@{\,}c@{\,}l}
\dis\sum_{m:\,m(x)\neq x}r_m
&\leq&\dis\sum_{j\in A}\sum_{m\in\Gi}1_{\Di(m)}(j)r_m
+\sum_{j\in\La\beh A}\sum_{i\in A}\sum_{m\in\Gi}1_{\Ri^\up_i(m)}(j)r_m\\[5pt]
&=&\dis\sum_{i\in A}\Big(\sum_{m\in\Gi}1_{\Di(m)}(i)r_m
+\sum_{m\in\Gi}|\Ri^\up_i(m)\beh A|r_m\Big),
\end{array}\]
which is finite by (\ref{upsum}).

It remains to prove that $G$ is nonexplosive. We apply
Theorem~\ref{T:Lyap} to the Lyapunov function
\begin{equation}\label{Lynorm}
L(x):=|x|\quad\big(x\in\Si_{\rm fin}(\La)\big).
\end{equation}
Note that $c(x):=-G(x,x)=\sum_{m:\,m(x)\neq x}r_m$, so our previous
calculation shows that the function $L$ satisfies condition~(i) of
Theorem~\ref{T:Lyap}. It remains to check condition~(ii).
Since each $m\in\Gi$ satisfies $m(\un 0)=\un 0$,
if $m(x)(j)\neq 0$ for some $j\in\La$, then there must be an $i\in\La$
such that $x(i)\neq 0$ and $j\in\Ri^\up_i(m)$, which allows us to
estimate
\[
|m(x)|=\sum_{i:\,x(i)\neq 0}\sum_j1_{\Ri^\up_i(m)}(j)=\sum_{i:\,x(i)\neq 0}\big|\Ri^\up_i(m)\big|.
\]
It follows that
\[\begin{array}{r@{\,}c@{\,}l}
\dis GL(x)
&=&\dis\sum_{m\in\Gi}r_m\big\{L\big(m(x)\big)-L\big(x\big)\big\}\\[5pt]
&\leq&\dis\sum_{m\in\Gi}r_m\sum_{i:\,x(i)\neq 0}\big(|\Ri^\up_i(m)|-1\big)
=\sum_{i:\,x(i)\neq 0}\sum_{m\in\Gi}r_m\big(|\Ri^\up_i(m)|-1\big)\\[5pt]
&\leq&\dis L(x)\sup_{i\in\La}\sum_{m\in\Gi}r_m\big(|\Ri^\up_i(m)|-1\big)\qquad(x\in\Si_{\rm fin}(\La)).
\end{array}\]
This shows that condition~(ii) of Theorem~\ref{T:Lyap} is satisfied
with $\la=K_\up$, so $G$ is nonexplosive. Formula (\ref{expsiz}) now
follows from the exponential bound in Theorem~\ref{T:Lyap}.
\end{Proof}

It is instructive to see some concrete examples of interacting
particle systems to which Theorem~\ref{T:finIPS} is
applicable. Generalizing (\ref{Gvot}), if $\la\cn\La^2\to\half$ is a
function, then we can define a \emph{voter model}\index{voter model}
generator by
\begin{equation}\label{Gvot2}
G_{\rm vot}f(x)
:=\sum_{i,j\in\La^2}
\la(i,j)\big\{f\big({\tt vot}_{ij}(x)\big)-f\big(x\big)\big\}
\qquad(x\in S^\La),
\end{equation}
where $\la(i,j)\geq 0$ is the Poisson rate at which site $j$ adopts
the type of site~$i$.

\begin{Exercise}
Show\label{E:finvot} that the generator $G_{\rm vot}$ satisfies the
assumptions of Theorem~\ref{T:finIPS} if
\[
\sup_{i\in\La}\big[\sum_{j\in\La}\la(j,i)+\sum_{j\in\La}\la(i,j)\big]<\infty.
\]
\end{Exercise}

Similarly, generalizing (\ref{Gcontact}), we can define a
\emph{contact process}\index{contact process} generator by
\begin{equation}\begin{array}{r@{\,}c@{\,}l}\label{Gcont}
\dis G_{\rm cont}f(x)
&:=&\dis\sum_{i,j\in\La}\la(i,j)
\big\{f\big({\tt bra}_{ij}(x)\big)-f\big(x\big)\big\}\\[5pt]
&&\dis+\de\sum_{i\in\La}\big\{f\big({\tt death}_i(x)\big)-f\big(x\big)\big\},
\end{array}\end{equation}
where $\la(i,j)\geq 0$ is the \emph{infection rate}\index{infection rate}
from $i$ to $j$ and $\de\geq 0$ is the \emph{death rate}.\index{death rate}

\begin{Exercise}
Show\label{E:fincont} that the generator $G_{\rm cont}$ satisfies the
assumptions of Theorem~\ref{T:finIPS} if
\[
\sup_{i\in\La}\big[\sum_{j\in\La}\la(j,i)+\sum_{j\in\La}\la(i,j)\big]<\infty.
\]
\end{Exercise}

For processes with a sort of translation invariant
structure\footnote{To formalize this, let us call a bijection
$\psi\cn\La\to\La$ such that $\la\big(\psi(i),\psi(j)\big)=\la(i,j)$
for all $i,j\in\La$ an \emph{automorphism}\label{autfoot} of $\la$.
In analogy with the terminology for graphs, we can define $\la$ to be
\emph{vertex transitive}\index{transitivity!of infection rates}
\index{vertex transitive!infection rates} if for each $i,j\in\La$,
there exists an automorphism $\psi$ such that $\psi(i)=j$. If $\la$ is
vertex transitive, then $\sum_{j\in\La}\la(j,i)+\sum_{j\in\La}\la(i,j)$
does not depend on $i\in\La$.} the expression
$\sum_{j\in\La}\la(j,i)+\sum_{j\in\La}\la(i,j)$ does not depend on
$i\in\La$. Using this, one can check that for translation invariant
voter models, the condition in Exercise~\ref{E:finvot} is
optimal. Indeed, if we start the process with a single one at $i$,
then $\sum_{j\in\La}\la(j,i)$ is the rate at which this one becomes a
zero while $\sum_{j\in\La}\la(i,j)$ is the rate at which this one
produces another one somewhere. In case of the contact process, we can
actually do a bit better than Exercise~\ref{E:fincont}.

\begin{proposition}[Finite contact processes]
Assume\label{P:fincont} that
\begin{equation}\label{outrate}
r:=\sup_{i\in\La}\sum_{j\in\La}\la(i,j)<\infty.
\end{equation}
Then $G_{\rm cont}$, defined in (\ref{Gcont}), is the generator of a
nonexplosive continuous-time Markov chain with state space $\Si_{\rm fin}(\La)$.
Moreover,
\begin{equation}\label{contexp}
\E^x\big[|X_t|\big]\leq\ex{(r-\de)t}|x|
\qquad\big(t\geq 0,\ x\in\Si_{\rm fin}(\La)\big).
\end{equation}
\end{proposition}

\begin{Proof}
For any $x\in\Si_{\rm fin}(\La)$, we can estimate the quantity $c(x)$
from (\ref{locfin}) by
\[
\sum_{m:\,m(x)\neq x}r_m=\de|x|+\sum_{i:\,x(i)=1}\sum_{j:\,x(j)=0}\la(i,j)\leq(\de+r)|x|,
\]
where $r$ is the quantity in (\ref{outrate}). Since this is finite for
each $x\in\Si_{\rm fin}(\La)$, $G_{\rm cont}$ is the generator of a
(possibly explosive) continuous-time Markov chain with state space
$\Si_{\rm fin}(\La)$.

To see that $G_{\rm cont}$ is nonexplosive we apply
Theorem~\ref{T:Lyap} to the Lyapunov function in (\ref{Lynorm}). Our
previous calculation shows that $L$ satisfies condition~(i) of
Theorem~\ref{T:Lyap} so it remains to check condition~(ii). We
estimate
\[\begin{array}{r@{\,}c@{\,}l}
\dis GL(x)
&=&\dis\sum_{m\in\Gi}r_m\big\{L\big(m(x)\big)-L\big(x\big)\big\}\\[5pt]
&=&\dis\sum_{i,j\in\La}\la(i,j)1_{\{x(i)=1,\ x(j)=0\}}
-\de\sum_{i\in\La}1_{\{x(i)=1\}}\leq(r-\de)|x|,
\end{array}\]
from which we see that condition~(ii) of Theorem~\ref{T:Lyap} is
satisfied with $\la=r-\de$. Theorem~\ref{T:Lyap} now tells us that
$G_{\rm cont}$ is nonexplosive and (\ref{contexp}) holds.
\end{Proof}

In particular, Proposition~\ref{P:fincont} tells us that if $r<\de$,
then the contact process \emph{dies out} in the sense that
\[
\P^x\big[X_t=\un 0]\asto{t}1\qquad\forall x\in\Si_{\rm fin}(\La).
\]
This is quite natural since $r$, defined in (\ref{outrate}), is the
maximal reproduction rate of an individual (assuming all other sites
are vacant). If this is less than the death rate, then each individual
produces on average less than one offspring before it dies, leading to
an exponential decay of the population size.

\chapter{The mean-field limit}\label{C:meanfield}

\section{Processes on the complete graph}

In Chapter~\ref{C:intro}, we have made acquaintances with a number of
interacting particle systems. While some properties of these systems
turn out easy to prove, other seemingly elementary questions can sometimes be
remarkably difficult. A few examples of such hard problems have been mentioned
in Chapter~\ref{C:intro}. In view of this, interacting particle systems are
being studied by a range of different methods, from straightforward
numerical simulations as we have seen in Chapter~\ref{C:intro}, to nonrigorous
renormalization group techniques and rigorous mathematical methods. All these
approaches complement each other. In addition, when a given problem appears
too hard, one often looks for simpler models that (one hopes) still catch the
essence, or at least some essential features of the behavior that one is
interested in.

A standard way to turn a difficult model into an (often) much easier model is
to take the \emph{mean-field limit}, which we explain in the present chapter.
Basically, this means that one replaces the graph structure of the underlying
lattice that one is really interested in (in practice often $\Z^d$) by the
structure of the complete graph with $N$ vertices, and then takes the limit
$N\to\infty$. As we will see, many properties of ``real'' interacting particle
systems are already reflected in these mean-field models. In particular, phase
transitions can often already be observed and even the values of critical
exponents of high-dimensional models are correctly predicted by the mean-field
model. In view of this, studying the mean-field limit is a wise first step in
the study of any more complicated model that one may encounter.

Of course, not all phenomena can be captured by replacing the graph structure
that one is really interested in by the complete graph. Comparing the real
model with the mean-field model, one can learn which elements of the observed
behavior are a consequence of the specific spatial structure of the lattice,
and which are not. Also for this reason, studying the mean-field limit should
be part of a complete study of any interacting particle system.

\section{The mean-field limit of the Ising model}\label{S:Ismean}
\index{Ising model!mean-field}

In this section we study the mean-field Ising model, also known as the
\emph{Curie--Weiss model}\index{Curie--Weiss model}, with Glauber dynamics.
We recall from formulas (\ref{Mxi}) and (\ref{Glauber}) in
Chapter~\ref{C:intro} that the Ising model is an interacting particle system
with local state space $S=\{-1,+1\}$, where each site $i$ updates its spin
value $x(i)\in\{-1,+1\}$ at rate one. When a spin value is updated, the
probability that the new value is $+1$ respectively $-1$ is proportional to
$e^{\bet N_{x,i}(+1)}$ respectively $e^{\bet N_{x,i}(-1)}$, where
$N_{x,i}(\sig):=\sum_{j\in\Ni_i}1_{\{x(j)=\sig\}}$ denotes the number of
neighboring sites that have the spin value $\sig$.

For the aim of taking the mean-field model, it will be convenient to formulate
the model slightly differently. We let
\[
\ov N_{x,i}:=\frac{1}{|\Ni_i|}\sum_{j\in\Ni_i}1_{\{x(j)=\sig\}}
\]
denote the fraction of neighbors that have the spin value $\sig$, and consider
the model where (compare (\ref{Glauber}))
\begin{equation}\label{meanGlauber}
\mbox{site $i$ flips to the value $\sig$ with rate}\quad
\frac{e^{\bet\ov N_{x,i}(\sig)}}{\sum_{\tau\in\{-1,+1\}}e^{\bet\ov N_{x,i}(\tau)}}.
\end{equation}
Assuming that $|\Ni_i|$ is just a constant that does not depend on
$i\in\La$ (as is the case for the model on $\Z^d$ and more generally
on any vertex transitive graph), this is just a reparametrization of
the original model where the parameter $\bet$ is replaced by
$\bet/|\Ni_i|$.

We now wish to construct the mean-field model, that is, the model on a complete
graph $\La_N$ with $|\La_N|=N$ vertices (sites), where each site is a neighbor
of each other site. For mathematical simplicity, we even count a site as a
neighbor of itself, that is, we set
\[
\Ni_i:=\La_N\quand|\Ni_i|=N.
\]
A consequence of this choice is that the \emph{average magnetization}
\[
\ov X_t:=\frac{1}{N}\sum_{i\in\La_N}X_t(i)\qquad(t\geq 0)
\]
forms a Markov process $\ov X=(\ov X_t)_{t\geq 0}$. Indeed, $\ov X_t$ takes
values in the space
\[
\big\{-1,-1+\ffrac{2}{N},\ldots,1-\ffrac{2}{N},1\big\},
\]
and jumps
\[\begin{array}{rcl}
\dis\ov x\mapsto\ov x+\ffrac{2}{N}
&\dis\quad\mbox{with rate}\quad
&\dis N_x(-1)\frac{e^{\bet N_x(+1)/N}}{e^{\bet N_x(-1)/N}+e^{\bet N_x(+1)/N}},\\[10pt]
\dis\ov x\mapsto\ov x-\ffrac{2}{N}
&\dis\quad\mbox{with rate}\quad
&\dis N_x(+1)\frac{e^{\bet N_x(-1)/N}}{e^{\bet N_x(-1)/N}+e^{\bet N_x(+1)/N}},
\end{array}\]
where $N_x(\sig):=N_{x,i}(\sig)=\sum_{j\in\La_n}1_{\{x(j)=\sig\}}$ does not
depend on $i\in\La_N$. We observe that
\[
N_x(+1)/N=(1+\ov x)/2
\quand
N_x(-1)/N=(1-\ov x)/2.
\]
In view of this, we can rewrite the jump rates of $\ov X$ as
\[\begin{array}{rcl}
\dis\ov x\mapsto\ov x+\ffrac{2}{N}
&\dis\quad\mbox{with rate}\quad
&\dis N(1-\ov x)/2\,\frac{e^{\bet(1+\ov x)/2}}
{e^{\bet(1-\ov x)/2}+e^{\bet(1+\ov x)/2}},\\[10pt]
\dis\ov x\mapsto\ov x-\ffrac{2}{N}
&\dis\quad\mbox{with rate}\quad
&\dis N(1+\ov x)/2\,\frac{e^{\bet(1-\ov x)/2}}
{e^{\bet(1-\ov x)/2}+e^{\bet(1+\ov x)/2}}.
\end{array}\]
In particular, since these rates are a function of $\ov x$ only (and do not
depend on other functions of $x=(x(i))_{i\in\La_N}$), we see that
$\ov X=(\ov X_t)_{t\geq 0}$, on its own, is a Markov process.
(This argument will be made rigorous in Section~\ref{S:Marfun} below.)
Canceling a common factor $e^{\bet/2}$ in the nominator and denominator of
the rates, we can simplify our formulas a bit to
\begin{equation}\begin{array}{rcl}\label{CurieWeiss}
\dis\ov x\mapsto\ov x+\ffrac{2}{N}
&\dis\quad\mbox{with rate}\quad
&\dis r^N_+(\ov x):=N(1-\ov x)/2\,\frac{e^{\bet\ov x/2}}
{e^{-\bet\ov x/2}+e^{\bet\ov x/2}},\\[10pt]
\dis\ov x\mapsto\ov x-\ffrac{2}{N}
&\dis\quad\mbox{with rate}\quad
&\dis r^N_-(\ov x):=N(1+\ov x)/2\,\frac{e^{-\bet\ov x/2}}
{e^{-\bet\ov x/2}+e^{\bet\ov x/2}}.
\end{array}\end{equation}

\begin{figure}[htb]
\begin{center}
\inputtikz{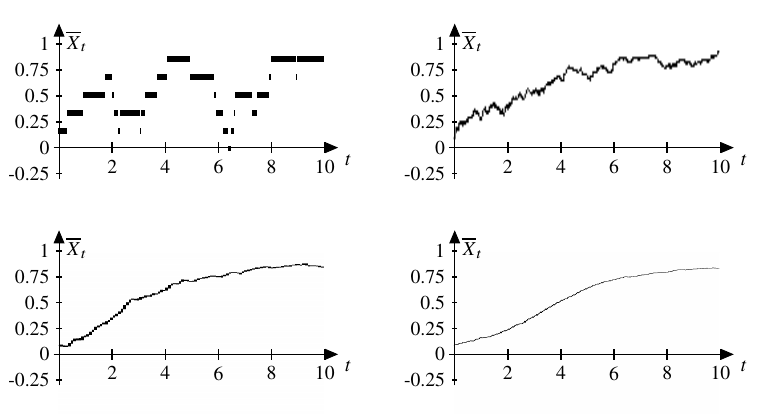}
\caption{The mean-field Ising model on lattice with $N=10$, $100$, $1000$, and
$10,000$ sites, respectively. In these simulations, the parameter is $\bet=3$,
and the initial state is $\ov X_0=0.1$, except in the first picture, where
$\ov X_0=0.2$.}
\label{fig:Ismean}
\commentAlt{Figure~\ref{fig:Ismean}}{Four graphs of functions. The
  first graph takes just five different values and jumps between
  them. The second looks like a Brownian motion with drift
  upwards. The third is more smooth and the last is very smooth.}
\end{center}
\end{figure}

In Figure~\ref{fig:Ismean} we can see simulations of the Markov process in
(\ref{CurieWeiss}) on a lattice with $N=10$, $100$, $1000$, and
$10,000$ sites, respectively. It appears that in the limit $N\to\infty$,
the process $\ov X_t$ is given by a smooth, deterministic function.

It is not hard to guess what this function is. Indeed, denoting the generator
of the process in (\ref{CurieWeiss}) by $\ov G_{N,\bet}$, we see that
\[
\E^{\ov x}[\ov X_t]=\ov x+tg_\bet(\ov x)+O(t^2)
\quad\mbox{where}\quad
g_\bet(\ov x):=\ov G_{N,\bet}f(\ov x)
\quad\mbox{with}\quad
f(\ov x):=\ov x.
\]
We call the function $g_\bet$ the \emph{local drift}\index{drift function}
of the process $\ov X$. We calculate
\begin{equation}\begin{array}{r@{\,}c@{\,}l}\label{gbet}
\dis g_\bet(\ov x)&=&\dis r^N_+(\ov x)\cdot\frac{2}{N}-r^N_-(\ov x)\cdot\frac{2}{N}
=\frac{(1-\ov x)e^{\bet\ov x/2}-(1+\ov x)e^{-\bet\ov x/2}}
{e^{\bet\ov x/2}+e^{-\bet\ov x/2}}\\[10pt]
&=&\dis\frac{e^{\bet\ov x/2}-e^{-\bet\ov x/2}}{e^{\bet\ov x/2}+e^{-\bet\ov x/2}}-\ov x
=\tanh(\ha\bet\ov x)-\ov x.
\end{array}\end{equation}
Note that the constant $N$ cancels out of this formula. When $N$ is
large, as long as the process is near the point $\ov x$, it locally
behaves as a rescaled random walk with drift $g_\bet(\ov x)$. In view
of this, by some law of large numbers (that will be made rigorous in
Theorem~\ref{T:toODE} below), we expect $(\ov X_t)_{t\geq 0}$ to
converge in distribution, as $N\to\infty$, to a solution of the
differential equation
\begin{equation}\label{difovX}
\dif{t}\ov X_t=g_\bet(\ov X_t)\qquad(t\geq 0).
\end{equation}

\section{Analysis of the mean-field model}

\begin{figure}[htb]
\begin{center}
\inputtikz{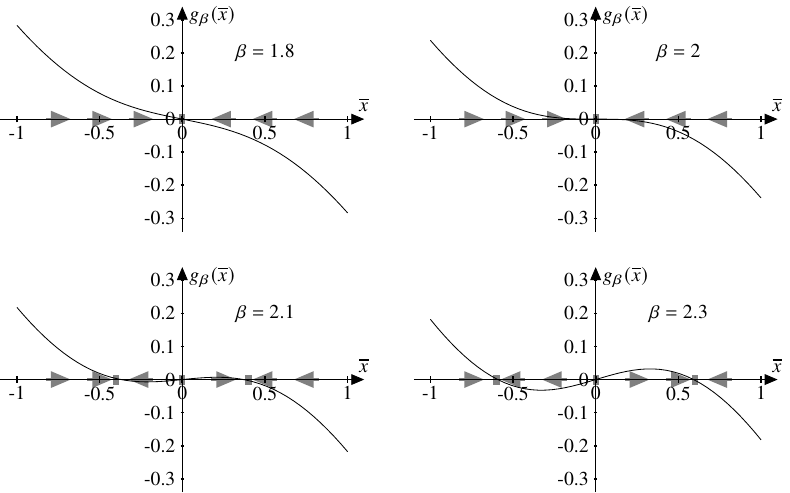}
\caption{The drift function $g_\bet$ for $\bet=1.8$, 2, 2.1, and 2.3,
  respectively. For $\bet>2$, the fixed point $\ov x=0$ becomes
  unstable and two new fixed points appear.}
\label{fig:gbet}
\commentAlt{Figure~\ref{fig:gbet}}{Four pictures of the graph of the
  function $g_\beta$ for different values of $\beta$, demonstrating
  how the zero fixed point becomes unstable and two new fixed points
  emerge. See long description.}
\commentLongAlt{Figure~\ref{fig:gbet}}{In the first picture the graph
  is strictly decreasing and passes zero in zero. The fixed point zero
  is attractive, as indicated by arrows on the horizontal axis
  pointing towards zero. The second picture is the same but the
  derivative of the function in zero is zero. In the third and fourth
  picture the derivative at zero is positive and the graph passes the
  horizontal axis three times, corresponding to three fixed points of
  which the middle one, at zero, is unstable and the other two are
  stable. This is indicated by arrows on the horizontal axis pointing
  towards the stable fixed points. The derivative at zero is larger in
  the fourth picture and the stable fixed points are further from
  zero.}
\end{center}
\end{figure}

Assuming the correctness of (\ref{difovX}) for the moment, we can study the
behavior of the mean-field Ising model $\ov X$ in the limit that we first send
$N\to\infty$, and then $t\to\infty$. A simple analysis of the function
$g_\bet$ (see Figure~\ref{fig:gbet}) reveals that the differential equation
(\ref{difovX}) has a single fixed point for $\bet\leq 2$, and three fixed
points for $\bet>2$. Here, with a \emph{fixed point} of the differential
equation, we mean a point $z$ such that $\ov x_0=z$ implies $\ov x_t=z$ for all
$t\geq 0$, that is, this is a point such that $g_\bet(z)=0$.

\begin{figure}[htb]
\begin{center}
\inputtikz{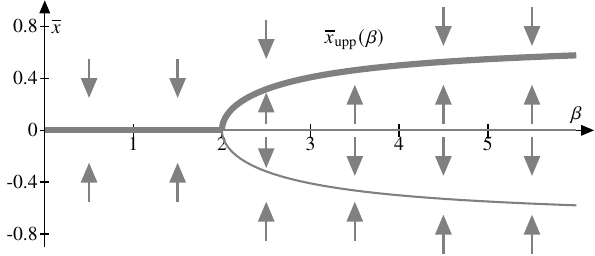}
\caption{Fixed points of the mean-field Ising model as a function of $\bet$,
  with their domains of attraction. The upper fixed point as a function of
  $\bet$ is indicated with a bold line.}
\label{fig:meanphase}
\commentAlt{Figure~\ref{fig:meanphase}}{Graph showing the values of
  the fixed points as a function of $\beta$, with arrows indicating
  the domains of attraction. One fixed point before the point 2, three
  thereafter.}
\end{center}
\end{figure}

Indeed, using the facts that $\tanh$ is an odd function that is concave on
$\half$ and satisfies $\dif{x}\tanh(x)|_{x=0}=1$, we see that:
\begin{itemize}
\item For $\bet\leq 2$, the equation $g_\bet(x)=0$ has the
 unique solution $x=0$.
\item For $\bet>2$, the equation $g_\bet(x)=0$ has three solutions $x_-<0<x_+$.
\end{itemize}

For $\bet\leq 2$, solutions to the differential equation (\ref{difovX})
converge to the unique fixed point $x=0$ as time tends to infinity. On the other
hand, for $\bet>2$, the fixed point $x=0$ becomes unstable. Solutions $\ov X$
to the differential equation (\ref{difovX}) starting in $\ov X_0>0$ converge
to $x_+$, while solutions starting in $\ov X_0<0$ converge to $x_-$.

In Figure~\ref{fig:meanphase}, we have plotted the three fixed points
$x_-<0<x_+$ as a function of $\bet$, and indicated their domains of
attraction. The function
\begin{equation}\index{0xupp@$x_{\rm upp}(\bet)$}\label{xupp}
x_{\rm upp}(\bet):=\left\{\begin{array}{ll}
0\quad&\mbox{if }\bet\leq 2,\\[5pt]
\mbox{the unique positive solution of }\tanh(\ha\bet x)=x\quad&\mbox{if }\bet>2
\end{array}\right.
\end{equation}
plays a similar role as the spontaneous magnetization $m_\ast(\bet)$ for the
Ising model on $\Z^d$ (see formula (\ref{mast})). More precisely, for
mean-field processes started in initial states $\ov X_0>0$, the quantity
$x_{\rm upp}$ describes the double limit
\begin{equation}\label{doublim}
\lim_{t\to\infty}\lim_{N\to\infty}\ov X_t=x_{\rm upp}.
\end{equation}

We see from (\ref{xupp}) that the mean-field Ising model (as formulated in
(\ref{meanGlauber})) exhibits a second-order (that is, continuous) phase
transition at the critical point $\bet_{\rm c}=2$. Since
\[
x_{\rm upp}(\bet)\propto(\bet-\bet_{\rm c})^{1/2}
\quad\mbox{as}\quad\bet\down\bet_{\rm c},
\]
the \emph{mean-field critical exponent} associated with the
magnetization\footnote{In general, for a given second-order phase transition,
  there are several quantities of interest that all show power-law behavior
  near the critical point, and hence there are also several critical exponents
  associated with a given phase transition.} is $c=1/2$, which is the same as
for the Ising model on $\Z^d$ in dimensions $d\geq 4$ (see
Section~\ref{S:phase}). Understanding why the mean-field model correctly
predicts the critical exponent in sufficiently high dimensions goes beyond the
scope of the present chapter.

To conclude the present section, we note that the two limits in
(\ref{doublim}) cannot be interchanged. Indeed, for each fixed $N$,
the Markov process $\ov X$ is irreducible, and hence, by
Theorem~\ref{T:finergo}, has a unique equilibrium law that is the
long-time limit of the law at time $t$, started from an arbitrary initial
state. In view of the symmetry of the problem, the magnetization in
equilibrium must be zero, so regardless of the initial state, we have,
for each fixed $N$,
\[
\lim_{t\to\infty}\E[\ov X_t]=0.
\]
The reason why this can be true while at the same time (\ref{doublim}) also
holds is that the speed of convergence to equilibrium of the Markov process
$\ov X$ becomes very slow as $N\to\infty$.

\begin{figure}[htb]
\begin{center}
\inputtikz{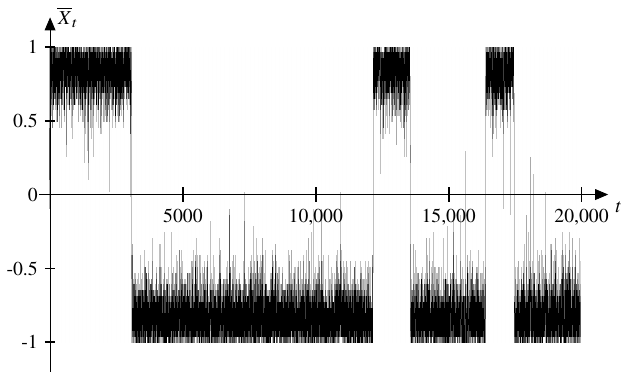}
\caption{Metastable behavior of a mean-field Ising model with $N=50$ and
  $\bet=3$. Note the different time scale compared to Figure~\ref{fig:Ismean}.}
\label{fig:meta}
\commentAlt{Figure~\ref{fig:meta}}{Graph of the value of the
  magnetization in one random run between time zero and 20,000. The
  value mostly stays near the values 0.9 and -0.9, with occasional
  quick transitions from one to the other. See long description.}
\commentLongAlt{Figure~\ref{fig:meta}}{Since time is very long the
  function moves very wildly. When it is positive, it mostly stays
  between 0.7 and 1, so this area is almost all black. There are many
  small spikes pointing downwards indicating that the function
  occasionally reached a low value for a very short time before moving
  back up. Around time 3,000, one such spike suddenly crosses the
  horizontal axis and then the functions stays between -1 and -0.7 for
  a long time, till around time 13,000 it moves back up. Thee more
  transitions (down, up, down again) can be observed.}
\end{center}
\end{figure}

In Figure~\ref{fig:meta}, we have plotted the time evolution of a mean-field
Ising model $\ov X$ on a lattice with $N=50$ sites, for a value of $\bet$
above the critical point (concretely $\bet=3$, which lies above $\bet_{\rm
  c}=2$). Although the average of $\ov X$ in the long run is $0$, we see that
the process spends most of its time around the values $x_{\rm upp}$ and
$-x_{\rm upp}$, with rare transitions between the two. This sort of behavior
is called \emph{metastable behavior}.\index{metastable behavior}

The value $N=50$ was near the highest possible value for which I could still
numerically observe this sort of behavior. For $N=100$ the transitions between
the two metastable states $x_{\rm upp}$ and $-x_{\rm upp}$ become so rare that
my program was no longer able to see them within a reasonable runtime.  With
the help of \emph{large deviations theory}, one can show that the time that
the system spends in one metastable state is approximately exponentially
distributed (with a large mean), and calculate the asymptotics of the mean
waiting time as $N\to\infty$. It turns out that the mean time one has to wait
for a transition grows exponentially fast in $N$.

\section{Functions of Markov processes}\label{S:Marfun}

In the present section we formulate a proposition and a theorem that we have
already implicitly used. Both are concerned with functions of Markov
processes. Let $X=(X_t)_{t\geq 0}$ be a Markov process with finite state space
$S$, generator $G$, and semigroup $(P_t)_{t\geq 0}$. Let $T$ be another finite
set and let $f\cn S\to T$ be a function. For each $x\in S$ and $y'\in T$ such
that $f(x)\neq y'$, let
\begin{equation}\label{Hxyy}
\Hi(x,y'):=\sum_{x'\in S:\, f(x')=y'}G(x,x')
\end{equation}
be the total rate at which $f(X_t)$ jumps to the state $y'$, when the present
state is $X_t=x$. The next proposition says that if these rates are a function
of $f(x)$ only, then the process $Y=(Y_t)_{t\geq 0}$ defined by
\begin{equation}\label{YfX}
Y_t:=f(X_t)\qquad(t\geq 0)
\end{equation}
is itself a Markov process.

\begin{proposition}[Autonomous Markov process]\label{P:auto}\hspace{-4.25pt}
Assume that the rates in (\ref{Hxyy}) are of the form
\begin{equation}\label{HiHa}
\Hi(x,y')=H\big(f(x),y'\big)\qquad(x\in S,\ y'\in T,\ f(x)\neq y')
\end{equation}
where $H$ is a Markov generator of some process in $T$.
Then the process $Y$ defined in (\ref{YfX}) is a Markov process
with generator $H$. Conversely, if for each initial law of the process $X$, it
is true that $Y$ is a Markov process with generator $H$, then (\ref{HiHa})
must hold.
\end{proposition}

\begin{Proof}[of Proposition~\ref{P:auto}]
Let us define $\Hi(x,y')$ as in (\ref{Hxyy}) also when $f(x)=y'$.
We start by noting that if (\ref{HiHa}) holds for all $x\in S$ and $y'\in T$
such that $f(x)\neq y'$, then it also holds when $f(x)=y'$. To see this, we
write
\[\begin{array}{l}
\dis H\big(f(x),f(x)\big)
=-\sum_{y':\,y'\neq f(x)}H(f(x),y')
=-\sum_{y':\,y'\neq f(x)}\Hi(x,y')\\[5pt]
\dis\quad=-\sum_{y':\,y'\neq f(x)}\sum_{x':\,f(x')=y'}G(x,x')
=-\sum_{x':\,f(x')\neq f(x)}G(x,x')
=\sum_{x':\,f(x')=f(x)}G(x,x'),
\end{array}\]
where we have used that since $H$ and $G$ are Markov generators, one has
$\sum_{y'\in T}H(f(x),y')=0$ and $\sum_{x'\in S}G(x,x')=0$.
We have thus shown that (\ref{HiHa}) is equivalent to
\begin{equation}\label{HiHa2}
H\big(f(x),y'\big)=\sum_{x':\,f(x')=y'}G(x,x')
\qquad\big(x\in S,\ y'\in T\big).
\end{equation}
We claim that this is equivalent to
\begin{equation}\label{QP}
Q_t\big(f(x),y'\big)=\sum_{x':\,f(x')=y'}P_t(x,x')
\qquad\big(t\geq 0,\ x\in S,\ y'\in T\big),
\end{equation}
where $(Q_t)_{t\geq 0}$ is the semigroup generated by $H$. To prove this, we
start by observing that for any function $g\cn T\to\R$,
\[\begin{array}{l}
\dis G(g\circ f)(x)=\sum_{x'}G(x,x')g(f(x'))
=\sum_{y'}\sum_{x':\,f(x')=y'}G(x,x')g(y'),\\[5pt]
\dis(Hg)\circ f(x)=\sum_{y'}H(f(x),y')g(y').
\end{array}\]
The right-hand sides of these equations are equal for all $g\cn T\to\R$ if and
only if (\ref{HiHa2}) holds, so (\ref{HiHa2}) is equivalent to the statement
that
\begin{equation}\label{HiHa3}
G(g\circ f)=(Hg)\circ f\qquad(g\cn T\to\R).
\end{equation}
By exactly the same argument with $G$ replaced by $P_t$ and $H$ replaced by
$Q_t$, we see that (\ref{QP}) is equivalent to
\begin{equation}\label{QP2}
P_t(g\circ f)=(Q_tg)\circ f\qquad(t\geq 0,\ g\cn T\to\R).
\end{equation}
To see that (\ref{HiHa3}) and (\ref{QP2}) are equivalent, we write
\begin{equation}\label{PQexp}
P_t=e^{Gt}=\sum_{n=0}^\infty\frac{1}{n!}t^nG^n
\quand
Q_t=e^{Ht}=\sum_{n=0}^\infty\frac{1}{n!}t^nH^n.
\end{equation}
We observe that (\ref{HiHa3}) implies
\[
G^2(g\circ f)=G\big((Hg)\circ f\big)=(H^2g)\circ f,
\]
and similarly, by induction, $G^n(g\circ f)=(H^ng)\circ f$ for all $n\geq 0$,
which by (\ref{PQexp}) implies (\ref{QP2}). Conversely, if (\ref{QP2}) holds
for all $t\geq 0$, then it must hold up to first order in $t$ as $t\down 0$,
which implies (\ref{HiHa3}). This completes the proof that (\ref{HiHa}) is
equivalent to (\ref{QP}).

If (\ref{QP}) holds, then by (\ref{Markfdd}), the finite dimensional
distributions of $Y$ are given by
\begin{equation}\begin{array}{l}\label{XfY}
\dis\P\big[Y_0=y_0,\ldots,Y_{t_n}=y_n\big]\\[5pt]
\dis\quad=\sum_{x_0:\,f(x_0)=y_0}\cdots\sum_{x_n:\,f(x_n)=y_n}
\P[X_0=x_0]P_{t_1-t_0}(x_0,x_1)\cdots P_{t_n-t_{n-1}}(x_{n-1},x_n)\\[5pt]
\dis\quad=\P[Y_0=y_0]
Q_{t_1-t_0}(y_0,y_1)\cdots Q_{t_n-t_{n-1}}(y_{n-1},y_n)
\end{array}\end{equation}
$(0=t_0<\cdots<t_n)$. Again by (\ref{Markfdd}), this implies that $Y$
is a Markov process with generator $H$.

Conversely, if $Y$ is a Markov process with generator $H$ for each initial
state of $X$, then for each $x_0\in S$, (\ref{XfY}) must hold when $X_0=x_0$
a.s.\ and for $n=1$, from which we see that (\ref{QP}) and hence (\ref{HiHa})
hold.
\end{Proof}

Summarizing, Proposition~\ref{P:auto} says that if $Y_t=f(X_t)$ is a function
of a Markov process, and the jump rates of $Y$ are a function of the present
state of $Y$ only (and do not otherwise depend on the state of $X$), then $Y$
is itself a Markov process. In such a situation, we will say that $Y$ is an
\emph{autonomous}\index{autonomous Markov process} Markov process.
We have already implicitly used Proposition~\ref{P:auto} in
Section~\ref{S:Ismean}, when we claimed that the process $\ov X$ is a Markov
process with jump rates as in (\ref{CurieWeiss}).\med

\noi
\textbf{Remark} For the final statement of the proposition, it is essential
that $Y$ is a Markov process for \emph{each} initial law $X$. There exist
interesting examples of functions of Markov processes that are not autonomous
Markov processes, but nonetheless are Markov processes for some \emph{special}
initial laws of the original Markov process. This is closely related to the
concept of intertwining of Markov processes that will briefly be mentioned
in Section~\ref{S:other} below.\med

Our next aim is to make the claim rigorous that for large $N$, the
process $\ov X$ can be approximated by solutions to the differential
equation (\ref{difovX}). The basic idea is that since the process
makes many small steps, as long as $\ov X_t\approx\ov x$, by some sort
of ``local'' law of large numbers, the process should
deterministically increase at speed $\approx g_\bet(\ov x)$. To make
this precise we will apply a theorem from \cite{DN08}. Although the
proof is not very complicated, it is a bit lengthy and would detract
from our main objects of interest here, so we only show how the
theorem below can be deduced from a theorem in \cite{DN08}. That paper
also treats the multi-dimensional case and gives explicit estimates on
probabilities of the form (\ref{toODE}) below. An alternative, more
probabilistic approach to mean-field equations is described in
\cite{MSS20}.

For each $N\geq 1$, let $X^N=(X^N_t)_{t\geq 0}$ be a Markov process
with finite state space $S_N$, generator $G_N$, and semigroup
$(P^N_t)_{t\geq 0}$, and let $f_N\cn S_N\to\R$ be functions. We will
be interested in conditions under which the processes
$(f_N(X^N_t))_{t\geq 0}$ approximate the solution $(y_t)_{t\geq 0}$ of
a differential equation, in the limit $N\to\infty$.  Note that we do
not require that $f_N(X^N_t)$ is an autonomous Markov process.  To
ease notation, we will sometimes drop the super- and subscripts $N$
when no confusion arises.

We define two functions $\al=\al_N$ and $\bet=\bet_N$ that describe the
quadratic variation and drift, respectively, of the process $f(X_t)$. More
precisely, these functions are given by
\[\begin{array}{r@{\,}c@{\,}l}
\dis\al(x)&:=&\dis\sum_{x'\in S}G(x,x')\big(f(x')-f(x)\big)^2,\\[5pt]
\dis\bet(x)&:=&\dis\sum_{x'\in S}G(x,x')\big(f(x')-f(x)\big).
\end{array}\]
The idea is that if $\al$ tends to zero and $\bet$ approximates a nice,
Lipschitz continuous function of $f(X_t)$, then $f(X_t)$ should in the limit
be given by the solution of a differential equation.

We assume that the functions $f_N$ all take values in a closed interval
$I\sub\R$ with left and right boundaries $I_-:=\inf I$ and $I_+:=\sup I$,
which may be finite or infinite. We also assume that there exists a 
globally Lipschitz function $b\cn I\to\R$ such that
\begin{equation}\label{betb}
\sup_{x\in S_N}\big|\bet_N(x)-b\big(f_N(x)\big)\big|\asto{N}0,
\end{equation}
that is, the drift function $\bet$ is uniformly approximated by $b\circ f_N$.
Assuming also that
\begin{equation}
b(I_-)\geq 0\quad\mbox{if }I_->-\infty
\quand
b(I_+)\leq 0\quad\mbox{if }I_+<-\infty,
\end{equation}
the differential equation
\[
\dif{t}y_t=b(y_t)\qquad(t\geq 0)
\]
has a unique $I$-valued solution $(y_t)_{t\geq 0}$ for each initial state
$y_0\in I$. The following theorem gives sufficient conditions for the
$I$-valued processes $(f_N(X^N_t))_{t\geq 0}$ to approximate a solution of the
differential equation.

\begin{theorem}[Limiting differential equation]\label{T:toODE}\hspace{5pt}
Assume that $f_N(X^N_0)$ converges in probability to $y_0$ and that as well as
(\ref{betb}), one moreover has
\begin{equation}\label{ala}
\sup_{x\in S_N}\al_N(x)\asto{N}0.
\end{equation}
Then, for each $T<\infty$ and $\eps>0$,
\begin{equation}\label{toODE}
\P\big[|f_N(X^N_t)-y_t|\leq\eps\ \forall t\in[0,T]\big]\asto{N}1.
\end{equation}
\end{theorem}

\begin{Proof}
We apply \cite[Thm~4.1]{DN08}. Fix $T<\infty$ and $\eps>0$ and also fix
$y_0\in I$. Let $L$ denote the Lipschitz constant of $b$. The assumptions of
\cite[Thm~4.1]{DN08} allow for the case that $f_N$ does not in general take
values in $I$, but only under the additional condition that $f_N(x)$ is not
further than $\eps$ from a possible value the solution of the differential
equation. In our case, these more general assumptions are automatically
satisfied. Set $\de:=\frac{1}{3}\eps e^{-LT}$. We consider the events
\[
\Om_0:=\big\{|f(X_0)-y_0|\leq\de\big\}
\quand
\Om_1:=\big\{\int_0^T\!|\bet(X_t)-b\big(f(X_t)\big)|\,\di t\leq\de\big\}.
\]
For $K>0$, we also define
\[
\Om_{K,2}:=\big\{\int_0^T\!\al(X_t)\,\di t\leq KT\big\}.
\]
Then \cite[Thm~4.1]{DN08} tells us that
\begin{equation}\label{DN08}
\P\big[\sup_{t\in[0,T]}|f(X_t)-y_t|>\eps\big]
\leq 4KT\de^{-2}+\P\big(\Om_0^{\rm c}\cup\Om_1^{\rm c}\cup\Om_{K,2}^{\rm c}\big).
\end{equation}
Our assumption that $f_N(X^N_0)\to y_0$ in probability implies that
$\P(\Om_0^{\rm c})\to 0$ as $N\to\infty$. Set
\[
A_N:=\sup_{x\in S_N}\al_N(x)
\quand
B_N:=\sup_{x\in S_N}\big|\bet_N(x)-b\big(f_N(x)\big)\big|
\]
Then $A_N\to 0$ by (\ref{ala}) and $B_N\to 0$ by (\ref{betb}).
Since
\[
\int_0^T|\bet(X_t)-b\big(f(X_t)\big)|\,\di t\leq B_NT\leq\de
\]
for $N$ sufficiently large, we see that $\P(\Om_1^{\rm c})=0$
for $N$ sufficiently large. Also, since
\[
\int_0^T\al(X_t)\,\di t\leq A_NT,
\]
we see that $\P(\Om_{A_N,2}^{\rm c})=0$ for all $N$. Inserting $K=A_N$ in
(\ref{DN08}), we see that the right-hand side tends to zero as $N\to\infty$.
\end{Proof}

Using Theorem~\ref{T:toODE}, we can make the approximation of the mean-field
Ising model by the differential equation (\ref{difovX}) rigorous. Let
$X^N=(X^N_t)_{t\geq 0}$ denote the Markov process with state space
$\{-1,+1\}^{\La_N}$, where $\La_N$ is a set containing $N$ elements and the
jump rates of $X^N$ are given in (\ref{meanGlauber}). By
Proposition~\ref{P:auto}, the process
$\ov X^N_t:=\frac{1}{N}\sum_{i\in\La_N}X_t(i)$
is itself a Markov process with jump rates as in
(\ref{CurieWeiss}). We can either apply Theorem~\ref{T:toODE} directly
to the Markov processes $X^N$ and the functions
$f_N(x):=\frac{1}{N}\sum_{i\in\La_N}x(i)$, or we can apply
Theorem~\ref{T:toODE} to the Markov processes $\ov X^N$ and choose for
$f_N$ the identity function $f_N(\ov x)=\ov x$. In either case, the
assumption (\ref{betb}) has already been verified in (\ref{gbet}). To check
also (\ref{ala}), we calculate
\[
\al_N(x)=r^N_+(\ov x)\Big(\frac{2}{N}\Big)^2+r^N_-(\ov x)\Big(\frac{2}{N}\Big)^2
=\frac{2}{N}\Big(1+\ov x\frac{e^{-\bet\ov x/2}-e^{\bet\ov x/2}}
{e^{-\bet\ov x/2}+e^{\bet\ov x/2}}\Big),
\]
which clearly tends uniformly to zero as $N\to\infty$.

\section{The mean-field contact process}\index{contact process!mean-field}
\label{S:contmean}

Recall the definition of the generator of the contact process from
(\ref{Gcontact}). We slightly reformulate this as
\begin{equation}\begin{array}{r@{\,}c@{\,}l}\label{Gcontact2}
\dis G_{\rm cont}f(x)
&:=&\dis\la\sum_{i\in\Z^d}\frac{1}{|\Ni_i|}\sum_{j\in\Ni_i}
\big\{f(\big({\tt bra}_{ij}(x))-f\big(x\big)\big\}\\[5pt]
&&\dis+\sum_{i\in\Z^d}\big\{f(\big({\tt death}_i(x))-f\big(x\big)\big\}
\qquad(x\in\{0,1\}^\La),
\end{array}\end{equation}
where as customary we have set the death rate to $\de=1$, and we have also
reparametrized the infection rate so that $\la$ denotes the total rate of all
outgoing infections from a given site, instead of the infection rate per
neighbor.

We will be interested in the contact process on the complete graph,
which means that we take for $\La=\La_N$ a set with $N$ elements,
which we equip with the structure of a complete graph with
(undirected) edge set $E=E_N:=\{\{i,j\}:i,j\in\La_N\}$ and
corresponding set of directed edges $\Ei=\Ei_N$. We view $i$ as a
neighbor of itself, but since ${\tt bra}_{ii}$ is the identity map
this has no effect. We will be interested in the fraction of infected
sites
\[
\ov X_t=\ov X^N_t:=\frac{1}{N}\sum_{i\in\La_N}X_t(i)\qquad(t\geq 0),
\]
which jumps with the following rates
\begin{equation}\begin{array}{rcl}\label{meancontact}
\dis\ov x\mapsto\ov x+\ffrac{1}{N}
&\dis\quad\mbox{with rate}\quad
&\dis r^N_+(\ov x):=\la N\ov x(1-\ov x),\\[10pt]
\dis\ov x\mapsto\ov x-\ffrac{1}{N}
&\dis\quad\mbox{with rate}\quad
&\dis r^N_-(\ov x):=N\ov x.
\end{array}\end{equation}
Here $N(1-\ov x)$ is the number of healthy sites, each of which gets infected
with rate $\la\ov x$, and $N\ov x$ is the number of infected sites, each of
which recovers with rate one. Note that since these rates are a function of
$\ov x$ only, by Proposition~\ref{P:auto}, the process $(\ov X_t)_{t\geq 0}$
is an autonomous Markov process.

\begin{figure}[htb]
\begin{center}
\inputtikz{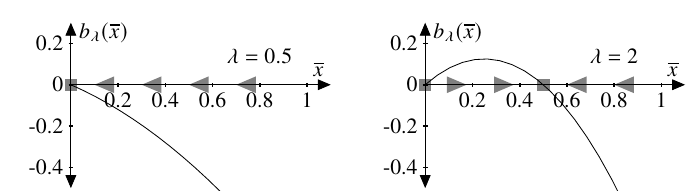}
\caption{The function $b_\la$ for two values of $\la$.}
\label{fig:contdrift}
\commentAlt{Figure~\ref{fig:contmean}}{Two pictures showing the graph
  of the function $b_\lambda$ for two values of $\lambda$,
  demonstrating how the zero fixed point becomes unstable and a new
  fixed point emerges.}
\end{center}
\end{figure}

We wish to apply Theorem~\ref{T:toODE} to conclude that $\ov X$ can for large
$N$ be approximated by the solution of a differential equation. To this aim,
we calculate the drift $\bet$ and quadratic variation function $\al$.
\[\begin{array}{r@{\,}c@{\,}l}
\dis\al_N(x)&=&\dis r^N_+(\ov x)\ffrac{1}{N^2}+r^N_-(\ov x)\ffrac{1}{N^2}
=\frac{1}{N}\big(\la\ov x(1-\ov x)+\ov x\big),\\[8pt]
\dis\bet_N(x)&=&\dis r^N_+(\ov x)\ffrac{1}{N}-r^N_-(\ov x)\ffrac{1}{N}
=\la\ov x(1-\ov x)-\ov x.
\end{array}\]
By Theorem~\ref{T:toODE}, it follows that in the mean-field limit
$N\to\infty$, the fraction of infected sites can be approximated by solutions
of the differential equation
\begin{equation}\label{meancont}
\dif{t}\ov X_t=b_\la(\ov X_t)\quad(t\geq 0),
\quad\mbox{where}\quad
b_\la(\ov x):=\la\ov x(1-\ov x)-\ov x.
\end{equation}
The equation $b_\la(\ov x)=0$ has the solutions
\[
\ov x=0\quand\ov x=1-\la^{-1}.
\]
The second solution lies inside the interval $[0,1]$ of possible values of
$\ov X_t$ if and only if $\la\geq 1$. Plotting the function $b_\la$ for
$\la<1$ and $\la>1$ yields Figure~\ref{fig:contdrift}.

\begin{figure}[htb]
\begin{center}
\inputtikz{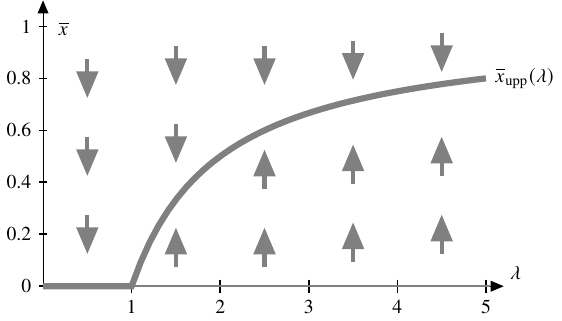}
\caption{Mean-field analysis of the contact process.}
\label{fig:contmean}
\commentAlt{Figure~\ref{fig:contmean}}{Graph showing the values of the
  fixed points as a function of $\lambda$, with arrows indicating the
  domains of attraction.}
\end{center}
\end{figure}

We see from this that the fixed point $\ov x=0$ is stable for $\la\leq 1$ but
becomes unstable for $\la>1$, in which case $\ov x=1-\la^{-1}$ is the only
stable fixed point that attracts all solutions started in a nonzero initial
state. The situation is summarized in Figure~\ref{fig:contmean}.

Letting $\ov x_{\rm upp}(\la):=0\vee(1-\la^{-1})$ denote the stable fixed
point, we see that the mean-field contact process exhibits a second-order
phase transition at the critical point $\la_{\rm c}=1$. Since
\[
\ov x_{\rm upp}(\la)\propto(\la-\la_{\rm c})\quad\mbox{as }\la\down\la_{\rm c},
\]
the associated critical exponent is $c=1$, in line with what we know for contact
processes in dimensions $d\geq 4$ (see the discussion in
Section~\ref{S:phase}).

\section{The mean-field voter model}

Recall the definition of the generator of the voter model from (\ref{Gvot}).
For simplicity, we will only consider the two-type model and as the local
state space we will choose $S=\{0,1\}$. Specializing to the complete graph
$\La=\La_N$ with $N$ vertices, the generator becomes
\begin{equation}\label{meanvoter}
G_{\rm vot}f(x)
=\frac{1}{|\La|}\sum_{(i,j)\in\Ei}
\big\{f(\big({\tt vot}_{ij}(x))-f\big(x\big)\big\}
\qquad(x\in\{0,1\}^\La).
\end{equation}
Note that the factor $|\La|^{-1}$ says that each site $i$ updates its type with
rate one, and at such an event chooses a new type from a uniformly chosen site
$j$ (allowing for the case $i=j$, which has no effect).

We are interested in the fraction of sites of type $1$,
\[
\ov X_t=\ov X^N_t:=\frac{1}{N}\sum_{i\in\La_N}X_t(i)\qquad(t\geq 0),
\]
which jumps as (compare (\ref{meancontact}))
\[\begin{array}{rcl}
\dis\ov x\mapsto\ov x+\ffrac{1}{N}
&\dis\quad\mbox{with rate}\quad
&\dis r^N_+(\ov x):=N\ov x(1-\ov x),\\[10pt]
\dis\ov x\mapsto\ov x-\ffrac{1}{N}
&\dis\quad\mbox{with rate}\quad
&\dis r^N_-(\ov x):=N\ov x(1-\ov x).
\end{array}\]
Note that $N(1-\ov x)$ is the number of sites of type 0, and that each
such site adopts the type 1 with rate $\ov x$. The derivation of
$r^N_-(\ov x)$ is similar.  We calculate the drift $\bet$ and
quadratic variation function $\al$.
\[\begin{array}{r@{\,}c@{\,}l}
\dis\al_N(x)&=&\dis r^N_+(\ov x)\ffrac{1}{N^2}+r^N_-(\ov x)\ffrac{1}{N^2}
=\frac{2}{N}\ov x(1-\ov x),\\[8pt]
\dis\bet_N(x)&=&\dis r^N_+(\ov x)\ffrac{1}{N}-r^N_-(\ov x)\ffrac{1}{N}=0.
\end{array}\]
Applying Theorem~\ref{T:toODE}, we see that in the limit $N\to\infty$, the
process $(\ov X_t)_{t\geq 0}$ is well approximated by solutions to the
differential equation
\[
\dif{t}\ov X_t=0\qquad(t\geq 0),
\]
that is, $\ov X_t$ is approximately constant as a function of $t$.

\begin{figure}[htb]
\begin{center}
\inputtikz{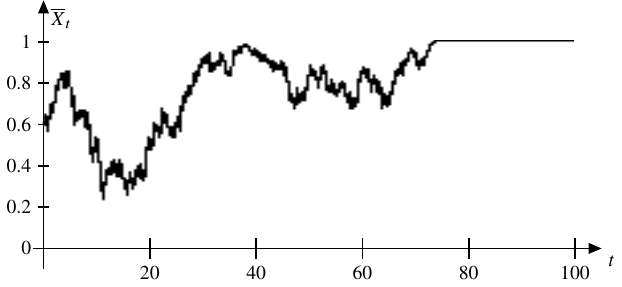}
\caption{The fraction of type 1 individuals in the mean-field voter model
  from (\ref{speedvoter}) on $N=100$ sites. This process approximates the
  Wright--Fisher diffusion.}
\label{fig:meanvoter}
\commentAlt{Figure~\ref{fig:meanvoter}}{Graph of a random realization
  of a Wright-Fisher diffusion. The process locally looks like a
  Brownian motion but moves a bit less fast near 1 and after some time
  gets absorbed in 1.}
\end{center}
\end{figure}

Of course, if we go to larger time scales, then $\ov X_t$ will no longer be
constant; compare Figure~\ref{fig:meta}. In fact, we can determine the time
scale at which $\ov X_t$ fluctuates quite precisely. Speeding up time by a
factor $|\La|=N$ is the same as multiplying all rates by a factor $|\La|$.
If we repeat our previous calculations for the process with generator
\begin{equation}\label{speedvoter}
G_{\rm vot}f(x)
=\sum_{(i,j)\in\Ei}
\big\{f(\big({\tt vot}_{ij}(x))-f\big(x\big)\big\}
\qquad(x\in\{0,1\}^\La),
\end{equation}
then the drift and quadratic variation are given by
\[\begin{array}{r@{\,}c@{\,}l}
\dis\al_N(x)&=&\dis 2\ov x(1-\ov x),\\[8pt]
\dis\bet_N(x)&=&\dis 0.
\end{array}\]
In this case, the quadratic variation does not go to zero, so
Theorem~\ref{T:toODE} is no longer applicable. One can show, however,
that in the limit $N\to\infty$ the new, sped-up process is well approximated
by solutions to the (It\^o) stochastic differential equation (SDE)
\[
\di\ov X_t=\sqrt{2\ov X_t(1-\ov X_t)}\,\di B_t\qquad(t\geq 0),
\]
where $2\ov X_t(1-\ov X_t)=\al(X_t)$ is of course the quadratic
variation function we have just calculated. Solutions to this SDE are
\emph{Wright--Fisher diffusions},\index{Wright--Fisher diffusion} that
is, Markov processes with continuous sample paths and generator
\begin{equation}\label{WF}
Gf(\ov x)=\ov x(1-\ov x)\diff{\ov x}f(\ov x).
\end{equation}
These calculations can be made rigorous using methods from the theory of
convergence of Markov processes; see, for example, the book \cite{EK86}.
See Figure~\ref{fig:meanvoter} for a simulation of the process $\ov X$
when $X$ has the generator in (\ref{speedvoter}) and $N=100$.

\section{Exercises}

\begin{figure}[htb]
\begin{center}
\inputtikz{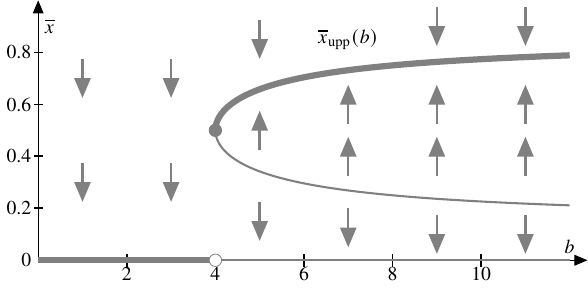}
\caption{Mean-field analysis of a model with cooperative branching and deaths.}
\label{fig:coopmean}
\commentAlt{Figure~\ref{fig:coopmean}}{Graph showing the values of
    the fixed points as a function of $b$, with arrows indicating the
    domains of attraction. Initially zero is the only fixed point but
    at b=4 a second fixed point appears at 0.5 that then splits into
    two fixed points.}
\end{center}
\end{figure}

\begin{Exercise}
Do a mean-field analysis of the process with generator
\[
Gf(x)
=b|\La|^{-2}\sum_{ii'j}
\big\{f\big({\tt coop}_{ii'j}x\big)-f\big(x\big)\big\}+\sum_i
\big\{f\big({\tt death}_ix\big)-f\big(x\big)\big\},
\]
where the maps ${\tt coop}_{ii'j}$ and ${\tt death}_i$ are defined in
(\ref{coopmap}) and (\ref{deathmap}), respectively.
Do you observe a phase transition? Is it first- or second order?
\emph{Hint:} Figure~\ref{fig:coopmean}.
\end{Exercise}

\begin{Exercise}
Same as above for the model with generator
\[
Gf(x)
=b|\La|^{-2}\sum_{ii'j}
\big\{f\big({\tt coop}_{ii'j}x\big)-f\big(x\big)\big\}+|\La|^{-1}\sum_{ij}
\big\{f\big({\tt rw}_{ij}x\big)-f\big(x\big)\big\}.
\]
\end{Exercise}

\begin{Exercise}
Derive an SDE in the limit $|\La|\to\infty$ for the density of the mean-field
voter model with small bias and death rates, with generator
\[\begin{array}{r@{\,}c@{\,}l}
\dis Gf(x)
&=&\dis\sum_{ij\in\La}
\big\{f\big({\tt vot}_{ij}x\big)-f\big(x\big)\big\}
+s|\La|^{-1}\sum_{ij\in\La}
\big\{f\big({\tt bra}_{ij}x\big)-f\big(x\big)\big\}\\[5pt]
&&\dis+d\sum_{i\in\La}
\big\{f\big({\tt death}_ix\big)-f\big(x\big)\big\}.
\end{array}\]
\emph{Hint:} You should find expressions of the form
\[\begin{array}{r@{\,}c@{\,}l}
\dis\E^{\ov x}\big[(\ov X_t-\ov x)\big]&=&\dis b(\ov x)\cdot t+O(t^2),\\[5pt]
\dis\E^{\ov x}\big[(\ov X_t-\ov x)^2\big]&=&\dis a(\ov x)\cdot t+O(t^2),
\end{array}\]
which leads to a limiting generator of the form
\[
Gf(\ov x)=\ha a(\ov x)\diff{\ov x}f(\ov x)+b(\ov x)\dif{x}f(\ov x).
\]
\end{Exercise}

\begin{Exercise}
Do a mean-field analysis of the following more general version of the
Neuhauser-Pacala model \cite{NP99} from (\ref{NP99}). In the more
general model, the site $i$ flips
\[\begin{array}{l}
0\mapsto 1\quad\mbox{with rate}\quad\big(f_0+\al_{01}f_1\big)f_1,\\[5pt]
1\mapsto 0\quad\mbox{with rate}\quad\big(f_1+\al_{10}f_0\big)f_0,
\end{array}\]
where $\al_{01},\al_{10}>0$ and
$f_\tau=|\Ni_i|^{-1}\sum_{j\in\Ni_i}1_{\{x(j)=\tau\}}$ is the relative
frequency of type $\tau$ in the neighborhood of $i$.
Find all stable and unstable fixed points of the mean-field model in the
regimes: I.\ $\al_{01},\al_{10}<1$, II.\ $\al_{01}<1<\al_{10}$,
III.\ $\al_{10}<1<\al_{01}$, IV.\ $1<\al_{01},\al_{10}$.
\end{Exercise}

\begin{Exercise}
Do a mean-field analysis of the cycle conform model introduced in
Section~\ref{S:period}. \emph{Warning:} this is considerably more complicated
than the previous exercises. Working hard enough, it should be
possible to establish the following facts. Let $F_t(s)$ denote the
fraction of vertices that at time $t$ are in the local state
$s\in\{0,1,2\}$. Then in the mean-field limit, one has
\[
F_t(s)=\ffrac{1}{3}+\ffrac{2}{3}\Re(\ex{is\phi}f_t)
\qquad(s\in\{0,1,2\},\ t\geq 0),
\]
where $\phi:=2\pi/3$ and $t\mapsto f_t\in\C$ solves the differential equation
\[
\dif{t}f_t=\ffrac{1}{3}\al\big\{f_t+f_t^\ast f_t^\ast-2f_tf_tf_t^\ast\big\}
+(1-\al)\big(\ex{-i\phi}-1\big)f_t\qquad(t\geq 0).
\]
Here $\Re(z)$ denotes the real part of a complex number $z$ and
$z^\ast$ denotes its complex conjugate. For $0\leq\al<6/7$ the only
fixed point of this equation is $f_t=0$. This fixed point is stable
for $\al<9/11$ but unstable for $\al>9/11$. In the intermediate regime
$9/11<\al<6/7$, there are periodic solutions.
\end{Exercise}

\chapter{Construction and ergodicity}\label{C:construct}

\section{Introduction}\label{S:conintro}

Throughout this chapter, $S$ is a finite set called the \emph{local
state space}, $\La$ is a countable set called the \emph{lattice},
$\Gi$ is a countable collection of continuous maps $m\cn S^\La\to
S^\La$, and $(r_m)_{m\in\Gi}$ are nonnegative rates. Our aim is to
construct a Markov process with state space $S^\La$ and formal
generator of the form
\begin{equation}\label{Gmap}
Gf(x)=\sum_{m\in\Gi}r_m\big\{f\big(m(x)\big)-f\big(x\big)\big\}
\qquad(x\in S^\La).
\end{equation}
We will usually be interested in the case that all $m\in\Gi$ are local
maps, that is, the set $\Di(m)$ of lattice points whose values can be
changed by $m$ as defined in Section~\ref{S:local} is finite, but our
proofs do not need the finiteness of $\Di(m)$. Under the assumption
\[
\sum_{m\in\Gi}r_m1_{\Di(m)}(i)<\infty\quad(i\in\La)
\]
one can show that $Gf$ is well-defined for all functions $f\cn
S^\La\to\R$ that depend on finitely many coordinates. In general, we
will need stronger conditions on the rates $(r_m)_{m\in\Gi}$ to ensure
that $G$ generates a Markov process. Most of the interacting particle
systems introduced in Chapter~\ref{C:intro} have a generator that can
naturally be written in the form (\ref{Gmap}). The only processes for
which this is not so obvious are the stochastic Ising and Potts
models. Nevertheless, we will show in Section~\ref{S:Isap} below that
also the generator of the stochastic Ising model with Glauber dynamics
can be written in the form (\ref{Gmap}), and we will use this random
mapping representation of $G$ to prove ergodicity for small values
of~$\bet$.

The space $S^\La$ is uncountable except in the trivial case that $S$
has cardinality one. This means that we cannot use the theory of
continuous-time Markov chains. Instead, we will rely on the theory of
Feller processes. In Section~\ref{S:Feller}, we will collect some
general facts about Feller processes, which are a class of Markov
processes with compact, metrizable state spaces, that are uniquely
characterized by their generators. Since this is rather functional
analytic material, which is moreover well-known, we will state the
main facts without proof, but give references to places where proofs
can be found.

In Section~\ref{S:Poiscon}, we then present a Poisson construction of
interacting particle systems (including proofs) that is similar to the
Poisson construction of continuous-time Markov chains. To some degree,
this provides a probabilistic alternative to the functional analytic
approach via Feller processes. To get the full picture, however, one
needs both approaches, so in Section~\ref{S:Gencon}, we link the
Feller process we have constructed in Section~\ref{S:Poiscon} to the
generator defined in (\ref{Gmap}).

Luckily, all this abstract theory gives us more than just the
information that the systems we are interested in are well defined. In
Section~\ref{S:ergod}, we will see that as a side-result of our
proofs, we can derive sufficient conditions for an interacting
particle system to be ergodic, that is, to have a unique invariant law
that is the long-time limit starting from any initial
state.\footnote{We use the word ergodic in a different meaning than in
ergodic theory, see the discussion on page~\pageref{ergdisc}.} In
Section~\ref{S:Isap} we apply this to derive lower bounds on the
critical points of the Ising model. The methods developed in
Section~\ref{S:Poiscon} will also prove useful in Chapter~\ref{C:dual}
where we discuss duality.

\section{Feller processes}\label{S:Feller}

In Section~\ref{S:finMark}, we showed how the semigroup of a Markov
process on a finite state space can be characterized in terms of its
generator and in Section~\ref{S:genchain} we generalized this to
countable state spaces. In the present section, we will treat a class
of Markov processes with compact metrizable state spaces. The basic
assumption we will make is that the transition probabilities
$(P_t)_{t\geq 0}$ are continuous, which means that we will be
discussing \emph{Feller processes}. We will later apply the theory of
Feller processes to state spaces of the form $S^\La$ equipped with the
product topology, which are compact by Tychonoff's theorem. It is easy
to see that the product topology on $S^\La$ is metrizable. For
example, if $(a_i)_{i\in\La}$ are strictly positive constants such
that $\sum_{i\in\La}a_i<\infty$, then
\[
d(x,y):=\sum_{i\in\La}a_i1_{\{x(i)\neq y(i)\}}
\]
defines a metric that generates the product topology.

Let $E$ be a compact metrizable space.\footnote{Such spaces are always
  separable and complete in any metric that generates the topology; in
  particular, they are Polish spaces.} We use the notation
\[\begin{array}{r@{\,}c@{\,}l}\index{0M1E@${\cal M}_1(E)$}
\index{0BE@${\cal B}(E)$}\index{0BE@$B(E)$}\index{0CE@${\cal C}(E)$}
\dis\Bi(E)&:=&\;\mbox{the Borel-\si-field on }E,\\
\dis B(E)&:=&\;\mbox{the space of bounded, Borel-measurable
 functions }f\cn E\to\R,\\
\dis\Ci(E)&:=&\;\mbox{the space of continuous functions }f\cn E\to\R,\\
\dis\Mi_1(E)&:=&\;\mbox{the space of probability measures $\mu$ on $E$.}
\end{array}\]
We equip $\Ci(E)$ with the supremum-norm\index{supremum-norm}
\[\index{00normi@$\Vert\,\cdot\,\Vert_\infty$}
\|f\|_\infty:=\sup_{x\in E}|f(x)|\qquad(f\in\Ci(E)),
\]
making $\Ci(E)$ into a Banach space. We equip $\Mi_1(E)$ with the topology of
weak convergence, \index{weak convergence} where by definition,\footnote{More
  precisely, the topology of weak convergence is the unique \emph{metrizable}
  topology with this property. Since in metrizable spaces, convergent
  subsequences uniquely characterize the topology, our definition is
  unambiguous.} $\mu_n$ converges weakly to $\mu$, denoted
$\mu_n\Rightarrow\mu$, if $\int f\,\di\mu_n\to\int f\,\di\mu$ for all
$f\in\Ci(E)$. With this topology, $\Mi_1(E)$ is a compact metrizable space.
The following lemma is often convenient.

\begin{lemma}[Convergence criterion]
For\label{L:fco} $f_n,f\in\Ci(E)$, the following conditions are equivalent:
\begin{enumerate}[(ii)]
\item $\dis\|f_n-f\|_\infty\asto{n}0$,
\item $\dis f_n(x_n)\asto{n}f(x)$ for all $x_n,x\in E$ with $x_n\to x$.
\end{enumerate}
\end{lemma}

\begin{Proof}
Since
\[
\big|f_n(x_n)-f(x)\big|\leq\|f_n-f\|_\infty+\big|f(x_n)-f(x)\big|,
\]
(i) and the continuity of $f$ imply (ii). On the other hand, if (i)
does not hold, then we can choose $x_n\in E$ such that
$\limsup_{n\to\infty}|f_n(x_n)-f(x_n)|>0$. By the compactness of $E$,
going to a subsequence if necessary, we can assume that $x_n\to x$ for
some $x\in E$. Since
\[
\big|f_n(x_n)-f(x)\big|\geq\big|f_n(x_n)-f(x_n)\big|-\big|f(x_n)-f(x)\big|,
\]
using the continuity of $f$, we see that $f_n(x_n)\not\to f(x)$.
\end{Proof}

A \emph{probability kernel}\index{probability kernel} on $E$ is a function
$K\cn E\times\Bi(E)\to\R$ such that
\begin{enumerate}[(ii)]
\item $K(x,\,\cdot\,)$ is a probability measure on $E$ for each $x\in E$,
\item $K(\,\cdot\,,A)$ is a real measurable function on $E$ for each
$A\in\Bi(E)$.
\end{enumerate}
This is equivalent to the statement that $x\mapsto K(x,\,\cdot\,)$ is
a measurable map from $E$ to $\Mi_1(E)$ (where the latter is equipped
with the topology of weak convergence and the associated
Borel-\si-field). By definition, a probability kernel is
\emph{continuous}\index{probability kernel!continuous} if the map
$x\mapsto K(x,\,\cdot\,)$ is continuous (with respect to the topologies
with which we have equipped these spaces). A probability kernel is
\emph{deterministic}\index{probability kernel!deterministic} if it is
of the form $K(x,\,\cdot\,)=\de_{m(x)}$ for some measurable map $m\cn
E\to E$, where $\de_{m(x)}$ denotes the delta-measure at $m(x)$. It is
easy to see that a deterministic kernel is continuous if and only if
$m$ is a continuous map. A \emph{random mapping representation}
\index{random mapping representation!of a probability kernel} of a
probability kernel $K$ is a random measurable map\footnote{More
formally, this means that $M\cn\Om\times E\to E$ is measurable with
respect to the product-\si-field $\Fi\otimes\Bi(E)$, where
$(\Om,\Fi,\P)$ is the underlying probability space.} $M\cn E\to E$
such that $K(x,\,\cdot\,)=\P[M(x)\in\,\cdot\,]$ $(x\in
E)$.\footnote{For infinite spaces, it is not so clear if every
probability kernel has a random mapping representation. One could also
ask if every continuous probability kernel has a representation in
terms of continuous maps. Although these questions are interesting, we
will neglect them here.}

If $K(x,\di y)$ is a probability kernel on a Polish space $E$, then setting
\[\index{0Ptf@$P_tf$}\index{0Kz@$Kf$}
Kf(x):=\int_EK(x,\di y)f(y)\qquad\big(x\in E,\ f\in B(E)\big)
\]
defines a linear operator $K\cn B(E)\to B(E)$. We define the \emph{composition}
of two probability kernels $K,L$ as
\[
(KL)(x,A):=\int_EK(x,\di y)L(y,A)\qquad\big(x\in E,\ A\in\Bi(E)\big).
\]
Then $KL$ is again a probability kernel on $E$ and the linear operator
$(KL)\cn B(E)\to B(E)$ associated with this kernel is the composition of the
linear operators $K$ and $L$. It follows from the definition of weak
convergence that a kernel $K$ is continuous if and only if its
associated linear operator maps the space $\Ci(E)$ into itself.
If $\mu$ is a probability measure and $K$ is a probability kernel,
then
\[\index{0zzmmuK@$\mu K$}\index{0zzmmuPt@$\mu P_t$}
(\mu K)(A):=\int\mu(\di x)K(x,A)\qquad\big(A\in\Bi(E)\big)
\]
defines another probability measure $\mu K$. Introducing the notation $\mu
f:=\int\!f\,\di\mu$,\index{0zzmmuf@$\mu f$} one has $(\mu K)f=\mu(Kf)$ for all
$f\in B(E)$.

By definition, a \emph{continuous transition probability}
\index{transition probability, continuous} on $E$ is a collection
$(P_t)_{t\geq 0}$ of probability kernels on $E$, such that
\begin{enumerate}[(ii)]
\item $(x,t)\mapsto P_t(x,\,\cdot\,)$ is a continuous map from
  $E\times\half$ into $\Mi_1(E)$,
\item $P_0=1$ and $P_sP_t=P_{s+t}\quad(s,t\geq 0)$.
\end{enumerate}
In particular, (i) implies that each $P_t$ is a continuous probability kernel,
so each $P_t$ maps the space $\Ci(E)$ into itself. One has
\[\begin{array}{rl}
{\rm(i)}&\lim_{t\to 0}P_tf=P_0f=f\qquad(f\in\Ci(E)),\\[4pt]
{\rm(ii)}&P_sP_tf=P_{s+t}f\qquad(s,t\geq 0),\\[4pt]
{\rm(iii)}&f\geq0\mbox{ implies }P_tf\geq 0,\\[4pt]
{\rm(iv)}&P_t1=1,
\end{array}\]
and conversely, each collection of linear operators $P_t\cn\Ci(E)\to\Ci(E)$ with
these properties corresponds to a unique continuous transition probability on
$E$. Such a collection of linear operators $P_t\cn\Ci(E)\to\Ci(E)$ is called a
\emph{Feller semigroup}.\index{Feller semigroup} We note that in (i), the
limit is (of course) with respect to the topology we have chosen on $\Ci(E)$,
that is, with respect to the supremum-norm.

By definition, a function $w\cn\half\to E$ is \emph{cadlag}\index{cadlag}
if it is right-continuous with left limits,\footnote{The word cadlag is an
  abbreviation of the French \textbf{c}ontinue \textbf{\`a} \textbf{d}roit, 
\textbf{l}imite \textbf{\`a} \textbf{g}auche.} that is,
\[\begin{array}{rll}
{\rm(i)}&\dis\lim_{t\down s}w_t=w_s\qquad&(s\geq 0),\\[5pt]
{\rm(ii)}&\dis\lim_{t\up s}w_t=:w_{s-}\mbox{ exists}\qquad&(s>0).
\end{array}\]
Let $(P_t)_{t\geq 0}$ be a Feller semigroup. By definition a
\emph{Feller process} with semigroup $(P_t)_{t\geq 0}$ is a stochastic
process $X=(X_t)_{t\geq 0}$ with cadlag sample paths\footnote{It is
possible to equip the space $\Di_E\half$ of cadlag functions
$w\cn\half\to E$ with a (rather natural) topology, called the
\emph{Skorohod topology}, \index{Skorohod topology} such that
$\Di_E\half$ is a Polish space and the Borel-\si-field on $\Di_E\half$
is generated by the coordinate projections $w\mapsto w_t$ $(t\geq
0)$. As a result, we can view a stochastic process $X=(X_t)_{t\geq 0}$
with cadlag sample paths as a single random variable $X$ taking values
in the space $\Di_E\half$. The law of such a random variable is then
uniquely determined by the finite dimensional distributions of
$(X_t)_{t\geq 0}$.} such that
\begin{equation}\label{Mark2}
\P\big[X_u\in\,\cdot\,\big|\,(X_s)_{0\leq s\leq t}\big]
=P_{u-t}(X_t,\,\cdot\,)\quad{\rm a.s.}
\qquad(0\leq t\leq u).
\end{equation}
Here we condition on the \si-field generated by the
random variables $(X_s)_{0\leq s\leq t}$. Formula (\ref{Mark2}) is equivalent
to the statement that the finite dimensional distributions of $X$ are
given by
\begin{equation}\begin{array}{l}\label{Markfdd2}
\dis\P\big[X_0\in\di x_0,\ldots,X_{t_n}\in\di x_n\big]\\[5pt]
\dis\quad=\P[X_0\in\di x_0]P_{t_1-t_0}(x_0,\di x_1)\cdots
P_{t_n-t_{n-1}}(x_{n-1},\di x_n)
\end{array}\end{equation}
$(0<t_1<\cdots<t_n)$. Formula (\ref{Markfdd2}) is symbolic notation,
which means that
\[\begin{array}{l}
\dis\E\big[f(X_0,\ldots,X_{t_n})\big]\\[5pt]
\dis\ =\!\int\!\P[X_0\in\di x_0]\int\! P_{t_1-t_0}(x_0,\di x_1)\cdots
\int\! P_{t_n-t_{n-1}}(x_{n-1},\di x_n)f(x_0,\ldots,x_n)
\end{array}\]
for all $f\in B(E^{n+1})$. By (\ref{Markfdd2}), the law of a Feller
process $X$ is uniquely determined by its initial law
$\P[X_0\in\,\cdot\,]$ and its transition probabilities $(P_t)_{t\geq 0}$.
Existence is less obvious than uniqueness, but the next theorem
says that this holds in full generality.

\begin{theorem}[Construction of Feller processes]\label{T:Fellexist}
Let $E$ be a compact metrizable space, let $\mu$ be a probability
measure on $E$, and let $(P_t)_{t\geq 0}$ be a Feller semigroup. Then
there exists a Feller process $X=(X_t)_{t\geq 0}$ with initial law
$\P[X_0\in\,\cdot\,]=\mu$, and such a process is unique in
distribution.
\end{theorem}

Just as in the case for finite state space, we would like to characterize a
Feller semigroup by its generator. This is somewhat more complicated than in
the finite setting since in general, it is not possible to make sense of the
exponential formula $P_t=e^{tG}:=\sum_{n=0}^\infty\frac{1}{n!}(tG)^n$. This is
related to the fact that if $G$ is the generator of a Feller semigroup, then
in general it is not possible to define $Gf$ for all $f\in\Ci(E)$, as we now
explain.

Let $\Vi$ be a Banach space. (In our case, the only Banach spaces that we will
need are spaces of the form $\Ci(E)$, equipped with the supremum-norm.) By
definition, a \emph{linear operator} on $\Vi$ is a pair $(A,\Di(A))$ where
$\Di(A)$ is a linear subspace of $\Vi$, called the \emph{domain}
\index{domain of a linear operator} and $A$ is a linear map
$A\cn\Di(A)\to\Vi$. Even though a linear operator is really a pair $(A,\Di(A))$,
one often writes sentences such as ``let $A$ be a linear operator'' without
explicitly mentioning the domain. This is similar to phrases like: ``let $\Vi$
be a Banach space'' (without mentioning the norm) or ``let $M$ be a measurable
space'' (without mentioning the \si-field).

We say that a linear operator $A$ (with domain $\Di(A)$) on a Banach space
$\Vi$ is \emph{closed}\index{closed linear operator} if and only if its
graph\index{graph of an operator} $\{(f,Af):f\in\Di(A)\}$ is a closed subset
of $\Vi\times\Vi$. By definition, a linear operator $A$ (with domain $\Di(A)$)
on a Banach space $\Vi$ is \emph{closable}\index{closable linear operator} if
the closure of its graph (as a subset of $\Vi\times\Vi$) is the graph of a
linear operator $\ov A$ with domain $\Di(\ov A)$. This operator is then called
the \emph{closure} of $A$.\index{closure of linear operator} We mention the
following theorem.

\begin{theorem}[Closed graph theorem]\label{T:closgra}
Let $\Vi$ be a Banach space and let $A$ be a linear operator that is
everywhere defined, that is, $\Di(A)=\Vi$. Then the following statements are
equivalent.
\begin{enumerate}[(iii)]
\item $A$ is continuous as a map from $\Vi$ into itself.
\item $A$ is \emph{bounded},\index{bounded linear operator} that is,
  there exists a constant $C<\infty$ such that $\|Af\|\leq C\|f\|$
  $(f\in\Vi)$.
\item $A$ is closed.
\end{enumerate}
\end{theorem}

Theorem~\ref{T:closgra} shows in particular that if $A$ is an unbounded
operator (that is, there exists $0\neq f_n\in\Di(A)$ such that
$\|Af_n\|/\|f_n\|\to\infty$) and $A$ is closable, then its closure $\ov A$
will not be everywhere defined. Closed (but possibly unbounded) linear
operators are in a sense ``the next good thing'' after bounded operators.

As before, let $E$ be a compact metrizable space and let $(P_t)_{t\geq 0}$ be
a continuous transition probability (or equivalently Feller semigroup) on $E$.
By definition, the \emph{generator}\index{generator!Feller semigroup} of
$(P_t)_{t\geq 0}$ is the linear operator
\[
Gf:=\lim_{t\to 0}t^{-1}\big(P_tf-f),
\]
with domain
\[
\Di(G):=\big\{f\in\Ci(E):\mbox{the limit }\lim_{t\to
  0}t^{-1}\big(P_tf-f)\mbox{ exists}\big\}.
\]
Here, when we say that the limit exists, we mean (of course) with
respect to the topology on $\Ci(E)$, that is, w.r.t.\ the
supremum-norm. The following lemma says that generators are closed,
densely defined operators.

\begin{lemma}[Elementary properties of generators]
Let\label{L:Gclos} $G$ be the generator of a Feller semigroup $(P_t)_{t\geq 0}$.
Then $G$ is closed and $\Di(G)$ is a dense subspace of $\Ci(E)$.
\end{lemma}

Since $G$ is closed, Theorem~\ref{T:closgra} tells us that $G$ is
everywhere defined (that is, $\Di(G)=\Ci(E)$) if and only if $G$ is
bounded. For bounded generators, it is not hard to show that the
exponential formula $e^{tG}f:=\sum_{n=0}^\infty\frac{1}{n!}(tG)^nf$ $(f\in\Ci(E))$
converges in the norm on $\Ci(E)$ and that the Feller semigroup with
generator $G$ is given by $P_t=e^{tG}$. On the other hand, if $G$ is
unbounded, then it is in general not possible to make sense of the
exponential formula.\footnote{In order for
$\sum_{n=0}^\infty\frac{1}{n!}t^nG^nf$ to make sense, we need that
$G^nf$ is well-defined for all $n\geq 0$. For $n=1$ this already
requires that $f\in\Di(G)$ but for higher $n$ we need even more since
it is in general not true that $G$ maps $\Di(G)$ into itself. Thus, it
is not even clear for which class of functions we can make sense of
each term in the expansion separately.} In the context of interacting
particle systems, it is not hard to show that a generator of the form
(\ref{Gmap}) is bounded if $\sum_{m\in\Gi}r_m<\infty$. For the
particle systems we will be interested in, this sum will usually be
infinite and the generator will be unbounded.

Since we cannot use the exponential formula $P_t=e^{tG}$, we need
another way to characterize $(P_t)_{t\geq 0}$ in terms of $G$. Similar
to what we did in Section~\ref{S:genchain}, we will use the backward
equation instead. Let $A$ be a linear operator on $\Ci(E)$. By
definition, we say that a function $\half\ni t\mapsto u_t\in\Ci(E)$
solves the \emph{Cauchy equation} \index{Cauchy equation}
\begin{equation}\label{Cauchy}
\dif{t}u_t=Au_t\qquad(t\geq 0)
\end{equation}
if $u_t\in\Di(A)$ for all $t\geq 0$, the maps $t\mapsto u_t$ and $t\mapsto
Au_t$ are continuous (w.r.t.\ the topology on $\Ci(E)$), the limit
$\dif{t}u_t:=\lim_{s\to 0}s^{-1}(u_{t+s}-u_s)$ exists (w.r.t.\ the topology on
$\Ci(E)$) for all $t\geq 0$, and (\ref{Cauchy}) holds. The following
proposition shows that a Feller semigroup is uniquely characterized by its
generator.

\begin{proposition}[Cauchy problem]\label{P:Cauchy}
Let $G$ be the generator of a Feller semigroup $(P_t)_{t\geq 0}$. Then, for
each $f\in\Di(G)$, the Cauchy equation $\dif{t}u_t=Gu_t$ $(t\geq 0)$ has a
unique solution $(u_t)_{t\geq 0}$ with initial state $u_0=f$. Denoting this
solution by $U_tf:=u_t$ defines for each $t\geq 0$ a linear operator $U_t$
with domain $\Di(G)$, of which $P_t=\ov U_t$ is the closure.
\end{proposition}

We need a way to check that (the closure of) a given operator is the generator
of a Feller semigroup. For a given linear operator $A$, constant $\la>0$, and
$f\in\Ci(E)$, we say that a function $p\in\Ci(E)$ solves the
\emph{Laplace equation}\index{Laplace equation}
\begin{equation}\label{Laplace}
(\la-A)p=f
\end{equation}
if $p\in\Di(A)$ and (\ref{Laplace}) holds. The following lemma shows how
solutions to Laplace equations typically arise. 

\begin{lemma}[Laplace equation]\label{L:Laplace}
Let $G$ be the generator of a Feller semigroup $(P_t)_{t\geq 0}$
on $\Ci(E)$, let $\la>0$ and $f\in\Ci(E)$. Then the Laplace equation
$(\la-G)p=f$ has a unique solution, that is given by
\[
p=\int_0^\infty P_tf\,e^{-\la t}\di t.
\]
\end{lemma}

We say that an operator $A$ on $\Ci(E)$ with domain $\Di(A)$ satisfies the
\emph{positive maximum principle}\index{positive maximum principle}
\index{maximum principle} if, whenever a function
$f\in\Di(A)$ assumes its maximum over $E$ in a point $x\in E$ and
$f(x)\geq 0$, we have $Af(x)\leq 0$. The following proposition
gives necessary and sufficient conditions for a linear operator $G$ to be the
generator of a Feller semigroup.


\begin{theorem}[Generators of Feller semigroups]\label{T:HY1}
A linear operator $G$ on $\Ci(E)$ is the generator of a Feller semigroup
$(P_t)_{t\geq 0}$ if and only if
\begin{enumerate}[(iii)]
\item $1\in\Di(G)$ and $G1=0$.
\item $G$ satisfies the positive maximum principle.
\item $\Di(G)$ is dense in $\Ci(E)$.
\item For every $f\in\Ci(E)$ and $\la>0$, the Laplace equation
  $(\la-G)p=f$ has a solution.
\end{enumerate}
\end{theorem}

In practice, it is rarely possible to give an explicit description of
the (full) domain of a Feller generator. Rather, one often starts with
an operator that is defined on a smaller domain of ``nice'' functions
and then takes its closure. In general, if $G$ is a closed linear
operator and $\Di'\sub\Di(G)$ is a linear subspace of $\Di(G)$, then
we let $G|_{\Di'}$ denote the restriction of $G$ to $\Di'$, that is,
$G|_{\Di'}$ is the linear operator with domain $\Di(G|_{\Di'}):=\Di'$
defined as $G|_{\Di'}f:=Gf$ for all $f\in\Di'$. We say that $\Di'$ is
a \emph{core} \index{core of closed operator} for $G$ if
$\ov{G|_{\Di'}}=G$.

\begin{lemma}[Core of a generator]
Let\label{L:core} $G$ be the generator of a Feller semigroup and let
$\Di'$ be a linear subspace of $\Di(G)$. Assume that $\Di'$ is dense
in $\Ci(E)$. Then the following conditions are equivalent:
\begin{enumerate}[(iii)]
\item $\Di'$ is a core for $G$,
\item the set $\{(\la-A)p:p\in\Di'\}$ is dense in $\Ci(E)$ for some $\la>0$,
\item the set $\{(\la-A)p:p\in\Di'\}$ is dense in $\Ci(E)$ for all $\la>0$.
\end{enumerate}
\end{lemma}

Note that by condition~(ii) of Lemma~\ref{L:core}, to check that a
dense set $\Di'\sub\Ci(E)$ is a core for $G$, it suffices to show that
for some $\la>0$, there exists a dense subspace $\Ri\sub\Ci(E)$ such
that for every $f\in\Ri$, the Laplace equation $(\la-A)p=f$ has a
solution $p\in\Di'$. Using Lemma~\ref{L:core}, one can prove the
following version of the Hille--Yosida theorem.\index{Hille--Yosida}

\begin{theorem}[Hille--Yosida]\label{T:HY2}
A linear operator $A$ on $\Ci(E)$ with domain $\Di(A)$ is closable and its
closure $G:=\ov A$ is the generator of a Feller semigroup if
and only if
\begin{enumerate}[(iii)]
\item There exist $f_n\in\Di(A)$ such that $f_n\to 1$ and $Af_n\to 0$.
\item $A$ satisfies the positive maximum principle.
\item $\Di(A)$ is dense in $\Ci(E)$.
\item For some (and hence for all) $\la\in(0,\infty)$, there exists
  a dense subspace $\Ri\sub\Ci(E)$ such that for every $f\in\Ri$, the Laplace
  equation $(\la-A)p=f$ has a solution $p$.
\end{enumerate}
\end{theorem}

Conditions~(i)--(iii) are usually easy to verify for a given operator $A$, but
condition~(iv) is the ``hard'' condition since this means that one has to
prove existence of solutions to the Laplace equation $(\la-G)p=f$ for a dense
set of functions $f$.

If $K$ is a probability kernel on $E$ and $r>0$, then
\begin{equation}\label{bdedG}
Gf:=r(Kf-f)\qquad\big(f\in\Ci(E)\big)
\end{equation}
defines a Feller generator that is everywhere defined (that is, $\Di(G)=\Ci(E)$)
and hence, in view of Theorem~\ref{T:closgra}, a bounded operator. For
generators of this simple form, one can construct the corresponding semigroup
by the exponential formula
\[
P_tf=\ex{tG}f:=\sum_{n=0}^\infty\frac{1}{n!}(tG)^nf,
\]
where the infinite sum converges in $\Ci(E)$. The corresponding Markov process
has a simple description: with rate $r$, the process jumps from its current
position $x$ to a new position chosen according to the probability law
$K(x\,\,\cdot\,)$.

As soon as Feller processes get more complicated in the sense that ``the total
rate of all things that can happen'' is infinite (as will be the case for
interacting particle systems), one needs the more complicated Hille--Yosida
theory. To demonstrate the strength of Theorem~\ref{T:HY2}, consider
$E:=[0,1]$ and the linear operator $A$ defined by $\Di(A):=\Ci^2[0,1]$ (the
space of twice continuously differentiable functions on $[0,1]$) and
\begin{equation}\label{WF2}
Af(x):=x(1-x)\diff{x}f(x)\qquad\big(x\in[0,1]\big).
\end{equation}
One can show that $A$ satisfies the conditions of Theorem~\ref{T:HY2} and
hence $\ov A$ generates a Feller semigroup. The corresponding Markov process
turns out to have continuous sample paths and is indeed the
\emph{Wright--Fisher diffusion}\index{Wright--Fisher diffusion} that we met
before in formula (\ref{WF}).

\begin{Exercise}[Brownian motion]
Let $(P_t)_{t\geq 0}$ denote the transition kernels of Brownian motion
on $\R^d$. Let $E:=\R^d\cup\{\infty\}$ denote the one-point
compactification of $\R^d$ and extend $P_t$ $(t\geq 0)$ to probability
kernels on $E$ by setting $P_t(\infty,\,\cdot\,):=\de_\infty$. Show
that $(P_t)_{t\geq 0}$ is a Feller semigroup.
\end{Exercise}

In Chapter~\ref{C:Markov}, we viewed (possibly explosive)
continuous-time Markov chains with a countable state space $S$ as
Markov processes on the extended state space $S_\infty$, where
$S_\infty$ is the one-point compactification of $S$. It is natural to
ask if they are in fact Feller processes on $S_\infty$. The answer is,
in general, negative. The reason is that the extended transition
kernels $(\ov P_t)_{t\geq 0}$ on $S_\infty$ may fail to be continuous
at $\infty$, that is, $\ov P_t(x_n,\,\cdot\,)$ may fail to converge to
$\ov P_t(\infty,\,\cdot\,)$ if $x_n\to\infty$. In many cases where
this problem occurs, it can be solved by choosing another
compactification of $S$ (that is, by adding more points at
infinity). Whether this can be done in general I don't know.

\begin{Exercise}[Wright--Fisher diffusion]
Show\label{E:WF} that the operator $A$ defined in (\ref{WF2}) satisfies the
conditions of Theorem~\ref{T:HY2}. \emph{Hint:} show that if $f$ is a
polynomial of order $n$, then so is $Af$. Use this to show that the
Cauchy equation $\dif{t}u_t=Au_t$ has a solution for each initial
state $u_0=f$ that is a polynomial. Then show that $p:=\int_0^\infty
u_t\,e^{-\la t}\,\di t$ solves the Laplace equation $(\la-A)p=f$.
\end{Exercise}

\subsection*{Some notes on the proofs}

In the remainder of this section, we indicate where proofs of the stated
theorems can be found. Readers who are more interested in interacting particle
systems than in functional analysis may skip from here to the next section.

The fact that there is a one-to-one correspondence between continuous
transition probabilities and collections $(P_t)_{t\geq 0}$ of linear operators
satisfying the assumptions (i)--(iv) of a Feller semigroup follows from
\cite[Prop.~17.14]{Kal97}.

Theorem~\ref{T:Fellexist} (including a proof) can be found in
\cite[Thm~17.15]{Kal97} and \cite[Thm~4.2.7]{EK86}.
Theorem~\ref{T:closgra} (the closed graph theorem and characterization of
continuous linear maps) can be found in many places (including Wikipedia).

Lemma~\ref{L:Gclos} follows from \cite[Corollary~I.1.6]{EK86}. The
statements of this lemma can also easily be derived from the
Hille--Yosida theorem (see below). Proposition~\ref{P:Cauchy}
summarizes a number of well-known facts.  The fact that $u_t:=P_tf$
solves the Cauchy equation if $f\in\Di(G)$ is proved in
\cite[Prop~1.1.5~(b)]{EK86}, \cite[Thm~17.6]{Kal97}, and
\cite[Thm~3.16~(b)]{Lig10}. To see that solutions to the Cauchy
equation are unique, we use the following fact.

\begin{lemma}[Positive maximum principle]\label{L:maxprinc}
Let $A$ be a linear operator on $\Ci(E)$ and let $u=(u_t)_{t\geq 0}$ be a
solution to the Cauchy equation $\dif{t}u_t=Au_t$ $(t\geq 0)$.
Assume that $A$ satisfies the positive maximum principle and $u_0\geq 0$.
Then $u_t\geq 0$ for all $t\geq 0$.\index{positive maximum principle}
\end{lemma}

\begin{Proof}
By linearity, we may equivalently show that $u_0\leq 0$ implies $u_t\leq 0$.
Assume that $u_t(x)>0$ for some $x\in E$. By the compactness of $E$, the
function $(x,t)\mapsto e^{-t}u_t(x)$ must assume its maximum over
$E\times[0,t]$ in some point $(y,s)$. Since $u$ is positive somewhere
on $E\times[0,t]$ we have $e^{-s}u_s(y)>0$ and hence $s>0$ by the fact that
$u_0\leq 0$. But now, since $A$ satisfies the positive maximum principle,
\[
0\leq\dif{s}\big(e^{-s}u_s(y)\big)=-e^{-s}u_s(y)+e^{-s}Au_s(y)\leq-e^{-s}u_s(y)<0,
\]
so we arrive at a contradiction.
\end{Proof}

By linearity, Lemma~\ref{L:maxprinc} implies that if $u,v$ are two solutions
to the same Cauchy equation and $u_0\leq v_0$, then $u_t\leq v_t$ for all
$t\geq 0$. In particular, since by Theorem~\ref{T:HY1}, Feller generators
satisfy the positive maximum principle, this implies uniqueness of solutions
of the Cauchy equation in Proposition~\ref{P:Cauchy}. Again by
Theorem~\ref{T:HY1}, the domain of a Feller semigroup is a dense subspace of
of $\Ci(E)$, so the final statement of Proposition~\ref{P:Cauchy} follows from
the following simple lemma and the fact that $\|P_tf\|_\infty\leq\|f\|_\infty$.

\begin{lemma}[Closure of bounded operators]
Let $(\Vi,\|\,\cdot\,\|)$ be a Banach space and let $A$ be a linear operator on
$\Vi$ such that $\Di(A)$ is dense and $\|Af\|\leq C\|f\|$ $(f\in\Di(A))$ for
some $C<\infty$. Then $A$ is closable, $\Di(\ov A)=\Vi$, and $\|\ov Af\|\leq
C\|f\|$ $(f\in\Vi)$.
\end{lemma}

\begin{Proof}[(sketch)]
Since $\Di(A)$ is dense, for each $f\in\Vi$ we can choose $\Di(A)\ni f_n\to
f$.  Using the fact that $A$ is bounded, it is easy to check that if
$(f_n)_{n\geq 0}$ is a Cauchy sequence and $f_n\in\Di(A)$ for all $n$, then
$(Af_n)_{n\geq 0}$ is also a Cauchy sequence. By the completeness of $\Vi$, it
follows that the limit $\ov Af:=\lim_{n\to\infty}Af_n$ exists for all
$f\in\Vi$. To see that this defines $\ov A$ unambiguously, assume that $f_n\to
f$ and $g_n\to f$ and observe that $\|Af_n-Ag_n\|\leq C\|f_n-g_n\|\to 0$.  The
fact that $\|\ov Af\|\leq C\|f\|$ $(f\in\Vi)$ follows from the continuity of
the norm.
\end{Proof}

Lemma~\ref{L:Laplace} follows from \cite[Prop~1.2.1]{EK86}.
Theorems~\ref{T:HY1} and \ref{T:HY2} both go under the name of the
Hille--Yosida theorem. Often, they are stated in a more general form without
condition~(i). In this generality, the operator $G$ generates a semigroup of
\emph{subprobability kernels}\index{subprobability kernel} $(P_t)_{t\geq 0}$,
that is, $P_t(x,\,\cdot\,)$ is a measure with total mass $P_t(x,E)\leq 1$.  In
this context, a Feller semigroup with $P_t(x,E)=1$ for all $t,x$ is called
\emph{conservative}. \index{conservative Feller semigroup} It is clear from
Proposition~\ref{P:Cauchy} that condition~(i) in Theorems~\ref{T:HY1} and
\ref{T:HY2} is necessary and sufficient for the Feller group to be
conservative.

The versions of the Hille--Yosida theorem stated in \cite{EK86,Kal97} are more
general than Theorems~\ref{T:HY1} and \ref{T:HY2} since they allow for the
case that $E$ is not compact but only locally compact. This is not really more
general, however, since what these books basically do if $E$ is not compact is
the following. First, they construct the one-point compactification $\ov
E=E\cup\{\infty\}$ of $E$. Next, they extend the transition probabilities to
$\ov E$ by putting $P_t(\infty,\,\cdot\,):=\de_\infty$ for all $t\geq 0$.
Having proved that they generate a conservative Feller semigroup on $\ov E$ of
this form, they then still need to prove that the associated Markov process
does not explode in the sense that $\P^x[X_t\in E\ \forall t\geq 0]=1$ for all
$x\in E$. In practical situations (such as when constructing Markov processes
with state space $\R^d$) it is usually better to explicitly work with the
one-point compactification of $\R^d$ instead of trying to formulate
theorems for locally compact spaces that try to hide this
compactification in the background.

Theorems~\ref{T:HY1} and \ref{T:HY2} are special cases of more general
theorems (also called Hille--Yosida theorem) for strongly continuous
contraction semigroups taking values in a general Banach space. In this
context, the positive maximum principle is replaced by the assumption that the
operator under consideration is \emph{dissipative}. In this more general
setting, Theorems~\ref{T:HY1} and \ref{T:HY2} correspond to \cite[Thms~1.2.6
  and 1.2.12]{EK86}. Lemma~\ref{L:core} follows from \cite[Lemma~1.2.11 and
  Prop~1.3.1]{EK86}. In the more specific set-up of Feller semigroups,
versions of Theorem~\ref{T:HY2} can be found in \cite[Thm~4.2.2]{EK86} and
\cite[Thm~17.11]{Kal97}. There is also an account of Hille--Yosida theory for
Feller semigroups in \cite[Chap~3]{Lig10}, but this reference does not mention
the positive maximum principle (using a dissipativity assumption instead).

Feller semigroups with bounded generators such as in (\ref{bdedG}) are treated
in \cite[Sect~4.2]{EK86} and \cite[Prop~17.2]{Kal97}.

\section{Poisson construction}\label{S:Poiscon}

We briefly recall the set-up introduced in
Section~\ref{S:conintro}. $S$ is a finite set, called the \emph{local
state space}, \index{local state space} and $\La$ is a countable set,
called the \emph{lattice}.\index{lattice} We equip the product space
$S^\La$ with the product topology, making it into a compact metrizable
space. Elements of $S^\La$ are denoted $x=(x(i))_{i\in\La}$. We fix a
countable set $\Gi$ whose elements are continuous maps $m\cn S^\La\to
S^\La$ as well as nonnegative rates $(r_m)_{m\in\Gi}$. Our aim is to
construct a Markov process with formal generator of the form
(\ref{Gmap}), that is,
\[
Gf(x)=\sum_{m\in\Gi}r_m\big\{f\big(m(x)\big)-f\big(x\big)\big\}
\qquad(x\in S^\La).
\]
Our approach is very similar to the Poisson construction of
continuous-time Markov chains described in
Section~\ref{S:Poichain}. We equip the space $\Gi\times\R$ with the
measure
\begin{equation}\label{rhoips}
\rho\big(\{m\}\times[s,t]\big):=r_m(t-s)\qquad\big(m\in\Gi,\ s\leq t\big).
\end{equation}
Let $\om$ be a Poisson point set with intensity $\rho$. We call $\om$
the \emph{graphical representation} associated with the random mapping
representation (\ref{Grep}). Since $\Gi$ is countable, by the argument
used in Section~\ref{S:Poichain}, the time coordinates of points
$(m,t)\in\om$ are all different. Therefore, as we did in the case of
continuous-time Markov chains, we can unambiguously define a random
function $\R\ni t\mapsto\mk^\om_t\in\Gi$ by setting
\begin{equation}\index{0mt@$\mk^\om_t$}\label{mkom}
\mk^\om_t:=\left\{\begin{array}{ll}
m\quad&\mbox{if }(m,t)\in\om,\\[5pt]
1\quad&\mbox{otherwise,}
\end{array}\right.\end{equation}
where we write $1$ to denote the identity map. By definition, we say
that a random function $X\cn[s,\infty)\to S^\La$ solves the evolution
  equation
\begin{equation}\label{Xevolve}
X_t=\mk^\om_t(X_{t-})\qquad(t>s),
\end{equation}
if $[s,\infty)\ni t\mapsto X_t\in S^\La$ is cadlag and (\ref{Xevolve})
holds. We recall that for any continuous map $m$ and site $i\in\La$,
the set $\Ri^\down_i(m)$ has been defined in Section~\ref{S:finsys}.
Here is the main result of the present section.

\begin{theorem}[Poisson construction]
Assume\label{T:Poispart} that the rates $(r_m)_{m\in\Gi}$ satisfy
\begin{equation}\label{downsum}
{\rm(i)}\ \sup_{i\in\La}\sum_{m\in\Gi}r_m1_{\Di(m)}(i)<\infty,\quad
{\rm(ii)}\ \sup_{i\in\La}\sum_{m\in\Gi}r_m\big|\Ri^\down_i(m)\beh\{i\}\big|<\infty.
\end{equation}
Then almost surely, for each $s\in\R$ and $x\in S^\La$, there exists a
unique solution $(X^{s,x}_t)_{t\geq s}$ to the evolution equation
(\ref{Xevolve}) with initial state $X^{s,x}_s=x$. Setting
\begin{equation}\label{partflow}\index{0Xb@$\Xb_{s,t}$}
\Xb_{s,t}(x):=X^{s,x}_t\qquad(s\leq t,\ x\in S^\La)
\end{equation}
defines a collection of continuous maps $(\Xb_{s,t})_{s\leq t}$ from
$S^\La$ into itself such that
\begin{equation}\label{flow}
\Xb_{s,s}=1
\quand\Xb_{t,u}\circ\Xb_{s,t}=\Xb_{s,u}\qquad(s\leq t\leq u).
\end{equation}
Setting
\begin{equation}\label{Ppart}
P_t(x,\,\cdot\,):=\P\big[\Xb_{0,t}(x)\in\,\cdot\,\big]\qquad(t\geq 0,\ x\in S^\La)
\end{equation}
defines the semigroup of a Feller process with state space $S^\La$. If
$s\in\R$ and $X_0$ is an $S^\La$-valued random variable with law
$\mu$, independent of $\om$, then the process $(X_t)_{t\geq 0}$
defined as
\begin{equation}\label{flowMark}
X_t:=\Xb_{s,s+t}(X_0)\qquad(t\geq 0)
\end{equation}
is distributed as the Feller process with semigroup $(P_t)_{t\geq 0}$
and initial law~$\mu$.
\end{theorem}

One may notice the similarity between condition (\ref{downsum}) and
the condition (\ref{upsum}) from Chapter~\ref{C:Markov}, the only
difference being that $\Ri^\up_i(m)$ in condition (\ref{upsum})~(ii)
is replaced by $\Ri^\down_i(m)$, which looks ``downwards'' in
time. (Here we use our usual convention of plotting time upwards in
pictures of graphical representations so that downwards means back in
time.) In Proposition~\ref{P:finpert} of Section~\ref{S:finax} below,
we will see that the ``upward'' condition (\ref{upsum})~(ii) in
general guarantees that finite perturbations of the initial state have
finite consequences at later times. For later use, we introduce three
constants whose finiteness is guaranteed by conditions
(\ref{downsum})~(ii) and (\ref{upsum})~(ii), and by part~(i) of either
of these equations.
\begin{equation}\begin{array}{c}\label{KKK}\index{0K0@$K_0$}
\index{0Kup@$K_\up$}\index{0Kdw@$K_\down$}
\dis K_\down:=\sup_{i\in\La}\sum_{m\in\Gi}r_m\big(|\Ri^\down_i(m)|-1\big),
\qquad
K_\up:=\sup_{i\in\La}\sum_{m\in\Gi}r_m\big(|\Ri^\up_i(m)|-1\big),\\[5pt]
\dis K_0:=\sup_{i\in\La}\sum_{m\in\Gi}r_m1_{\Di(m)}(i).
\end{array}\end{equation}

The proof of Theorem~\ref{T:Poispart} will take up the rest of this
section. At first, it may be surprising that solutions of the
evolution equation (\ref{Xevolve}) with a given initial state are
unique. After all, if we replace the compact set $S^\La$ by $[0,1]$,
then there are many cadlag functions $X\cn[s,\infty)\to[0,1]$ with a
given initial state that make no jumps at all. The following exercise
shows that at least in the case when $\om=\emptyset$, the equation
(\ref{Xevolve}) has a unique solution.

\begin{Exercise}[Total disconnectedness]
A\label{E:discon} topological space $E$ is \emph{totally disconnected}
if for each $x_1,x_2\in E$ with $x_1\neq x_2$, there exist open sets
$O_1\ni x_1$ and $O_2\ni x_2$ such that $O_1\cap O_2=\emptyset$ and
$O_1\cup O_2=E$. Prove that $S^\La$ is totally disconnected. Prove
that if $E$ is a totally disconnected space, then each continuous
function $f\cn\half\to E$ is constant.
\end{Exercise}

The difficulty with proving that for each $s\in\R$ and $x\in S^\La$,
the equation (\ref{Xevolve}) has a unique solution $X^{s,x}$ is that
typically condition (\ref{locfin}) will be violated. As a result,
$\{t>s:(m,t)\in\om,\ m(x)\neq x\}$ is a dense subset of $[s,\infty)$
and solutions to (\ref{Xevolve}) will not be piecewise
constant. However, most of the jumps of $(X^{s,x}_t)_{t\geq s}$ will
involve sites that are far away, and $t\mapsto X^{s,x}_t(i)$ will
still be piecewise constant for each fixed $i\in\La$.

The trick for proving uniqueness of solutions of (\ref{Xevolve}) is to
look backwards in time. We recall from Lemma~\ref{L:contprod} that if
$T$ is a finite set and $f\cn S^\La\to T$ is a continuous function,
then $f$ depends on finitely many coordinates. A consequence of this
is that the space $\Ci(S^\La,T)$ of continuous functions $\phi\cn
S^\La\to T$ is countable. It turns out that if we fix
$\phi\in\Ci(S^\La,T)$ and a time $u\in\R$, then the stochastic flow
$(\Xb_{s,t})_{s\leq t}$ that we are about to construct has the
property that the process
\[
\Phi_t:=\phi\circ\Xb_{u-t,u}\qquad(t\geq 0)
\]
is a nonexplosive continuous-time Markov chain with countable state
space $\Ci(S^\La,T)$. We call $(\Phi_t)_{t\geq 0}$ the \emph{backward
in time process.}\index{backward!in time process} The first step
towards proving Theorem~\ref{T:Poispart} is showing that this
continuous-time Markov chain is nonexplosive.

\begin{proposition}[Backward in time process]
Assume\label{P:back} (\ref{downsum}). Let $T$ be a finite set. Then setting
\begin{equation}\label{Hgen}
Hf(\phi):=\sum_{m\in\Gi}r_m\big\{f(\phi\circ m)-f(\phi)\big\}
\end{equation}
for all bounded $f\cn\Ci(S^\La,T)\to\half$ defines the generator of a
nonexplosive continuous-time Markov chain $(\Phi_t)_{t\geq 0}$ with
state space $\Ci(S^\La,T)$. This Markov chain satisfies
\begin{equation}\label{Phiexp}
\E^\phi\big[\big|\Ri(\Phi_t)\big|\big]\leq|\Ri(\phi)|\ex{K_\down t}
\qquad\big(t\geq 0,\ \phi\in\Ci(S^\La,T)\big),
\end{equation}
where $K_\down$ is defined in (\ref{KKK}).
\end{proposition}

\begin{figure}[htb]
\begin{center}
\inputtikz{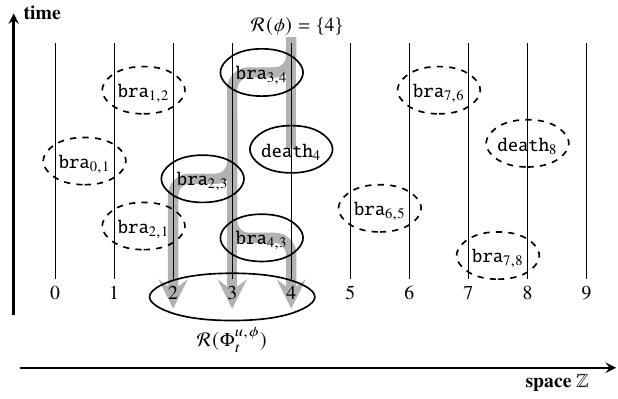}
\caption{Graphical representation of a one-dimensional contact
  process, with the backwards in time process
  $(\Phi^{u,\phi}_t)_{t\in(-\infty,u]}$. In this example
  $T=S=\{0,1\}$, $\phi\in\Ci(S^\La,T)$ is defined by $\phi(x):=x(4)$
  $(x\in S^\La)$, and one has
$\Phi^{u,\phi}_t=\phi\circ{\tt bra}_{3,4}\circ{\tt death}_4\circ{\tt bra}_{2,3}\circ{\tt bra}_{4,3}$.
  The gray arrows pointing downwards show the evolution backwards in
  time of the set $\Ri(\Phi^{u,\phi}_t)$ of relevant lattice points
  for the map $\Phi^{u,\phi}_t$.}
\label{fig:contgraph}
\commentAlt{Figure~\ref{fig:contgraph}}{Picture of the graphical
  representation omega. Time is plotted upwards. Branching maps and
  death maps are written on the places in space and time where they
  are applied. See long description.}
\commentLongAlt{Figure~\ref{fig:contgraph}}{Gray arrows indicate how
  the set of relevant points for the backward in time process
  evolves. Initially, on top of the picture, the only relevant point
  is the point 4. Branching maps lead to the creation of new relevant
  points while a death map causes one point to be removed from the set
  of relevant points. Finally, at the bottom of the picture one needs
  to know the local state in the points 2, 3 and 4 in order to
  determine the value of the process in the point 4 on the top of the
  picture.}
\end{center}
\end{figure}

Before we prove Proposition~\ref{P:back}, we first discuss its
consequences. Combining Proposition~\ref{P:back} with
Theorem~\ref{T:backflow}, one obtains that almost surely, for each
$u\in\R$ and $\phi\in\Ci(S^\La,T)$, there exists a unique cadlag
function $\Phi^{u,\phi}\cn(-\infty,u]\to\Ci(S^\La,T)$ such that
$\Phi^{u,\phi}_u=\phi$ and
\begin{equation}\label{backevol}
\Phi^{u,\phi}_{t-}=\left\{\begin{array}{ll}
\dis\Phi^{u,\phi}_t\circ m\quad&\dis\mbox{if }(m,t)\in\om,\\[5pt]
\dis\Phi^{u,\phi}_t\quad&\dis\mbox{otherwise}
\end{array}\right.\qquad(t\leq u).
\end{equation}
Setting
\begin{equation}\label{Fback}\index{0Fb@$\Fb_{t,s}$}
\Fb_{u,t}(\phi):=\Phi^{u,\phi}_t\qquad\big(t\leq u,\ \phi\in\Ci(S^\La,T)\big)
\end{equation}
defines a collection of maps $(\Fb_{u,t})_{u\geq t}$ from
$\Ci(S^\La,T)$ into itself such that
\[
\Fb_{u,u}=1
\quand\Fb_{t,s}\circ\Fb_{u,t}=\Fb_{u,s}\qquad(u\geq t\geq s).
\]
If $u\in\R$ and $\Phi_0$ is an $\Ci(S^\La,T)$-valued random variable
with law $\mu$, independent of $\om$, then the process
$(\Phi_t)_{t\geq 0}$ defined as
\[
\Phi_t:=\Fb_{u,u-t}(\Phi_0)\qquad(t\geq 0)
\]
is distributed as the left-continuous modification of the
continuous-time Mar\-kov chain with generator $H$ from (\ref{Hgen})
and initial law~$\mu$. These facts are illustrated in
Figure~\ref{fig:contgraph}.\med

The proof of Proposition~\ref{P:back} needs one preparatory lemma.

\begin{lemma}[Well-defined generator]
Assume\label{L:Hwell} that the rates $(r_m)_{m\in\Gi}$ satisfy
\begin{equation}\label{Hwell}
\sum_{m\in\Gi}1_{\Di(m)}(i)r_m<\infty\quad\forall i\in\La.
\end{equation}
Then formula (\ref{Hgen}) defines the generator of a (possibly
explosive) continu\-ous-time Markov chain $(\Phi_t)_{t\geq 0}$ with
state space $\Ci(S^\La,T)$.
\end{lemma}

\begin{Proof}
We must check condition (\ref{locfin}), which in the present context
says that
\[
c(\phi):=\sum_{m:\,\phi\circ m\neq\phi}r_m<\infty\quad\mbox{for all }\phi\in\Ci(S^\La,T).
\]
We observe that $\phi(m(x))\neq\phi(x)$ for some $x\in S^\La$ implies
that $\Di(m)\cap\Ri(\phi)\neq\emptyset$, so we can estimate
\begin{equation}\label{Hcalc}
\sum_{m:\,\phi\circ m\neq\phi}r_m\leq\sum_{m:\,\Di(m)\cap\Ri(\phi)\neq\emptyset}r_m
\leq\sum_{i\in\Ri(\phi)}\sum_m1_{\Di(m)}(i)r_m,
\end{equation}
which is finite by (\ref{Hwell}) and the finiteness of $\Ri(\phi)$.
\end{Proof}

\begin{Proof}[of Proposition~\ref{P:back}]
Condition (\ref{downsum})~(i) clearly implies (\ref{Hwell}) so by
Lemma~\ref{L:Hwell} $H$ is the generator of a (possibly explosive)
continuous-time Markov chain with state space $\Ci(S^\La,T)$. To prove
that $H$ is nonexplosive we apply Theorem~\ref{T:Lyap} to the Lyapunov
function
\[
L(\phi):=|\Ri(\phi)|\qquad\big(\phi\in\Ci(S^\La,T)\big).
\]
Formula (\ref{Hcalc}) shows that
\[
\sup\big\{c(\phi):\phi\in\Ci(S^\La,T),\ L(\phi)<C\big\}
\leq C\sup_{i\in\La}\sum_m1_{\Di(m)}(i)r_m,
\]
which by (\ref{downsum})~(i) implies that $L$ satisfies condition~(i)
of Theorem~\ref{T:Lyap}. It remains to check condition~(ii). We will
show that $HL\leq K_\down L$ where $K_\down$ is the constant defined
in (\ref{KKK}), which is finite by (\ref{downsum})~(ii). We observe
that
\[
HL(\phi)
=\sum_{m\in\Gi}r_m\big\{L(\phi\circ m)-L(\phi)\big\}
=\sum_{m\in\Gi}r_m\big\{|\Ri(\phi\circ m)|-|\Ri(\phi)|\big\}.
\]
Since
\[
\Ri(\phi\circ m)\sub\bigcup_{i\in\Ri(\phi)}\Ri^\down_i(m),
\]
we can estimate
\[
|\Ri(\phi\circ m)|-|\Ri(\phi)|\leq\sum_{i\in\Ri(\phi)}\big(|\Ri^\down_i(m)|-1\big).
\]
It follows that
\[
HL(\phi)
\leq\sum_{m\in\Gi}\sum_{i\in\Ri(\phi)}r_m\big(|\Ri^\down_i(m)|-1\big)
\leq K_\down\big|\Ri(\phi)\big|=K_\down L(\phi),
\]
so condition~(ii) of Theorem~\ref{T:Lyap} is satisfied with
$\la=K_\down$. It follows that $H$ is nonexplosive and (\ref{Phiexp})
holds.
\end{Proof}

Proposition~\ref{P:back} is the cornerstone of the proof of
Theorem~\ref{T:Poispart}. We will prove the following result, that is
sometimes applicable even when condition (\ref{downsum}) of
Theorem~\ref{T:Poispart} is not satisfied. In
Exercise~\ref{E:backexpl} below we will see that if the
continuous-time Markov chain with generator $H$ from (\ref{Hgen}) is
explosive, then solutions to the the evolution equation
(\ref{Xevolve}) may fail to be unique.

\begin{theorem}[Graphical construction]
Assume\label{T:graph} that the rates $(r_m)_{m\in\Gi}$ satisfy
(\ref{Hwell}) and that the continuous-time Markov chain with generator
$H$ from (\ref{Hgen}) is nonexplosive. Then the conclusions of
Theorem~\ref{T:Poispart} remain true.
\end{theorem}

It turns out that the condition (\ref{downsum}) is more or less
optimal in a translation invariant setting, but not necessarily for
inhomogeneous systems, as we now explain. For each bijection
$\psi\cn\La\to\La$ and map $m\cn S^\La\to S^\La$, we define a
translated map $T_\psi m\cn S^\La\to S^\La$ by
\[
T_\psi m(x):=m\big(x\circ\psi^{-1}\big)\qquad(x\in S^\La),
\]
where $x\circ\psi^{-1}$ denotes the concatenation of the functions
$\psi^{-1}\cn\La\to\La$ and $x\cn\La\to S$. Let $\Gi$ be a collection
of continuous maps and let $(r_m)_{m\in\Gi}$ be nonnegative rates. By
definition, an \emph{automorphism} of $(r_m)_{m\in\Gi}$ is a bijection
$\psi\cn\La\to\La$ such that $T_\psi m\in\Gi$ for all $m\in\Gi$ and
\[
r_{T_\psi m}=r_m\qquad(m\in\Gi).
\]
We say that the rates $(r_m)_{m\in\Gi}$ are \emph{transitive} if for
each $i,j\in\La$, there exists an automorphism $\psi$ of
$(r_m)_{m\in\Gi}$ such that $\psi(i)=\psi(j)$. In such a situation, we
also say that the associated random mapping representation of a
generator as in (\ref{Gmap}) is
\emph{transitive}.\index{transitivity!of random mapping representations}
For transitive random mapping representations, the expressions
\[
{\rm(i)}\ \sum_{m\in\Gi}r_m1_{\Di(m)}(i)\quand
{\rm(ii)}\ \sum_{m\in\Gi}r_m\big|\Ri^\down_i(m)\beh\{i\}\big|
\]
do not depend on $i\in\La$ and the suprema in (\ref{downsum}) (i) and
(ii) can be dropped. In such situations, one can show that the
condition (\ref{downsum}) is more or less optimal. If $\La$ is a
transitive graph, then the random mapping representations of the voter
model, the contact process and other systems on $\La$ that we have
already seen, are transitive. On the other hand, it is nowadays common
to study interacting particle systems in a random environment, for
example on random graphs. The random mapping representations of such
systems are clearly not transitive. In such situations, the supremum
in (\ref{downsum}) is a nuisance and Theorem~\ref{T:graph} may be
applicable even when Theorem~\ref{T:Poispart} is not applicable.

We now set out to prove Theorem~\ref{T:graph}, which by
Proposition~\ref{P:back} implies Theorem~\ref{T:Poispart}.

\begin{lemma}[Evolution equation]
Under\label{L:Xevolve} the assumptions of Theorem \ref{T:graph},
almost surely, for each $s\in\R$ and $x\in S^\La$, the evolution
equation (\ref{Xevolve}) has a unique solution $(X^{s,x}_t)_{t\geq s}$
with initial state $X^{s,x}_s=x$. For any finite set $T$ and
$\phi\in\Ci(S^\La,T)$, this solution satisfies
\begin{equation}\label{FXdu}
\phi(X^{s,x}_t)=\Fb_{t,s}(\phi)(x)\qquad(t\geq 0),
\end{equation}
where $(\Fb_{t,s})_{t\geq s}$ is the backward stochastic flow defined
in (\ref{Fback}).
\end{lemma}

\begin{Proof}
By Theorem~\ref{T:backflow}, the assumptions of Theorem~\ref{T:graph}
guarantee that the backward stochastic flow in (\ref{Fback}) is
well-defined. For each $i\in\La$, we define $\phi_i\in\Ci(S^\La,S)$ by
$\phi_i(x):=x(i)$ $(x\in S^\La)$. We now fix $s\in\R$ and $x\in S^\La$
and define $(X_t)_{t\geq s}$ by
\[
X_t(i):=\Fb_{t,s}(\phi_i)(x)\qquad(t\geq s,\ i\in\La),
\]
where $(\Fb_{t,s})_{t\geq s}$ is the backward stochastic flow defined
in (\ref{Fback}). Then $X_s=x$. Moreover $t\mapsto X_t(i)$ is
piecewise constant, right-continuous, and
\[
X_t(i)=\left\{\begin{array}{ll}
\dis m(X_{t-})(i)\quad&\mbox{if }\exists(m,t)\in\om
\mbox{ s.t.\ }\phi_i\circ m\neq\phi_i,\\[5pt]
X_{t-}(i)\quad&\mbox{otherwise.}
\end{array}\right.
\]
Using (\ref{Hwell}) and the fact that $\phi_i\circ m\neq\phi_i$
implies $i\in\Di(m)$, we see that $(X_t)_{t\geq s}$ is cadlag and
solves (\ref{Xevolve}). This establishes existence of solutions. To
prove uniqueness, assume that $(X_t)_{t\geq s}$ is cadlag and solves
(\ref{Xevolve}) with $X_s=x$. Let $T$ be a finite set and fix
$\phi\in\Ci(S^\La,T)$. We claim that for each $u>s$, the function
\begin{equation}\label{Fpath}
[s,u]\ni t\mapsto\Fb_{u,t}(\phi)(X_t)
\end{equation}
is constant. Indeed, $t\mapsto\Fb_{u,t}(\phi)$ is piecewise constant
and right-continuous and takes values in $\Ci(S^\La,T)$, so by the
fact that functions in $\Ci(S^\La,T)$ depend on finitely many
coordinates and $t\mapsto X_t(i)$ is piecewise constant and
right-continuous for each $i\in\La$, we see that also
$t\mapsto\Fb_{u,t}(\phi)(X_t)$ is piecewise constant and
right-continuous. For each $(m,t)\in\om$ with $t\in(s,u]$, one has
\[
\Fb_{u,t-}(\phi)(X_{t-})=\Fb_{u,t}(\phi)\circ m(X_{t-})=\Fb_{u,t}(\phi)(X_t),
\]
while at the remaining times in $(s,u]$ trivially
$\Fb_{u,t-}(\phi)(X_{t-})=\Fb_{u,t}(\phi)(X_t)$. This proves that
$t\mapsto\Fb_{u,t}(\phi)(X_t)$ is constant on $(s,u]$ and hence, by
right-continuity, also on $[s,u]$. Since the function in
(\ref{Fpath}) is constant,
\[
\phi(X_u)=\Fb_{u,u}(\phi)(X_u)=\Fb_{u,s}(\phi)(X_s)=\Fb_{u,s}(\phi)(x).
\]
Since this holds for arbitrary $T$ and $\phi\in\Ci(S^\La,T)$, we conclude
that $(X_t)_{t\geq s}$ is unique and that (\ref{FXdu}) holds.
\end{Proof}

By Lemma~\ref{L:Xevolve}, under the assumptions of
Theorem~\ref{T:graph}, almost surely, for each $s\in\R$ and $x\in
S^\La$, the evolution equation (\ref{Xevolve}) has a unique solution
$(X^{s,x}_t)_{t\geq s}$ with initial state $X^{s,x}_s=x$. We use this
to define random maps $(\Xb_{s,t})_{s\leq t}$ as in
(\ref{partflow}). Then (\ref{FXdu}) implies that for each finite set
$T$, one has
\begin{equation}\label{FXrel}\index{0Fb@$\Fb_{t,s}$}
\Fb_{t,s}(\phi)=\phi\circ\Xb_{s,t}\qquad\big(s\leq t,\ \phi\in\Ci(S^\La,T)\big).
\end{equation}
It is straightforward from the definition that these maps satisfy
(\ref{flow}). The the stochastic flow $(\Xb_{s,t})_{s\leq t}$ is
clearly stationary. Using the fact that restrictions of a Poisson
point set to disjoint parts of the space are independent, we also see
that $(\Xb_{s,t})_{s\leq t}$ has independent increments.

\begin{lemma}[Continuity of the flow]
Under\label{L:Xcont} the assumptions of Theorem~\ref{T:graph}, almost
surely, the maps $\Xb_{s,t}\cn S^\La\to S^\La$ are continuous for all
$s\leq t$.
\end{lemma}

\begin{Proof}
Since $S^\La$ is equipped with the product topology, it suffices to
show that $x\mapsto\Xb_{s,t}(x)(i)$ is continuous for all $s\leq t$
and $i\in\La$. Using notation as in the proof of
Lemma~\ref{L:Xevolve}, we have by (\ref{FXrel}) that
\begin{equation}\label{XFi}
\Xb_{s,t}(x)(i)=\phi_i\big(\Xb_{s,t}(x)\big)=\Fb_{t,s}(\phi_i)(x).
\end{equation}
Since $\Fb_{t,s}(\phi_i)\in\Ci(S^\La,S)$, the map
$x\mapsto\Fb_{t,s}(\phi_i)(x)$ is continuous.
\end{Proof}

\begin{lemma}[Almost sure continuity]
Assume\label{L:ascont} that $x_n,x\in S^\La$ satisfy $x_n\to x$ in the
product topology and that $t_n,t\geq 0$ satisfy $t_n\to t$. Then,
under the assumptions of Theorem~\ref{T:graph},
\[
\Xb_{0,t_n}(x_n)\asto{n}\Xb_{0,t}(x)\quad{\rm a.s.}
\]
\end{lemma}

\begin{Proof}
In line with notation introduced in Section~\ref{S:local}, let
$\Xb_{s,u}[i]\cn S^\La\to S$\index{0Xbi@$\Xb_{s,t}[i]$} be defined as
$\Xb_{s,u}[i](x):=\Xb_{s,u}(x)(i)$.  Since $S^\La$ is equipped with
the product topology, it suffices to show that $\Xb_{0,t_n}(x_n)(i)$
converges to $\Xb_{0,t}(x)(i)$ for each $i\in\La$. Since $t$ is
deterministic, by (\ref{Hwell}) there a.s.\ exists an $\eps>0$ such
that
\[
i\not\in\Di(m)\quad\forall(m,r)\in\om\quad\mbox{with}\quad t-\eps<r<t+\eps.
\]
It follows that $\Xb_{0,t_n}[i]=\Xb_{0,t}[i]$ for all $n$
large enough such that $t-\eps<t_n<t+\eps$. Since $x_n\to x$,
Lemma~\ref{L:Xcont} now tells us that
$\Xb_{0,t}(x_n)(i)=\Xb_{0,t}(x)(i)$ for all $n$ large enough.
\end{Proof}

\begin{Proof}[of Theorem~\ref{T:graph}]
By Lemma~\ref{L:Xevolve}, almost surely, for each $s\in\R$ and $x\in
S^\La$, the evolution equation (\ref{Xevolve}) has a unique solution
$(X^{s,x}_t)_{t\geq s}$ with initial state $X^{s,x}_s=x$, which allows
us to define random maps $(\Xb_{s,t})_{s\leq t}$ satisfying
(\ref{flow}) as in (\ref{partflow}). By Lemma~\ref{L:Xcont} the maps
$\Xb_{s,t}\cn S^\La\to S^\La$ are continuous. To see that
(\ref{Ppart}) defines a Feller semigroup, we need to check that
\begin{enumerate}[(ii)]
\item $(x,t)\mapsto P_t(x,\,\cdot\,)$ is a continuous map from
  $E\times\half$ into $\Mi_1(E)$,
\item $P_0=1\mbox{ and }P_sP_t=P_{s+t}\quad(s,t\geq 0)$.
\end{enumerate}
Property~(i) follows from Lemma~\ref{L:ascont} and the fact that
a.s.\ convergence implies weak convergence in law. To prove (ii) we
observe that for each bounded measurable $f\cn S^\La\to\R$,
\[
P_0f(x)=\E\big[f(\Xb_{0,0}(x))\big]=f(x)
\]
and, for every measurable $A\sub S^\La$,
\[\begin{array}{l}
\dis P_{s+t}(x,A)=\P\big[\Xb_{s,s+t}\circ\Xb_{0,s}(x)\in A\big]\\[5pt]
\dis\quad=\int\P\big[\Xb_{0,s}(x)\in\di y\big]
\,\P\big[\Xb_{s,s+t}(y)\in A\,|\,\Xb_{0,s}(x)=y\big]\\[5pt]
\dis\quad=\int\P\big[\Xb_{0,s}(x)\in\di y\big]
\,\P\big[\Xb_{s,s+t}(y)\in A\big]\\[5pt]
\dis\quad=\int P_s(x,\di y)P_t(y,A)=P_sP_t(x,A),
\end{array}\]
where we have used the stationarity and independent increments of the
stochastic flow.

Using the fact that the stochastic flow $(\Xb_{s,t})_{s\leq t}$ is
stationary and has independent increments, it is now straightforward
to check that the process in (\ref{flowMark}) is distributed as a
Markov process with semigroup $(P_t)_{t\geq 0}$. Note that
$(X_t)_{t\geq 0}$ has cadlag sample paths since solutions to the
evolution equation (\ref{Xevolve}) are cadlag. This completes the
proof.
\end{Proof}

\section{Finite approximation and perturbations}\label{S:finax}

In this section we use the methods of the previous section to prove
two additional results that are sometimes useful. In (\ref{partflow}),
we defined a stochastic flow $(\Xb_{s,t})_{s\leq t}$ in terms of a
graphical representation $\om$ by means of the evolution equation
(\ref{Xevolve}). Our first aim in this section is to show that
$(\Xb_{s,t})_{s\leq t}$ can alternatively be defined by approximating
$\om$ with finite sets.

We first need a simple fact about continuous-time Markov chains. Let
$S$ be a countable set. Recall that a random mapping representation of
the generator $G$ of a continuous-time Markov chain with state space $S$
is an expression of the form (\ref{Grep}) where the rates
$(r_m)_{m\in\Gi}$ satisfy (\ref{locfin}). Let $\om$ be the graphical
representation associated with such a random mapping representation
and let $\om'$ be a finite subset of $\om$. Then for each $s\leq u$ we
define a map $\Xb^{\om'}_{s,u}\cn S\to S$ by
\begin{equation}\begin{array}{l}\label{Xom}
\dis\quad\Xb^{\om'}_{s,u}:=m_n\circ\cdots\circ m_1,\quad\mbox{where}\\[5pt]
\dis\big\{(m,t)\in\om':t\in(s,u]\big\}
=\big\{(m_1,t_1),\ldots,(m_n,t_n)\big\}\quad\mbox{with}\quad t_1<\cdots<t_n.
\end{array}\end{equation}
In words, this says that $\Xb^{\om'}_{s,u}$ is the concatenation of
the maps $m$ for which $(m,t)\in\om'$ with $t\in(s,u]$, ordered by
their times. Below, we equip $S$ with the discrete topology, so
(\ref{Mfincon}) simply says that $\Xb^{\om_n}_{s,u}(x)=\Xb_{s,u}(x)$
for all $n$ large enough.

\begin{lemma}[Finite approximation of Markov chains]
Let\label{L:Mfincon} $G$ be the generator of a nonexplosive
con\-tinu\-ous-time Markov chain with countable state space $S$, let
$\om$ be the graphical representation associated with a random mapping
representation of $G$, and let $(\Xb_{s,t})_{s\leq t}$ be the
stochastic flow defined in terms of $\om$. Then almost surely, for all
$s\leq u$, for all $x\in S$, and for each sequence $\om_n$ of finite
subsets of $\om$ such that $\om_n\up\om$, one has that\footnote{For
those who like abstract mathematics, another way to formulate
(\ref{Mfincon}) is as follows. Denoting by $F$ the set of finite
subsets of $\om$, equipped with the partial order of set inclusion, we
have that $(\Xb^\om_{s,u})^{\om\in F}$ is a \emph{net} in the
topological meaning of the word, and (\ref{Mfincon}) says that this
net converges pointwise to $\Xb_{s,u}$.}
\begin{equation}\label{Mfincon}
\Xb^{\om_n}_{s,u}(x)\asto{n}\Xb_{s,u}(x).
\end{equation}
An analogue statement holds for the backward stochastic flow defined
in terms of $\om$ as in Theorem~\ref{T:backflow}, where now
$\Xb^{\om_n}_{u,s}$ is defined by ordering the elements of
$\{(m,t)\in\om_n:t\in(s,u]\}$ in decreasing order of time.
\end{lemma}

\begin{Proof}
Fix $s\leq u$ and $x\in S$. Then $[s,u]\ni t\mapsto\Xb_{s,t}(x)$ is
piecewise constant and right-continuous. We set
\[
\om^\ast:=\big\{(m,t)\in\om:
t\in(s,u],\ m\big(\Xb_{s,t-}(x)\big)\neq\Xb_{s,t-}(x)\big\},
\]
and write
\[
\om^\ast=\big\{(m_1,t_1),\ldots,(m_k,t_k)\big\}
\quad\mbox{with}\quad t_1<\cdots<t_k.
\]
Then
\begin{equation}\label{concat}
\Xb_{s,u}(x)=m_k\circ\cdots\circ m_1(x).
\end{equation}
It follows from the definition of $\om^\ast$ that if $(m,t)\in\om\beh\om^\ast$
for some $t\in(s,u]$, then $m\big(\Xb_{s,t-}(x)\big)=\Xb_{s,t-}(x)$.
This implies that the right-hand side of (\ref{concat}) does not
change if in the concatenation of maps we add additional maps from
$\om\beh\om^\ast$ at their appropriate times. It follows that for all $n$
large enough such that $\om_n\supset\om^\ast$, one has
\[
\Xb^{\om_n}_{s,u}(x)=m_k\circ\cdots\circ m_1(x)=\Xb_{s,u}(x).
\]
\end{Proof}

If $\om$ is the graphical representation of an interacting particle
system and $\om'$ is a finite subset of $\om$, then we define
$\Xb^{\om'}_{s,t}$ in the same way as in (\ref{Xom}).

\begin{proposition}[Finite approximation]
Under\label{P:fincon} the assumptions of Theorem~\ref{T:graph}, the
stochastic flow $(\Xb_{s,t})_{s\leq t}$ has the property that almost
surely, for all $s\leq u$, for all $x\in S^\La$, and for each sequence
$\om_n$ of finite subsets of $\om$ such that $\om_n\up\om$, one has
that
\begin{equation}\label{fincon}
\Xb^{\om_n}_{s,u}(x)\asto{n}\Xb_{s,u}(x)
\end{equation}
with respect to the product topology on $S^\La$.
\end{proposition}

\begin{Proof}
We continue to use the notation $\Xb_{s,u}[i](x):=\Xb_{s,u}(x)(i)$
\index{0Xbi@$\Xb_{s,t}[i]$} $(i\in\La,\ x\in S^\La)$ and as we did before
define $\phi_i\in\Ci(S^\La,S)$ by $\phi_i(x):=x(i)$ $(i\in\La)$.
Fix $s\leq u$ and $i\in\La$. Then (\ref{FXrel}) tells us that
\[
\Xb_{s,u}[i]=\Fb_{u,s}(\phi_i).
\]
Similarly, for each finite $\om'\sub\om$, if
$(m_1,t_t),\ldots,(m_k,t_k)$ are the elements of
$\{(m,t)\in\om':t\in(s,u]\}$, ordered according to their times, then
\[
\Xb^{\om'}_{s,u}[i]=\phi_i\circ m_k\circ\cdots\circ m_1=\Fb^{\om'}_{u,s}(\phi_i).
\]
Applying Lemma~\ref{L:Mfincon} to the backward stochastic flow
$(\Fb_{u,s})_{u\geq s}$ we see that if $\om_n\up\om$, then for all $n$
large enough
\[
\Xb^{\om_n}_{s,u}[i]=\Fb^{\om_n}_{u,s}(\phi_i)=\Fb_{u,s}(\phi_i)=\Xb_{s,u}[i].
\]
Since this holds for each $i\in\La$, it follows that
$\Xb^{\om_n}_{s,u}(x)$ converges to $\Xb_{s,u}(x)$ in the product
topology for each $x\in S^\La$.
\end{Proof}

We have already noticed the similarity between condition (\ref{upsum})
of Theorem~\ref{T:finIPS} and condition (\ref{downsum}) of
Theorem~\ref{T:Poispart}, the only difference being that
$\Ri^\up_i(m)$ in (\ref{upsum})~(ii) is replaced by $\Ri^\down_i(m)$
in (\ref{downsum})~(ii). In Theorem~\ref{T:finIPS}, we assumed that
$m(\un 0)=\un 0$ for all $m\in\Gi$, which implies that the all-zero
configuration is a trivial fixed point of the evolution equation
(\ref{Xevolve}). The theorem then showed that under the condition
(\ref{upsum}), finite perturbations of the all-zero initial state have
finite consequences at later times. We show here that this statement
generalizes to arbitrary initial states.

\begin{proposition}[Finite perturbations]
Assume\label{P:finpert} that the conditions of Theorem~\ref{T:graph}
are satisfied and that (\ref{upsum}) holds, that is,
\begin{equation}\label{upsum2}
{\rm(i)}\ \sup_{i\in\La}\sum_{m\in\Gi}r_m1_{\Di(m)}(i)<\infty,\quad
{\rm(ii)}\ \sup_{i\in\La}\sum_{m\in\Gi}r_m\big|\Ri^\up_i(m)\beh\{i\}\big|<\infty.
\end{equation}
Then almost surely, for each $x,y\in S^\La$ such that
$A:=\{i\in\La:x(i)\neq y(i)\}$ is finite, one has that
\begin{equation}\label{Ast}
A_{s,t}:=\big\{i\in\La:\Xb_{s,t}(x)(i)\neq\Xb_{s,t}(y)(i)\big\}
\end{equation}
is finite for all $s\leq t$. If $s,t,x,y$ are deterministic, then
\begin{equation}\label{expert}
\E\big[|A_{s,t}|\big]\leq|A|\ex{K_\up(t-s)},
\end{equation}
where $K_\up$ is defined in (\ref{KKK}).
\end{proposition}

\begin{Proof}
We assume without loss of generality that $|S|\geq 2$. Indeed, if
$|S|=1$, then $A_{s,t}=\emptyset$ for all $s,t,x,y$ and the statements
are trivial. It follows from (\ref{upsum2}) that if
$|\Ri^\up_i(m)|=\infty$ for some $i\in\La$ and $m\in\Gi$, then
$r_m=0$, so we can without loss of generality assume that
$\Ri^\up_i(m)$ is finite for all $i\in\La$ and $m\in\Gi$. Let $\Pc$
denote the set of all subsets of $\La$ and set
$\Pc_{\rm fin}:=\{A\in\Pc:|A|<\infty\}$. For each $m\in\Gi$, we define a map
$\ov m\cn\Pc_{\rm fin}\to\Pc_{\rm fin}$ by
\[
\ov m(A):=\bigcup_{i\in A}\Ri^\up_i(m)\qquad(A\in\Pc_{\rm fin}).
\]
We claim that setting
\[
\ov Gf(A):=\sum_{m\in\Gi}r_m\big\{f\big(\ov m(A)\big)-f\big(A\big)\big\}
\]
defines the generator of a nonexplosive continuous-time Markov chain
$(Z_t)_{t\geq 0}$ with state space $\Pc_{\rm fin}$, and that this
Markov chain satisfies
\begin{equation}\label{Zex}
\E^A\big[|Z_t|\big]\leq|A|\ex{K_\up t}\qquad(t\geq 0,\ A\in\Pc_{\rm fin}).
\end{equation}
We start by checking (\ref{locfin}), which in the present context says that
\[
\sum_{m:\,\ov m(A)\neq A}r_m<\infty\quad\forall A\in\Pc_{\rm fin}.
\]
If $\ov m(A)\neq A$, then either there exists an $i\in A$ such that
$i\not\in\ov m(A)$, or there exists a $j\not\in A$ such that $j\in\ov
m(A)$. In the first case $i\not\in\Ri^\up_i(m)$ which by the fact that
$|S|\geq 2$ implies $i\in\Di(m)$. In the second case there must be an
$i\in A$ such that $j\in\Ri^\up_i(m)$. It follows that
\begin{equation}\begin{array}{r@{\,}c@{\,}l}\label{cruc}
\dis\sum_{m:\,\ov m(A)\neq A}r_m
&\leq&\dis\sum_{m\in\Gi}\sum_{i\in A}1_{\Di(m)}(i)r_m
+\sum_{m\in\Gi}\sum_{i\in A}\sum_{j\in\Ri^\up_i(m)\beh A}r_m\\[5pt]
&\leq&\dis|A|\sup_{i\in\La}\sum_{m\in\Gi}1_{\Di(m)}(i)r_m
+|A|\sup_{i\in\La}\sum_{m\in\Gi}\big|\Ri^\up_i(m)\beh\{i\}\big|r_m,
\end{array}\end{equation}
which is finite by (\ref{upsum2}). Lemma~\ref{L:randmap} now implies
that $\ov G$ is the generator of a (possibly explosive)
continuous-time Markov chain with state space $\Pc_{\rm fin}$. To
prove nonexplosiveness and the exponential bound (\ref{Zex}) we apply
Theorem~\ref{T:Lyap} to the Lyapunov function $L(A):=|A|$. Formula
(\ref{cruc}) shows that condition~(i) of Theorem~\ref{T:Lyap} is
satisfied so it remains to check that $\ov GL\leq K_\up L$. This is
very similar to the last steps of the proof of
Proposition~\ref{P:back}. Indeed,
\[
\big|\ov m(A)\big|-\big|A\big|\leq\sum_{i\in A}\big(|\Ri^\up_i(m)|-1\big),
\]
which implies that
\[
\ov GL(A)\leq\sum_{m\in\Gi}\sum_{i\in A}r_m\big(|\Ri^\up_i(m)|-1\big)
\leq K_\up|A|=K_\up L(A).
\]

Let $(\Zb_{s,t})_{s\leq t}$ be the stochastic flow constructed from
the Poisson point set $\ov\om:=\{(\ov m,t):(m,t)\in\om\}$. Since $\ov
G$ is nonexplosive, Theorem~\ref{T:stochflow} tells us that almost
surely, $\Zb_{s,t}$ maps $\Pc_{\rm fin}$ into itself for each $s\leq t$.
To complete the proof, we will show that almost surely, for each
$s\leq t$ and $x,y\in S^\La$ such that $A:=\{i\in\La:x(i)\neq y(i)\}$
is finite, the set $A_{s,t}$ defined in (\ref{Ast}) is contained in
$\Zb_{s,t}(A)$. The basic observation we need is that for each $m\in\Gi$
and $x,y\in S^\La$ such that $\{i\in\La:x(i)\neq y(i)\}$ is finite
\[
\big\{i\in\La:x(i)\neq y(i)\big\}\sub A\quad\mbox{implies}\quad
\big\{i\in\La:m(x)(i)\neq m(y)(i)\big\}\sub\ov m(A).
\]
Using this, we see that for each $x,y\in S^\La$ such that
$\{i\in\La:x(i)\neq y(i)\}$ is finite, for each $s\leq u$, and for
each finite $\om_n\sub\om$,
\[
\big\{i\in\La:x(i)\neq y(i)\big\}\sub A\quad\mbox{implies}\quad
\big\{i\in\La:\Xb^{\om_n}_{s,u}(x)(i)\neq\Xb^{\om_n}_{s,u}(y)(i)\big\}
\sub\Zb^{\om_n}_{s,u}(A).
\]
The claim now follows by letting $\om_n\up\om$, using
Lemma~\ref{L:Mfincon} and Proposition~\ref{P:fincon}.
\end{Proof}

\section{Generator construction}\label{S:Gencon}

Although Theorem~\ref{T:Poispart} gives us an explicit way how to construct the
Feller semigroup associated with an interacting particle system, it does not
tell us very much about its generator. To fill this gap, we need a bit more
theory. For any continuous function $f\cn S^\La\to\R$ and $i\in\La$, we define
\[
\de f(i):=\sup\big\{|f(x)-f(y)|:x,y\in S^\La,\ x(j)=y(j)\ \forall j\neq i\big\}.
\]
Note that $\de f(i)$ measures how much $f(x)$ can change if we change
$x$ only in the point $i$. We call $\de f$ the
\emph{variation}\index{variation of a function on $S^\La$} of
$f$.\footnote{This definition is similar to, but different from the
more usual definition of the (total) variation of a function of one
real variable. With functions of one real variable, the total
variation is the maximal sum of all changes in the value of the
function as one gradually increases the real variable. For functions
on $S^\La$, the idea is similar but instead of increasing a real
variable we will gradually change a configuration $x$ by modifying its
coordinates one by one.}

\begin{lemma}[Variation of a function]
Let\label{L:varsum} $f\in\Ci(S^\La)$. Then
\begin{equation}\label{varsum}
\big|f(x)-f(y)\big|\leq\sum_{i:\,x(i)\neq y(i)}\de f(i)
\qquad\big(f\in\Ci(S^\La),\ x,y\in S^\La\big).
\end{equation}
\end{lemma}

\begin{Proof}
Let $n$ be the number of sites $i$ where $x$ and $y$
differ. Enumerate these sites as $\{i:x(i)\neq y(i)\}=\{i_1,\ldots,i_n\}$ or
$=\{i_1,i_2,\ldots\}$ depending on whether $n$ is finite or not.
For $0\leq k<n+1$, set
\[
z_k(i):=\left\{\begin{array}{ll}
y(i)\quad\mbox{if }i\in\{i_1,\ldots,i_k\},\\[5pt]
x(i)\quad\mbox{otherwise.}
\end{array}\right.
\]
If $n$ is finite, then
\[
\big|f(x)-f(y)\big|\leq\sum_{k=1}^n\big|f(z_k)-f(z_{k-1})\big|
\leq\sum_{k=1}^n\de f(i_k)
\]
and we are done. If $n$ is infinite, then the same argument gives
\[
\big|f(x)-f(z_m)\big|\leq\sum_{k=1}^m\de f(i_k)\qquad(m\geq 1).
\]
Since $z_m\to y$ pointwise and $f$ is continuous, (\ref{varsum}) now follows by
letting $m\to\infty$.
\end{Proof}

We define spaces of functions by
\[\begin{array}{r@{\,}c@{\,}l}\index{0Csum@$\Ci_{\rm sum}$}
\index{0Cfin@$\Ci_{\rm fin}$}
\dis\Ci_{\rm sum}=\Ci_{\rm sum}(S^\La)
&:=&\dis\big\{f\in\Ci(S^\La):\sum_i\de f(i)<\infty\big\},\\[5pt]
\dis\Ci_{\rm fin}=\Ci_{\rm fin}(S^\La)
&:=&\dis\big\{f\in\Ci(S^\La):\de f(i)=0
\mbox{ for all but finitely many }i\big\}.
\end{array}\]
We say that functions in $\Ci_{\rm sum}$ are of \emph{summable
variation}.\index{summable variation} The next exercise shows that
functions in $\Ci_{\rm fin}$ depend on finitely many coordinates only.

\begin{Exercise}\label{E:findens}
Let us say that a function $f\cn S^\La\to\R$ \emph{depends on finitely
  many coordinates} if there exists a finite set $A\sub\La$ and a function
$f'\cn S^A\to\R$ such that
\[
f\big((x(i))_{i\in\La}\big)=f'\big((x(i))_{i\in A}\big)
\qquad\big(x\in S^\La\big).
\]
Show that each function that depends on finitely many coordinates is
continuous, that
\[
\Ci_{\rm fin}(S^\La)=\big\{f\in\Ci(S^\La):f
\mbox{ depends on finitely many coordinates}\big\},
\]
and that $\Ci_{\rm fin}(S^\La)$ is a dense linear subspace of the
Banach space $\Ci(S^\La)$ of all continuous real functions on
$S^\La$, equipped with the supremum-norm.
\end{Exercise}

\begin{Exercise}
Define $f\cn\{0,1\}^\N\to\R$ by
\[
f(x):=\frac{1}{1+r}\quad\mbox{with}\quad r:=\inf\{i\geq 0:x(i)=1\}.
\]
Show that $f\in\Ci(\{0,1\}^\N)$ but $f\not\in\Ci_{\rm sum}(\{0,1\}^\N)$.
\end{Exercise}

In what follows, we assume that $\Gi$ is a countable collection of
continuous maps $m\cn S^\La\to S^\La$ and that $(r_m)_{m\in\Gi}$ are
nonnegative rates.

\begin{lemma}[Domain of pregenerator]
Assume\label{L:welldef} (\ref{downsum})~(i) and let $K_0$ be the
constant defined in (\ref{KKK}). Then, for each $f\in\Ci_{\rm sum}(S^\La)$,
\[
\sum_{m\in\Gi}r_m\big|f(m(x))-f(x)|\leq K_0\sum_{i\in\La}\de f(i).
\]
In particular, for each $f\in\Ci_{\rm sum}(S^\La)$ and $x\in S^\La$,
the right-hand side of (\ref{Gmap}) is absolutely summable and $Gf$ is
well-defined.
\end{lemma}

\begin{Proof}
This follows by writing
\[
\sum_{m\in\Gi}r_m\big|f(m(x))-f(x)|
\leq\sum_{m\in\Gi}r_m\sum_{i\in\Di(m)}\de f(i)
=\sum_{i\in\La}\de f(i)\!\!\!\sum\subb{m\in\Gi}{\Di(m)\ni i}\!\!r_m
\leq K_0\sum_{i\in\La}\de f(i).
\vspace{-15pt}
\]
\end{Proof}

The following theorem is the main result of the present section.

\begin{theorem}[Generator construction of particle systems]\label{T:genpart}
Assume that the rates $(r_m)_{m\in\Gi}$ satisfy
(\ref{downsum}), let $(P_t)_{t\geq 0}$ be the Feller semigroup defined
in (\ref{Ppart}) and let $G$ be the linear operator with domain
$\Di(G):=\Ci_{\rm sum}$ defined by (\ref{Gmap}). Then $G$ is closable
and its closure $\ov G$ is the generator of $(P_t)_{t\geq 0}$. Moreover,
$\Ci_{\rm fin}$ is a core for $G$.
\end{theorem}

To prepare for the proof of Theorem~\ref{T:genpart} we need a few lemmas.

\begin{lemma}[Generator on local functions]\label{L:locgen} 
Under the assumptions of Theorem~\ref{T:genpart}, one has $\lim_{t\down 0}t^{-1}(P_tf-f)=Gf$ for all $f\in\Ci_{\rm fin}$, where the limit holds with respect to the topology on $\Ci(S^\La)$.
\end{lemma}

\begin{Proof}
Since $f\in\Ci_{\rm fin}$, by Exercise~\ref{E:findens}, there exists
some finite $A\sub\La$ such that $f$ depends only on the coordinates
in $A$.  Let $\Gi^A:=\{m\in\Gi:\Di(m)\cap A\neq\emptyset\}$ denote the
set of maps $m\in\Gi$ that can potentially change the state in $A$. We
introduce the notation
\[\index{0zzzomsu@$\om_{s,t}$}
\om_{s,t}:=\big\{(m,r)\in\om:r\in(s,t]\big\}\qquad(s\leq t)
\]
and we let $\om^A_{s,t}$ denote the set of Poisson points
$(m,r)\in\om_{s,t}$ with $m\in\Gi^A$. If $\om^A_{0,t}=\emptyset$, then
$f(\Xb_{0,t}(x))=f(x)$. Also, if $\om^A_{0,t}$ contains a single
element $(m,s)$, then $f(\Xb_{0,t}(x))=f(m(x))$. Therefore
\[\begin{array}{r@{\,}l}
\dis P_tf(x)
=&\dis\E\big[f\big(\Xb_{0,t}(x)\big)\big]=f(x)\P[\om^A_{0,t}=\emptyset]\\[5pt]
&\dis+\sum_{m\in\Gi^A}f\big(m(x)\big)
\P\big[\om^A_{0,t}=\{(m,s)\}\mbox{ for some }0<s\leq t\big]\\[5pt]
&+\E\big[f\big(\Xb_{0,t}(x)\big)1_{\txt\{|\om^A_{0,t}|\geq 2\}}\big].
\end{array}\]
Here, setting $R:=\sum_{m\in\Gi^A}r_m$, which is finite by the
finiteness of $A$ and (\ref{downsum})~(i), we have
\[\begin{array}{l}
\dis\P[\om^A_{0,t}=\emptyset]=e^{-Rt},\\[5pt]
\dis\P\big[\om^A_{0,t}=\{(m,s)\}\mbox{ for some }0<s\leq t\big]
=tr_me^{-Rt}\quad(m\in\Gi^A),\\[5pt]
\dis\P[|\om^A_{0,t}|\geq 2]=1-e^{-Rt}-tRe^{-Rt}.
\end{array}\]
It follows that
\begin{equation}\begin{array}{r@{\,}c@{\,}l}\label{Pord}
\dis P_tf(x)
&=&\dis e^{-Rt}f(x)+te^{-Rt}\sum_{m\in\Gi^A}r_m f\big(m(x)\big)+O_x(t^2),\\[5pt]
&=&\dis f(x)+t\sum_{m\in\Gi^A}r_m\big\{f\big(m(x)\big)-f\big(x\big)\big\}+O_x(t^2),
\end{array}\end{equation}
where $O_x(t^2)$ denotes a function, which may differ from line to
line, that has the property that
$\limsup_{t\to\infty}t^{-2}\sup_{x\in S^\La}|O_x(t^2)|<\infty$. Indeed, in the first
line of (\ref{Pord}),
\[
O_x(t^2)=\E\big[f\big(\Xb_{0,t}(x)\big)1_{\txt\{|\om^A_{0,t}|\geq 2\}}\big].
\]
This can be estimated as
\[
\sup_{x\in S^\La}|O_x(t^2)|\leq\|f\|_\infty\P[|\om^A_{0,t}|\geq 2]=\|f\|_\infty\big(1-e^{-Rt}-tRe^{-Rt}\big),
\]
which using the fact that $e^{-Rt}=1-Rt+O(t^2)$ as $t\to 0$ shows that
the error term in the first line of (\ref{Pord}) is of order $t^2$
uniformly in $x\in S^\La$. The second line of (\ref{Pord}) now also
follows readily, with a somewhat different definition of $O_x(t^2)$.

Since $f(m(x))=f(x)$ if $m\not\in\Gi^A$, formula (\ref{Pord}) implies that
\[
t^{-1}\big(P_tf(x)-f(x)\big)=Gf(x)+O_x(t),
\]
where $O_x(t)$ denotes a term that is of order $t$ as $t\to 0$,
uniformly in $x\in S^\La$. This shows that
\[
\lim_{t\to\infty}\big\|t^{-1}(P_tf-f)-Gf\big\|_\infty=0,
\]
as claimed.
\end{Proof}

\begin{lemma}[Approximation by local functions]
Assume\label{L:locapp} (\ref{downsum})~(i). Then for all $f\in\Ci_{\rm sum}$
there exist $f_n\in\Ci_{\rm fin}$ such that $\|f_n-f\|_\infty\to 0$
and $\|Gf_n-Gf\|_\infty\to 0$.
\end{lemma}

\begin{Proof}
Choose finite $\La_n\up\La$, set $\Ga_n:=\La\beh\La_n$, fix $z\in
S^\La$, and for each $x\in S^\La$ define $x^n\to x$ by
\[
x^n(i):=\left\{\begin{array}{ll}x(i)\quad&\mbox{if }i\in\La_n,\\
z(i)\quad&\mbox{if }i\in\Ga_n.\end{array}\right.
\]
Fix $f\in\Ci_{\rm sum}$ and define $f_n(x):=f(x^n)$ $(x\in
S^\La)$. Then $f_n$ depends only on the coordinates in $\La_n$, hence
$f_n\in\Ci_{\rm fin}$. If $x_n,x\in S^\La$ satisfy $x_n\to x$, then by
the continuity of $f$ we have $f_n(x_n)=f(x^n_n)\to f$, so applying
Lemma~\ref{L:fco} we see that
\[
\|f_n-f\|_\infty\asto{n}0.
\]
To prove that also $\|Gf_n-Gf\|_\infty\to 0$ we observe that
\begin{equation}\begin{array}{l}\label{GGest}
\dis|Gf_n(x)-Gf(x)|\\[5pt]
\dis\quad=\big|\sum_{m\in\Gi}r_m\big(f_n(m(x))-f_n(x)\big)
-\sum_{m\in\Gi}r_m\big(f(m(x))-f(x)\big)\big|\\[5pt]
\dis\quad\leq\sum_{m\in\Gi}r_m\big|f(m(x)^n)-f(x^n)-f(m(x))+f(x)\big|.
\end{array}\end{equation}
On the one hand, we have
\[\begin{array}{l}
\dis\big|f(m(x)^n)-f(x^n)-f(m(x))+f(x)\big|\\[5pt]
\dis\quad\leq\big|f(m(x)^n)-f(x^n)\big|+\big|f(m(x))-f(x)\big|
\leq 2\sum_{i\in\Di(m)}\de f(i),
\end{array}\]
while on the other hand, we can estimate the same quantity as
\[
\quad\leq\big|f(m(x)^n)-f(m(x))\big|+\big|f(x^n)-f(x)\big|
\leq 2\sum_{i\in\Ga_n}\de f(i).
\]
Let $A\sub\La$ be finite. Inserting either of our two estimates into
(\ref{GGest}), depending on whether $\Di(m)\cap A\neq\emptyset$ or not, we
find that
\[\begin{array}{r@{\,}c@{\,}l}
\dis\|Gf_n-Gf\|_\infty
&\leq&\dis 2\!\!\!\!\!\!\sum\subb{m\in\Gi}{\Di(m)\cap A\neq\emptyset}
\!\!\!\!\!\!r_m\sum_{i\in\Ga_n}\de f(i)
+2\!\!\!\!\!\!\sum\subb{m\in\Gi}{\Di(m)\cap A=\emptyset}
\!\!\!\!\!\!r_m\sum_{i\in\Di(m)}\de f(i)\\
&\leq&\dis 2K_0|A|\sum_{i\in\Ga_n}\de f(i)
+2\sum_{i\in\La}\de f(i)
\!\!\!\!\!\!\sum\subbb{m\in\Gi}{\Di(m)\cap A=\emptyset}{\Di(m)\ni i}
\!\!\!\!\!\!r_m,
\end{array}\]
where $K_0$ is the constant defined in (\ref{KKK}). It follows that
\[
\limsup_{n\to\infty}\|Gf_n-Gf\|_\infty
\leq 2\sum_{i\in\La\beh A}\de f(i)
\!\!\!\!\!\!\sum\subb{m\in\Gi}{\Di(m)\ni i}
\!\!\!\!\!\!r_m
\leq 2K_0\sum_{i\in\La\beh A}\de f(i).
\]
Since $A$ is arbitrary, letting $A\up\La$, we see that
$\limsup_n\|Gf_n-Gf\|_\infty=0$.
\end{Proof}

\begin{lemma}[Functions of summable variation]
Under\label{L:Pfsum} the assumptions of Theorem~\ref{T:genpart}, one has
\[
\sum_{i\in\La}\de P_tf(i)\leq\ex{K_\down t}\sum_{i\in\La}\de f(i)
\qquad\big(t\geq 0,\ f\in\Ci_{\rm sum}(S^\La)\big),
\]
where $K_\down$ is the constant from (\ref{KKK}). In particular, for
each $t\geq 0$, $P_t$ maps $\Ci_{\rm sum}(S^\La)$ into itself.
\end{lemma}

\begin{Proof}
Fix $i\in\La$. As we have done before, for $j\in\La$ we define
$\phi_j\in\Ci(S^\La,S)$ by $\phi_j(x):=x(j)$ $(x\in S^\La)$. Then for
each $x,y\in S^\La$ such that $x(k)=y(k)$ for all $k\neq i$, we can
estimate using (\ref{varsum}) and (\ref{FXrel})
\[\begin{array}{l}
\dis|P_tf(x)-P_tf(y)|
=\big|\E[f(\Xb_{0,t}(x))]-\E[f(\Xb_{0,t}(y))]\big|\\[5pt]
\dis\quad\leq\E\big[|f(\Xb_{0,t}(x))-f(\Xb_{0,t}(y))|\big]
\leq\E\big[{\txt\sum_{j:\,\Xb_{0,t}(x)(j)\neq\Xb_{0,t}(y)(j)}}\de f(j)\big]\\[5pt]
\dis\quad=\sum_j\P\big[\Xb_{0,t}(x)(j)\neq\Xb_{0,t}(y)(j)\big]\de f(j)
=\sum_j\P\big[\Fb_{t,0}(\phi_j)(x)\neq\Fb_{t,0}(\phi_j)(y)\big]\de f(j)\\[5pt]
\dis\quad\leq\sum_j\P\big[i\in\Ri(\Fb_{t,0}(\phi_j))\big]\de f(j).
\end{array}\]
By formula (\ref{Phiexp}) of Proposition~\ref{P:back}, it follows that
\[\begin{array}{l}
\dis\sum_i\de P_tf(i)\leq\sum_{ij}\P\big[i\in\Ri(\Fb_{t,0}(\phi_j))\big]\de f(j)\\[5pt]
\dis\quad=\sum_j\E\big[|\Ri(\Fb_{t,0}(\phi_j))|\big]\de f(j)\leq\ex{K_\down t}\sum_j\de f(j).
\end{array}\]
\end{Proof}

\begin{Proof}[of Theorem~\ref{T:genpart}]
Let $H$ be the full generator of $(P_t)_{t\geq 0}$ and let $\Di(H)$
denote it domain. Then Lemma~\ref{L:locgen} shows that $\Ci_{\rm fin}\sub\Di(H)$
and $Gf=Hf$ for all $f\in\Ci_{\rm fin}$. By Lemma~\ref{L:locapp}, it
follows that $\Ci_{\rm sum}\sub\Di(H)$ and $Gf=Hf$ for all
$f\in\Ci_{\rm sum}$. To complete the proof, it suffices to show that
$\Ci_{\rm fin}$, and hence also the larger $\Ci_{\rm sum}$, is a core
for $H$.

We first prove that $\Ci_{\rm sum}$ is a core for $H$. We will apply
Lemma~\ref{L:core}. We will show that for each $r>K_\down$, where
$K_\down$ is the constant from (\ref{KKK}), and for each
$f\in\Ci_{\rm sum}(S^\La)$, there exists a $p_r\in\Ci_{\rm sum}(S^\La)$ that
solves the Laplace equation $(r-G)p_r=f$. Since $\Ci_{\rm sum}(S^\La)$
is dense in $\Ci(S^\La)$ by Exercise~\ref{E:findens}, the claim then
follows from the equivalence of (i) and (ii) of Lemma~\ref{L:core}.

Fix $r>K_\down$ and $f\in\Ci_{\rm sum}(S^\La)$. We need to find a
$p_r\in\Ci_{\rm sum}(S^\La)$ that solves the Laplace equation
$(r-G)p_r=f$. In the light of Lemma~\ref{L:Laplace} a natural
candidate for such a function is
\[
p_r:=\int_0^\infty e^{-rt}P_tf\,\di t
\]
and we will show that this $p_r$ indeed satisfies $p_r\in\Ci_{\rm sum}(S^\La)$
and $(r-G)p_r=f$. It follows from Theorem~\ref{T:HY2} that $p_r\in\Di(H)$ and
$(r-H)p_r=f$. Thus, it suffices to show that $p_r\in\Ci_{\rm sum}$.
To see this, note that if $x(j)=y(j)$ for all $j\neq i$, then
\[\begin{array}{l}
\dis|p_r(x)-p_r(y)|=\Big|\int_0^\infty e^{-rt}P_tf(x)\,\di t
-\int_0^\infty e^{-rt}P_tf(y)\di t\Big|\\[5pt]
\dis\quad\leq\int_0^\infty e^{-rt}\big|P_tf(x)-P_tf(y)\big|\,\di t
\leq\int_0^\infty e^{-rt}\de P_tf(i)\,\di t,
\end{array}\]
and therefore, by Lemma~\ref{L:Pfsum}, and our assumption that $r>K_\down$
\[
\sum_i\de p(i)\leq\int_0^\infty e^{-rt}\sum_i\de P_tf(i)\,\di t
\leq\big(\sum_i\de f(i)\big)\int_0^\infty e^{-rt}e^{K_\down t}\,\di t<\infty,
\]
which proves that $p_r\in\Ci_{\rm sum}$. This completes the proof that
$\Ci_{\rm sum}$ is a core for $H$, that is, the closure of
$G|_{\Ci_{\rm sum}}$ is $H$. By Lemma~\ref{L:locapp}, the closure of
$G|_{\Ci_{\rm fin}}$ contains $G|_{\Ci_{\rm sum}}$, so we see that
$\Ci_{\rm fin}$ is also a core for $H$.
\end{Proof}

The following lemma is sometimes useful.

\begin{lemma}[Differentiation of semigroup]
Assume\label{L:difP} that the rates $(r_m)_{m\in\Gi}$ satisfy (\ref{downsum}), let
$(P_t)_{t\geq 0}$ be the Feller semigroup defined in (\ref{Ppart}) and
let $G$ be the linear operator with domain $\Di(G):=\Ci_{\rm sum}(S^\La)$
defined by (\ref{Gmap}). Then, for each $f\in\Ci_{\rm sum}(S^\La)$,
$t\mapsto P_tf$ is a continuously differentiable function from $\half$
to $\Ci(S^\La)$ satisfying $P_0f=f$, $P_tf\in\Ci_{\rm sum}(S^\La)$,
and $\dif{t}P_tf=GP_tf=P_tGf$ for each $t\geq 0$.
\end{lemma}

\begin{Proof}
The statement that $\dif{t}P_tf=P_tGf$ holds for any Feller semigroup
and $f$ in the domain of its generator, see
\cite[Prop~1.1.5]{EK86}. The remaining statements are a direct
consequence of Proposition~\ref{P:Cauchy}, Lemma~\ref{L:Pfsum}, and
Theorem~\ref{T:genpart}. A direct proof based on our definition of
$(P_t)_{t\geq 0}$ (not using Hille--Yosida theory) is also possible,
but quite long and technical.
\end{Proof}

We conclude this section by proving an analogue of
Proposition~\ref{P:cov} for interacting particle systems. We continue
to assume that $S$ is a finite set, $\La$ is countable, $\Gi$ is a
collection of continuous maps $m\cn S^\La\to S^\La$, and
$(r_m)_{m\in\Gi}$ are nonnegative rates satisfying (\ref{downsum}). By
Theorem~\ref{T:genpart}, the linear operator $G$ with domain
$\Di(G):=\Ci_{\rm sum}$ defined in (\ref{Gmap}) is closable and its
closure $\ov G$ generates a Feller semigroup $(P_t)_{t\geq 0}$. We
need the following simple lemma.

\begin{lemma}[Closedness under multiplication]
For\label{L:prodclos} each $f,g\in\Ci_{\rm sum}$, the pointwise
product $fg$ is an element of $\Ci_{\rm sum}$.
\end{lemma}

\begin{Proof}
For each $i\in\La$ and $x,y\in S^\La$ such that $x(j)=y(j)$ for all
$j\neq i$, we can estimate
\[\begin{array}{r@{\,}c@{\,}l}
\dis\big|f(x)g(x)-f(y)g(y)\big|
&\leq&\dis\big|f(x)g(x)-f(y)g(x)\big|+\big|f(y)g(x)-f(y)g(y)\big|\\[5pt]
&\leq&\dis\|g\|_\infty\cdot\big|f(x)-f(y)\big|+\|f\|_\infty\big|g(x)-g(y)\big|,
\end{array}\]
which tells us that
\[
\de(fg)(i)\leq\|g\|_\infty\de f(i)+\|f\|_\infty\de g(i).
\]
Summing over $i$ yields the claim.
\end{Proof}

Lemma~\ref{L:prodclos} allows us to define
$\Ga_G\cn\Ci_{\rm sum}\times\Ci_{\rm sum}\to\Ci$ by
\[\index{0zzC@$\Ga_G(f,g)$}
\Ga_G(f,g):=G(fg)-(Gf)g-f(Gg)\qquad(f,g\in\Ci_{\rm sum}).
\]
A calculation similar to the one below (\ref{Gadef}) shows that
\begin{equation}\label{Gadef2}
\Ga_G(f,g)(x)=\sum_{m\in\Gi}r_m\big\{f\big(m(x)\big)-f\big(x\big)\big\}\big\{g\big(m(x)\big)-g\big(x\big)\big\}.
\end{equation}
The following proposition generalizes Proposition~\ref{P:cov} to
interacting particle systems.

\begin{proposition}[Covariance formula]
Assume\label{P:cov2} that the rates $(r_m)_{m\in\Gi}$ satisfy
(\ref{downsum}), let $(P_t)_{t\geq 0}$ be the Feller semigroup defined
in (\ref{Ppart}) and let $G$ be the linear operator with domain
$\Di(G):=\Ci_{\rm sum}$ defined in (\ref{Gmap}). Then for each
probability measure $\mu$ on $S^\La$, one has
\[
\cov_{\mu P_t}(f,g)=\cov_\mu(P_tf,P_tg)+\int_0^t\!\di s\,\mu P_{t-s}\Ga_G(P_sf,P_sg)\quad(f,g\in\Ci_{\rm sum}).
\]
\end{proposition}

\begin{Proof}
The proof is essentially the same as in the finite case
(Proposition~\ref{P:cov}). Using Lemmas \ref{L:difP} and
\ref{L:prodclos}, and the continuity of $P_{t_1}$, we obtain that
\[\begin{array}{r@{\,}c@{\,}l}
\dis\dif{t_1}P_{t_1}\big((P_{t_2}f)(P_{t_3}g)\big)&=&\dis P_{t_1}G\big((P_{t_2}f)(P_{t_3}g)\big),\\[5pt]
\dis\dif{t_2}P_{t_1}\big((P_{t_2}f)(P_{t_3}g)\big)&=&\dis P_{t_1}\big((GP_{t_2}f)(P_{t_3}g)\big),\\[5pt]
\dis\dif{t_3}P_{t_1}\big((P_{t_2}f)(P_{t_3}g)\big)&=&\dis P_{t_1}\big((P_{t_2}f)(GP_{t_3}g)\big).
\end{array}\]
The rest of the proof is the same.
\end{Proof}

\subsection*{Some bibliographical remarks}

In 1972, several authors published results of various degree of
generality showing that interacting particle systems on infinite
lattices are well-defined. Harris \cite{Har72} used the Poisson
approach. His result applies only to finite range interactions on
$\Z^d$. Instead of using the backward in time process he argued
forwards in time, using percolation theory to show that if $t$ is
small enough, then the lattice can randomly be divided into finite
pieces that mutually do not interact with each other during the time
interval $(0,t]$.

Liggett\index{Liggett} \cite{Lig72}, on the other hand, gave a direct
proof that the closure of $G$ generates a Feller semigroup
$(P_t)_{t\geq 0}$, and then invoked the abstract result
Theorem~\ref{T:Fellexist} about Feller processes to prove the
existence of a corresponding Markov process with cadlag sample
paths. This result is more widely applicable than Harris' result and
made it to Liggett's famous book \cite[Theorem~I.3.9]{Lig85}.
Liggett's conditions are similar to condition (\ref{downsum}) of our
Theorem~\ref{T:Poispart} but there are also some differences. Liggett
does not write his generators in terms of local maps, but in terms of
local probability kernels. This way of writing the generator is more
general and sometimes (for example for stochastic Ising models) more
natural than our approach using local maps. It is worth noting that
Liggett's construction, like ours, depends on a clever way of writing
the generator that is in general not unique.

Liggett's book \cite{Lig85} does not treat graphical representations
in the generality of our Theorem~\ref{T:Poispart} but he does use
explicit Poisson constructions for some specific systems, such as the
contact process. He does not actually prove that these Poisson
constructions yield the same process as the generator construction,
but apparently finds this self-evident. (Equivalence of the two
constructions follows from our Theorem~\ref{T:genpart} but
alternatively can also be proved by approximation with finite systems,
using approximation results such as \cite[Cor.~I.3.14]{Lig85}.)

Liggett's \cite[Theorem~I.3.9]{Lig85} allows for the case that the
local state space $S$ is a (not necessarily finite) compact metrizable
space. This is occasionally convenient. For example, this allows one
to construct voter models with infinitely many types, where at time
zero, the types $(X_0(i))_{i\in\La}$ are i.i.d.\ and uniformly
distributed on $S=[0,1]$. We have made essential use of the finiteness
of $S$ in several places. For example, the state space $\Ci(S^\La,T)$
of the backward in time process is no longer countable if $S$ is not
finite and, as explained above Exercise~\ref{E:discon}, solutions to
the evolution equation (\ref{Xevolve}) may no longer be unique if $S$
is allowed to be a general compact metrizable space. With some extra
work, these difficulties can presumably be overcome (for example by
requiring that solutions to (\ref{Xevolve}) are cadlag with respect to
the discrete topology on $S$) but for simplicity we restrict ourselves
to finite local state spaces. Non-compact local state spaces are more
tricky, see \cite[Chapter~IX]{Lig85}. An alternative treatment of
non-compact local state spaces, that works only for processes with
finite range interactions, is given in \cite{Pen08}.

The backward in time process of Proposition~\ref{P:back} will come
back in Chapter~\ref{C:dual} when we discuss duality of interacting
particle systems. It is also interesting to look at the mean-field
limit of this process. One can show that in the mean-field limit, the
process
\[
\big(\Ri(\Phi_t)\big)_{t\geq 0}\quad\mbox{with}\quad
\Phi_t:=\Fb_{u,u-t}(\Phi_0)\quad(t\geq 0)
\]
behaves as a branching process. As a result, solutions to the
mean-field equation can be represented in terms of a stochastic
process on the genealogical tree of a branching process. This is
explained in \cite{MSS20}.

\section{Ergodicity}\label{S:ergod}

The proofs of Theorems \ref{T:Poispart} and \ref{T:genpart} were quite
long. Luckily, they yield more information than just the fact that the
interacting particle systems we are interested in are
well-defined. The basic phenomenon that motivates the study of
interacting particle systems is collective behavior. The general
picture is that for weak strengths of the interaction, different parts
of space behave essentially independently, but for sufficiently strong
interaction it may happen that all sites start to coordinate their
behavior, giving rise to multiple invariant laws or even more exotic
phenomena such as periodic behavior.

As a result of the methods of the previous sections, we will be able
to prove results that confirm the ``easy'' part of this picture,
namely the absence of collective behavior for weak strengths of the
interaction.

If $X$ is a Markov process with state space $E$ and transition
probabilities $(P_t)_{t\geq 0}$, then by definition, an
\emph{invariant law}\index{invariant law} of $X$ is a probability
measure $\nu$ on $E$ such that
\[
\nu P_t=\nu\qquad(t\geq 0).
\]
This says that if we start the process in the initial law
$\P[X_0\in\,\cdot\,]=\nu$, then $\P[X_t\in\,\cdot\,]=\nu$ for all
$t\geq 0$. As a consequence, one can construct a stationary process
$(X_t)_{t\in\R}$ such that (compare (\ref{Mark2}))
\begin{equation}\label{statMark}
\P\big[X_u\in\,\cdot\,\big|\,(X_s)_{-\infty<s\leq t}\big]
=P_{u-t}(X_t,\,\cdot\,)\quad{\rm a.s.}
\qquad(t\leq u),
\end{equation}
and $\P[X_t\in\,\cdot\,]=\nu$ for all $t\in\R$. Conversely, the
existence of such a stationary Markov process implies that the law at
any time $\nu:=\P[X_t\in\,\cdot\,]$ must be an invariant law.

\begin{theorem}[Ergodicity]
Let\label{T:ergo} $X$ be an interacting particle system with state
space of the form $S^\La$ and generator $G$ of the form (\ref{Gmap}),
and assume that the rates $(r_m)_{m\in\Gi}$ satisfy
(\ref{downsum}). Let $T$ be a finite set with at least two
elements.\med

\noi
\textbf{\rm(a)} Assume that the constant $K_\down$ from (\ref{KKK}) satisfies
$K_\down<0$. Then the backward in time process satisfies
\begin{equation}\label{dudie}
\lim_{t\to-\infty}\big|\Ri(\Fb_{u,u-t}(\phi)\big|=0\quad{\rm a.s.}
\quad\big(u\in\R,\ \phi\in\Ci(S^\La,T)\big).
\end{equation}
\textbf{\rm(b)} Assume that the backward in time process satisfies
(\ref{dudie}). Then the interacting particle system $X$ has a unique
invariant law $\nu$, and
\begin{equation}\label{ergo}
\P^x\big[X_t\in\,\cdot\,\big]\Asto{t}\nu\qquad(x\in S^\La).
\end{equation}
Moreover, there exists an a.s.\ unique cadlag process $(X_t)_{t\in\R}$
such that
\begin{equation}\label{infevolve}
X_t=\mk^\om_t(X_{t-})\qquad(t\in\R),
\end{equation}
and $(X_t)_{t\in\R}$ is distributed as the stationary Markov process
corresponding to the invariant law $\nu$.
\end{theorem}

\begin{Proof}
Part~(a) is immediate from formula (\ref{Phiexp}) of
Proposition~\ref{P:back}. Let $\Psi_i$ denote the set of functions
$\phi\cn S^\La\to T$ that depend only on $x(i)$. Note that this set is
finite. By (\ref{FXrel}) and the assumption that $T$ has at least two
elements
\[
\Ri\big(\Xb_{s,t}[i]\big)\sub\bigcup_{\phi\in\Psi_i}\Ri\big(\Fb_{t,s}(\phi)\big)
\quad(s\leq t,\ i\in\La),
\]
so (\ref{dudie}) implies that
\begin{equation}\label{dudie2}
\lim_{s\to-\infty}\big|\Ri(\Xb_{s,t}[i])\big|=0\quad{\rm a.s.}
\quad\big(t\in\R,\ i\in\La\big).
\end{equation}
It follows from the definition of $(\Xb_{s,t})_{s\leq t}$ in
(\ref{partflow}) that the function $t\mapsto\Xb_{s,t}[i]$ jumps only
at times for which there exists a $(m,t)\in\om$ such that
$i\in\Di(m)$. Since this set is locally finite by (\ref{downsum}), we
can replace the order of the ``almost sure'' and ``for all $t$''
statements, that is, (\ref{dudie2}) holds almost surely for all $t$
simultaneously. Formula (\ref{dudie2}) says that for low enough $s$,
the function $\Xb_{s,t}[i]$ is constant, which implies that for each
$z\in S^\La$ the a.s.\ limit
\begin{equation}\label{statdef}
X_t(i):=\lim_{s\to-\infty}\Xb_{s,t}(z)(i)\qquad(i\in\La,\ t\in\R)
\end{equation}
exists and does not depend on the choice of the configuration $z\in
S^\La$. Using the continuity of $\Xb_{s,u}$ (which is proved in
Theorem~\ref{T:Poispart}) and the flow property, we see that a.s.
\[
\Xb_{t,u}(X_t)=\lim_{s\to-\infty}\Xb_{t,u}\circ\Xb_{s,t}(z)=X_u\qquad(t\leq u),
\]
which implies that $(X_t)_{t\in\R}$ solves (\ref{infevolve}). If
$(X'_t)_{t\in\R}$ is another solution, then for all $s$ low enough so
that $\Ri(\Xb_{s,t}[i])=\emptyset$,
\[
X'_t(i)=\Xb_{s,t}(X'_s)(i)=\Xb_{s,t}(z)(i)=X_t(i),
\]
which shows that solutions to (\ref{infevolve}) are unique.

We claim that $X=(X_t)_{t\in\R}$ is Markov with respect to the
transition probabilities $(P_t)_{t\geq 0}$ in the sense of
(\ref{statMark}). Indeed, for almost every trajectory
$(x_s)_{-\infty<s\leq t}$ with respect to the law of
$(X_s)_{-\infty<s\leq t}$, we have
\[\begin{array}{l}
\dis\P\big[X_u\in\,\cdot\,\big|
\,(X_s)_{-\infty<s\leq t}=(x_s)_{-\infty<s\leq t}\big]\\[5pt]
\dis\quad=\P\big[\lim_{s\to-\infty}\Xb_{t,u}\circ\Xb_{s,t}(z)\in\,\cdot\,\big|
\,(X_s)_{-\infty<s\leq t}=(x_s)_{-\infty<s\leq t}\big]\\[5pt]
\dis\quad\stackrel{1}{=}\P\big[\Xb_{t,u}(X_t)\in\,\cdot\,\big|
\,(X_s)_{-\infty<s\leq t}=(x_s)_{-\infty<s\leq t}\big]\\[5pt]
\dis\quad\stackrel{2}{=}\P\big[\Xb_{t,u}(x_t)\in\,\cdot\,\big|
\,(X_s)_{-\infty<s\leq t}=(x_s)_{-\infty<s\leq t}\big]\\[5pt]
\dis\quad\stackrel{3}{=}\P\big[\Xb_{t,u}(x_t)\in\,\cdot\,\big]
=P_{u-t}(x_t,\,\cdot\,),
\end{array}\]
where in step~1 we have used the continuity of the map $\Xb_{t,u}$, in
step~2 we have replaced $X_t$ by $x_t$, and in step~3 we have used
that the random variables $\Xb_{t,u}$ and $(X_s)_{-\infty<s\leq t}$
are independent, since they are functions of the restriction of the
Poisson set $\om$ to the disjoint sets $\Gi\times(t,u]$ and
$\Gi\times(-\infty,t]$, respectively. By the stationarity of the
stochastic flow,
\[
\nu:=\P\big[X_t\in\,\cdot\,]\qquad(t\in\R)
\]
does not depend on $t\in\R$, and since $X$ is Markov this defines an invariant
law $\nu$. Since
\[
\P^x\big[X_t\in\,\cdot\,\big]=\P\big[\Xb_{-t,0}(x)\in\,\cdot\,\big]
\]
and since by (\ref{statdef}), we have
\[
\Xb_{-t,0}(x)\asto{t}X_0\quad{\rm a.s.}\qquad(x\in S^\La)
\]
with respect to the topology of pointwise convergence, we conclude that
(\ref{ergo}) holds.
\end{Proof}

\noi
\textbf{Remark} It is possible for an interacting particle systems to
be ergodic in the sense of (\ref{ergo}) while (\ref{dudie}) does not
hold. In such a situation, it is not clear if solutions to
(\ref{infevolve}) are a.s.\ unique. Even if there are multiple
invariant laws, one can ask if (\ref{infevolve}) has an a.s.\ unique
solution subject to the condition that $(X_t)_{t\in\R}$ is stationary
with a given invariant law. Not much is known about this, but these
questions are related to the concept of endogeny of recursive tree
processes \cite{AB05,MSS20}.\med

We\label{ergdisc} note that (\ref{ergo}) says that if we start the
process in an arbitrary initial state $x$, then the law at time $t$
converges weakly\footnote{Here weak convergence is of course
w.r.t.\ our topology on $S^\La$, that is, w.r.t.\ the product
topology.} as $t\to\infty$ to the invariant law $\nu$. This property
is often described by saying that the interacting particle system is
\emph{ergodic}.\index{ergodic interacting particle system} Indeed,
this implies that the corresponding stationary process
$(X_t)_{t\in\R}$ is ergodic in the usual sense of that word, that is,
the \si-field of events that are invariant under translations in time
is trivial. The converse conclusion cannot be drawn, however, so the
traditional way of describing (\ref{ergo}) as ``ergodicity'' is a bit
of a bad habit.

We have split Theorem~\ref{T:ergo} into a part~(a) and (b) since the
condition (\ref{dudie}) is sometimes satisfied even when the constant
$K_\down$ from (\ref{KKK}) is positive. Indeed, we will later see that
for the contact process, the condition (\ref{dudie}) is sharp but the
condition $K_\down<0$ is not. In Exercise~\ref{E:context} below, we
will calculate the constant $K_\down$ for the contact process and
deduce that this process is ergodic for small values of the infection
rate.

Theorem~\ref{T:ergo} is similar, but not identical to
\cite[Thm~I.4.1]{Lig85}. For Theorem~\ref{T:ergo}~(a) and (b) to be
applicable, one needs to be able to express the generator in terms of local
maps such that the constant $K_\down$ from (\ref{KKK}) is negative.
For \cite[Thm~I.4.1]{Lig85}, one needs to express the generator in a
convenient way in terms of local transition kernels. For certain problems, the
latter approach is more natural and \cite[Thm~I.4.1]{Lig85} yields sharper
estimates for the regime where ergodicity holds.

\section{Application to the Ising model}\label{S:Isap}

The Ising model with Glauber dynamics has been introduced in
Section~\ref{S:IsingPotts}. So far, we have not shown how to represent
the generator of this interacting particle system in terms of local
maps. In the present section, we will fill this gap. We willl only
consider the ferromagnetic case $\bet\geq 0$. As an application of the
theory developed so far, we will then show that the Ising model with
Glauber dynamics is well-defined for all values of its parameter, and
ergodic for $\bet$ sufficiently small. Our construction will also
prepare for the next chapter, where we discuss monotone interacting
particle systems, by showing that the Ising model with Glauber
dynamics can be represented in monotone maps.

We recall from Section~\ref{S:IsingPotts} that the Ising model with Glauber
dynamics on a graph $(\La,E)$ is the interacting particle system with state
space $\{-1,+1\}^\La$ and dynamics such that
\[
\mbox{site $i$ flips to the value $\sig$ with rate}\quad
r^\sig_i(x):=\frac{e^{\bet N_{x,i}(\sig)}}{e^{\bet N_{x,i}(+1)}+e^{\bet N_{x,i}(-1)}},
\]
where
\[
N_{x,i}(\sig):=\sum_{j\in\Ni_i}1_{\txt\{x(j)=\sig\}}\qquad\big(\sig\in\{-1,+1\}\big)
\]
denotes the number of neighbors of $i$ that have the spin value
$\sig$. For each $i\in\La$, let $K^\bet_i$ denote the
probability kernel on $\{-1,+1\}^{\La}$ defined as
\[
K^\bet_i(x,y):=\left\{\begin{array}{ll}
\dis r^\sig_i(x)\quad&\dis\mbox{if }
y=m^\sig_i(x)\quad\big(\sig\in\{-1,+1\}\big),\\[5pt]
\dis 0\quad&\dis\mbox{otherwise,}
\end{array}\right.\]
where $m^\sig_i(x)$, defined in (\ref{msig}), denotes the
configuration $x$ with the spin at $i$ flipped to $\sig$. Then the
generator (\ref{GPotts}) of the Ising model takes the form
\begin{equation}\label{Isingkern}
G_{\rm Ising}f=\sum_{i\in\La}\big\{K^\bet_if-f\big\},
\end{equation}
which is an expression of the form (\ref{GK}) but not a random mapping
representation of the form (\ref{Gmap}). To find a random mapping
representation for $G_{\rm Ising}$ in terms of local maps as in
(\ref{Gmap}), it suffices to find a random mapping representation for
the kernels $K^\bet_i$. This needs some preparations. Let
\[\index{0Mxi@$M_{x,i}$}
M_{x,i}:=N_{x,i}(+)-N_{x,i}(-)=\sum_{j\in\Ni_i}x(j)
\]
denote the \emph{local magnetization} in the neighborhood $\Ni_i$ of $i$.
Since $N_{x,i}(+)+N_{x,i}(-)=|\Ni_i|$, we can rewrite the
probability under $K^\bet_i$ of flipping to the spin value $+1$ as
\[\begin{array}{l}
\dis r^+_i(x)=\frac{e^{\bet N_{x,i}(+1)}}{e^{\bet N_{x,i}(+1)}+e^{\bet N_{x,i}(-1)}}
=\frac{e^{\bet(|\Ni_i|+M_{x,i})/2}}
{e^{\bet(|\Ni_i|+M_{x,i})/2}+e^{\bet(|\Ni_i|-M_{x,i})/2}}\\[10pt]
\dis\quad=\frac{e^{\haa\bet M_{x,i}}}{e^{\haa\bet M_{x,i}}+e^{-\haa\bet M_{x,i}}}
=\ha\Big(1+\frac{e^{\haa\bet M_{x,i}}-e^{-\haa\bet M_{x,i}}}
{e^{\haa\bet M_{x,i}}+e^{-\haa\bet M_{x,i}}}\Big)\\[12pt]
\dis\quad=\ha\big(1+\tanh(\ha\bet M_{x,i})\big).
\end{array}\]
Similarly, the probability of flipping to $-1$ is
$r^-_i(x)=\ha(1-\tanh(\ha\bet M_{x,i}))=1-r^+_i(x)$.

\begin{figure}[htb]
\begin{center}
\inputtikz{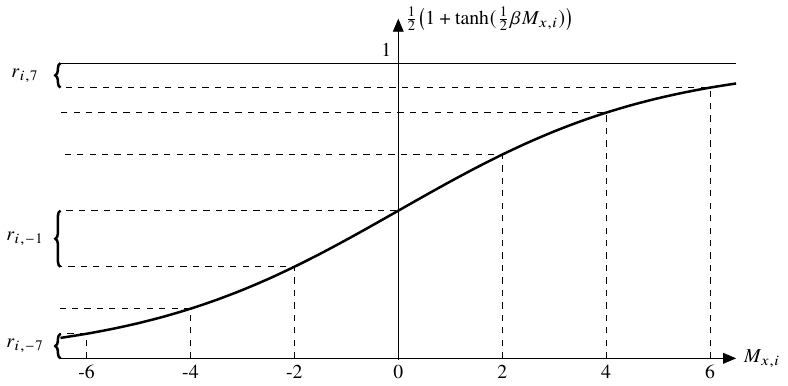}
\caption{Definition of the rates $r_{i,L}$ from (\ref{ril}). In this
  example $N=6$ and $\bet=0.4$.}
\label{fig:ril}
\commentAlt{Figure~\ref{fig:ril}}{Graph of the hyperbolic
  tangent. Dashed upward lines starting at the points $-6,-4,...,6$
  meet the graph at different heights, as indicated by horizontal
  lines. The rates are the distances between these lines.}
\end{center}
\end{figure}

For (mainly notational) simplicity, let us assume that each site $i$
has the same number of neighbors in the graph $(\La,E)$, so that the
size of the neighborhood
\[
N:=|\Ni_i|\qquad(i\in\La)
\]
does not depend on $i\in\La$. Then $M_{x,i}$ takes values in
$\{-N,-N+2,\ldots,N\}$. We observe that for $\bet>0$ the function
$z\mapsto\ha(1+\tanh(\ha\bet z))$ is increasing (see
Figure~\ref{fig:ril}). Inspired by this, for $L=-N-1,-N+1,\ldots,N+1$, we
define local maps $m_{i,L}$ by
\begin{equation}\label{mpmiL}
m_{i,L}(x)(j):=\left\{\begin{array}{ll}
+1\quad&\mbox{if }j=i\mbox{ and }M_{x,i}>L,\\[5pt]
-1\quad&\mbox{if }j=i\mbox{ and }M_{x,i}<L,\\[5pt]
x(j)\quad&\mbox{if }j\neq i.\end{array}\right.
\end{equation}
We try a generator of the form
\begin{equation}\label{Isingrep}
G_{\rm Ising}f(x)=\sum_{i\in\La}\sum_{L=-N-1}^{N+1}
r_{i,L}\big\{f\big(m_{i,L}(x)\big)-f\big(x\big)\big\},
\end{equation}
where we sum only over odd $L$ and the constants $r_{i,L}\geq 0$ are
probabilities that need to be chosen in such a way that
\begin{equation}\label{Krarep}
K^\bet_i(x,y)=\sum_{L=-N-1}^{N+1}r_{i,L}1_{\txt\{m_{i,L}(x)=y\}}
\end{equation}
is a random mapping representation of the kernel $K^\bet_i$. Consider
$x,y$ such that $x(i)=-1$, $y(i)=+1$, and $x(j)=y(j)$ for all $j\neq i$.
For such $x,y$, (\ref{Krarep}) yields the equation
\[
\ha\big(1+\tanh(\ha\bet M_{x,i})\big)=r^+_i(x)=\sum_{L=-N-1}^{M_{x,i}-1}r_{i,L}
\]
Similarly, for $x,y$ such that $x(i)=+1$, $y(i)=-1$, and $x(j)=y(j)$ for all
$j\neq i$, formula (\ref{Krarep}) yields
\[
\ha\big(1-\tanh(\ha\bet M_{x,i})\big)=r^-_i(x)=\sum_{L=M_{x,i}+1}^{N+1}r_{i,L}.
\]
From this, we see that (\ref{Krarep}) is satisfied for
(see Figure~\ref{fig:ril})
\begin{equation}\label{ril}
r_{i,L}:=\left\{\begin{array}{ll}
\ha\big(1+\tanh(-\ha\bet N)\big)\quad&\mbox{if }L=-N-1,\\[5pt]
\ha\tanh(\ha\bet(L+1))-\ha\tanh(\ha\bet(L-1))\quad
&\mbox{if }-N+1\leq L\leq N-1,\\[5pt]
\ha\big(1-\tanh(\ha\bet N)\big)\quad&\mbox{if }L=N+1,
\end{array}\right.
\end{equation}
which has the effect that the generator in (\ref{Isingrep}) equals the
one in (\ref{Isingkern}). We observe that this is even true for
$\bet=0$: in this case, $r_{i,-N-1}=\ha=r_{i,N+1}$ and all other
probabilities are zero.

\begin{theorem}[Existence and ergodicity of the Ising model]
Consider\label{T:Iserg} an Ising model with Glauber dynamics on a
countable graph $\La$ in which each lattice point $i$ has exactly
$|\Ni_i|=N\geq 2$ neighbors, that is, the Markov process $X$ with state
space $\{-1,+1\}^\La$ and generator $G_{\rm Ising}$ given by
(\ref{Isingkern}).  Then, for each $\bet\geq 0$, the closure of
$G_{\rm Ising}$ generates a Feller semigroup. Moreover, for each
\begin{equation}\label{ergcond}
0\leq\bet<N^{-1}\big(\log(N+1)-\log(N-1)\big),
\end{equation}
the Markov process with generator $\ov G_{\rm Ising}$ has a unique
invariant law $\nu$, and the process started in an arbitrary initial
state $x$ satisfies
\[
\P^x\big[X_t\in\,\cdot\,\big]\Asto{t}\nu\qquad\big(x\in\{-1,+1\}^\La\big).
\]
\end{theorem}

For the critical point of the Ising model on $\Z^2$,
Theorem~\ref{T:Iserg} yields the bound
\[
0.1277\approx\ffrac{1}{4}(\log 5-\log 3)\leq\bet_{\rm c}
\]
which should be compared with the known value
$\bet_{\rm c}=\log(1+\sqrt 2)\approx 0.8814$ from (\ref{Onsager}).\med

\begin{Proof}[of Theorem~\ref{T:Iserg}]
We use the representation (\ref{Isingrep}). We observe that $\Di(m_{i,L})=\{i\}$
is the set of lattice points whose spin value can be changed by the
map $m_{i,L}$. The set of lattice points that are
$m_{i,L}$-relevant for $i$ is given by
\[
\Ri^\down_i(m_{i,L})=\left\{\begin{array}{ll}
\emptyset\quad&\mbox{if }L=-N-1\quad\mbox{or}\quad L=N+1,\\[5pt]
\Ni_i\quad&\mbox{otherwise.}
\end{array}\right.
\]
Here we have used that $-N-1<M_{x,i}\leq N+1$ holds always, so
$m_{-N-1}(x)(i)=+1$ and $m_{N+1}(x)(i)=-1$ regardless of what $x$ is. On
the other hand, in all other cases, the value of each lattice point
$j\in\Ni_i$ can potentially make a difference for the outcome
$m_{i,L}(x)(i)$.

By Theorem~\ref{T:genpart}, to conclude that the closure of
$G_{\rm Ising}$ generates a Feller semigroup, it suffices to check
that the expressions in (\ref{downsum}) (i) and (ii) are finite. In
our case, these expressions are
\[
{\rm(i)}\ \sup_{i\in\La}\sum_{L=-N-1}^{N+1}r_{i,L}=1
\qquad
{\rm(ii)}\ \sup_{i\in\La}\sum_{L=-N+1}^{N-1}r_{i,L}\big|\Ni_i\big|\leq N.
\]
To prove ergodicity for $\bet$ small enough, we apply
Theorem~\ref{T:ergo}. We calculate the constant $K_\down$ from
(\ref{KKK}).
\[\begin{array}{r@{\,}c@{\,}l}
\dis K_\down
&=&\dis\sum_{L=-N-1}^{N+1}r_{i,L}\big(|\Ri^\down_i(m_{i,L})|-1\big)
=-r_{i,-N-1}-r_{i,N+1}+\sum_{L=-N+1}^{N-1}r_{i,L}\big(N-1\big)\\[5pt]
&=&\dis-1+N\sum_{L=-N+1}^{N-1}r_{i,L}=N\tanh(\ha\bet N)-1.
\end{array}\]
This is negative if and only if
\[\begin{array}{l}
\dis
N\frac{e^{\ha\bet N}-e^{-\ha\bet N}}{e^{\ha\bet N}+e^{-\ha\bet N}}<1\quad\desd\quad
N\big(e^{\ha\bet N}-e^{-\ha\bet N}\big)<e^{\ha\bet N}+e^{-\ha\bet N}\\[5pt]
\dis\quad\desd\quad
(N-1)e^{\ha\bet N}<(N+1)e^{-\ha\bet N}\quad\desd\quad
e^{\bet N}<\frac{N+1}{N-1},
\end{array}\]
which yields condition (\ref{ergcond}).
\end{Proof}

\begin{Exercise}
Show that the antiferromagnetic Ising model is ergodic if
\[
-N^{-1}\big(\log(N+1)-\log(N-1)\big)<\bet\leq 0.
\]
\end{Exercise}


\section{Further results}

In the present section we collect a number of technical results of a general
nature that will be needed in later chapters. On a first reading, readers are
advised to skip the present section and refer back to specific results when
the need arises. The only result of the present section that is perhaps of
some intrinsic value is Theorem~\ref{T:Fellim} which together with
Corollary~\ref{C:Partlim} below implies that the transition probabilities of
interacting particle systems on infinite lattices can be approximated by those
on finite lattices, something that we have been using implicitly when doing
simulations. An alternative way to see this is to use
Proposition~\ref{P:fincon} that we have already proved.

Let $E$ be a compact metrizable space. By definition, a collection of
functions $\Hi\sub\Ci(E)$ is \emph{distribution determining}
\index{distribution determining} if for probability measures $\mu,\nu$ on $E$
\[
\mu f=\nu f\ \forall f\in\Hi\quad\mbox{implies}\quad\mu=\nu.
\]
We say that $\Hi$ \emph{separates points}\index{separation of points} if for
all $x,y\in E$ such that $x\neq y$, there exists an $f\in\Hi$ such that
$f(x)\neq f(y)$. We say that $\Hi$ is \emph{closed under products} if
$f,g\in\Hi$ implies $fg\in\Hi$.

\begin{lemma}[Application of Stone--Weierstrass]\label{L:SW}
Let $E$ be a compact metrizable space.\index{Stone--Weierstrass theorem}
Assume that $\Hi\sub\Ci(E)$ separates points and is closed under
products. Then $\Hi$ is distribution determining.
\end{lemma}

\begin{Proof}
If $\mu f=\nu f$ for all $f\in\Hi$, then we can add the constant function $1$
to $\Hi$ and retain this property. In a next step, we can add all linear
combinations of functions in $\Hi$ to the set $\Hi$; by the linearity of the
integral, it will then still be true that $\mu f=\nu f$ for all $f\in\Hi$.
But now $\Hi$ is an algebra that separates points and vanishes nowhere, so by
the Stone--Weierstrass theorem, $\Hi$ is dense in $\Ci(E)$. If $f_n\in\Hi$,
$f\in\Ci(E)$, and $\|f_n-f\|_\infty\to 0$, then $\mu f_n\to\mu f$ and likewise
for $\nu$, so we conclude that $\mu f=\nu f$ for all $f\in\Ci(E)$. If $A\sub
E$ is a closed set, then the function $f(x):=d(x,A)$ is continuous, where
$d$ is a metric generating the topology on $E$ and $d(x,A):=\inf_{y\in
  A}d(x,y)$ denotes the distance of $x$ to $A$. Now the functions
$f_n:=1\wedge nf$ are also continuous and $f_n\up 1_{A^{\rm c}}$, so by the
continuity of the integral with respect to increasing sequences we see that
$\mu(O)=\nu(O)$ for every open set $O\sub E$. Since the open sets are closed
under intersections, it follows that $\mu(A)=\nu(A)$ for every element $A$ of
the \si-algebra generated by the open sets, that is, the Borel-\si-field
$\Bi(E)$.
\end{Proof}

\begin{lemma}[Weak convergence]\label{L:weaksuf}
Let $E$ be a compact metrizable space. Assume that $\mu_n\in\Mi_1(E)$ have the
property that $\lim_{n\to\infty}\mu_n f$ exists for all $f\in\Hi$, where
$\Hi\sub\Ci(E)$ is distribution determining. Then there exists a
$\mu\in\Mi_1(E)$ such that $\mu_n\Rightarrow\mu$.
\end{lemma}

\begin{Proof}
By Prohorov's theorem, the space $\Mi_1(E)$, equipped with the topology of
weak convergence, is compact. Therefore, to prove the statement, it suffices
to show that the sequence $\mu_n$ has not more than one cluster point, that is,
it suffices to show that if $\mu,\mu'$ are subsequential limits, then
$\mu'=\mu$. Clearly, $\mu,\mu'$ must satisfy $\mu'f=\mu f$ for all $f\in\Hi$,
so the claim follows from the assumption that $\Hi$ is distribution
determining.
\end{Proof}

\begin{lemma}[Continuous probability kernels]
Let\label{L:muKf} $E$ be a compact metrizable space and let $K$ be a
continuous probability kernel on $E$. Then, for any
$\mu_n,\mu\in\Mi_1(E)$ and $f_n,f\in\Ci(E)$,
\[\begin{array}{l}
\dis\phantom{\quand}
\mu_n\Asto{n}\mu\quad\mbox{implies}\quad\mu_nK\Asto{n}\mu K\\[5pt]
\dis\quand
\|f_n-f\|_\infty\asto{n}0\quad\mbox{implies}\quad\|Kf_n-Kf\|_\infty\asto{n}0.
\end{array}\]
\end{lemma}

\begin{Proof}
Since $K$ is a continuous probability kernel, its associated operator maps the
space $\Ci(E)$ into itself, so $\mu_n\Rightarrow\mu$ implies that
$\mu_n(Kf)\Rightarrow\mu(Kf)$ for all $f\in\Ci(E)$, or equivalently
$(\mu_nK)f\Rightarrow(\mu K)f$ for all $f\in\Ci(E)$, that is, the measures
$\mu_nK$ converge weakly to $\mu K$.

The second statement follows from the linearity and monotonicity of $K$ and
the fact that $K1=1$, which together imply that
$\|Kf_n-Kf\|_\infty\leq\|f_n-f\|_\infty$.
\end{Proof}

\begin{lemma}[Long-time limits]\label{L:longtime}
Let $E$ be a compact metrizable space and let $(P_t)_{t\geq 0}$ be the
transition probabilities of a Feller process in $E$. Let $\mu\in\Mi_1(E)$ and
assume that
\[
\mu P_t\Asto{t}\nu
\]
for some $\nu\in\Mi_1(E)$. Then $\nu$ is an invariant law of the Feller
process with transition probabilities $(P_t)_{t\geq 0}$.
\end{lemma}

\begin{Proof}
Using Lemma~\ref{L:muKf}, this follows by writing
\[
\nu P_t=(\lim_{s\to\infty}\mu P_s)P_t
=\lim_{s\to\infty}\mu P_sP_t=\lim_{s\to\infty}\mu P_{s+t}=\nu.
\]
\end{Proof}

The following theorem follows from \cite[Thm~17.25]{Kal97}, where it is
moreover shown that the condition (\ref{Gconv}) implies convergence in
distribution of the associated Feller processes, viewed as random variables
taking values in the space $\Di_E\half$ of cadlag paths with values in $E$.
Note that in (\ref{Gconv}) below, $\to$ (of course) means convergence in the
topology we have defined on $\Ci(E)$, that is, convergence w.r.t.\ the
supremum-norm.

\begin{theorem}[Limits of semigroups]
Let\label{T:Fellim} $E$ be a compact metrizable space and let $G_n,G$
be generators of Feller processes in $E$. Assume that there exists a
linear operator $A$ on $\Ci(E)$ such that $\ov A=G$ and
\begin{equation}\label{Gconv}
\forall f\in\Di(A)\ \exists f_n\in\Di(G_n)\mbox{ such that }f_n\to f
\quand G_nf_n\to Af.
\end{equation}
Then the Feller semigroups $(P^n_t)_{t\geq 0}$ and $(P_t)_{t\geq 0}$ with
generators $G_n$ and $G$, respectively, satisfy
\[
\sup_{t\in[0,T]}\|P^n_tf-P_tf\|_\infty\asto{n}0\quad\big(f\in\Ci(E),\ T<\infty\big).
\]
Moreover, if $\mu_n,\mu\in\Mi_1(E)$, then
\[
\mu_n\Asto{n}\mu
\quad\mbox{implies}\quad
\mu_nP^n_t\Asto{n}\mu P_t\qquad(t\geq 0).
\]
\end{theorem}

We note that in the case of interacting particle systems,
Theorem~\ref{T:genpart} implies the following.

\begin{corollary}[Convergence of particle systems]
Let\label{C:Partlim} $S$ be a finite set and let $\La$ be
countable. Let $G_n,G$ be generators of interacting particle systems
in $S^\La$ and assume that $G_n,G$ can be written in the form
(\ref{Gmap}) with rates satisfying (\ref{downsum}). Assume moreover
that
\[
\|G_nf-Gf\|_\infty\asto{n}0\qquad\big(f\in\Ci_{\rm fin}(S^\La)\big).
\]
Then the generators $G_n,G$ satisfy (\ref{Gconv}) with $A$ the
restriction of $G$ to $\Ci_{\rm fin}(S^\La)$.
\end{corollary}

Theorem~\ref{T:Fellim} has the following useful consequence.

\begin{proposition}[Limits of invariant laws]
Let\label{P:invlim} $E$ be a compact metrizable space and let $G_n,G$
be generators of Feller processes in $E$ satisfying (\ref{Gconv}). Let
$\nu_n,\nu\in\Mi_1(E)$ and assume that for each $n$, the measure
$\nu_n$ is an invariant law of the Feller process with generator
$G_n$. Then $\nu_n\Rightarrow\nu$ implies that $\nu$ is an invariant
law of the Feller process with generator $G$.
\end{proposition}

\begin{Proof}
Using Theorem~\ref{T:Fellim}, this follows simply by observing that
\[
\nu P_t=\lim_{n\to\infty}\nu_nP^n_t=\lim_{n\to\infty}\nu_n=\nu
\]
for each $t\geq 0$.
\end{Proof}

\chapter{Monotonicity}\label{C:monot}

\section{The stochastic order}

If the local state space $S$ of an interacting particle system is
partially ordered, then we equip the product space $S^\La$ with the
\emph{product order}
\[
x\leq y\quad\mbox{iff}\quad x(i)\leq y(i)\ \forall i\in\La.
\]
Many well-known interacting particle systems use the local state space
$S=\{0,1\}$, which is of course equipped with a natural order $0\leq
1$. Often, it is often useful to prove comparison results, that say
that two interacting particle systems $X$ and $Y$ can be coupled in
such a way that $X_t\leq Y_t$ for all $t\geq 0$. Here $X$ and $Y$ may
be different systems, started in the same initial state, or also two
copies of the same interacting particle system, started in initial
states such that $X_0\leq Y_0$. A useful tool in such comparison
arguments is the stochastic order, which is the subject of the present
section. We will come back to interacting particle systems in the next
section.

We recall that if $S$ and $T$ are partially ordered sets, then a
function $f\cn S\to T$ is called
\emph{monotone}\index{monotone!function} iff $x\leq y$ implies
$f(x)\leq f(y)$. In particular, this definition also applies to
real-valued functions (where we equip $\R$ with the well-known
order). Throughout this section, $E$ is a compact metrizable space
that is equipped with a partial order $\leq$ that is \emph{compatible
with the topology}\index{compatible with the topology, partial order}
in the sense that
\[
\big\{(x,y)\in E^2:x\leq y\big\}\mbox{ is closed in the product topology on }E^2.
\]
We recall that $B(E)$ and $\Ci(E)$ denote the spaces of Borel
measurable bounded functions and continuous functions $f\cn E\to\R$,
respectively. We set
\[\index{0BpE@$B^+(E)$}\index{0CaE@${\cal C}^+(E)$}
B^+(E):=\big\{f\in B(E):f\mbox{ is monotone}\big\}
\quand
\Ci^+(E):=B^+(E)\cap\Ci(E).
\]
We need the following technical result.

\begin{lemma}[Distribution determining property]
If\label{L:mondense} $\mu,\nu$ are probability measures on $E$ such
that $\int\mu(\di x)f(x)=\int\nu(\di x)f(x)$ for all $f\in\Ci^+(E)$,
then $\mu=\nu$.
\end{lemma}

\begin{Proof}
Let $\Fi:=\{f\in\Ci^+(\Xc):f\geq 0\}$. By Lemma~\ref{L:SW} it suffices
to show that $\Fi$ is closed under products in the sense that
$f,g\in\Fi$ imply $fg\in\Fi$, and separates points in the sense that
for each $x,y\in\Xc$ with $x\neq y$, there exists an $f\in\Fi$ such
that $f(x)\neq f(y)$. Closedness under products is trivial. Showing
that $\Fi$ separates points is considerably more work. Assume that
$x,y\in E$ satisfy $x\neq y$. Then either $x\not\leq y$ or $y\not\leq
x$. By symmetry we may assume hat we are in the second case. Let
$E_0:=\{z\in E:z\leq x\}$ and $E_1:=\{z\in E:z\geq y\}$. Then $E_0\cap
E_1=\emptyset$. Since the partial order is compatible with the
topology, $E_0$ and $E_1$ are closed subsets of $E$. Our assumptions
on $E$ imply that it is a ``compact ordered space'' as defined in
\cite[Section~1.4]{Nac65}, which by the corollary to
\cite[Theorem~I.3.4]{Nac65} implies that $E$ is a ``normally ordered
space'' as defined in \cite[Section~1.2]{Nac65}. We can then apply
\cite[Theorem~I.2.1]{Nac65} to conclude that there exists a continuous
monotone function $f\cn E\to[0,1]$ such that $f=0$ on $E_0$ and $f=1$
on $E_1$.
\end{Proof}

The following theorem gives necessary and sufficient conditions for it
to be possible to couple two random variables $X$ and $Y$ with values
in $E$ such that $X\leq Y$. A \emph{coupling} \index{coupling} of two
random variables $X$ and $Y$, in the most general sense of the word,
is \emph{a way to construct $X$ and $Y$ together on one underlying
probability space $(\Om,\Fi,\P)$}. More precisely, if $X$ and $Y$ are
random variables defined on different underlying probability spaces,
then a coupling of $X$ and $Y$ is a pair of random variables $(X',Y')$
defined on one underlying probability space $(\Om,\Fi,\P)$, such that
$X'$ is equally distributed with $X$ and $Y'$ is equally distributed
with $Y$. Equivalently, since the laws of $X$ and $Y$ are all we
really care about, we may say that a \emph{coupling} of two
probability laws $\mu,\nu$ defined on measurable spaces $(E,\Ei)$ and
$(F,\Fi)$, respectively, is a probability measure $\rho$ on the
product space $(E\times F,\Ei\otimes\Fi)$ such that the first marginal
of $\rho$ is $\mu$ and its second marginal is $\nu$. If two
probability laws $\mu,\nu$ satisfy the equivalent conditions of the
following theorem, then we say that $\mu$ and $\nu$
are \emph{stochastically ordered}\index{stochastic order} and we
write\footnote{This notation may look a bit confusing at first sight,
since, if $\mu,\nu$ are probability measures on a measurable space
$(\Om,\Fi)$, then one might interpret $\mu\leq\nu$ in a pointwise
sense, that is, in the sense that $\mu(A)\leq\nu(A)$ for all
$A\in\Fi$. In practice, this does not lead to confusion, since a
pointwise inequality for probability measures is a very uninteresting
property. Indeed, it is easy to check that probability measures
$\mu,\nu$ satisfy $\mu\leq\nu$ in a pointwise sense if and only if
$\mu=\nu$.} $\mu\leq\nu$.

\begin{theorem}[Stochastic order]
Let\label{T:stochord} $E$ be a compact metrizable space that is
equipped with a partial order that is compatible with the topology,
and let $\mu,\nu$ be probability laws on $E$. Then the following
statements are equivalent:
\begin{enumerate}[(iii)]
\item $\int\mu(\di x)f(x)\leq\int\nu(\di x)f(x)\ \forall f\in\Ci^+(E)$,
\item $\int\mu(\di x)f(x)\leq\int\nu(\di x)f(x)\ \forall f\in B^+(E)$,
\item It is possible to couple random variables $X,Y$ with laws
 $\mu=P[X\in\cdot\,]$ and $\nu=P[Y\in\cdot\,]$ in such a way that $X\leq Y$.
\end{enumerate}
Moreover, setting $\mu\leq\nu$ if and only if these conditions are
satisfied defines a partial order on the space of probability measures
on $E$.
\end{theorem}

\begin{Proof}
The implication (iii)$\volgt$(ii) is easy: if $X$ and $Y$ are coupled
such that $X\leq Y$ and $f$ is monotone, then
\[
\E\big[f(Y)\big]-\E\big[f(X)\big]
=\E\big[f(Y)-f(X)\big]\geq 0,
\]
since $f(Y)-f(X)\geq 0$ a.s. The implication (ii)$\volgt$(i) is
trivial. For the nontrivial implication (i)$\volgt$(iii) we refer to
\cite[Theorem~II.2.4]{Lig85}. For finite spaces, a nice intuitive
proof based on the max flow min cut theorem can be found in
\cite{Pre74}. To see that (i)--(iii) defines a partial order on the
space of probability measures on $E$ we must check that
1.\ $\mu\leq\mu$, 2.\ $\mu\leq\nu$ and $\nu\leq\mu$ imply $\mu=\nu$,
and 3.\ $\mu\leq\nu\leq\rho$ implies $\mu\leq\rho$. Properties 1 and 3
are immediate from condition~(i). Property~2 follows by combining (i)
with Lemma~\ref{L:mondense}.
\end{Proof}

Sometimes it is more convenient (or intuitive) to work with events
than with real functions. A set $A$ is called
\emph{increasing}\index{increasing event} if its indicator function
$1_A$ is monotone.

\begin{lemma}[Increasing events]
Let\label{L:incr} $E$ be a compact metrizable space that is equipped
with a partial order that is compatible with the topology and let
$\mu,\nu$ be probability measures on $E$. Then $\mu\leq\nu$ if and
only if
\begin{equation}\label{incr}
\mu(A)\leq\nu(A)\quad\mbox{for all closed increasing }A\sub E.
\end{equation}
\end{lemma}

\begin{Proof}
Condition~(ii) of Theorem~\ref{T:stochord} clearly implies
(\ref{incr}). To prove the converse, by condition~(i) of
Theorem~\ref{T:stochord}, it suffices to show that (\ref{incr})
implies that $\int\mu(\di x)f(x)\leq\int\nu(\di x)f(x)$ for all
$f\in\Ci^+(E)$. Fix $f\in\Ci^+(E)$. By adding a constant and
multiplying with a positive constant we can without loss of generality
assume that $f$ takes values in $[0,1]$. Define sets of dyadic
rationals by $D_n:=\{k2^{-n}:0\leq k\leq 2^n\}$ $(n\geq 1)$ and set
$f_n(x):=\sup\{d\in D_n:d\leq f(x)\}$. Then
\[
f_n=2^{-n}\sum_{k=1}^{2^n}1_{A_{n,k}}\quad\mbox{with}\quad
A_{n,k}:=\big\{x\in E:k2^{-n}\leq f(x)\big\}.
\]
The sets $A_{n,k}$ are closed and increasing, so (\ref{incr}) implies
that $\int\mu(\di x)f_n(x)\leq\int\nu(\di x)f_n(x)$ for all $n$. Since
$f_n\up f$, the claim follows.
\end{Proof}

\begin{Exercise}
Let $\Ci^+(\R^d)$ denote the space of bounded continuous functions
$f\cn\R^d\to\R$ that are monotone with respect to the product order on
$\R^d$. It is well-known \cite{KKO77} that for two probability
measures $\mu,\nu$ on $\R^d$ the following statements are equivalent:
(i) $\int\mu(\di x)f(x)\leq\int\nu(\di x)f(x)$ for all $f\in\Ci^+(E)$,
(ii) it is possible to couple random variables $X,Y$ with laws
$\mu=P[X\in\cdot\,]$ and $\nu=P[Y\in\cdot\,]$ in such a way that
$X\leq Y$. Give a proof of this fact using
Theorem~\ref{T:stochord}. \emph{Hint:} compactify.
\end{Exercise}

\section{Monotone interacting particle systems}\label{S:monsys}

In this section we specialize to spaces of the form $S^\La$ where $S$
is a finite partially ordered set and $\La$ is countable. In
particular, since $\La$ can be a set with only one element, this
includes arbitrary finite partially ordered sets. We equip $S^\La$
with the product topology and product partial order. We start with a
simple observation.

\begin{lemma}[Compatibility of the product order]
Let\label{L:prodcomp} $S$ be a finite partially ordered set and let
$\La$ be countable. Then the product partial order on $S^\La$ is
compatible with the product topology on $S^\La$ .
\end{lemma}

\begin{Proof}
Assume that $x_n,y_n,x,y\in S^\La$ satisfy $x_n\leq y_n$ for all $n$
and $x_n\to x$ and $y_n\to y$ in the product topology. Then for each
$i\in\La$, there exists an $N$ such that $x_n(i)=x(i)$ and
$y_n(i)=y(i)$ for all $n\geq N$ and hence $x(i)\leq y(i)$.
\end{Proof}

Because of Lemma~\ref{L:prodcomp}, Theorem~\ref{T:stochord} is
applicable with $E=S^\La$. The following lemma shows that for $E$ of
this form, in condition~(i) of Theorem~\ref{T:stochord}, we can
replace $\Ci^+(S^\La)$ by
$\Ci^+_{\rm fin}(S^\La):=\Ci^+(S^\La)\cap\Ci_{\rm fin}(S^\La)$, the space of
monotone functions $f\cn S^\La\to\R$ that depend on finitely many
coordinates.

\begin{lemma}[Local monotone functions]
The\label{L:locmon} space $\Ci^+_{\rm fin}(S^\La)$ is dense in $\Ci^+(S^\La)$.
\end{lemma}

\begin{Proof}
Fix $z\in S^\La$, choose finite $\La_n\up\La$, and for each $z\in
S^\La$ define $x^n(i):=x(i)$ if $i\in\La_n$ and $:=z(i)$
otherwise. Fix $f\in\Ci^+(S^\La)$ and define $f_n(x):=f(x^n)$. Then
clearly $f_n\in\Ci^+_{\rm fin}(S^\La)$ for all $n$. The argument used
at the beginning of the proof of Lemma~\ref{L:locapp} shows that
$\|f_n-f\|_\infty\to 0$, completing the proof.
\end{Proof}

By definition, a probability kernel $K$ on $S^\La$ is
\emph{monotone}\index{monotone!probability kernel} if it satisfies the
following equivalent conditions. Note that in (i) below, $\leq$
denotes the stochastic order. The equivalence of (i)--(iii) is a
trivial consequence of Theorem~\ref{T:stochord}.
\begin{enumerate}[(iii)]
\item $K(x,\,\cdot\,)\leq K(y,\,\cdot\,)$ for all $x\leq y$.
\item $Kf\in B^+(S^\La)$ for all $f\in\Ci^+(S^\La)$.
\item $Kf\in B^+(S^\La)$ for all $f\in B^+(S^\La)$.
\end{enumerate}
We note that if $K$ is monotone, then
\begin{equation}\label{ordpres}
\mu\leq\nu\quad\mbox{implies}\quad\mu K\leq\nu K.
\end{equation}
Indeed, this follows from (iii) since $f\in B^+(S^\La)$ implies $Kf\in
B^+(S^\La)$ and hence $\mu Kf\leq\nu Kf$ since $\mu\leq\nu$.

Recall from (\ref{Krep}) that a \emph{random mapping
representation}\index{random mapping representation!of a probability kernel}
of a probability kernel $K$ is a random map $M$ such that
\begin{equation}\label{Krep2}
K(x,\,\cdot\,)=\P[M(x)\in\,\cdot\,]\quad\forall x.
\end{equation}
We say that $K$ can be represented in the class of monotone maps, or
that $K$ is \emph{monotonically representable},
\index{monotone representability} if there exists a random monotone map $M$ such
that (\ref{Krep2}) holds. In Chapter~\ref{C:construct} we based our
construction of an interacting particle system on a random mapping
representation of its generator $G$ in terms of continuous maps, of the
form
\begin{equation}\label{Gdef4}
Gf(x)=\sum_{m\in\Gi}r_m\big\{f\big(m(x)\big)-f\big(x\big)\big\},
\end{equation}
where the rates satisfy (\ref{downsum}) or possibly the weaker
conditions from Theorem~\ref{T:graph}. If there exists such a random
mapping representation for which all maps $m\in\Gi$ are monotone, then
we say that $G$ is \emph{monotonically representable}.

\begin{lemma}[Monotone representability]
Each monotonically represent\-able probability kernel is monotone.
If the generator of an interacting particle system is monotonically
representable, then, for each $t\geq 0$, the transition probability $P_t$ is a
monotonically representable probability kernel.
\end{lemma}

\begin{Proof}
If a probability kernel $K$ can be written in the form (\ref{Krep2}) with $M$ a
random monotone map, then for each $x\leq y$, the random variables $M(x)$ and
$M(y)$ are coupled such that $M(x)\leq M(y)$ a.s., so their laws are
stochastically ordered as $K(x,\,\cdot\,)\leq K(y,\,\cdot\,)$. Since this
holds for all $x\leq y$, the kernel $K$ is monotone.

Given a random mapping representation of the form (\ref{Gdef4}) of the
generator $G$ of an interacting particle system, we can construct a
stochastic flow $(\Xb_{s,t})_{s\leq t}$ as in Theorem~\ref{T:Poispart}
based on a graphical representation $\om$. If all maps $m\in\Gi$ are
monotone, then for each finite $\om'\sub\om$, the maps
$(\Xb^{\om'}_{s,t})_{s\leq t}$ defined in (\ref{Xom}) are also
monotone, since they are the concatenation of finitely many maps from
$\Gi$. By Proposition~\ref{P:fincon}, this implies that the maps
$(\Xb_{s,t})_{s\leq t}$ are also monotone. It follows that
\[
P_t(x,\,\cdot\,)=\P\big[\Xb_{0,t}(x)\in\cdot\,\big]
\]
is a representation of $P_t$ in terms of the random monotone map $\Xb_{0,t}$,
so $P_t$ is monotonically representable.
\end{Proof}

We say that an interacting particle system is \emph{monotone}
\index{monotone!interacting particle system} if its transition kernels are
monotone probability kernels, and we say that it is \emph{monotonically
  representable} \index{monotone representability} if its generator is
monotonically representable. Somewhat surprisingly, it turns out that for
probability kernels, ``monotonically representable'' is a strictly stronger
concept than being ``monotone''. See \cite{FM01} for an example of a
probability kernel on $\{0,1\}^2$ that is monotone but not monotonically
representable. Nevertheless, it turns out that (almost) all monotone
interacting particle systems that one encounters in practice are also
monotonically representable.

The following maps are examples of monotone maps:
\begin{itemize}
\item The voter map ${\tt vot}_{ij}$ defined in (\ref{votmap}).
\item The branching map ${\tt bra}_{ij}$ defined in (\ref{bramap}).
\item The death map ${\tt death}_i$ defined in (\ref{deathmap}).
\item The asymmetric exclusion map ${\tt asep}_{ij}$ defined in (\ref{asep}).
\item The exclusion map ${\tt excl}_{ij}$ defined in (\ref{exclmap}).
\item The coalescing random walk map ${\tt rw}_{ij}$ defined in (\ref{rwmap}).
\item The cooperative branching map ${\tt coop}_{ij}$ defined in (\ref{coopmap}).
\item The maps $m_{i,L}$ defined in (\ref{mpmiL}) to construct the Ising model with Glauber dynamics.
\end{itemize}
As a result, the following interacting particle systems are monotonically
representable (and hence, in particular, monotone):
\begin{itemize}
\item The voter model with generator as in (\ref{Gvot}).
\item The contact process with generator as in (\ref{Gcontact}).
\item The ferromagnetic Ising model with Glauber dynamics, since its generator
  can be written as in (\ref{Isingrep}).
\item The biased voter model with generator as in (\ref{Gbiasvot}).
\item The exclusion process with generator as in (\ref{Gexcl}).
\item Systems of coalescing random walks with generator as in (\ref{Grw}).
\item Systems with cooperative branching and coalescence as in
  Figure~\ref{fig:coop}.
\end{itemize}
On the other hand, the following maps are \emph{not} monotone:
\begin{itemize}
\item The annihilating random walk map ${\tt arw}_{ij}$ defined
in (\ref{annmap}).
\item The killing map ${\tt kill}_{ij}$ defined
in (\ref{killmap}).
\end{itemize}
Examples of interacting particle systems that are not
monotone\footnote{Note that the fact that a given interacting particle
system is represented in maps that are not monotone does not prove
that the system is not monotone. Indeed, it is conceivable that the
same system can also be monotonically represented. See Exercises
\ref{E:threshold1} and \ref{E:threshold2} for an interacting particle
system that has a monotone and a non-monotone random mapping
representation, both of which are useful.} are:
\begin{itemize}
\item The antiferromagnetic Ising model with Glauber dynamics.
\item The Neuhauser--Pacala model from (\ref{NP99}) for small values of $\al$.
\item Systems of annihilating random walks.
\item The biased annihilating branching process of (\ref{BABP}).
\end{itemize}


\section{Positive correlations}\label{S:poscor}

In this section, we study positive correlations. Positive correlations
play an important role in many more involved arguments but
unfortunately we will not see an example of this so readers may skip
the present section at an initial reading. Let $S$ be a finite
partially ordered set and let $\La$ be countable. A probability
measure $\mu$ on $S^\La$ has \emph{positive correlations}
\index{positive correlations} if it satisfies the
equivalent conditions of the following lemma. Recall from (\ref{cov})
that $\cov_\mu(f,g)$ denotes the covariance of $f$ and $g$ under
$\mu$.

\begin{lemma}[Positive correlations]
Let\label{L:poscor} $S$ be a finite partially ordered set, let $\La$
be countable, and let $\mu$ be a probability measure on $S^\La$. Then
the following conditions are equivalent:
\begin{enumerate}[(ii)]
\item $\cov_\mu(f,g)\geq 0\quad\forall f,g\in\Ci^+_{\rm fin}(S^\La)$,
\item $\cov_\mu(f,g)\geq 0\quad\forall f,g\in B^+(S^\La)$.
\end{enumerate}
\end{lemma}

\begin{Proof}
We introduce the following notation. For any $f\geq 0$ with $\mu f>0$,
we let $\mu_f$ denote the probability measure on $S^\La$ defined as
$\mu_f(g):=\mu(fg)/\mu f$ $(g\in B(S^\La))$. In particular, if $f$ is
the indicator function of an event, then $\mu_f$ is the law obtained
from $\mu$ by conditioning on this event. If $f,g\in B(S^\La)$, $f\geq
0$, and $\mu f>0$, then
\[
\cov_\mu(f,g)\geq 0
\quad\desd\quad\mu(fg)\geq(\mu f)(\mu g)
\quad\desd\quad\mu_fg\geq\mu g.
\]
We now prove the equivalence of (i) and (ii). Trivially (ii) implies
(i). Conversely, if (i) holds, then for each $f,g\in\Ci^+_{\rm fin}(S^\La)$
with $f\geq 1$ one has $\mu_fg\geq\mu g$ which using
Lemma~\ref{L:locmon} implies that $\mu_f\geq\mu$ in the stochastic
order. By Theorem~\ref{T:stochord} this implies that $\mu_fg\geq\mu g$
for all $g\in B^+(S^\La)$ and hence $\cov_\mu(f,g)\geq 0$ for all
$f\in\Ci^+_{\rm fin}(S^\La)$ with $f\geq 1$ and $g\in
B^+(S^\La)$. Since adding a constant to $f$ does not change the
covariance, we can remove the condition $f\geq 1$. Repeating the
argument with the roles of $f$ and $g$ interchanged then yields (ii).
\end{Proof}

\begin{Exercise}
In the context of Lemma~\ref{L:poscor}, show that a probability
measure $\mu$ on $S^\La$ has positive correlations if and only if
$\mu(A|B)\geq\mu(A)$ for all measurable increasing events $A,B\sub
S^\La$ such that $\mu(B)>0$. Show that it suffices to check this
condition for events that depend on finitely many
coordinates. \emph{Hint:} Lemma~\ref{L:incr}.
\end{Exercise}

\begin{Exercise}
Show that a probability measure $\mu$ on $\{0,1\}^2$ has positive
correlations if and only if $\mu(00)\mu(11)\geq\mu(01)\mu(10)$. For
probability measures on $\{0,1\}^n$ with $n\geq 3$, it is known that
the FKG condition\index{FKG condition} $\mu(x\wedge y)\mu(x\vee y)
\geq\mu(x)\mu(y)$ implies (but is not equivalent to) positive
correlations, see \cite{FKG71}.
\end{Exercise}

The following proposition gives sufficient conditions for the time
evolution of an interacting particle system to preserve the space of
probability measures with positive correlations. Condition~(i) says
that the interacting particle system is monotone. Note that we do not
assume monotone representability. Condition~(ii) says that the system
only jumps between comparable states. Note that this condition is
actually independent of the graphical representation (assuming all
rates $r_m$ are strictly positive). Variations of this result can be
found in \cite{Har77,Cox84} and \cite[Thm~II.2.14]{Lig85}.

\begin{proposition}[Preservation of positive correlations]
Let\label{P:pospres} $S$ be a finite partially ordered set, let $\La$
be countable, and let $(P_t)_{t\geq 0}$ be the semigroup of an
interacting particle system with generator of the form (\ref{Gmap})
with the rates satisfying (\ref{downsum}). Assume that
\begin{enumerate}[(ii)]
\item $P_t$ is monotone for each $t\geq 0$,
\item for each $x\in S^\La$ and $m\in\Gi$, either $x\leq m(x)$ or $x\geq m(x)$.
\end{enumerate}
Assume that $\mu$ is a probability measure on $S^\La$ with positive
correlations. Then $\mu P_t$ has positive correlations for each $t\geq 0$.
\end{proposition}

\begin{Proof}
Let $\Ci^+_{\rm sum}=\Ci^+_{\rm sum}(S^\La):=\Ci^+(S^\La)\cap\Ci_{\rm sum}(S^\La)$.
By Lemma~\ref{L:difP}, $P_t$ maps $\Ci_{\rm sum}$ into itself for each
$t\geq 0$. Condition~(i) then implies that $P_t$ maps $\Ci^+_{\rm sum}$
into itself. Condition~(ii) and formula (\ref{Gadef2}) imply
that
\[
\Ga_G(f,g)=\sum_{m\in\Gi}r_m\big\{f\big(m(x)\big)-f\big(x\big)\big\}
\{g\big(m(x)\big)-g\big(x\big)\big\}\geq 0
\]
for all $f,g\in\Ci^+_{\rm sum}$. Proposition~\ref{P:cov2} tells us
that for $f,g\in\Ci^+_{\rm sum}$,
\[
\cov_{\mu P_t}(f,g)=\cov_\mu(P_tf,P_tg)+\int_0^t\!\di s\,\mu P_{t-s}\Ga_G(P_sf,P_sg).
\]
Using the fact that $\mu$ has positive correlations and our previous
observations, we see that the right-hand side of this equation is
nonnegative, proving that $\mu P_t$ has positive correlations.
\end{Proof}

\begin{Exercise}
Let $\La$ be a countable set and let $n\geq 1$. Show that product
measures on $\{0,\ldots,n\}^\La$ have positive
correlations. \emph{Hint:} construct an interacting particle system
that has the desired product measure as its invariant law.
\end{Exercise}

\section{The upper and lower invariant laws}

In the present section, we assume that the local state space is
$S=\{0,1\}$, which covers all examples of monotone interacting
particle systems mentioned in Section~\ref{S:monsys}. We will use the
phrase ``an interacting particle system with state space
$\{0,1\}^\La$'' as a shorthand for any interacting particle that can
be constructed from a graphical representation with rates that satisfy
(\ref{downsum}) or possibly the weaker conditions from
Theorem~\ref{T:graph}. We use the symbols $\un 0$ and $\un 1$
\index{00yyy@$\un 0,\un 1$} to denote the states in $S^\La$ that are
identically $0$ or $1$, respectively. Below, $\de_{\un 0}$ denotes the
delta measure at the configuration that is identically $0$, so
$\de_{\un 0}P_t$ denotes the law at time $t$ of the process started in
$X_0(i)=0$ a.s. $(i\in\La)$.

\begin{theorem}[Upper and lower invariant laws]\label{T:upinv}
Let $X$ be an interacting particle system with state space of the form
$\{0,1\}^\La$ and transition probabilities $(P_t)_{t\geq 0}$. Assume
that $X$ is monotone. Then there exist invariant laws $\un\nu$ and
$\ov\nu$ such that
\[
\de_{\un 0}P_t\Asto{t}\un\nu
\quand
\de_{\un 1}P_t\Asto{t}\ov\nu.
\]
If $\nu$ is any other invariant law, then $\un\nu\leq\nu\leq\ov\nu$.
\end{theorem}

The invariant laws $\un\nu$ and $\ov\nu$ from Theorem~\ref{T:upinv} are called
\emph{lower} \index{lower invariant law} and \emph{upper invariant law},
\index{upper invariant law} respectively. Before we give the proof of
Theorem~\ref{T:upinv}, we start with two preparatory lemmas.

\begin{lemma}[Equal mean]\label{L:meanord}
Let $\mu,\nu$ be probability laws on $\{0,1\}^\La$ such that $\mu\leq\nu$ and
\[
\int\!\mu(\di x)\,x(i)\geq\int\!\nu(\di x)\,x(i)\qquad(i\in\La).
\]
Then $\mu=\nu$.
\end{lemma}

\begin{Proof}
By Theorem~\ref{T:stochord}, we can couple random variables with laws
$\P[X\in\,\cdot\,]=\mu$ and $\P[Y\in\,\cdot\,]=\nu$ in such a way that $X\leq
Y$. Now $\E[X(i)]\geq\E[Y(i)]$ implies $\E[Y(i)-X(i)]\leq 0$. Since
$Y(i)-X(i)\geq 0$ a.s., it follows that $X(i)=Y(i)$ a.s. In particular, if this
holds for all $i\in\La$, then $\mu=\nu$.
\end{Proof}

\begin{lemma}[Monotone convergence of probability laws]
Let\label{L:monlaw} $(\nu_n)_{n\geq 0}$ be a sequence of probability
laws on $\{0,1\}^\La$ that are stochastically ordered as
$\nu_k\leq\nu_{k+1}$ $(k\geq 0)$. Then there exists a probability law
$\nu$ on $\{0,1\}^\La$ such that $\nu_n\Rightarrow\nu$, that is, the
$\nu_n$ converge weakly to $\nu$.
\end{lemma}

\begin{Proof}
Since $\nu_nf$ increases to a finite limit for each
$f\in\Ci^+(\{0,1\}^\La)$, this is an immediate consequence of
Lemmas~\ref{L:mondense} and \ref{L:weaksuf}.
\end{Proof}

\begin{Proof}[of Theorem~\ref{T:upinv}]
By symmetry, it suffices to prove the statement for $\un\nu$.  Since
$\un 0$ is the lowest possible state, for each $t\geq 0$, we trivially
have
\[
\de_{\un 0}\leq\de_{\un 0}P_t
\]
By (\ref{ordpres}), this implies that
\[
\de_{\un 0}P_s\leq\de_{\un 0}P_tP_s=\de_{\un 0}P_{t+s}\qquad(s,t\geq 0),
\]
which shows that $t\mapsto\de_{\un 0}P_t$ is nondecreasing with respect to the
stochastic order. By Lemma~\ref{L:monlaw}, each monotone sequence of
probability laws has a weak limit, so there exists a probability law $\un\nu$
on $\{0,1\}^\La$ such that
\[
\de_{\un 0}P_t\Asto{t}\un\nu.
\]
It follows from Lemma~\ref{L:longtime} that $\un\nu$ is an invariant law.

To complete the proof of the theorem, we observe that if $\nu$ is any other
invariant law, then, by (\ref{ordpres}),
\[
\de_{\un 0}\leq\nu
\quad\volgt\quad
\de_{\un 0}P_t\leq\nu P_t=\nu\quad(t\geq 0).
\]
Since $\de_{\un 0}P_t\Rightarrow\nu$ as $t\to\infty$, if follows that $\un\nu f\leq\nu f$ for all $f\in\Ci^+(\{0,1\}^\La)$, which by Theorem~\ref{T:stochord} implies that
$\un\nu\leq\nu$.
\end{Proof}

\begin{theorem}[Ergodicity of monotone systems]\label{T:lowup}
Let $X$ be a monotone interacting particle system with state space
$\{0,1\}^\La$ and lower and upper invariant laws $\un\nu$ and $\ov\nu$.
If
\begin{equation}\label{lowisup}
\int\!\un\nu(\di x)x(i)=\int\!\ov\nu(\di x)x(i)\qquad\forall i\in\La,
\end{equation}
then $X$ has a unique invariant law $\nu:=\un\nu=\ov\nu$ and is ergodic in the
sense that
\[
\P^x\big[X_t\in\,\cdot\,\big]\Asto{t}\nu\qquad(x\in\{0,1\}^\La).
\]
On the other hand, if (\ref{lowisup}) does not hold, then $X$ has at least two
invariant laws.
\end{theorem}

\begin{Proof}
By Lemma~\ref{L:meanord}, (\ref{lowisup}) is equivalent to the condition that
$\un\nu=\ov\nu$. It is clear that if $\un\nu\neq\ov\nu$, then $X$ has at least
two invariant laws and ergodicity cannot hold. On the other hand, by
Theorem~\ref{T:upinv}, any invariant law $\nu$ satisfies
$\un\nu\leq\nu\leq\ov\nu$, so if $\un\nu=\ov\nu$, then $\nu=\un\nu=\ov\nu$.

To complete the proof, we must show that $\un\nu=\ov\nu=:\nu$ implies
$\de_xP_t\Rightarrow\nu$ as $t\to\infty$ for all $x\in\{0,1\}^\La$.
Since
\[
\de_{\un 0}P_tf\leq\de_xP_tf\leq\de_{\un 1}P_tf
\]
for all $f\in\Ci^+(\{0,1\}^\La)$, we see that
\[
\un\nu f\leq\liminf_{t\to\infty}P_tf\leq\limsup_{t\to\infty}P_tf\leq\ov\nu f
\]
for all $f\in\Ci^+(\{0,1\}^\La)$. The claim now follows from
Lemmas~\ref{L:weaksuf} and \ref{L:mondense}.
\end{Proof}

To state the final result of this section, we need a bit of theory. We
observe that for any interacting particle system, the set
$\Ii$\index{0Ii@$\Ii$} of all invariant laws is a compact, convex
subset of the space $\Mi_1(S^\La)$ of probability measures on $S^\La$,
equipped with the topology of weak convergence. Indeed, if $\mu$ and
$\nu$ are invariant laws and $p\in[0,1]$, then clearly
\[
\big(p\mu+(1-p)\nu\big)P_t=p\mu P_t+(1-p)\nu P_t=p\mu+(1-p)\nu\quad(t\geq 0),
\]
proving that $p\mu+(1-p)\nu$ is an invariant law. The fact that $\Ii$ is
closed follows from Proposition~\ref{P:invlim}. Since $\Mi_1(S^\La)$ is
compact, $\Ii$ is also compact.

By definition, an element $\nu\in\Ii$ is called \emph{extremal}
\index{extremal invariant law} if it cannot be written as a nontrivial convex
combination of other elements of $\Ii$, that is,
\[
\nu=p\nu_1+(1-p)\nu_2\qquad(0<p<1,\ \nu_1,\nu_2\in\Ii)
\quad\mbox{implies}\quad\nu_1=\nu_2=\nu.
\]
We let
\[\index{0Iie@$\Ii_{\rm e}$}
\Ii_{\rm e}:=\{\nu\in\Ii:\nu\mbox{ is an extremal element of }\Ii\}.
\]
Since $\Ii$ is compact and convex, Choquet's theorem implies that each
invariant law $\nu$ can be written as
\[
\nu=\int_{\Ii_{\rm e}}\rho_\nu(\di\mu)\mu,
\]
where $\rho_\nu$ is a probability measure on $\Ii_{\rm e}$. In practice, it
happens quite often\footnote{The the voter model in dimensions $d\geq 3$ is a
  counterexample. The Ising model in dimensions $d\geq 3$ is also a
  counterexample, although for the Ising model, it is still true that $\un\nu$
  and $\ov\nu$ are the only extremal invariant measures that are moreover
  translation invariant.} that $\Ii_{\rm e}$ is a finite set.\footnote{This
  may, however, be quite difficult to prove!} In this case, Choquet's theorem
simply says that each invariant law is a convex combination of the extremal
invariant laws, that is, each invariant law is of the form
\[
\nu=\sum_{\mu\in\Ii_{\rm e}}p(\mu)\mu,
\] 
where $(p(\mu))_{\mu\in\Ii_{\rm e}}$ are nonnegative constants, summing up to one.
In view of this, we are naturally interested in finding all extremal invariant
laws of a given interacting particle system.

\begin{lemma}[The lower and upper invariant law are extremal]
Let\label{L:upextr} $X$ be a monotone interacting particle system with
state space $\{0,1\}^\La$ and lower and upper invariant laws $\un\nu$
and $\ov\nu$. Then $\un\nu$ and $\ov\nu$ are extremal invariant laws
of $X$.
\end{lemma}

\begin{Proof}
By symmetry, it suffices to prove the statement for $\ov\nu$. 
Imagine that
\[
\ov\nu=p\nu_1+(1-p)\nu_2\quad\mbox{for some}\quad 0<p<1,\ \nu_1,\nu_2\in\Ii.
\]
By Theorem~\ref{T:upinv}, for each $f\in B^+(\{0,1\}^\La)$,
one has $\nu_1f\leq\ov\nu f$ and $\nu_2f\leq\ov\nu f$. Since
\[
p(\ov\nu f-\nu_1f)+(1-p)(\ov\nu f-\nu_2f)=0,
\]
it follows that $\ov\nu f=\nu_1f=\nu_2f$. Since this holds for each
monotone $f$, we conclude (by Lemma~\ref{L:mondense}) that $\ov\nu=\nu_1=\nu_2$.
\end{Proof}

\begin{Exercise}
Let $X$ be an interacting particle system with state space $\{0,1\}^\La$ and
generator $G$. Assume that $G$ has a random mapping representation in terms of
monotone maps and let $(\Xb_{s,t})_{s\leq t}$ be the corresponding stochastic
flow as in Theorem~\ref{T:Poispart}. Show that the a.s.\ limits
\[\left.\begin{array}{r@{\,}c@{\,}l}
\dis\un X_t&:=&\dis\lim_{s\to-\infty}\Xb_{s,t}(\un 0),\\[5pt]
\dis\ov X_t&:=&\dis\lim_{s\to-\infty}\Xb_{s,t}(\un 1)
\end{array}\ \right\}\quad(t\in\R)
\]
define stationary Markov processes $(\un X_t)_{t\in\R}$ and $(\ov X_t)_{t\in\R}$
whose invariant laws
\[
\un\nu=\P[\un X_t\in\,\cdot\,]
\quand
\ov\nu=\P[\ov X_t\in\,\cdot\,]\qquad(t\in\R)
\]
are the lower and upper invariant law of $X$, respectively.
Show that (\ref{lowisup}) implies that
\[
\lim_{s\to-\infty}\Xb_{s,t}(x)=\un X_t=\ov X_t
\quad{\rm a.s.}\quad(x\in\{0,1\}^\La,\ t\in\R).
\]
\end{Exercise}

\section{The contact process}\label{S:contmon}\index{contact process}

The contact process has been defined on $\Z^d$ in (\ref{Gcontact}) and
on a very general class of lattices in (\ref{Gcont}). In the present
section, we will look at a class of contact processes that are more
general than those in (\ref{Gcontact}) but a bit less general than
those in (\ref{Gcont}). Throughout this section, $\La$ will be a
countable set and $p$ will be a probability kernel on $\La$ that is
\emph{symmetric} in the sense that $p(i,j)=p(j,i)$ $(i,j\in\La)$
and satisfies $p(i,i)=0$ $(i\in\La)$. By definition, an \emph{automorphism}
\index{automorphism!of probability kernel} of $(\La,p)$ is a bijection
$\psi\cn\La\to\La$ such that $p\big(\psi(i),\psi(j)\big)=p(i,j)$ ($i,j\in\La)$
(compare the footnote on page~\pageref{autfoot}). We will assume that
$(\La,p)$ is \emph{vertex transitive} \index{transitivity!of probability kernel}
\index{vertex transitive!probability kernel} in the sense that
\begin{equation}\label{ptrans}
\forall i,j\in\La\ \exists\mbox{ automorphism $\psi$ of }(\La,p)
\mbox{ s.t.\ }\psi(i)=j.
\end{equation}
We will be interested in contact processes with generator of the form
\begin{equation}\label{Gcont3}
G_{\rm cont}f(x):=\la\sum_{i,j\in\La}p(i,j)
\big\{f\big({\tt bra}_{ij}(x)\big)-f\big(x\big)\big\}
+\sum_{i\in\La}\big\{f\big({\tt death}_i(x)\big)-f\big(x\big)\big\},
\end{equation}
where $\la\geq 0$ is the \emph{infection rate}\index{infection rate}
and the death rate is one. Note that our present
definition differs a bit from the classical definition of the contact
process on $\Z^d$ in (\ref{Gcontact}) in the sense that in
(\ref{Gcont3}) the total rate of all infections out of a site $i$ is
$\la$, while in (\ref{Gcontact}) it is $\la|\Ni_i|$, where $|\Ni_i|$
is the number of neighbors of $i$. Already when we studied the
mean-field limit of the contact process, we discovered that the
normalization in (\ref{Gcont3}) is often more natural, compare
(\ref{Gcontact2}).

Since both the branching and death map are monotone, the contact
process is a monotonically representable interacting particle system,
so by Theorem~\ref{T:upinv}, it has a lower and upper invariant law
$\un\nu$ and $\ov\nu$. Since ${\tt bra}_{ij}(\un 0)=\un 0$ and
${\tt death}_i(\un 0)=\un 0$ for each $i,j\in\La$, the all-zero
configuration $\un 0$ is a trap for the contact process, so
$\de_{\un 0}P_t=\de_{\un 0}$ for all $t\geq 0$ and hence
\[
\un\nu=\de_{\un 0}.
\]
Therefore, by Theorem~\ref{T:lowup}, the contact process is ergodic if and
only if the function
\begin{equation}\label{tetdef}\index{0zzhthetalambda@$\tet(\la)$}
\tet(\la):=\int\!\ov\nu_\la(\di x)\,x(i)\qquad(i\in\Z^d)
\end{equation}
satisfies $\tet(\la)=0$. Here $\ov\nu_\la$ denotes the upper invariant
law of the contact process with infection rate $\la$ and the
right-hand side of (\ref{tetdef}) does not depend on $i\in\La$ by our
assumption that $(\La,p)$ is vertex transitive. For reasons that will
become clear in the next chapter (see Lemma~\ref{L:theta}),
$\tet(\la)$ is actually the same as the survival probability started
from a single occupied site, that is, this is the function in
Figure~\ref{fig:surprob}.

By definition, we say that a probability law $\mu$ on $\{0,1\}^\La$ is
\emph{nontrivial}\index{nontrivial law} if
\[
\mu(\{\un 0\})=0,
\]
that is, if $\mu$ gives zero probability to the all-zero configuration.

\begin{lemma}[Nontriviality of the upper invariant law]\label{L:nontriv}
For the contact process, if $\ov\nu\neq\de_{\un 0}$, then $\ov\nu$ is nontrivial.
\end{lemma}

\begin{Proof}
We can always write $\ov\nu=(1-p)\de_{\un 0}+p\mu$ where $p\in[0,1]$
and $\mu$ is a nontrivial law. By assumption, $\ov\nu\neq\de_{\un 0}$,
so $p>0$. Since $\ov\nu$ and $\de_{\un 0}$ are invariant laws, $\mu$
must be an invariant law too. By Lemma~\ref{L:upextr}, $\ov\nu$ cannot
be written as a nontrivial convex combination of other invariant laws,
so we conclude that $p=1$.
\end{Proof}

\begin{proposition}[Monotonicity in the infection rate]\label{P:moninla}
Let $\ov\nu_\la$ denote the upper invariant law of the contact process with
infection rate $\la$. Then $\la\leq\la'$ implies $\ov\nu_\la\leq\ov\nu_{\la'}$.
In particular, the function $\la\mapsto\tet(\la)$ is nondecreasing.
\end{proposition}

\begin{Proof}
Let $X$ and $X'$ be contact processes started in the initial state
$X_0=\un 1=X'_0$ and with infection rates $\la$ and $\la'$. It suffices to prove
that $X$ and $X'$ can be coupled such that $X_t\leq X'_t$ for all $t\geq 0$.

We will couple the graphical representations of the processes with
infection rates $\la$ and $\la'$. We write $\Gi=\Gi_{\rm bra}\cup\Gi_{\rm death}$
where
\[
\Gi_{\rm bra}:=\big\{{\tt bra}_{ij}:i,j\in\La\big\}
\quand
\Gi_{\rm death}:=\big\{{\tt death}_i:i\in\La\big\}.
\]
Then $X$ can be constructed as in Theorem~\ref{T:Poispart} from a Poisson point
set $\om$ on $\Gi\times\R$ with intensity measure $\rho_\la$ given by
\[\left.\begin{array}{r@{\,}c@{\,}l}
\dis\rho_\la(\{{\tt bra}_{ij}\}\times[s,t])&:=&\la p(i,j)(t-s),\\[5pt]
\dis\rho_\la(\{{\tt death}_i\}\times[s,t])&:=&(t-s),
\end{array}\right\}\quad(i,j\in\La,\ s\leq t).
\]
Likewise, $X'$ can be constructed from a Poisson point set $\om'$ with intensity
$\rho_{\la'}$. We claim that we can couple $\om$ and $\om'$ in such a
way that the latter has more branching maps, and the same death
maps as $\om$. This can be done as follows. Let $\om''$ be a
Poisson point set on $\Gi\times\R$, independent of $\om$, with
intensity measure $\rho'':=\rho_{\la'}-\rho_\la$, that is,
\[\left.\begin{array}{r@{\,}c@{\,}l}
\dis\rho''(\{{\tt bra}_{ij}\}\times[s,t])&:=&(\la'-\la)p(i,j)(t-s),\\[5pt]
\dis\rho''(\{{\tt death}_i\}\times[s,t])&:=&0,
\end{array}\right\}\quad(i,j\in\La,\ s\leq t).
\]
Since the sum of two independent Poisson sets yields another Poisson set,
setting
\[
\om':=\om+\om''
\]
defines a Poisson point set with intensity $\rho_{\la'}$.
We observe that
\[\begin{array}{l}
\dis x\leq x'\quad\mbox{implies}\quad
{\tt bra}_{ij}(x)\leq{\tt bra}_{ij}(x'),\\[5pt]
\dis x\leq x'\quad\mbox{implies}\quad
{\tt death}_i(x)\leq{\tt death}_i(x'),\\[5pt]
\dis x\leq x'\quad\mbox{implies}\quad
x\leq{\tt bra}_{ij}(x').
\end{array}\]
The first two statements just say that the maps ${\tt bra}_{ij}$ and
${\tt death}_i$ are monotone. The third statement says that if we apply a
branching map only to the larger configuration $x'$, then the order between
$x$ and $x'$ is preserved.

Since $\om'$ has the same branching and death maps as $\om$, plus
some extra branching maps, using Proposition~\ref{P:fincon} we
conclude that the stochastic flows $(\Xb_{s,t})_{s\leq t}$ and
$(\Xb'_{s,t})_{s\leq t}$ constructed from $\om$ and $\om'$ satisfy
\[
x\leq x'\quad\mbox{implies}\quad\Xb_{s,t}(x)\leq\Xb'_{s,t}(x')
\qquad(s\leq t).
\]
In particular, setting $X_t:=\Xb_{0,t}(\un 1)$ and $X'_t:=\Xb'_{0,t}(\un 1)$
yields the desired coupling between $X$ and $X'$.
\end{Proof}

\begin{Exercise}\label{E:context}
For contact processes with generator of the form (\ref{Gcont3}),
calculate the constant $K_\down$ from (\ref{KKK}) and apply
Theorem~\ref{T:ergo} to conclude that
\[
\la<1\quad\mbox{implies}\quad\ov\nu=\de_{\un 0}.
\]
\end{Exercise}

In Chapter~\ref{C:percol}, we will prove that $\tet(\la)>0$ for $\la$
sufficiently large.

\section{Other examples}\label{S:oth}

\subsubsection*{The Ising model with Glauber dynamics}

We have seen in (\ref{Isingrep}) that the generator of the Ising model
with Glauber dynamics is monotonically representable, so by
Theorem~\ref{T:upinv},\footnote{The difference between the local state
space $\{-1,1\}$ of the Ising model and $\{0,1\}$ of
Theorem~\ref{T:upinv} is of course entirely notational.} it has a
lower and upper invariant law $\un\nu$ and $\ov\nu$.  We let
\[
m_\ast(\bet):=\int\!\ov\nu(\di x)\,x(i),
\]
which is independent of $i$ if the processes has some translation invariant
structure (like the nearest neighbor or range $R$ processes on $\Z^d$). For
reasons that cannot be explained here, this function is actually the same as
the one defined in (\ref{mast}), that is, this is the \emph{spontaneous
  magnetization}\index{spontaneous magnetization} of the Ising model, see
Figure~\ref{fig:Ismag}. By the symmetry between $-1$ and $+1$ spins, we
clearly have
\[
\int\!\un\nu(\di x)\,x(i)=-m_\ast(\bet).
\]
By Theorem~\ref{T:Iserg}, we have
\[
\bet<N^{-1}\big(\log(N+1)-\log(N-1)\big)
\quad\mbox{implies}\quad\un\nu=\ov\nu,
\]
from which we conclude that $m_\ast(\bet)=0$ for $\bet$ sufficiently small,

The function $\bet\mapsto m_\ast(\bet)$ is nondecreasing, but
this cannot be proved with the sort of techniques used in
Proposition~\ref{P:moninla}. The lower and upper invariant laws of the Ising
model with Glauber dynamics are infinite volume Gibbs measures, and much of
the analysis of the Ising model is based on this fact. In fact, the Ising
model with Glauber dynamics is just one example of an interacting particle
system that has these Gibbs measures as its invariant laws. In general,
interacting particle systems with this property are called stochastic Ising
models, and the Gibbs measures themselves are simply called the Ising model.
We refer to \cite[Chapter~IV]{Lig85} for an exposition of this material.
In particular, in \cite[Thm~IV.3.14]{Lig85}, it is shown that for the
nearest-neighbor model on $\Z^2$, one has $m_\ast(\bet)>0$ for $\bet$
sufficiently large.

\subsubsection*{The voter model}

Consider a voter model with local state space $S=\{0,1\}$.  Since the voter
maps ${\tt vot}_{ij}$ from (\ref{votmap}) are monotone, the voter model is
monotonically representable. Since both the constant configurations $\un 0$ and
$\un 1$ are traps,
\[
\un\nu=\de_{\un 0}\quand\ov\nu=\de_{\un 1},
\]
so we conclude (recall Theorem~\ref{T:lowup}) that the voter model is
never ergodic. For the model on $\Z^d$, it is proved in
\cite[Thm~V.1.8]{Lig85} that if $d=1,2$, then $\de_{\un 0}$ and
$\de_{\un 1}$ are the only extremal invariant laws.  On the other
hand, in dimensions $d\geq 3$, the set $\Ii_{\rm e}$ of extremal
invariant laws is of the form $\{\nu_p:p\in[0,1]\}$ where the
invariant measure $\nu_p$ has intensity $\int\!\nu_p(\di x)\,x(i)=p$.
We will give a partial proof of these statements in Section~\ref{S:votinv}.

\section{Exercises}

\begin{Exercise}
Give an example of two probability measures $\mu,\nu$ on a
set of the form $\{0,1\}^\La$ that satisfy
\[
\int\mu(\di x)x(i)
\leq\int\nu(\di x)x(i)\qquad(i\in\La),
\]
but that are \emph{not} stochastically ordered as $\mu\leq\nu$.
\end{Exercise}

\begin{Exercise}\label{E:tetright}
Let $(X^\la_t)_{t\geq 0}$ denote the contact process with
infection rate $\la$ (and death rate one), started in $X^\la_0=1$.
Apply Corollary~\ref{C:Partlim} to prove that for each fixed $t\geq 0$,
the function
\begin{equation}
\tet_t(\la):=\P[\Xb^\la_{0,t}(1)(i)=1]
\end{equation}
depends continuously on $\la$. Use this to conclude that the function
$\tet(\la)$ from (\ref{tetdef}) is right-continuous. \emph{Hint:} Use
that the decreasing limit of continuous functions is upper semi-continuous.
\end{Exercise}

For the next exercise, let us define a \emph{double death}
\index{double death map} map
\begin{equation}\label{doubdeath}\index{0deathz@${\tt death}_{ij}$}
{\tt death}_{ij}(x)(k):=\left\{\begin{array}{cl}
0&\mbox{if }k\in\{i,j\},\\
x(k)&\mbox{otherwise.}
\end{array}\right.
\end{equation}
Recall the branching map ${\tt bra}_{ij}$ defined in (\ref{bramap}),
the death map ${\tt death}_i$ defined in (\ref{deathmap}), and the
cooperative branching map ${\tt coop}_{ij}$ defined in (\ref{coopmap}).
Consider the cooperative branching process $X$ with values in $\{0,1\}^\Z$
with generator
\[
G_Xf(x)=\la\sum_{i\in\Z}\sum_{\sig\in\{-1,+1\}}
\big\{f\big({\tt coop}_{i+2\sig,i+\sig,i}(x)\big)-f\big(x\big)\big\}
+\sum_{i\in\Z}\big\{f\big({\tt death}_i(x)\big)-f\big(x\big)\big\},
\]
and the contact process with double deaths $Y$ with generator
\[
G_Yf(y)
=\la\sum_{i\in\Z}\sum_{\sig\in\{-1,+1\}}
\big\{f\big({\tt bra}_{i+\sig,i}(y)\big)-f\big(y\big)\big\}+\sum_{i\in\Z}
\big\{f\big({\tt death}_{i,i+1}(y)\big)-f\big(y\big)\big\}.
\]

\begin{Exercise}
Let\label{E:doubdeath} $X$ be the process with cooperative branching defined
above and set
\[
X^{(2)}_t(i):=1_{\{X_t(i)=1=X_t(i+1)\}}\qquad(i\in\Z,\ t\geq 0).
\]
Show that $X$ can be coupled to a contact process with double deaths $Y$ (with
the same parameter $\la$) in such a way that
\[
Y_0\leq X^{(2)}_0\quad\mbox{implies}\quad
Y_t\leq X^{(2)}_t\quad(t\geq 0).
\]
\end{Exercise}

\begin{Exercise}
Show that a system $(X_t)_{t\geq 0}$ of annihilating random
walks can be coupled to a system $(Y_t)_{t\geq 0}$ of coalescing random
walks such that
\[
X_0\leq Y_0\quad\mbox{implies}\quad
X_t\leq Y_t\quad(t\geq 0).
\]
Note that the annihilating random walks are not a monotone particle system.
\end{Exercise}

\begin{Exercise}
Let $X$ be a system of branching and coalescing random
walks with generator
\[\begin{array}{r@{\,}c@{\,}l}
\dis G_Xf(x)
&=&\dis\ha b\sum_{i\in\Z}\sum_{\sig\in\{-1,+1\}}
\big\{f\big({\tt bra}_{i,i+\sig}x\big)-f\big(x\big)\big\}\\[5pt]
&&\dis+\ha\sum_{i\in\Z}\sum_{\sig\in\{-1,+1\}}
\big\{f\big({\tt rw}_{i,i+\sig}x\big)-f\big(x\big)\big\},
\end{array}\]
and let $Y$ be a system of coalescing random
walks with positive drift, with generator
\[\begin{array}{r@{\,}c@{\,}l}
\dis G_Yf(y)
&=&\dis\ha(1+b)\sum_{i\in\Z}
\big\{f\big({\tt rw}_{i,i+1}y\big)-f\big(y\big)\big\}\\[5pt]
&&\dis+\ha\sum_{i\in\Z}
\big\{f\big({\tt rw}_{i,i-1}y\big)-f\big(y\big)\big\}.
\end{array}\]
Show that $X$ and $Y$ can be coupled such that
\[
Y_0\leq X_0\quad\mbox{implies}\quad
Y_t\leq X_t\quad(t\geq 0).
\]
\end{Exercise}

\begin{Exercise}\label{E:dimcomp}
Let $d<d'$ and identify $\Z^d$ with the subset of $\,\Z^{d'}$
consisting of all $(i_1,\ldots,i_{d'})$ with
$(i_{d+1},\ldots,i_{d'})=(0,\ldots,0)$. Let $X$ and $X'$ denote the
nearest-neighbor contact processes on $\Z^d$ and $\Z^{d'}$,
respectively, with generator as in (\ref{Gcontact}), with the same
infection rate $\la$ and death rate $\de$. Show that $X$ and $X'$ can
be coupled such that
\[
X_0(i)\leq X'_0(i)\quad(i\in\Z^d)
\qquad\mbox{implies}\qquad
X_t(i)\leq X'_t(i)\quad(t\geq 0,\ i\in\Z^d).
\]
Prove the same when $X$ is the nearest-neighbor process and $X'$ is
the range $R$ process (both on $\Z^d$). (Note that for these
comparison arguments, the normalization in (\ref{Gcontact}) is more
convenient than the normalization in (\ref{Gcont3}).)
\end{Exercise}

\chapter{Duality}\label{C:dual}

\section{Basic definitions}\label{S:dubasic}

Let $S$ be a finite set, let $\La$ be countable, let $\Gi$ be a
collection of local maps $m\cn S^\La\to S^\La$, and let
$(r_m)_{m\in\Gi}$ be nonnegative rates satisfying
(\ref{downsum}). Then Theorem~\ref{T:Poispart} tells us how the
interacting particle system with generator
\begin{equation}\label{Gmap2}
Gf(x):=\sum_{m\in\Gi}r_m\big\{f\big(m(x)\big)-f\big(x\big)\big\}
\qquad(x\in S^\La).
\end{equation}
can be constructed from a graphical representation $\om$, which is a
Poisson point set on $\Gi\times\R$ with intensity as in
(\ref{rhoips}). More precisely, in (\ref{partflow}), we have seen how
in terms of $\om$ it is possible to define a stationary stochastic
flow $(\Xb_{s,t})_{s\leq t}$ with independent increments, so that if
$s\in\R$ and $X_0$ is an $S^\La$-valued random variable, independent
of $\om$, then setting
\begin{equation}\label{forpr}
X_t:=\Xb_{s,s+t}(X_0)\qquad(t\geq 0)
\end{equation}
defines a Feller process $(X_t)_{t\geq 0}$ whose generator is (the
closure of) $G$ defined in (\ref{Gmap}). We call this the interacting
particle system with generator $G$.

Key to the proof of Theorem~\ref{T:Poispart} was the backward in time
process. For any finite set $T$, by (\ref{FXrel}), setting
\begin{equation}\label{backF}\index{0Fb@$\Fb_{t,s}$}
\Fb_{t,s}(\phi):=\phi\circ\Xb_{s,t}\qquad\big(s\leq t,\ \phi\in\Ci(S^\La,T)\big),
\end{equation}
defines a backward stochastic flow $(\Fb_{t,s})_{t\geq s}$ on the
countable set $\Ci(S^\La,T)$. Alternatively\footnote{In fact, in
Chapter~\ref{C:construct}, we first defined $(\Fb_{t,s})_{t\geq s}$ by
(\ref{Fback}) and then used it to prove the existence of a forward
stochastic flow $(\Xb_{s,t})_{s\leq t}$ as in Theorem~\ref{T:Poispart}
such that (\ref{backF}) holds.} $(\Fb_{t,s})_{t\geq s}$ can be defined
directly in terms of the graphical representation $\om$ as in
(\ref{Fback}). If $u\in\R$ and $\Phi_0$ is a random variable with
values in $\Ci(S^\La,T)$, independent of $\om$, then setting
\begin{equation}\label{backpr}
\Phi_t:=\Fb_{u,u-t}(\Phi_0)\qquad(t\geq 0)
\end{equation}
defines a nonexplosive continuous-time Markov chain $(\Phi_t)_{t\geq
  0}$ with values in $\Ci(S^\La,T)$ and generator $H$ as in
(\ref{Hgen}). We called this the \emph{backward in time
process}.\index{backward!in time process} Due to the reversal of time,
this Markov process has, somewhat unusually, left-continuous sample
paths.

As we will see in the present chapter, the backward in time process is
not just a useful tool in the construction of interacting particle
systems but also in their further study. In fact, many of the most
tractable and most studied interacting particle systems (such as the
voter model and the contact process) are tractable precisely because
their backward in time process is of a simple nature. The study of the
backward in time process naturally leads to Markov process duality,
which is the topic of the present chapter.

Let $S,R$, and $T$ be sets and let $\psi\cn S\times R\to T$ be a
function. Then we say that two maps $m\cn S\to S$ and $\hat m\cn R\to
R$ are \emph{dual}\index{duality!of maps} to each other with respect
to the \emph{duality function}\index{duality!function} $\psi$ if
\[\index{0mh@$\hat m$}
\psi\big(m(x),y\big)=\psi\big(x,\hat m(y)\big)\qquad(x\in S,\ y\in R).
\]
If $(\Xb_{s,t})_{s\leq t}$ is a stochastic flow (in the sense of
(\ref{stochflow})) on $S$ and $(\Yb_{t,s})_{t\geq s}$ is a backward
stochastic flow (in the sense of (\ref{backflow})) on $R$, then we say
that $(\Xb_{s,t})_{s\leq t}$ and $(\Yb_{t,s})_{t\geq s}$ are
\emph{dual}\index{duality!of stochastic flows} to each other with
respect to the duality function $\psi$ if
\begin{equation}\label{flowdual}\index{0Yb@$\Yb_{t,s}$}
\psi\big(\Xb_{s,t}(x),y\big)=\psi\big(x,\Yb_{t,s}(y)\big)\qquad(s\leq t,\ x\in S,\ y\in R).
\end{equation}
Fix $s<u$, let $X_0$ and $Y_0$ be independent of each other and of the
stochastic flows $(\Xb_{s,t})_{s\leq t}$ and $(\Yb_{t,s})_{t\geq s}$,
and let $(X_t)_{t\geq 0}$ and $(Y_t)_{t\geq 0}$ be the Markov
processes defined by
\[
X_t:=\Xb_{s,s+t}(X_0)
\quand
Y_t:=\Yb_{u,u-t}(Y_0)
\qquad(t\geq 0).
\]
Then we claim that
\begin{equation}\label{pathdual}
\mbox{the function }[s,u]\ni t\mapsto\psi(X_{t-s},Y_{u-t})\mbox{ is constant.}
\end{equation}
Indeed, the duality of $(\Xb_{s,t})_{s\leq t}$ and $(\Yb_{t,s})_{t\geq s}$ implies
\[\begin{array}{l}
\dis\psi(X_{t-s},Y_{u-t})=\psi\big(\Xb_{s,t}(X_0),\Yb_{u,t}(Y_0)\big)\\[5pt]
\dis\quad=\psi\big(\Xb_{t,u}\circ\Xb_{s,t}(X_0),Y_0\big)=\psi\big(\Xb_{s,u}(X_0),Y_0\big),
\end{array}\]
which clearly does not depend on $t$. A relation of the form
(\ref{pathdual}) is called a \index{duality!pathwise}\emph{pathwise
duality}.\footnote{This terminology was first introduced in
\cite{JK14b}.} In particular, setting $t=s,u$ in (\ref{pathdual}) we
see that
\[
\psi(X_u,Y_0)=\psi(X_0,Y_u).
\]
In the special case that $\psi$ takes values in a linear space such as
$\R$ or $\C$, we can take expectations and conclude that
\begin{equation}\label{dual}
\E\big[\psi(X_u,Y_0)\big]=\E\big[\psi(X_0,Y_u)\big]\qquad(u\geq 0),
\end{equation}
whenever $X_u$ is independent of $Y_0$ and $X_0$ is independent of
$Y_u$. A relation of the form (\ref{dual}) is called a
\emph{duality}\index{duality!of Markov processes} between the Markov
processes $(X_t)_{t\geq 0}$ and $(Y_t)_{t\geq 0}$.

We claim that any interacting particle system of the type described at
the beginning of this section trivially has at least one pathwise
dual, which is the backward in time process. To see this, fix a finite
set $T$ and let $\psi\cn S^\La\times\Ci(S^\La,T)\to T$ be the duality
function defined by
\begin{equation}
\psi(x,\phi):=\phi(x)\qquad\big(x\in S^\La,\ \phi\in\Ci(S^\La,T)\big).
\end{equation}
Then it is straightforward to check that the stochastic flow
$(\Xb_{s,t})_{s\leq t}$ and the backward stochastic flow
$(\Fb_{t,s})_{t\geq s}$ are dual with respect to the duality function
$\psi$, and hence the interacting particle system $(X_t)_{t\geq 0}$
and the backward in time process $(\Phi_t)_{t\geq 0}$ are pathwise
dual.

Although formally, we have now found a pathwise dual for each
interacting particle system, this dual is of little use in practice
since the backward in time process is in general very complicated and
the space $\Ci(S^\La,T)$ is very large. As we will see in the coming
sections, however, it sometimes happens that $\Ci(S^\La,T)$ contains a
subspace of ``nice'' functions that is mapped into itself under the
backward stochastic flow $(\Fb_{t,s})_{t\geq s}$, and this then leads
to a more useful pathwise duality.

\section{Additive systems}\label{S:add}

There exists a useful duality theory for additive systems. We first
discuss additive systems and then in the next section discuss their
duality. In line with notation introduced in Section~\ref{S:finsys},
but specialized to $S=\{0,1\}$, we set
\[\index{0SiLa@$\Si(\La),\Si_{\rm fin}(\La)$}
\Si(\La):=\{0,1\}^\La
\quand
\Si_{\rm fin}(\La):=\big\{x\in\Si(\La):|x|<\infty\big\},
\]
where
\[\index{00normx@$\vert x\vert$}
|x|:=\sum_{i\in\La}x(i)\qquad\big(x\in\Si(\La)\big).
\]
We equip $\Si(\La):=\{0,1\}^\La$ with the product topology and
$\Si_{\rm fin}(\La)$, which is countable, with the discrete
topology. As we have done before, we let $\un 0$ denote the
configuration that is identically zero. We let $(x\vee y)(i):=x(i)\vee y(i)$
\index{00normy@$x\vee y$} denote the pointwise maximum of two
configurations $x,y\in\Si(\La)$. Let $\La,\De$ be countable
sets. By definition, a map $m\cn\Si(\La)\to\Si(\De)$ is
\emph{additive}\index{additive!map} iff
\[\begin{array}{r@{\ }l}
{\rm(i)}&m(\un 0)=\un 0,\\[5pt]
{\rm(ii)}&m(x\vee y)=m(x)\vee m(y)\qquad\big(x,y\in\Si(\La)\big).
\end{array}\]
The same definition applies to maps $m\cn\Si_{\rm fin}(\La)\to\Si(\De)$,
where now (ii) needs to hold only for $x,y\in\Si_{\rm fin}(\La)$.
If $m$ is moreover continuous, then
\begin{equation}\label{infadd}
m\Big(\bigvee_{k=1}^\infty x_k\Big)=\bigvee_{k=1}^\infty m(x_k)
\qquad\big(x_k\in\Si(\La)\big),
\end{equation}
as follows by writing $\dis m\Big(\lim_{n\to\infty}\bigvee_{k=1}^n x_k\Big)
=\lim_{n\to\infty}m\Big(\bigvee_{k=1}^n x_k\Big)
=\lim_{n\to\infty}\bigvee_{k=1}^n m(x_k)$.
Each additive map is monotone, as follows by observing that $y\geq x$
implies $m(y)=m(x\vee y)=m(x)\vee m(y)\geq m(x)$. A lot of local maps
we have already seen are additive. Examples are:
\begin{itemize}
\item The voter map ${\tt vot}_{ij}$ defined in (\ref{votmap}).
\item The branching map ${\tt bra}_{ij}$ defined in (\ref{bramap}).
\item The death map ${\tt death}_i$ defined in (\ref{deathmap}).
\item The exclusion map ${\tt excl}_{ij}$ defined in (\ref{exclmap}).
\item The coalescing random walk map ${\tt rw}_{ij}$ defined in (\ref{rwmap}).
\end{itemize}
On the other hand, the following local maps are monotone, but not additive:
\begin{itemize}
\item The asymmetric exclusion map ${\tt asep}_{ij}$ defined in (\ref{asep}).
\item The cooperative branching map ${\tt coop}_{ijk}$ defined in
  (\ref{coopmap}).
\item The maps $m_{i,L}$ defined in (\ref{mpmiL}) to construct the
  Ising model with Glauber dynamics.
\end{itemize}
An interacting particle system is called
\emph{additive}\index{additive!interacting particle system} if its
generator can be represented in additive local maps. Examples of
additive particle systems are:
\begin{itemize}
\item The voter model with generator as in (\ref{Gvot}).
\item The contact process with generator as in (\ref{Gcontact}).
\item The biased voter model with generator as in (\ref{Gbiasvot}).
\item The symmetric exclusion process with generator as in (\ref{Gexcl}).
\item Systems of coalescing random walks with generator as in (\ref{Grw}).
\end{itemize}

We define $e_i\in\Si(\La)$ by $e_i(j):=1$ if $i=j$ and $:=0$
otherwise.\index{0ei@$e_i$} There is a useful graphical way to
describe a continuous additive map $m\cn\Si(\La)\to\Si(\La)$, that
works as follows:
\begin{itemize}
\item For each $i,j\in\La$ with $i\neq j$ such that $m(e_i)(j)=1$, we
  draw an arrow from $i$ to $j$.
\item For each $i\in\La$ such that $m(e_i)(i)=0$, we draw a blocking
  symbol \blok at $i$.
\end{itemize}
The following lemma says that continuous additive maps are fully
described by their arrows and blocking symbols.

\begin{lemma}[Graphical description]
Let\label{L:adloc} $m\cn\Si(\La)\to\Si(\La)$ be a continuous additive
map and let $x\in\{0,1\}^\La$. Then $m(x)(j)=1$ if and only if at
least one of the following conditions is satisfied:
\begin{enumerate}[(ii)]
\item for some $i\in\La\beh\{j\}$ with $x(i)=1$ there is an arrow from
  $i$ to $j$,
\item $x(j)=1$ and there is no blocking symbol at $j$.
\end{enumerate}
\end{lemma}

\begin{Proof}
This follows from (\ref{infadd}) by writing
\[
m(x)=m\big(\bigvee_{i:\,x(i)=1}e_i\big)=\bigvee_{i:\,x(i)=1}m(e_i).
\]
\end{Proof}

In terms of arrows and blocking symbols, the branching map ${\tt bra}_{ij}$,
the death map ${\tt death}_i$, the voter map ${\tt vot}_{ij}$, the
coalescing random walk map ${\tt rw}_{ij}$, and the exclusion map
${\tt excl}_{ij}$ look as follows:

\begin{equation}\label{admaps}
\inputtikz{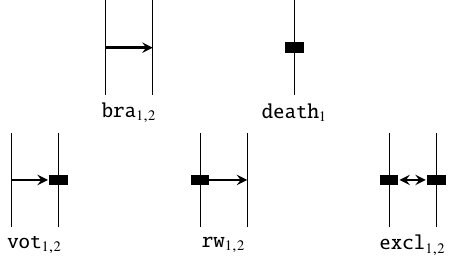}
\commentAlt{Formula (\ref{admaps})}{Drawing of the
  arrows and blocking symbols of the maps under discussion. See long
  description.}
\commentLongAlt{Formula (\ref{admaps})}{
The map ${bra}_{1,2}$ has an arrow from 1 to 2.
The map ${death}_1$ has a blocking symbol at 1.
The map ${vot}_{1,2}$ has an arrow from 1 to 2 and a blocking symbol at 2.
The map ${rw}_{1,2}$ has an arrow from 1 to 2 and a blocking symbol at 1.
The map ${exclr}_{1,2}$ has arrows between 1 and 2 in both directions and
blocking symbols at 1 and 2.}
\end{equation}

\begin{figure}[htb]
\begin{center}
\inputtikz{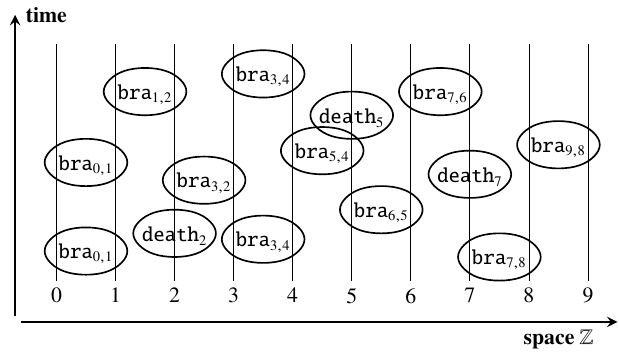}
\caption{Graphical representation of a contact process with maps.}
\label{fig:oldgrap}
\commentAlt{Figure~\ref{fig:oldgrap}}{Picture of the graphical
  representation omega. Time is plotted upwards. Branching maps and
  death maps are written on the places in space and time where they
  are applied.}
\end{center}
\end{figure}

We use our conventions of representing additive maps in terms of
arrows and blocking symbols to depict the graphical representation of
an additive interacting particle system in a more suggestive way. In
Figure~\ref{fig:contgraph}, we drew the graphical representation of a
contact process as in Figure~\ref{fig:oldgrap}.

\newpage

\begin{figure}[htb]
\begin{center}
\inputtikz{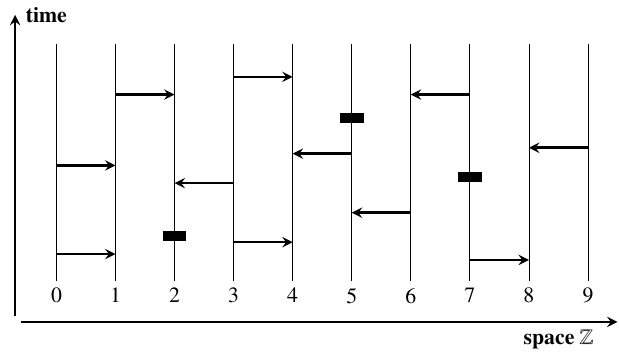}
\caption{Graphical representation of a contact process with arrows and blocking
  symbols.}
\label{fig:newgrap}
\commentAlt{Figure~\ref{fig:newgrap}}{Same as Figure~\ref{fig:oldgrap}
  except that maps are now represented by arrows and blocking
  symbols.}
\end{center}
\end{figure}

\noi
With our new conventions, the same graphical representation looks as in Figure~\ref{fig:newgrap}.

It is easy to see that the concatenation of two additive maps is again
additive. As a result, using Proposition~\ref{P:fincon}, we see that
if $(\Xb_{s,t})_{s\leq t}$ is the stochastic flow associated with the
graphical representation of an additive particle system, then the
functions $\Xb_{s,t}\cn\{0,1\}^\La\to\{0,1\}^\La$ are additive
maps. By Theorem~\ref{T:Poispart} they are also continuous. We claim
that $\Xb_{s,t}(x)(j)=1$ if and only if there is an $i\in\La$ such
that $x(i)=1$ and it is possible to walk through the graphical
representation from the space-time point $(i,s)$ to the space time
point $(j,t)$ along an upward path that may use arrows, but must avoid
the blocking symbols. We now make this claim more precise.

\begin{figure}[htb]
\begin{center}
\inputtikz{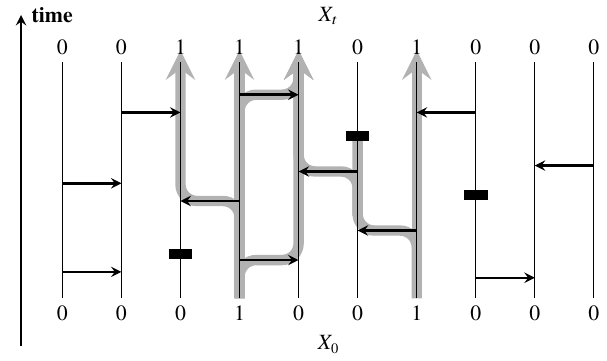}
\caption{Contact process defined by open paths.}
\label{fig:openpath}
\commentAlt{Figure~\ref{fig:openpath}}{Same as the previous figure
  except that with grey arrows it is indicated how the infection
  spreads through the graphical representation. Branching arrows yield
  new infections, at blocking symbols a site recovers.}
\end{center}
\end{figure}

For any $i,j\in\La$ and $s<u$, by definition, an \emph{open path}
\index{open path!in graphical representation} from $(i,s)$ to $(j,u)$
is a cadlag function $\ga\cn[s,u]\to\La$ such that $\ga_s=i$,
$\ga_u=j$, and
\begin{equation}\begin{array}{rl}\label{open}
{\rm(i)}&\mbox{if $\ga_{t-}\neq\ga_t$ for some $t\in(s,u]$, then there is an}\\
&\mbox{arrow from $(\ga_{t-},t)$ to $(\ga_t,t)$,}\\[5pt]
{\rm(ii)}&\mbox{there exist no $t\in(s,u]$ such that $\ga_{t-}=\ga_t$}\\
&\mbox{while there is a blocking symbol at $(\ga_t,t)$.}
\end{array}\end{equation}
We write $(i,s)\leadsto(j,u)$ \index{00pijl@$\leadsto$} if there
exists an open path from $(i,s)$ to $(j,u)$. With these definitions,
we can make our earlier claim precise. We claim that:
\begin{equation}\label{addflow}
\Xb_{s,t}(x)(j)=1\quad\mbox{iff}\quad
\exists i\in\La\mbox{ s.t.\ }x(i)=1\mbox{ and }(i,s)\leadsto(j,t).
\end{equation}
To prove (\ref{addflow}), it suffices to observe that if we define
$X_t(j):=1$ iff the condition on the right-hand side of
(\ref{addflow}) is satisfied, then the function $(X_t)_{t\geq s}$
solves the evolution equation (\ref{Xevolve}). For example, for the
graphical representation of the contact process that we earlier used
as an example, the time evolution of the process $X_t:=\Xb_{0,t}(X_0)$
$(t\geq 0)$ might look as in Figure~\ref{fig:openpath}.

Recall that if $E$ is any topological space, then a function $f\cn E\to\R$ is
called \emph{lower semi-continuous}\index{lower semi-continuous function}
if one (and hence both) of the following equivalent conditions are
satisfied:
\begin{enumerate}[(ii)]
\item $\dis\liminf_{n\to\infty}f(x_n)\geq f(x)\mbox{ whenever }x_n\to x$,
\item the level set $\{x\in E:f(x)\leq a\}$ is closed for each $a\in\R$.
\end{enumerate}

\begin{Exercise}
Show that (\ref{infadd}) remains true if $m$ is additive and lower
semi-continuous.
\end{Exercise}

\section{Additive duality}\label{S:addual}

Recall the definition of $\Si_{\rm fin}(\La)$ from the beginning of
the previous section. For countable sets $\La,\De$, we let
\[\index{0Cadd@$\Ci_{\rm add}$}
\Ci_{\rm add}\big(\Si(\La),\Si(\De)\big)
\quand\Ci_{\rm add}\big(\Si_{\rm fin}(\La),\Si_{\rm fin}(\De)\big)
\]
denote the space of continuous additive maps $m\cn\Si(\La)\to\Si(\De)$
and the space of additive maps $m\cn\Si_{\rm fin}(\La)\to\Si_{\rm fin}(\De)$,
respectively. Since $\Si_{\rm fin}(\La)$ is equipped with the discrete topology,
the latter are automatically continuous. We define a function
$\psi_{\rm add}\cn\Si(\La)\times\Si(\La)\to\{0,1\}$ by
\begin{equation}\label{psiadd}\index{0zzypsiadd@$\psi_{\rm add}$}
\psi_{\rm add}(x,y):=1_{\txt\{x\wedge y\neq\un 0\}}\qquad\big(x,y\in\Si(\La)\big),
\end{equation}
where $x\wedge y$\index{00normyy@$x\wedge y$} denotes the pointwise minimum of
$x$ and $y$ and $\un 0$ denotes the all zero configuration. We adopt the
following notation. For each $y\in\Si(\La)$, we define
$\psi_{\rm add}(\,\cdot\,,y)\cn\Si(\La)\to\{0,1\}$ by
\[
\psi_{\rm add}(\,\cdot\,,y)(x):=\psi_{\rm add}(x,y)\qquad\big(x\in\Si(\La)\big),
\]
and we let $\psi_{\rm add}(\,\ast\,,y)$ denote the restriction of
$\psi_{\rm add}(\,\cdot\,,y)$ to $\Si_{\rm fin}(\La)$. We define
$\psi_{\rm add}(x,\,\cdot\,)$ and $\psi_{\rm add}(x,\,\ast\,)$ in the
same way.

\begin{lemma}[Additive duality function]
One\label{L:psim} has
\begin{equation}\begin{array}{rr@{\,}c@{\,}l}\label{psim}
{\rm(i)}&\dis\Ci_{\rm add}\big(\Si(\La),\{0,1\}\big)&=&\dis\big\{\psi_{\rm add}(\,\cdot\,,y):y\in\Si_{\rm fin}(\La)\big\},\\[5pt]
{\rm(ii)}&\dis\Ci_{\rm add}\big(\Si_{\rm fin}(\La),\{0,1\}\big)&=&\dis\big\{\psi_{\rm add}(x,\,\ast\,):x\in\Si(\La)\big\}.
\end{array}\end{equation}
Moreover, $y\mapsto\psi_{\rm add}(\,\cdot\,,y)$ is a bijection from
$\Si_{\rm fin}(\La)$ to $\Ci_{\rm add}\big(\Si(\La),\{0,1\}\big)$ and
$x\mapsto\psi_{\rm add}(x,\,\ast\,)$ is a bijection from $\Si(\La)$
to $\Ci_{\rm add}\big(\Si_{\rm fin}(\La),\{0,1\}\big)$.
\end{lemma}

\begin{Proof}
It is straightforward to check that
$\Si(\La)\ni x\mapsto\psi_{\rm add}(x,y)\in\{0,1\}$ is additive for
each $y\in\Si(\La)$ and by symmetry an analogue statement holds for
$y\mapsto\psi_{\rm add}(x,y)$. It follows from Lemma~\ref{L:contprod}
that $\psi_{\rm add}(\,\cdot\,,y)$ is continuous if $y\in\Si_{\rm fin}(\La)$.
This proves the inclusions $\supset$ in (\ref{psim}) (i) and (ii).

To prove the converse inclusion in (\ref{psim}) (i), assume that
$\phi\cn\Si(\La)\to\{0,1\}$ is continuous and additive. Define
$y\in\Si(\La)$ by $y(i):=1$ if $\phi(e_i)=1$ and $:=0$
otherwise. Since $\phi(\un 0)=0$ we have $\phi(e_i)=0$ for all
$i\in\La\beh\Ri(\phi)$ and hence $y\in\Si_{\rm fin}(\La)$ by the
continuity of $\phi$. Now (\ref{infadd}) gives
\[
\phi(x)=\phi\big(\bigvee_{i:\,x(i)=1}e_i\big)=\bigvee_{i:\,x(i)=1}\phi(e_i)
=\psi_{\rm add}(x,y).
\]
The proof of the inclusion $\sub$ in (\ref{psim}) (ii) is similar. In
this case we can't use (\ref{infadd}) but we don't need to since
$y\in\Si_{\rm fin}(\La)$. To see that $y\mapsto\psi_{\rm add}(\,\cdot\,,y)$
and $x\mapsto\psi_{\rm add}(x,\,\ast\,)$ are
bijections, it suffices to note that if $x(i)\neq x'(i)$, then
$\psi_{\rm add}(x,e_i)\neq\psi_{\rm add}(x',e_i)$.
\end{Proof}

We now consider an interacting particle system whose generator $G$ has
a random mapping representation of the form (\ref{Gmap2}). We assume
that $S=\{0,1\}$ and all maps $m\in\Gi$ are additive. We also assume
that the rates satisfy (\ref{downsum}) so that the stochastic flow
$(\Xb_{s,t})_{s\leq t}$ and the backward stochastic flow
$(\Fb_{t,s})_{t\geq s}$ are well-defined. Since the concatenation of
two additive functions is again additive, we have that
\[
\phi\in\Ci_{\rm add}(S^\La,\{0,1\})\quad\mbox{implies}\quad
\Fb_{t,s}(\phi)\in\Ci_{\rm add}(S^\La,\{0,1\})\qquad(t\geq s).
\]
By Lemma~\ref{L:psim}, there is a one-to-one correspondence between
functions $\phi\in\Ci_{\rm add}(S^\La,\{0,1\})$ and configurations
$y\in\Si_{\rm fin}(\La)$. It follows that we can define a backward
stochastic flow $(\Yb_{t,s})_{t\geq s}$ on $\Si_{\rm fin}(\La)$ by
\begin{equation}\label{FY}
\Fb_{t,s}(\psi_{\rm add}(\,\cdot\,,y))=:\psi_{\rm add}\big(\,\cdot\,,\Yb_{t,s}(y)\big)
\qquad\big(t\geq s,\ y\in\Si_{\rm fin}(\La)\big).
\end{equation}
We will show that $(\Yb_{t,s})_{t\geq s}$ can be used to define a
Markov process $(Y_t)_{t\geq 0}$ that is itself an additive particle
system, and that is pathwise dual to the system $(X_t)_{t\geq 0}$ with
generator $G$. We first state the main facts, and then prove them.

\begin{lemma}[Dual maps]
For each\label{L:dumas} local additive map $m\cn\Si(\La)\to\Si(\La)$,
there exists a unique map $\hat m\cn\Si(\La)\to\Si(\La)$ that is dual
to $m$ with respect to the duality function $\psi_{\rm add}$, in the
sense that
\begin{equation}\label{mhm}\index{0mh@$\hat m$}
\psi_{\rm add}\big(m(x),y\big)=\psi_{\rm add}\big(x,\hat m(y)\big)
\qquad\big(x,y\in\Si(\La)\big).
\end{equation}
This dual map is also local and additive and uniquely characterized by
\begin{equation}\label{mrev}
m(e_i)(j)=1\quad\desd\quad\hat m(e_j)(i)=1\qquad(i,j\in\La).
\end{equation}
\end{lemma}

We observe that in terms of our graphical way of depicting additive
maps, formula (\ref{mrev}) has the following interpretation:
\begin{equation}\begin{array}{l}\label{gradu}
\mbox{$\hat m$ is obtained from $m$ by keeping the blocking symbols}\\
\mbox{and reversing the direction of all arrows.}
\end{array}\end{equation}
This means that the duals of the maps depicted below
Lemma~\ref{L:adloc} are given by:

\begin{equation}\label{dumapfig}
\inputtikz{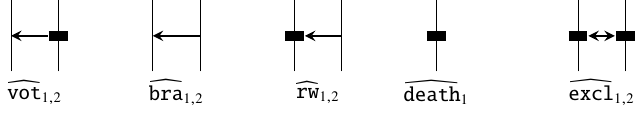}
\commentAlt{Formula (\ref{dumapfig})}{Drawing of the dual maps under discussion
 in terms of arrows and blocking symbols. See long
  description.}
\commentLongAlt{Formula (\ref{dumapfig})}{
The map $\widehat{vot}_{1,2}$ has an arrow from 2 to 1 and a blocking symbol at 2.
The map $\widehat{bra}_{1,2}$ has an arrow from 2 to 1.
The map $\widehat{rw}_{1,2}$ has an arrow from 2 to 1 and a blocking symbol at 1.
The map $\widehat{death}_1$ has a blocking symbol at 1.
The map $\widehat{exclr}_{1,2}$ has arrows between 1 and 2 in both directions and
blocking symbols at 1 and 2.}
\end{equation}

We see from this that:
\begin{equation}\begin{array}{c}\label{dumaps}
\widehat{\tt bra}_{ij}={\tt bra}_{ji},\quad\widehat{\tt death}_i={\tt death}_i,\\[5pt]
\widehat{\tt vot}_{ij}={\tt rw}_{ji},\quad\widehat{\tt rw}_{ij}={\tt vot}_{ji},\quad\widehat{\tt excl}_{ij}={\tt excl}_{ij}.
\end{array}\end{equation}

\begin{theorem}[Additive duality]
Let\label{T:addu} $G$ be the generator of an interacting particle
system $(X_t)_{t\geq 0}$ with state space $\Si(\La)$. Assume that $G$
has a random mapping representation of the form (\ref{Gmap2}) such
that all local maps $m\in\Gi$ are additive and the rates
$(r_m)_{m\in\Gi}$ satisfy (\ref{downsum}). Then
\begin{equation}\label{adduG}
\hat G:=\sum_{m\in\Gi}r_m\big\{f\big(\hat m(y)\big)-f\big(y\big)\big\}\qquad\big(y\in\Si_{\rm fin}(\La)\big)
\end{equation}
is the generator of a nonexplosive continuous-time Markov chain
$(Y_t)_{t\geq 0}$ with state space $\Si_{\rm fin}(\La)$. Let $\om$ be
a graphical representation associated with the random mapping
representation (\ref{Gmap2}) of $G$. Define a graphical representation
$\hat\om$ associated with the random mapping representation
(\ref{adduG}) of $\hat G$ by
\[
\hat\om:=\big\{(\hat m,t):(m,t)\in\om\big\}.
\]
Let $(\Xb_{s,t})_{s\leq t}$ be the stochastic flow on $\Si(\La)$
defined in terms of $\om$ as in Theorem~\ref{T:Poispart} and let
$(\Yb_{t,s})_{t\geq s}$ be the backward stochastic flow on $\Si_{\rm fin}(\La)$
defined in terms of $\hat\om$ as in Theorem~\ref{T:backflow}. Then almost surely
\begin{equation}\label{addu}
\psi_{\rm add}\big(\Xb_{s,t}(x),y\big)=\psi_{\rm add}\big(x,\Yb_{t,s}(y)\big)
\qquad\big(s\leq t,\ x\in\Si(\La),\ y\in\Si_{\rm fin}(\La)\big).
\end{equation}
If the random
mapping representation (\ref{adduG}) also satisfies (\ref{downsum}) so
that $(\Yb_{t,s})_{t\geq s}$ can be extended to $\Si(\La)$, then
(\ref{addu}) holds for all $x,y\in\Si(\La)$.
\end{theorem}

\begin{figure}[htb]
\begin{center}
\inputtikz{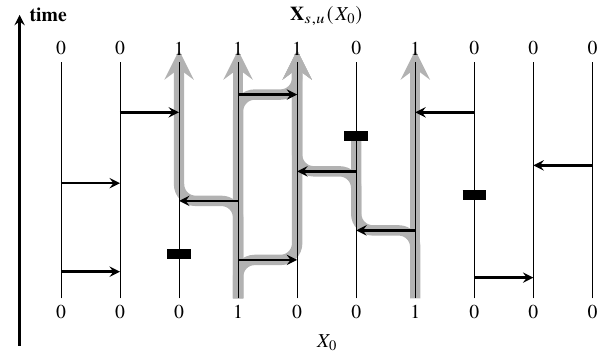}
\caption{Contact process defined by open paths.}
\label{fig:openpath2}
\commentAlt{Figure~\ref{fig:openpath2}}{Graphical representation of a
  contact process with grey arrows indicating how the infection
  spreads through the graphical representation, upwards in time.}
\end{center}
\end{figure}

Formula (\ref{addu}) says that the interacting particle systems with
generators $G$ and $\hat G$ are pathwise dual with duality function
$\psi_{\rm add}$, see the discussion in Section~\ref{S:dubasic}. This
duality has a simple graphical interpretation. We recall from
(\ref{gradu}) that we can obtain the dual of an additive local map by
reversing the direction of all arrows and keeping all blocking
symbols. In Section~\ref{S:add}, we constructed a contact process
$(X_t)_{t\geq 0}$ from its graphical representation in terms of open
paths, as in Figure~\ref{fig:openpath2}.

\begin{figure}[htb]
\begin{center}
\inputtikz{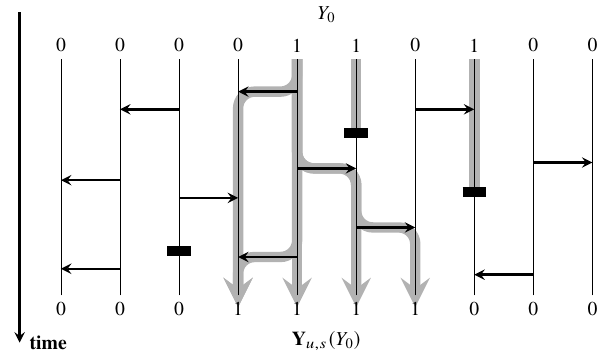}
\caption{Dual contact process defined by open paths.}
\label{fig:duflow}
\commentAlt{Figure~\ref{fig:duflow}}{Same as
  Figure~\ref{fig:openpath2} but the direction of the arrows is
  reversed and gray arrows indicate how the dual process spreads
  through the graphical representation, downwards in time.}
\end{center}
\end{figure}

Using the recipe ``reverse the arrows, keep the blocking symbols'' we
can construct the dual process $(Y_t)_{t\geq 0}$ as in Figure~\ref{fig:duflow}.

The duality relation between the forward and backward stochastic flows
then follows from the observation that
\[\begin{array}{l}
\psi_{\rm add}\big(\Xb_{s,u}(x),y\big)=1\\[5pt]
\dis\quad\desd\exists i,j\in\La\mbox{ s.t.\ }
x(i)=1,\ y(j)=1,\ (i,s)\leadsto(j,u)\\[5pt]
\dis\quad\desd\psi_{\rm add}\big(x,\Yb_{u,s}(y)\big)=1.
\end{array}\]
In our previous example of the contact process, the dual process is
also a contact process, but in general, the dual process can have a
different dynamics from the forward in time process. For example, we
see from (\ref{dumaps}) that the additive dual of the voter model is a
system of coalescing random walks.

We still need to prove Lemma~\ref{L:dumas} and Theorem~\ref{T:addu}.\med

\begin{Proof}[of Lemma~\ref{L:dumas}]
If $m(e_i)(j)=1$ for some $i\neq j$, then by the fact that $m(\un
0)=0$ we see that $j\in\Di(m)$ and $i\in\Ri(m[j])$. Also, if
$m(e_i)(i)=0$, then $i\in\Di(m)$. It follows that the
graphical representation of an additive local map contains only
finitely many arrows and blocking symbols. Setting $x=e_i$ and $y=e_j$
in (\ref{mhm}) gives (\ref{mrev}) so the latter is clearly
necessary. Using the recipe ``reverse the arrows, keep the blocking
symbols'' we can find a local map $\hat m$ such that (\ref{mrev})
holds. Using (\ref{infadd}), which is applicable since both $m$ and
$\hat m$ are continuous, we observe that
\[
\psi_{\rm add}\big(m(x),y\big)
=\psi_{\rm add}\Big(\bigvee_{i:\,x(i)=1}m(e_i),y\Big)
=\bigvee_{i:\,x(i)=1}\bigvee_{j:\,y(j)=1}1_{\txt\{m(e_i)(j)=1\}},
\]
which by (\ref{mrev}) and the same argument backwards is equal to
$\psi_{\rm add}\big(x,\hat m(y)\big)$. Since
$\hat m(y)(i)=\psi_{\rm add}(e_i,\hat m(y))=\psi_{\rm add}(m(e_i),y)$,
we see that $\hat m$ is the unique map from $\Si(\La)$ into itself
that is dual to $m$ with respect to $\psi_{\rm add}$.
\end{Proof}

\begin{Proof}[of Theorem~\ref{T:addu}]
Under the condition (\ref{downsum}), it has been shown in
Proposition~\ref{P:back} that
\[
Hf(\phi):=\sum_{m\in\Gi}r_m\big\{f(\phi\circ m)-f(\phi)\big\}
\]
is the generator of a nonexplosive continuous-time Markov chain,
called the backward in time process, with state space
$\Ci(\Si(\La),\{0,1\})$. In (\ref{Fback}) we used the graphical
representation $\om$ to define a backward stochastic flow
$(\Fb_{t,s})_{t\geq s}$ associated with this backward in time
process. Since the composition of two additive maps is additive, this
backward stochastic flow maps the space $\Ci_{\rm add}(\Si(\La),\{0,1\})$
into itself. By Lemma~\ref{L:psim}, $y\mapsto\Ci_{\rm add}(\Si(\La),\{0,1\})$
is a bijection from $\Si_{\rm fin}(\La)$ to $\Ci_{\rm add}(\Si(\La),\{0,1\})$
which allows us to define $(\Yb_{t,s})_{t\geq s}$ in terms of
$(\Fb_{t,s})_{t\geq s}$ as in (\ref{FY}). In view of (\ref{backevol}),
this means that for each $u\in\R$ and $y\in\Si_{\rm fin}(\La)$, the
function $(Y_t)_{t\leq u}$ defined as $Y_t:=\Yb_{u,t}(y)$ $(t\leq u)$
is the unique piecewise constant, right-continuous solution of the
evolution equation
\[
Y_{t-}=\left\{\begin{array}{ll}
\dis Y'\mbox{ where }\psi_{\rm add}(\,\cdot\,,Y')
:=\psi_{\rm add}(\,\cdot\,,Y_t)\circ m\quad&\mbox{if }(m,t)\in\om,\\[5pt]
Y_t\quad&\mbox{otherwise.}
\end{array}\right.
\]
We observe that for any $x\in\Si(\La)$ and $y\in\Si_{\rm fin}(\La)$,
\[
\psi_{\rm add}(\,\cdot\,,y)\circ m(x)=\psi_{\rm add}\big(m(x),y\big)
=\psi_{\rm add}\big(x,\hat m(y)\big),
\]
so $(\Yb_{t,s})_{t\geq s}$ is the backward stochastic flow associated
with the continuous-time Markov chain with generator $\hat G$ as in
(\ref{adduG}). In particular, the latter is nonexplosive since the
backward in time process is. Now (\ref{backF}) and (\ref{FY}) imply
that
\[
\psi_{\rm add}\big(\Xb_{s,t}(x),y\big)
=\psi_{\rm add}(\,\cdot\,,y)\circ\Xb_{s,t}(x)
=\Fb_{t,s}\big(\psi_{\rm add}(\,\cdot\,,y)\big)(x)
=\psi_{\rm add}\big(x,\Yb_{t,s}(y)\big)
\]
for all $t\geq s$, $x\in\Si(\La)$, and $y\in\Si_{\rm fin}(\La)$,
proving (\ref{addu}).

Assume that the random mapping representation (\ref{adduG}) also
satisfies (\ref{downsum}) so that $(\Yb_{t,s})_{t\geq s}$ can be
extended to $\Si(\La)$. For each $y\in\Si(\La)$ we can find
$y_n\in\Si_{\rm fin}(\La)$ such that $y_n\up y$. Then
$\Yb_{t,s}(y_n)\up\Yb_{t,s}(y)$ by the continuity of $\Yb_{t,s}$
(proved in Theorem~\ref{T:Poispart}) and the monotonicity of
$\Yb_{t,s}$ (which follows from its additivity), so taking the limit
in (\ref{addu}) we see that the latter holds for all $x,y\in\Si(\La)$.
\end{Proof}

We conclude this section with a simple lemma that is important in
applications of Theorem~\ref{T:addu} and in particular of the duality
relation (\ref{addu}). It says that the values of $\E[\psi_{\rm add}(X_t,y)]$
for all $y\in\Si_{\rm fin}(\La)$ determine the law of $X_t$ uniquely.

\begin{lemma}[Distribution determining functions]
The\label{L:adddist} class of functions
$\{\psi_{\rm add}(\,\cdot\,,y):y\in\Si_{\rm fin}(\La)\}$ is distribution
determining on $\Si(\La)$.
\end{lemma}

\begin{Proof}
We may equivalently prove that the functions
$g_y(x):=1-\psi_{\rm add}(\,\cdot\,,y)=1_{\{x\wedge y=\un 0\}}$ are distribution
determining. Since $1_{\{x\wedge e_i=\un 0\}}=1-x(i)$, the class
$\{g_y:y\in\Si_{\rm fin}(\La)\}$ separates points, and since
$g_yg_{y'}=g_{y\vee y'}$, this class is closed under products. The
claim now follows from Lemma~\ref{L:SW}.
\end{Proof}

\begin{Exercise}
Let $\Li_{\rm add}(\Si(\La),\{0,1\})$ denote the space of lower
semi-continu\-ous additive maps $m\cn\Si(\La)\to\{0,1\}$. In analogy
with Lemma~\ref{L:psim}, show that
\[
\Li_{\rm add}\big(\Si(\La),\{0,1\}\big)
=\big\{\psi_{\rm add}(\,\cdot\,,y):y\in\Si(\La)\big\}.
\]
Show that under the assumptions of Theorem~\ref{T:addu}, the backward
stochastic flow $(\Fb_{t,s})_{t\geq s}$ maps the space
$\Li_{\rm add}(\Si(\La),\{0,1\})$ into itself. Use this to show that the
definition of the backward stochastic flow $\Yb_{t,s}(y)$ in
(\ref{FY}) can be extended to $y\in\Si(\La)$, even if the random
mapping representation (\ref{adduG}) does not satisfy the summability
condition~\ref{downsum}).
\end{Exercise}

\begin{Exercise}
Generalize Lemma~\ref{L:dumas} to lower semi-continuous additive maps,
by showing that each lower semi-continuous additive map
$m\cn\Si(\La)\to\Si(\La)$ has a unique dual map $\hat m$ with respect
to the duality function $\psi_{\rm add}$, and that this dual map $\hat m$
is also lower semi-continuous and additive.
\end{Exercise}

\begin{Exercise}
Give\label{E:Kupdo} an example of an additive particle system for
which the forward generator $G$ satisfies the summability condition
(\ref{downsum}) but the dual generator $\hat G$ does not. \emph{Hint:}
consider a contact process and its dual contact process on a binary
tree where for the forward process all infections point away from the
root and for the dual process all infections point in the direction of
the root.
\end{Exercise}

\begin{Exercise}
Let\label{E:backexpl} $(\beta_i)_{i\in\N}$ be positive constants and
consider the interacting particle system with lattice $\La:=\N$, local
state space $S:=\{0,1\}$, and generator
\[
Gf(x):=\sum_{i=0}^\infty\beta_i\big\{f\big({\tt bra}_{i+1,i}(x)\big)
-f\big(x\big)\big\},
\]
where ${\tt bra}_{j,i}$ is the branching map defined in
(\ref{bramap}). This process is additive and the generator of the dual
process is
\[
\hat Gf(y)=\sum_{i=0}^\infty\beta_i\big\{f\big({\tt bra}_{i,i+1}(y)\big)
-f\big(y\big)\big\}.
\]
If the rates $\bet_i$ are not bounded as a function of $i$, then
(\ref{downsum}) fails so Theorem~\ref{T:Poispart} is not
applicable. Nevertheless, condition (\ref{Hwell}) is still satisfied
so by Lemma~\ref{L:Hwell} the backward in time process is still a
well-defined continuous-time Markov chain even though it may be
explosive if $\bet$ grows too fast as a function of $i$. (Compare
Exercise~\ref{E:explos}.) In this case, by Theorem~\ref{T:stochflow},
we can still construct the backward stochastic flow
$(\Yb_{t,s})_{t\geq s}$ associated with the continuous-time Markov
chain with generator $\hat G$. Here $\Yb_{t,s}$ is a map from the
space $\Si_{\rm fin}(\N)\cup\{\infty\}$ into itself, where
$\Yb_{t,s}(y)=\infty$ means that the dual process started at time $t$
and run backwards to time $s<t$ has exploded. We can use this to
define a function $X\cn[0,\infty)\to\{0,1\}^\N$ by
\[
X_t(i):=\left\{\begin{array}{ll}
0\quad&\mbox{if }\ \Yb_{t,0}(e_i)\in\Si_{\rm fin}(\N),\\[5pt]
1\quad&\mbox{if }\ \Yb_{t,0}(e_i)=\infty.
\end{array}\right.\qquad(i\in\N,\ t\geq 0).
\]
Show that $(X_t)_{t\geq 0}$ solves the evolution equation
(\ref{Xevolve}) with initial state $X_0=\un 0$. Since the all zero
function is another solution, this shows that for this interacting
particle system, solutions to the the evolution equation
(\ref{Xevolve}) are not unique if the backward in time process is
explosive.
\end{Exercise}

\section{Cancellative systems and their duals}\label{S:canc}

Let $\oplus$ denote addition modulo two, that is,
\[\index{00normz@$x\oplus y$}
0\oplus 0:=0,\quad
0\oplus 1:=1,\quad
1\oplus 0:=1,\quand
1\oplus 1:=0.
\]
Let $\La$ and $\De$ be countable sets. For $x,y\in\Si(\La)$, we define
$(x\oplus y)(i):=x(i)\oplus y(i)$ $(i\in\La)$ in a pointwise way. By
definition, a map $m\cn\Si(\La)\to\Si(\De)$ is
\emph{cancellative}\index{cancellative!map} iff
\[\begin{array}{r@{\ }l}
{\rm(i)}&m(\un 0)=\un 0,\\[5pt]
{\rm(ii)}&m(x\oplus y)=m(x)\oplus m(y)\qquad\big(x,y\in\Si(\La)\big).
\end{array}\]
The same definition applies to maps $m\cn\Si_{\rm fin}(\La)\to\Si_{\rm fin}(\De)$,
where in this case (ii) only needs to hold for
$x,y\in\Si_{\rm fin}(\La)$. An interacting particle system is called
\emph{cancellative}\index{cancellative!interacting particle system} if
its generator can be represented in cancellative local maps. Examples
of cancellative maps are:
\begin{itemize}
\addtolength\itemsep{7pt}
\item The voter map ${\tt vot}_{ij}$ defined in (\ref{votmap}).
\item The death map ${\tt death}_i$ defined in (\ref{deathmap}).
\item The exclusion map ${\tt excl}_{ij}$ defined in (\ref{exclmap}).
\item The annihilating random walk map ${\tt arw}_{ij}$ defined in
  (\ref{annmap}).
\item The annihilating branching map ${\tt abra}_{ij}$ defined in
  (\ref{branmap}) below.
\end{itemize}
Here, we define an \emph{annihilating branching map}
\index{annihilating branching map}
${\tt abra}_{ij}\cn\{0,1\}^\La\to\{0,1\}^\La$ by
\begin{equation}\label{branmap}\index{0abra@${\tt abra}_{ij}$}
{\tt abra}_{ij}(x)(k):=\left\{\begin{array}{ll}
x(i)\oplus x(j)\quad&\mbox{if }k=j,\\[5pt]
x(k)\quad&\mbox{otherwise.}
\end{array}\right.
\end{equation}
For countable sets $\La,\De$, we let
\[\index{0Ccanc@$\Ci_{\rm canc}$}
\Ci_{\rm canc}\big(\Si(\La),\Si(\De)\big)
\quand\Ci_{\rm canc}\big(\Si_{\rm fin}(\La),\Si_{\rm fin}(\De)\big)
\]
denote the space of continuous cancellative maps
$m\cn\Si(\La)\to\Si(\De)$ and the space of cancellative maps
$m\cn\Si_{\rm fin}(\La)\to\Si_{\rm fin}(\De)$, respectively. Similar
to what we did in the additive case, we can describe a continuous
cancellative map $m\cn\Si(\La)\to\Si(\La)$ in terms of arrows and
blocking symbols:
\begin{itemize}
\item For each $i,j\in\La$ with $i\neq j$ such that $m(e_i)(j)=1$, we
  draw an arrow from $i$ to $j$.
\item For each $i\in\La$ such that $m(e_i)(i)=0$, we draw a blocking
  symbol \blok at $i$.
\end{itemize}
The following lemma says that continuous cancellative maps are fully
described by their arrows and blocking symbols.

\begin{lemma}[Graphical description]
Let\label{L:adcanc} $m\cn\Si(\La)\to\Si(\La)$ be a continuous
cancellative map and let $x\in\{0,1\}^\La$. For each $j\in\La$, let
$R_j$ denote the set of $i\in\La$ such that either $i=j$ and there is
no blocking symbol at $i$, or $i\neq j$ and there is an arrow from $i$
to $j$. Then $R_j$ is finite and
\[
m(x)(j)=1\quad\desd\quad\big|\big\{i\in R_j:x(i)=1\big\}\big|\mbox{ is odd.}
\]
\end{lemma}

\begin{Proof}
Since $m(\un 0)=0$, we see that $R_j\sub\Ri(m[j])$, where by the
continuity of $m$, the latter is finite for each $j\in\La$. Using the
fact that we can change $x$ outside $\Ri(m[j])$ without changing
$m(x)(j)$, we see that
\[
m(x)(j)=m[j]\big(\bigoplus_{i\in\Ri(m[j]):\,x(i)=1}e_i\big)
=\bigoplus_{i\in\Ri(m[j]):\,x(i)=1}m[j](e_i),
\]
which is one if and only if $\{i\in R_j:x(i)=1\big\}$ has an odd
number of elements.
\end{Proof}

Every graphical representation involving arrows and blocking symbols
that can be used to define an additive particle system can also be
used to define a cancellative particle system. The cancellative maps
mentioned at the beginning of this section have the following
representations in terms of arrows and blocking symbols:

\begin{equation}\label{cancmaps}
\inputtikz{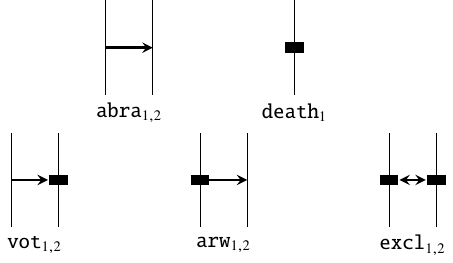}
\commentAlt{Formula (\ref{cancmaps})}{Drawing of the arrows and blocking
 symbols of the maps under discussion. See long
  description.}
\commentLongAlt{Formula (\ref{cancmaps})}{
The map ${abra}_{1,2}$ has an arrow from 1 to 2.
The map ${death}_1$ has a blocking symbol at 1.
The map ${vot}_{1,2}$ has an arrow from 1 to 2 and a blocking symbol at 2.
The map ${arw}_{1,2}$ has an arrow from 1 to 2 and a blocking symbol at 1.
The map ${exclr}_{1,2}$ has arrows between 1 and 2 in both directions and
blocking symbols at 1 and 2.}
\end{equation}

\begin{figure}[htb]
\begin{center}
\inputtikz{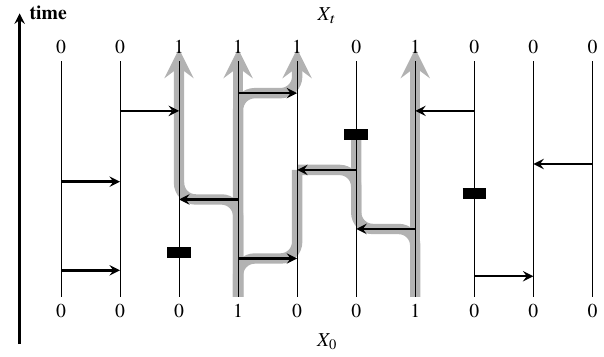}
\caption{Graphical representation of a cancellative version of the
  contact process.}
\label{fig:canccont}
\commentAlt{Figure~\ref{fig:canccont}}{The same arrows and blocking
  symbols that previously defined a contact process here define a very
  similar process, indicated by gray arrows. An infection of an
  already infected site leads to the site becoming vacant.}
\end{center}
\end{figure}

If we interpret the graphical representation of a contact process in a
cancellative way, then it becomes a graphical representation for an
interacting particle system involving the annihilating branching map
${\tt abra}_{ij}$ and the death map ${\tt death}_i$, see
Figure~\ref{fig:canccont}. This system has been studied in
\cite{BDD91}.

We define a \emph{cancellative duality function}
\index{cancellative!duality function}
$\psi_{\rm canc}\cn\Si(\La)\times\Si_{\rm fin}(\La)\to\{0,1\}$ by
\begin{equation}\label{psicanc}\index{0zzypsicanc@$\psi_{\rm canc}$}
\psi_{\rm canc}(x,y):=\bigoplus_{i\in\La}x(i)y(i)
\qquad\big(x\in\Si(\La),\ y\in\Si_{\rm fin}(\La)\big).
\end{equation}
Note that since $y\in\Si_{\rm fin}(\La)$, all but finitely many of the
summands are zero, so the infinite sum modulo two is
well-defined. Unlike in the additive case, there is no way to make
sense of $\psi_{\rm canc}(x,y)$ for general
$x,y\in\Si(\La)$.\footnote{For interacting particle systems on
$\La=\Z$, it is sometimes useful to consider the case that
$\sup\{i:x(i)=1\}<\infty$ and $\inf\{i:y(i)=1\}>-\infty$. Clearly,
$\psi_{\rm canc}(x,y)$ is well-defined for such $x,y$, even though
both may be infinite.} The following lemma is similar to
Lemma~\ref{L:psim}. Below, $\psi_{\rm add}(\,\cdot\,,y)$ denotes the
map $\Si(\La)\ni x\mapsto\psi_{\rm add}(x,y)$ and $\psi_{\rm add}(\,\ast\,,y)$
denotes its restriction to $\Si_{\rm fin}(\La)$.

\begin{lemma}[Cancellative duality function]
One\label{L:psimc} has
\begin{equation}\begin{array}{rr@{\,}c@{\,}l}\label{psimc}
{\rm(i)}&\dis\Ci_{\rm canc}\big(\Si(\La),\{0,1\}\big)&=&\dis\big\{\psi_{\rm canc}(\,\cdot\,,y):y\in\Si_{\rm fin}(\La)\big\},\\[5pt]
{\rm(ii)}&\dis\Ci_{\rm canc}\big(\Si_{\rm fin}(\La),\{0,1\}\big)&=&\dis\big\{\psi_{\rm canc}(x,\,\ast\,):x\in\Si(\La)\big\}.
\end{array}\end{equation}
Moreover, $y\mapsto\psi_{\rm canc}(\,\cdot\,,y)$ is a bijection from
$\Si_{\rm fin}(\La)$ to $\Ci_{\rm canc}\big(\Si(\La),\{0,1\}\big)$ and
$x\mapsto\psi_{\rm canc}(x,\,\ast\,)$ is a bijection from $\Si(\La)$
to $\Ci_{\rm canc}\big(\Si_{\rm fin}(\La),\{0,1\}\big)$.
\end{lemma}

\begin{Proof}
The proof is almost the same as the proof of Lemma~\ref{L:psim}, only
the proof of the inclusion $\sub$ in (\ref{psimc}) (i) is a bit more
complicated than in the additive case. To prove this inclusion, assume
that $\phi\cn\Si(\La)\to\{0,1\}$ is continuous and
cancellative. Define $y\in\Si(\La)$ by $y(i):=1$ if $\phi(e_i)=1$ and
$:=0$ otherwise. Since $\phi(\un 0)=0$ we have $\phi(e_i)=0$ for all
$i\in\La\beh\Ri(\phi)$ and hence $y\in\Si_{\rm fin}(\La)$ by the
continuity of $\phi$. Since we can change $x$ outside $\Ri(\phi)$
without changing $\phi(x)$, and $\phi(e_i)=0$ for
$i\in\La\beh\Ri(\phi)$, we now have
\[
\phi(x)=\phi\big(\bigoplus\subb{i\in\Ri(\phi)}{x(i)=1}e_i\big)=\bigoplus\subb{i\in\Ri(\phi)}{x(i)=1}\phi(e_i)=\bigoplus\subb{i\in\La}{x(i)=1}\phi(e_i)=\psi_{\rm canc}(x,y).
\]
The rest of the proof is the same as the proof of Lemma~\ref{L:psim}.
%
\end{Proof}

We now consider an interacting particle system whose generator $G$ has
a random mapping representation of the form (\ref{Gmap2}) such that
all local maps $m\in\Gi$ are cancellative. We also assume that the
rates satisfy (\ref{downsum}) so that the stochastic flow
$(\Xb_{s,t})_{s\leq t}$ and the backward stochastic flow
$(\Fb_{t,s})_{t\geq s}$ are well-defined. Since the concatenation of
two cancellative functions is again cancellative, $(\Fb_{t,s})_{t\geq s}$
maps the space $\Ci_{\rm canc}(S^\La,\{0,1\})$ into itself, so
using Lemma~\ref{L:psimc}, just as we did in the additive case (recall
(\ref{FY})), we can define a backward stochastic flow
$(\Yb_{t,s})_{t\geq s}$ on $\Si_{\rm fin}(\La)$ by
\begin{equation}\label{FYc}
\Fb_{t,s}(\psi_{\rm canc}(\,\cdot\,,y))
=:\psi_{\rm canc}\big(\,\cdot\,,\Yb_{t,s}(y)\big)
\qquad\big(t\geq s,\ y\in\Si_{\rm fin}(\La)\big).
\end{equation}
We will see that by the same arguments as in the additive case,
$(\Yb_{t,s})_{t\geq s}$ is the backward stochastic flow of a
cancellative particle system $(Y_t)_{t\geq 0}$ that is pathwise dual
to the system $(X_t)_{t\geq 0}$ with generator $G$. The following
lemma, that will be proved below, is similar to Lemma~\ref{L:dumas}.

\begin{lemma}[Dual maps]
For each\label{L:dumasc} local cancellative map
$m\cn\Si(\La)\to\Si(\La)$, there exists a unique map $\ti
m\cn\Si(\La)\to\Si(\La)$ that is dual to $m$ with respect to the
duality function $\psi_{\rm canc}$, in the sense that
\begin{equation}\label{mhmc}\index{0mg@$\ti m$}
\psi_{\rm canc}\big(m(x),y\big)=\psi_{\rm canc}\big(x,\ti m(y)\big)
\qquad\big(x,y\in\Si(\La),\ |x|\wedge|y|<\infty\big).
\end{equation}
This dual map is also local and cancellative and uniquely characterized by
\begin{equation}\label{mrevc}
m(e_i)(j)=1\quad\desd\quad\ti m(e_j)(i)=1\qquad(i,j\in\La).
\end{equation}
\end{lemma}

In terms of our graphical way of depicting cancellative maps, formula
(\ref{mrev}) can be described by saying that, just as in the additive
case,
\begin{equation}\begin{array}{l}\label{graduc}
\mbox{$\ti m$ is obtained from $m$ by keeping the blocking symbols}\\
\mbox{and reversing the direction of all arrows.}
\end{array}\end{equation}
Using this, we see that the cancellative duals of the maps listed at
the beginning of the section are given by:
\begin{equation}\begin{array}{c}\label{dumapsc}
\widetilde{\tt abra}_{ij}={\tt abra}_{ji},
\quad\widetilde{\tt death}_i={\tt death}_i,\\[5pt]
\widetilde{\tt vot}_{ij}={\tt arw}_{ji},
\quad\widetilde{\tt arw}_{ij}={\tt vot}_{ji},
\quad\widetilde{\tt excl}_{ij}={\tt excl}_{ij}.
\end{array}\end{equation}
Note that the voter map is both additive and cancellative, and has two
different dual maps depending on whether we are considering additive
or cancellative duality. The main result about cancellative duality is
the following analogue of Theorem~\ref{T:addu}.

\begin{theorem}[Cancellative duality]
Let\label{T:candu} $G$ be the generator of an interacting particle
system $(X_t)_{t\geq 0}$ with state space $\Si(\La)$. Assume that $G$
has a random mapping representation of the form (\ref{Gmap2}) such
that all local maps $m\in\Gi$ are cancellative and the rates
$(r_m)_{m\in\Gi}$ satisfy (\ref{downsum}). Then
\begin{equation}\label{canduG}
\ti G:=\sum_{m\in\Gi}r_m\big\{f\big(\ti m(y)\big)-f\big(y\big)\big\}\qquad\big(y\in\Si_{\rm fin}(\La)\big)
\end{equation}
is the generator of a nonexplosive continuous-time Markov chain
$(Y_t)_{t\geq 0}$ with state space $\Si_{\rm fin}(\La)$. Let $\om$ be
a graphical representation associated with the random mapping
representation (\ref{Gmap2}) of $G$. Define a graphical representation
$\ti\om$ associated with the random mapping representation
(\ref{canduG}) of $\ti G$ by
\[
\ti\om:=\big\{(\ti m,t):(m,t)\in\om\big\}.
\]
Let $(\Xb_{s,t})_{s\leq t}$ be the stochastic flow on $\Si(\La)$
defined in terms of $\om$ as in Theorem~\ref{T:Poispart} and let
$(\Yb_{t,s})_{t\geq s}$ be the backward stochastic flow on $\Si_{\rm fin}(\La)$
defined in terms of $\ti\om$ as in Theorem~\ref{T:backflow}. Then almost surely
\begin{equation}\label{candu}
\psi_{\rm canc}\big(\Xb_{s,t}(x),y\big)=\psi_{\rm canc}\big(x,\Yb_{t,s}(y)\big)
\qquad\big(s\leq t,\ x\in\Si(\La),\ y\in\Si_{\rm fin}(\La)\big).
\end{equation}
If the random mapping representation (\ref{adduG}) also satisfies
(\ref{downsum}), then $(\Yb_{t,s})_{t\geq s}$ can be extended to
$\Si(\La)$ and $(\Xb_{s,t})_{s\leq t}$ maps the space $\Si_{\rm fin}(\La)$
into itself. In this case (\ref{candu}) also holds for
$x\in\Si_{\rm fin}(\La)$ and $y\in\Si(\La)$.
\end{theorem}

\begin{Proof}[of Lemma~\ref{L:dumasc}]
As in the additive case, it is easy to see that a cancellative map is
local if and only if it is defined by finitely many arrows and
blocking symbols. Using the recipe ``reverse the arrows, keep the
blocking symbols'' we can find a local map $\ti m$ such that
(\ref{mrevc}) holds. To see that it is dual to $m$ in the sense of
(\ref{mhmc}), let $\De$ be the set of lattice points where a blocking
symbol is located or that are the starting point or endpoint of an
arrow. Let $x'$ be the restriction of $x$ to $\De$, that is,
$x'(i):=x(i)$ if $i\in\De$ and $:=0$ otherwise. Similarly, let $y'$
denote the restriction of $y$ to $\De$ and let $x''$ and $y''$ denote
the restrictions of $x$ and $y$ to $\La\beh\De$. Then
\[
\psi_{\rm canc}\big(m(x),y\big)
=\psi_{\rm canc}\big(m(x'),y'\big)\oplus\psi_{\rm canc}\big(x'',y''\big),
\]
where
\[
\psi_{\rm canc}\big(m(x'),y'\big)
=\psi_{\rm canc}\Big(\bigoplus_{i:\,x'(i)=1}m(e_i),y'\Big)
=\bigoplus_{i:\,x'(i)=1}\bigoplus_{j:\,y'(j)=1}1_{\txt\{m(e_i)(j)=1\}}.
\]
Rewriting $\psi_{\rm canc}\big(x,\ti m(y)\big)$ in the same way, using
(\ref{mrevc}), we see that $\ti m$ is dual to $m$ in the sense of
(\ref{mhmc}). As in the additive case, it is easy to see that $\ti m$
is uniquely determined by (\ref{mhmc}).
\end{Proof}

\begin{Proof}[of Theorem~\ref{T:candu}]
This is completely the same as the proof of Theorem~\ref{T:addu},
except at the very end. If the random mapping representation
(\ref{canduG}) also satisfies (\ref{downsum}), then by using what is
already proved with the roles of the forward and backward process
reversed, we see that $(\Xb_{s,t})_{s\leq t}$ maps the space $\Si_{\rm fin}(\La)$
into itself and (\ref{candu}) also holds for
$x\in\Si_{\rm fin}(\La)$ and $y\in\Si(\La)$.
\end{Proof}

The following lemma is similar to Lemma~\ref{L:adddist}.

\begin{lemma}[Distribution determining functions]
The\label{L:cancdist} class of functions
$\{\psi_{\rm canc}(\,\cdot\,,y):y\in\Si_{\rm fin}(\La)\}$ is distribution
determining on $\Si(\La)$.
\end{lemma}

\begin{Proof}
We may equivalently prove that the functions
\[
g_y(x):=1-2\psi_{\rm canc}(x,y)=(-1)^{\sum_ix(i)y(i)}
\]
are distribution determining. Since $g_yg_{y'}=g_{y\oplus y'}$, the
class $\{g_y:y\in\Si_{\rm fin}(\La)\}$ is closed under products and
since $g_{e_i}(x)=(-1)^{x(i)}$ this class separates points. The claim
now follows from Lemma~\ref{L:SW}.
\end{Proof}

Additive and cancellative duality are so similar that one wonders if
they can be treated in a unified way. This is indeed the case. For
local state spaces with three or more elements, an obvious thing one
can do is to replace the addition modulo two from cancellative systems
by addition modulo three or more. There are also less obvious
possibilities. The paper \cite{LS21} explores dualities where
$(\{0,1\},\vee)$ or $(\{0,1\},\oplus)$ are replaced by commutative
monoids or semirings.

Some models that a priori do not look like cancellative systems turn
out to be representable in cancellative maps. An example is the
Neuhauser--Pacala model,\index{Neuhauser--Pacala model} defined by its
transition rates in (\ref{NP99}). We define a \emph{rebellious
map}\index{rebellious map} by
\begin{equation}\label{rebelmap}\index{0rebel@${\tt rebel}_{ijk}$}
{\tt rebel}_{ijk}(x)(l):=\left\{\begin{array}{ll}
x(i)\oplus x(j)\oplus x(k)\quad&\mbox{if }l=k,\\[5pt]
x(l)\quad&\mbox{otherwise.}
\end{array}\right.
\end{equation}
In words, this says that $x(k)$ changes its state if $x(i)\neq x(j)$.

\begin{Exercise}
Show that the map ${\tt rebel}_{ijk}$ is cancellative. Show that the generator
of the Neuhauser--Pacala model defined in (\ref{NP99}) can be represented as
\[\begin{array}{r@{\,}c@{\,}l}
\dis G_{\rm NP}f(x)&=&\dis\frac{\al}{|\Ni_i|}\sum_i\sum_{j\in\Ni_i}
\big\{f\big({\tt vot}_{ji}(x)\big)-f\big(x\big)\big\}\\[5pt]
&=&\dis\frac{1-\al}{|\Ni_i|^2}\sum_i\sum_{k,j\in\Ni_i}
\big\{f\big({\tt rebel}_{kji}(x)\big)-f\big(x\big)\big\}.
\end{array}\]
\end{Exercise}

\begin{Exercise}
In\label{E:threshold1} the \emph{threshold voter
model},\index{threshold voter model} the site $i$ changes its type
$x(i)$ from 0 to 1 with rate one as long as at least one site in its
neighborhood $\Ni_i$ has type 1, and likewise, $i$ flips from 1 to 0
with rate one as long as at least one site in $\Ni_i$ has type 0. Show
that the generator of the threshold voter model can be written as
\[
G_{\rm thres}f(x)=2^{-|\Ni_i|+1}\;
\sum_i\!\!\sum\subb{\De\sub\Ni_i\cup\{i\}}{|\De|\mbox{ is even}}\!\!\!
\big\{f\big(m_{\De,i}(x)\big)-f\big(x\big)\big\},
\]
where $m_{\De,i}$ is the cancellative map defined by
\[
m_{\De,i}(x)(k):=\left\{\begin{array}{ll}
x(i)\oplus\bigoplus_{j\in\De}x(j)\quad&\mbox{if }k=i,\\[5pt]
x(k)\quad&\mbox{otherwise.}\end{array}\right.
\]
Cancellative duality for the threshold voter model is used extensively
in \cite{Han99}.
\end{Exercise}

\begin{Exercise}
Show\label{E:threshold2} that the threshold voter model is monotone.
\end{Exercise}

\section{Lloyd--Sudbury duality}\label{S:other}

The additive systems duality function (\ref{psiadd}) and cancellative
systems duality function (\ref{psicanc}) are not the only choices of
$\psi$ that lead to useful dualities. There are two approaches to
finding useful duality functions: the pathwise approach, that aims to
find dualities between stochastic flows in the sense of
(\ref{flowdual}), and the algebraic approach, that only aims to prove
distributional relations of the form (\ref{dual}). There has been a
lot of recent work on the algebraic approach, starting with the
paper \cite{GK+09}, linking dualities to representations of Lie
algebras. For an overview of this work, we refer to \cite{GR25}. We
will below present some older results, based on the algebraic
approach, due to Lloyd and Sudbury \cite{SL95,SL97,Sud99}.

The pathwise approach always depends on finding a clever stochastic
flow and then finding a suitable space of functions on $S^\La$ that is
mapped into itself by the stochastic flow of the backward in time
process. As we have seen, for additive and cancellative systems, the
spaces $\Ci_{\rm add}(\{0,1\}^\La,\{0,1\})$ and
$\Ci_{\rm canc}(\{0,1\}^\La,\{0,1\})$ are invariant, and this naturally leads
to additive and cancellative duality.

To explain a bit about the algebraic approach, which only aims to
prove relations of the form (\ref{dual}) without proving a duality of
stochastic flows, for technical simplicity, for the remainder of this
section we will restrict ourselves to finite state spaces. In general,
when trying to prove a duality for interacting particle systems on
infinite lattices, it is often a good idea to first prove the result
on finite lattices and then extend it to infinite lattices using
approximation results such as Theorem~\ref{T:Fellim} and
Corollary~\ref{C:Partlim}. We will demonstrate this method in
Section~\ref{S:convot} below.

As we have done before, we write $\E^x$ (respectively $\E^y$) to denote
expectation with respect to the law of the process $X$ started in
$X_0=x$ (respectively started in $Y_0=y$).

\begin{lemma}[Duality of finite Markov processes]
Let\label{L:dumat} $(X_t)_{t\geq 0}$ and $(Y_t)_{t\geq 0}$ be Markov
processes with finite state spaces $S$ and $R$, generators $G$ and
$H$, and Markov semigroups $(P_t)_{t\geq 0}$ and $(Q_t)_{t\geq 0}$. Then one has
\begin{equation}\label{dual2}
\E^x\big[\psi(X_t,y)\big]=\E^y\big[\psi(x,Y_t)\big]
\qquad(x\in S,\ y\in R,\ t\geq 0)
\end{equation}
if and only if
\begin{equation}\label{Gdual}
G\psi(\,\cdot\,,y)(x)=H\psi(x,\,\cdot\,)(y)\qquad(x\in S,\ y\in R).
\end{equation}
\end{lemma}

\begin{Proof}
The duality relation (\ref{dual2}) says that
\[
\sum_{x'\in S}P_t(x,x')\psi(x',y)=\sum_{y'\in R}\psi(x,y')Q_t(y,y')
\qquad(x\in S,\ y\in R,\ t\geq 0),
\]
which can in matrix form be written as
\begin{equation}\label{matrix}
P_t\psi=\psi Q^\dgg_t\qquad(t\geq 0),
\end{equation}
where $Q^\dgg_t(y',y):=Q_t(y,y')$ denotes the transpose of
$Q_t$. Differentiating with respect to $t$ and setting $t=0$, it
follows that
\[
G\psi=\psi H^\dgg
\]
which is just a more formal way of writing (\ref{Gdual}). Conversely,
if (\ref{Gdual}) holds, then $G^2\psi=G\psi H^\dgg=\psi(H^\dgg)^2$ and
by induction $G^n\psi=\psi(H^\dgg)^n$ for all $n\geq 0$. Using the
fact that
\[
P_t=\sum_{n=0}^\infty\frac{1}{n!}t^nG^n
\quand
Q_t=\sum_{n=0}^\infty\frac{1}{n!}t^nH^n,
\]
it follows that $P_t\psi=\psi Q^\dgg_t$ $(t\geq 0)$ and hence
(\ref{dual2}) holds.
\end{Proof}

Let $(X_t)_{t\geq 0}$ and $(Y_t)_{t\geq 0}$ be Markov processes with
finite state spaces $S$ and $R$, generators $G$ and $H$, and Markov
semigroups $(P_t)_{t\geq 0}$ and $(Q_t)_{t\geq 0}$. Let $K$ be a
probability kernel from $S$ to $R$. A relation of the form (compare
(\ref{matrix}))
\begin{equation}\label{twine}
P_tK=KQ_t\qquad(t\geq 0)
\end{equation}
is called an \emph{intertwining}\index{intertwining of Markov processes}
of Markov processes. Note that (\ref{twine}) says that
the following two procedures are equivalent for each $S$-valued random
variable $X_0$:
\begin{itemize}
\item Evolve the state $X_0$ for time $t$ under the evolution of the
  Markov process $(X_t)_{t\geq 0}$, then map the outcome $X_t$ into a
  random variable $Y_t$ using the kernel $K$.
\item Map $X_0$ into a random variable $Y_0$ using the kernel $K$,
  then evolve $Y_0$ for time $t$ under the evolution of the Markov
  process $(Y_t)_{t\geq 0}$.
\end{itemize}
We can summarize the situation in the following commutative diagram:

\begin{equation}\label{twinefig}
\inputtikz{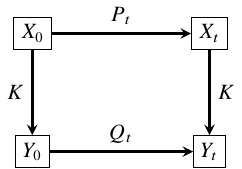}
\commentAlt{Formula (\ref{twinefig})}{Commutative diagram demonstrating
  the two possible routes from $X_0$ to $Y_t$: one via $P_t$ to $X_t$
  and then via $K$ to $Y_t$, the other via $K$ to $Y_0$ and then via
  $Q_t$ to $Y_t$.}
\end{equation}

\begin{lemma}[Intertwining of finite Markov processes]
The intertwining relation (\ref{twine}) is equivalent to
\begin{equation}\label{Gtwine}
GK=KH.
\end{equation}
\end{lemma}

\begin{Proof}
Analogue to the proof of Lemma~\ref{L:dumat}.
\end{Proof}

As one might guess, there is a close relationship between duality and
intertwining. If $(P_t)_{t\geq 0}$, $(Q_t)_{t\geq 0}$, and
$(R_t)_{t\geq 0}$ are Markov semigroups, $K$ is a probability kernel,
and $\psi$ a duality function such that
\[
P_tK=KQ_t\quand Q_t\psi=\psi R^\dgg_t\qquad(t\geq 0),
\]
then trivially
\begin{equation}\label{dualtwine}
P_t(K\psi)=KQ_t\psi=(K\psi)R^\dgg_t\qquad(t\geq 0),
\end{equation}
which says that the Markov processes with semigroups $(P_t)_{t\geq 0}$
and $(R_t)_{t\geq 0}$ are dual with duality function $K\psi$.

To see these general principles at work, let us look at interacting
particle systems with state space of the form $\{0,1\}^\La$ where
$\La$ is finite. For each $r>0$, we let $\psi_r$ denote the duality
function
\begin{equation}\label{LSpsi}
\psi_r(x,y):=\prod_{i\in\La}(1-r)^{\txt x(i)y(i)}\qquad\big(x,y\in\{0,1\}^\La\big).
\end{equation}
Using the fact that $0^n=1_{\{n=0\}}$, we observe that
\[\left.\begin{array}{r@{\,}c@{\,}l}
\dis\psi_1(x,y)&=&\dis 1-\psi_{\rm add}(x,y)\\[5pt]
\dis\psi_2(x,y)&=&\dis(-1)^{\txt\psi_{\rm canc}(x,y)}\end{array}\right\}\qquad\big(x,y\in\{0,1\}^\La\big).
\]
Therefore, duality with respect to duality functions of the form
(\ref{LSpsi}) includes additive and cancellative duality as special
cases. Duality functions of the form (\ref{LSpsi}) arose from the work
of Lloyd and Sudbury \cite{SL95,SL97}. Accordingly, we will call
$\psi_r$ the \emph{Lloyd--Sudbury duality function.}
\index{Lloyd--Sudbury duality}

There is a close connection between duality functions of the form
(\ref{LSpsi}) and \emph{thinning},\index{thinning} as we now
explain. Thinning has already been introduced in Section~\ref{S:Pois}
but for convenience we repeat the definition here. Let
$(\chi_p(i))_{i\in\La}$ be i.i.d.\ with $\P[\chi_p(i)=1]=p$ and
$\P[\chi_p(i)=0]=1-p$. Then
\begin{equation}\label{Kthin}
K_p(x,y):=\P\big[y(i)=\chi_p(i)x(i)\ \forall i\in\La\big]
\qquad\big(x,y\in\{0,1\}^\La\big)
\end{equation}
defines a \emph{thinning kernel}. Note that if we interpret sites $i$
with $x(i)=1$ as being occupied by a particle, then the effect of
$K_p$ is to independently throw away some of these particles, where
each particle has a probability $p$ to remain. We claim that
\begin{equation}\label{thinpsi}
K_pK_q=K_{pq}\quand K_p\psi_r=\psi_{pr}\qquad(0\leq p,q\leq 1,\ r>0).
\end{equation}
The first relation is clear from the interpretation in terms of
thinning, while the second relation follows by writing
\[\begin{array}{l}
\dis K_p\psi_r(x,z)=\sum_yK_p(x,y)\prod_{i\in\La}(1-r)^{\txt y(i)z(i)}
=\E\big[\prod_{i\in\La}(1-r)^{\txt\chi_p(i)x(i)z(i)}\big]\\[5pt]
\dis\quad=\prod_{i\in\La}\E\big[(1-r)^{\txt\chi_p(i)x(i)z(i)}\big]
=\prod_{i\in\La}(1-pr)^{\txt x(i)z(i)}=\psi_{pr}(x,z).
\end{array}\]
The following lemma says that if a particle system has two duals, one
with respect to the duality function $\psi_{r_1}$ and the other with
respect to the duality function $\psi_{r_2}$, then one of these duals
is a thinning of the other.

\begin{lemma}[Lloyd--Sudbury duals and thinning]
Let\label{L:LSthin} $G,H_1$ and $H_2$ be generators of Markov
processes with state space $\{0,1\}^\La$ where $\La$ is finite. Let
$0<r_1\leq r_2$ and set $p:=r_1/r_2$. Then of the relations
\[
{\rm (i)}\ H_1\psi_{r_1}=\psi_{r_1}G^\dgg,
\quad
{\rm (ii)}\ H_2\psi_{r_2}=\psi_{r_2}G^\dgg,
\quand
{\rm (iii)}\ H_1K_p=K_pH_2
\]
any two imply the third.
\end{lemma}

\begin{Proof}
Using all tree relations (i)--(iii) as well as (\ref{thinpsi}), we
have the ``circular'' sequence of equalities:
\[
H_1\psi_{r_1}
\stackrel{\rm(i)}{=}
\psi_{r_1}G^\dgg
=
K_p\psi_{r_2}G^\dgg
\stackrel{\rm(ii)}{=}
K_pH_2\psi_{r_2}
\stackrel{\rm(iii)}{=}
H_1K_p\psi_{r_2}
=
H_1\psi_{r_1}.
\]
From this, we immediately see that of the relations
\[
{\rm (i)}\ H_1\psi_{r_1}=\psi_{r_1}G^\dgg,
\quad
{\rm (ii)'}\ K_pH_2\psi_{r_2}=K_p\psi_{r_2}G^\dgg,
\quad
{\rm (iii)'}\ H_1K_p\psi_{r_2}=K_pH_2\psi_{r_2}
\]
any two imply the third. To complete the proof, it suffices to show
that $K_p$ and $\psi_r$ are invertible as matrices for all $p\in(0,1]$
and $r>0$, since we can then multiply (ii)' from the left with
$K_p^{-1}$ and (iii)' from the right with $\psi_{r_2}^{-1}$ to obtain
(ii) and (iii).

We can view the linear space of all functions $f\cn\{0,1\}^\La\to\R$
as the tensor product $\bigotimes_{i\in\La}\R^{\{0,1\}}$. In this
picture, the matrices $K_p$ and $\psi_r$ are the tensor product over
$\La$ of single-site matrices of the form
\[
\left(\begin{array}{cc}1&0\\ 1-p&p\end{array}\right)
\quand\left(\begin{array}{cc}1&1\\ 1&1-r\end{array}\right),
\]
respectively. These single-site matrices are invertible for all
$p\in(0,1]$ and $r\in(0,\infty)$ and hence the same is true for their
tensor products $K_p$ and $\psi_r$.
\end{Proof}

It is useful to look at a concrete example. Let $(\La,E)$ be a finite
graph, as in (\ref{Ni}) let $\Ni_i:=\big\{j\in\La:\{i,j\}\in E\big\}$
denote the neighborhood of a vertex $i\in\La$, let and assume that
$N:=|\Ni_i|$ does not depend on $i\in\La$. In line with notation
introduced in Section~\ref{S:setup}, we let
$\Ei:=\big\{(i,j)\in\La^2:\{i,j\}\in E\big\}$ denote the set of
directed edges associated with $E$. Let $G_{\rm vot},G_{\rm rw}$
$G_{\rm arw}$ be the Markov generators defined by
\[\begin{array}{r@{\,}c@{\,}l}
\dis G_{\rm vot}f(x)&:=&\dis N^{-1}\sum_{(i,j)\in\Ei}
\big\{f\big({\tt vot}_{ij}(x)\big)-f\big(x\big)\big\},\\[5pt]
\dis G_{\rm rw}f(x)&:=&\dis N^{-1}\sum_{(i,j)\in\Ei}
\big\{f\big({\tt rw}_{ij}(x)\big)-f\big(x\big)\big\},\\[5pt]
\dis G_{\rm arw}f(x)&:=&\dis N^{-1}\sum_{(i,j)\in\Ei}
\big\{f\big({\tt arw}_{ij}(x)\big)-f\big(x\big)\big\},
\end{array}\]
where the voter model map ${\tt vot}_{ij}$, the coalescing random walk
map ${\tt rw}_{ij}$, and the annihilating random walk map ${\tt arw}_{ij}$
are defined in (\ref{votmap}), (\ref{rwmap}), and
(\ref{annmap}), respectively. In words, $G_{\rm vot}$ is the generator
of a voter model in which each site $i\in\La$ adopts with rate one the
type of a randomly chosen neighbor. The processes with generators
$G_{\rm rw}$ and $G_{\rm arw}$ consist of coalescing and annihilating
particles that jump to a randomly chosen neighboring site with rate
one. We claim that
\[
{\rm(i)}\ G_{\rm rw}\psi_1=\psi_1G_{\rm vot}^\dgg,\quad
{\rm(ii)}\ G_{\rm arw}\psi_2=\psi_2G_{\rm vot}^\dgg,\quad
{\rm(iii)}\ G_{\rm rw}K_{1/2}=K_{1/2}G_{\rm arw}.
\]
Indeed, the voter model map is both additive and cancellative, so (i)
follows from Theorem~\ref{T:addu} since the coalescing random walk map
is the additive dual of the voter model map and likewise (ii) follows
from Theorem~\ref{T:candu} since the annihilating random walk map is
the cancellative dual of the voter model map. By Lemma~\ref{L:LSthin},
(i) and (ii) imply (iii), which says that annihilating random walks
are a $1/2$-thinning of coalescing random walks. In other words, for
each $t\geq 0$, the following two procedures are equivalent:
\begin{itemize}
\item Run coalescing random walk dynamics for time $t$ and then thin
  the resulting configuration with $1/2$.
\item Thin the initial configuration with $1/2$ and then run
  annihilating random walk dynamics for time $t$.
\end{itemize}
One can also verify this directly and use this to deduce (ii) from (i)
or vice versa. In Proposition~\ref{P:CVthin} below we will see a less
trivial example of a thinning relation between two interacting
particle systems.

Surprisingly, there exist many duality relations between interacting
particle systems with respect to the Lloyd--Sudbury duality function
$\psi_r$ for other values of $r$ than $r=1,2$. These dualities can
usually not be obtained as pathwise dualities.

Let $(\La,E)$ be a finite graph. The paper \cite{Sud00} considers
interacting particle systems on graphs where the configuration along
each edge makes the following transitions with the following
rates:\footnote{The meaning of the words ``annihilation'',
``branching'',\ldots here is a bit different from the way we have used
these words so far. In particular, the ``death'' rate $d$ refers only
to ``deaths while the neighboring site is empty'', while ``deaths
while the neighboring site is occupied'' are called ``coalescence''.}
\[\begin{array}{r@{\qquad}c@{\qquad}l}
\mbox{``annihilation''}&11\mapsto00&\mbox{at rate }a,\\[5pt]
\mbox{``branching''}&01\mapsto11\quand10\mapsto11&\mbox{each at rate }b,\\[5pt]
\mbox{``coalescence''}&11\mapsto01\quand11\mapsto10&\mbox{each at rate }c,\\[5pt]
\mbox{``death''}&01\mapsto00\quand10\mapsto00&\mbox{each at rate }d,\\[5pt]
\mbox{``exclusion''}&01\mapsto10\quand10\mapsto01&\mbox{each at rate }e.
\end{array}\]
More formally, for each $i,j\in\La$, we can define a map $m^{01\mapsto 11}_{ij}$
on $\{0,1\}^\La$ as follows:
\[
m^{01\mapsto 11}_{ij}(x)(k)=\left\{\begin{array}{ll}
1\quad&\mbox{if }k=i\mbox{ and }\big(x(i),x(j)\big)=(0,1),\\
1\quad&\mbox{if }k=j\mbox{ and }\big(x(i),x(j)\big)=(0,1),\\
x(k)\quad&\mbox{in all other cases.}
\end{array}\right.
\]
Defining $m^{11\mapsto 00}_{ij}$ etc.\ in a similar way, the generator
of the process we are interested in can be written as
\begin{equation}\begin{array}{r@{\,}c@{\,}l}\label{abcde}
\dis Gf(x)&=&\dis\sum_{\{i,j\}}a\big\{f\big(m^{11\mapsto 00}_{ij}(x)\big)-f\big(x\big)\big\}\\[5pt]
&&\dis+\sum_{(i,j)}\Big[b\big\{f\big(m^{01\mapsto 11}_{ij}(x)\big)-f\big(x\big)\big\}
+c\big\{f\big(m^{11\mapsto 01}_{ij}(x)\big)-f\big(x\big)\big\}\\[5pt]
&&\dis\phantom{\sum_{(i,j)}\Big[}
+d\big\{f\big(m^{01\mapsto 00}_{ij}(x)\big)-f\big(x\big)\big\}
+e\big\{f\big(m^{01\mapsto 10}_{ij}(x)\big)-f\big(x\big)\big\}\Big],
\end{array}\end{equation}
where the first sum runs over all (unordered) edges $\{i,j\}\in E$ and
the second sum runs over all ordered pairs $(i,j)$ such that
$\{i,j\}\in E$.

\begin{theorem}[Lloyd--Sudbury duality]
Let\label{T:LSdual} $G$ and $G'$ be defined as in (\ref{abcde}) in
terms of rates $a,b,c,d,e$ and $a',b',c',d',e'$, respectively, and let
$r>0$. Then one has
\begin{equation}\label{LSdual}
G\psi_r=\psi_r{G'}^\dgg
\end{equation}
if and only if $a'=a+2(1-r)\ga$, $b'=b+\ga$, $c'=c-(2-r)\ga$,
$d'=d+\ga$, and $e'=e-\ga$, where $\ga:=(a+c-d+(1-r)b)/r$.
\end{theorem}

\begin{Proof}
This follows from Lemma~\ref{L:dumat} by checking (\ref{Gdual}). The
calculations are a bit tedious, so we omit them here. They can be
found in \cite[formula (9)]{Sud00}, which is a simplification of
\cite[formula (21)]{SL95}.
\end{Proof}

We note that a generalization of Theorem~\ref{T:LSdual} to directed
graphs can be found in \cite[Prop.~10]{Swa06}.

\section{The contact-voter model}\label{S:convot}

As we have already seen, $\psi_1=1-\psi_{\rm add}$ and
$\psi_2=(-1)^{\psi_{\rm canc}}$ correspond to additive and
cancellative duality. It seems that for $r\neq 1,2$, dualities of the
form (\ref{LSdual}) are almost never\footnote{Except some very trivial
and pathological cases.} pathwise dualities. To give an example with
$r\neq 1,2$, consider an interacting particle system on a (possibly
infinite) graph $(\La,E)$ whose dynamics are a mixture of contact
process and voter model dynamics, with generator of the form:
\begin{equation}\begin{array}{r@{\,}c@{\,}l}\label{Gcovo}
\dis G_{\rm covo}f(x)
&:=&\dis\la\sum_{(i,j)\in\Ei}
\big\{f\big({\tt bra}_{ij}(x)\big)-f\big(x\big)\big\}
+\sum_{i\in\La}\big\{f\big({\tt death}_i(x)\big)-f\big(x\big)\big\}
\\[5pt]
&&\dis+\al\sum_{(i,j)\in\Ei}
\big\{f\big({\tt vot}_{ij}(x)\big)-f\big(x\big)\big\}
\qquad\qquad(x\in\{0,1\}^\La),
\end{array}\end{equation}
where $\Ei$ denotes the set of directed edges associated with
$E$. Letting $\Ni_i:=\big\{j\in\La:\{i,j\}\in E\big\}$ denote the set
of neighbors of $i$, we assume that $\La$ is countable and
\[
\sup_{i\in\La}|\Ni_i|<\infty,
\]
which implies that the generator in (\ref{Gcovo}) satisfies the
summability condition (\ref{downsum}) of Theorem~\ref{T:Poispart} and
hence corresponds to a well-defined interacting particle system. Such
systems are studied in \cite{DLZ14}, who are especially interested in
the fast-voting limit $\al\to\infty$. The contact-voter model is
additive (but not cancellative, because the branching map is not), and
by Theorem~\ref{T:addu} dual with respect to the duality function
$\psi_1=1-\psi_{\rm add}$ to the interacting particle system with
generator
\begin{equation}\begin{array}{r@{\,}c@{\,}l}\label{Gcorw}
\dis G_{\rm corw}f(y)
&:=&\dis\la\sum_{(i,j)\in\Ei}
\big\{f\big({\tt bra}_{ij}(y)\big)-f\big(y\big)\big\}
+\sum_{i\in\La}\big\{f\big({\tt death}_i(y)\big)-f\big(y\big)\big\}
\\[5pt]
&&\dis+\al\sum_{(i,j)\in\Ei}
\big\{f\big({\tt rw}_{ij}(y)\big)-f\big(y\big)\big\}\qquad(y\in\{0,1\}^\La),
\end{array}\end{equation}
which corresponds to a system of branching and coalescing random
walks. Perhaps surprisingly, the contact-voter model is also
self-dual.

\begin{proposition}[Self-duality of the contact-voter model]
Assume\label{P:covo} that $\la>0$. Then the con\-tact-voter model with
generator as in (\ref{Gcovo}) is self-dual\index{self-duality} with respect
to the duality function $\psi_r$ with $r:=\la/(\al+\la)$.
\end{proposition}

\begin{Proof}
We first consider the case that the graph $(\La,E)$ is finite. The
generator $G_{\rm covo}$ is a special case of the generators
considered in Theorem~\ref{T:LSdual} and corresponds to the choice of
parameters
\[
a=0,\quad b=\la+\al,\quad c=1,\quad d=1+\al,\quad e=0.
\]
We observe that setting $r:=\la/(\al+\la)$ makes the parameter $\ga$
from Theorem~\ref{T:LSdual} zero, which has the effect that $a'=a$,
$b'=b$, $c'=c$, $d'=d$, and $e'=e$, that is, we have found a
self-duality.

To extend the result to infinite graphs, we use an approximation
argument. We need to show that
\begin{equation}\label{infdu}
\E\big[\psi_r\big(\Xb_{0,t}(x),x'\big)\big]
=\E\big[\psi_r\big(x,\Xb_{0,t}(x')\big)\big]
\qquad\big(t\geq 0,\ x,x'\in\{0,1\}^\La\big),
\end{equation}
where $(\Xb_{s,u})_{s\leq u}$ denotes the stochastic flow defined by
the graphical representation of the contact-voter model, and
\begin{equation}\label{LSpsi2}
\psi_r(x,y):=\prod_{i\in\La}(1-r)^{\txt x(i)y(i)}\qquad\big(x,y\in\{0,1\}^\La\big).
\end{equation}
The argument will be a bit tricky since it is in general not true that
$x_n\to x$ and $y_n\to y$ pointwise imply that
$\psi_r(x_n,y_n)\to\psi_r(x,y)$. This prevents us from using general
approximation results like Corollary~\ref{C:Partlim}. We observe,
however, that $x_n\up x$ and $y_n\up x$ imply
$\psi_r(x_n,y_n)\down\psi_r(x,y)$. This is why we will base our argument
on monotone approximation.

Our first aim is to prove (\ref{infdu}) for $x,x'\in\Si_{\rm fin}(\La)$,
the set of finite configurations. Let $(\La_n,E_n)$ be
finite subgraphs of $(\La,E)$ that increase to the whole graph. Let
$(\Xb^n_{s,t})_{s\leq t}$ be the stochastic flow of a restricted
process defined by a graphical representation where we have removed
all death maps outside $\La_n$ and all branching and voter maps along
edges that are not in $E_n$. For all $n$ large enough so that $x$ and
$x'$ are zero outside $\La_n$, we have that also $\Xb^n_{0,t}(x)$ and
$\Xb^n_{0,t}(x')$ are zero outside $\La_n$ for all $t\geq 0$. Note
that if $x$ and $y$ are zero outside $\La_n$, then in (\ref{LSpsi2})
it does not matter if we take the product over $\La$ or
$\La_n$. Therefore, applying Theorem~\ref{T:LSdual} to the processes
on the finite graphs $(\La_n,E_n)$, we see that
\begin{equation}\label{findu}
\E\big[\psi_r\big(\Xb^n_{0,t}(x),x'\big)\big]
=\E\big[\psi_r\big(x,\Xb^n_{0,t}(x')\big)\big]
\qquad(t\geq 0)
\end{equation}
for all $n$ large enough. By Theorem~\ref{T:finIPS}, the contact-voter
model restricted to $\Si_{\rm fin}(\La)$ is a nonexplosive
continuous-time Markov chain. It follows that almost surely, there
exists a (random) $m<\infty$ such that for all $n\geq m$, the
unrestricted process $(\Xb_{0,s}(x))_{0\leq s\leq t}$ stays inside
$\La_n$ up to time $t$. But then $(\Xb^n_{0,s}(x))_{0\leq s\leq t}$
must be equal to $(\Xb_{0,s}(x))_{0\leq s\leq t}$ for all $n\geq m$ so
we see that almost surely $\Xb^n_{0,t}(x)\to\Xb_{0,t}(x)$ as
$n\to\infty$ with respect to the discrete topology on $\Si_{\rm fin}(\La)$,
and by the same argument also $\Xb^n_{0,t}(x')\to\Xb_{0,t}(x')$ with respect
to the discrete topology on $\Si_{\rm fin}(\La)$. Taking the limit
$n\to\infty$ in (\ref{findu}) it follows that (\ref{infdu}) holds for all
$x,x'\in\Si_{\rm fin}(\La)$.

For general $x,x'\in\{0,1\}^\La$ we can find $x_n,x'_n\in\Si_{\rm fin}(\La)$
such that $x_n\up x$ and $x'_n\up x'$. Then also
$\Xb_{0,t}(x_n)\up\Xb_{0,t}(x)$ and $\Xb_{0,t}(x'_n)\up\Xb_{0,t}(x')$
by the monotonicity and continuity of $\Xb_{0,t}$. Using the
continuity of $\psi_r$ with respect to increasing sequences, we obtain
(\ref{infdu}) in full generality.
\end{Proof}

We have already seen in Lemma~\ref{L:LSthin} that there is a close
connection between the Llyod-Sudbury duality functions $\psi_r$ and
thinning. The following proposition demonstrates this on our example
of the contact-voter model.

\begin{proposition}[Thinning of the contact-voter model]
Let\label{P:CVthin} $(P_t)_{t\geq 0}$ and $(Q_t)_{t\geq 0}$ denote the
semigroups of the contact-voter model with generator as in
(\ref{Gcovo}) and the system of branching and coalescing random walks
with generator as in (\ref{Gcorw}), respectively. Let $K_r$ denote the
thinning kernel defined in (\ref{Kthin}) with $p:=\la/(\al+\la)$. Then
\begin{equation}\label{CVthin}
P_tK_p=K_pQ_t\qquad(t\geq 0).
\end{equation}
\end{proposition}

\begin{Proof}
We first prove the statement for finite graphs. Additive duality tells
us that (i) $G_{\rm covo}\psi_1=\psi_1G_{\rm corw}^\dgg$, and
Proposition~\ref{P:covo} tells us that (ii)
$G_{\rm covo}\psi_p=\psi_pG_{\rm covo}^\dgg$. By Lemma~\ref{L:LSthin}, this
implies (iii) $G_{\rm covo}K_p=K_pG_{\rm corw}$, which implies
(\ref{CVthin}).

To also get the result for infinite graphs, we use approximation with
finite graphs. In this case, the argument is simpler than in the proof
of Proposition~\ref{P:covo}. We claim that thinning kernels are
continuous, that is, $f\in\Ci(\{0,1\}^\La)$ implies
$K_rf\in\Ci(\{0,1\}^\La)$. Indeed, if $f:\{0,1\}^\La\to\R$ is
continuous, $x_n\to x$ pointwise, and $(\chi_i(i))_{i\in\La}$ are
i.i.d.\ Bernoulli random variables with intensity $r$, then
\[
K_rf(x_n)=\E\big[f(\chi_r x_n)\big]\asto{n}\E\big[f(\chi_r x)\big]=K_rf(x),
\]
where $(\chi_r x)(i):=\chi_r(i)x(i)$ denotes the pointwise product of
$\chi_r$ and $x$. It follows that if $X^n$ are random variables with values in
$\{0,1\}^\La$ that converge weakly in law to $X$, and $Y^n$ and $Y$
are obtained from $X^n$ and $X$ by thinning with the kernel $K_r$,
then the $Y_n$ converge weakly in law to $Y$. As a result, we can use
Corollary~\ref{C:Partlim} to approximate infinite systems with finite
systems and take the limit to get the result for infinite systems.
\end{Proof}

We will continue our study of the contact-voter model in Section~\ref{S:eqcrit}.

\section{Invariant laws of the voter model}\label{S:votinv}

\begin{figure}[htb]
\begin{center}
\inputtikz{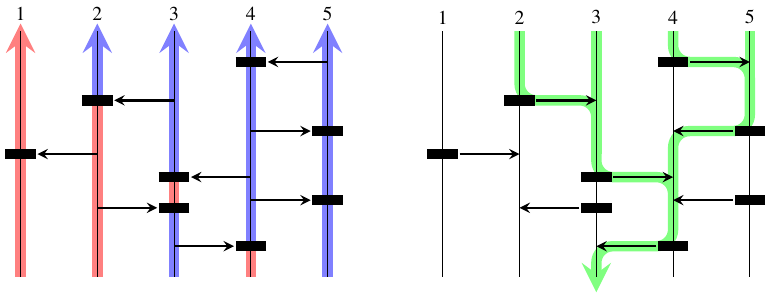}
\caption{Graphical representation of a one-dimensional voter model and
  its dual system of coalescing random walks. At the final time, the
  points 2,3,4, and 5 have the same type, because they descend from
  the same ancestor.}
\label{fig:votdual}
\commentAlt{Figure~\ref{fig:votdual}}{Graphical representation of a
  voter model and of its dual system of coalescing random
  walks. Space-time points of the voter model are colored according to
  their type. The coalescing random walks trace back where these types
  originate.}
\end{center}
\end{figure}

By Theorem~\ref{T:addu}, the voter model $X$ is pathwise dual, with
respect to the additive duality function $\psi_{\rm add}$ from
(\ref{psiadd}), to a collection $Y$ of coalescing random walks. Due to
the fact that $|Y_t|$ is a nonincreasing function of $t$ (that is, the
number of walkers can only decrease), it is much easier to work with
this dual system than with the voter model itself, so duality is
really the key to understanding the voter model.

\begin{proposition}[Clustering in low dimensions]\label{P:votclust}
Let $X$ be a nearest-neighbor or range $R$ voter model on $\Z^d$.
Assume that $d=1,2$. Then, regardless of the initial law,
\[
\P[X_t(i)=X_t(j)]\asto{t}1\qquad\forall i,j\in\Z^d.
\]
Moreover, the delta measures $\de_{\un 0}$ and $\de_{\un 1}$ on the constant
configurations are the only extremal invariant laws.
\end{proposition}

\begin{Proof}
In the graphical representation of the voter model, for each
$(i,t)\in\Z^d\times\R$ and $s\geq 0$, there is a unique site
\[\index{0zzoxiits@$\xi^{(i,t)}_s$}
j=:\xi^{(i,t)}_s\in\Z^d\mbox{ such that }(j,t-s)\leadsto(i,t).
\]
Here $(\xi^{(i,t)}_s)_{s\geq 0}$ is the path of a random walk starting at
$\xi^{(i,t)}_0=i$ and ``running downwards in the graphical representation''.
Two such random walks started from different space-time points $(i,t)$ and
$(i',t')$ are independent up to the first time they meet, and coalesce as soon
as they meet. Moreover, if $X_t=\Xb_{0,t}(X_0)$, then, as demonstrated in
Figure~\ref{fig:votdual},
\[
X_t(i)=X_{t-s}(\xi^{(i,t)}_s)\qquad(0\leq s\leq t),
\]
that is, $\xi^{(i,t)}_s$ traces back where the site $i$ at time $t$ got its type
from.\footnote{This construction works in fact generally for multitype voter
  models, where the local state space $S$ can be any finite set, and which are
  in general of course not additive systems. For simplicity, we will focus on
  the two-type voter model here.} 

Since the difference $\xi^{(i,t)}_s-\xi^{(j,t)}_s$ of two such random walks is
a random walk with absorption in the origin, and since random walk on $\Z^d$ in
dimensions $d=1,2$ is recurrent, we observe that
\[
\P[X_t(i)=X_t(j)]
\geq\P[\xi^{(i,t)}_t=\xi^{(j,t)}_t]
=\P[\xi^{(i,0)}_t=\xi^{(j,0)}_t]
\asto{t}1\qquad\forall i,j\in\Z^d.
\]
This clearly implies that all invariant laws must be concentrated on constant
configurations, that is, a general invariant law is of the form $p\de_{\un 0}+(1-p)\de_{\un 1}$
with $p\in[0,1]$.
\end{Proof}

For product initial laws we can be more precise. Although we state the
following theorem for two-type processes only, it is clear from the proof that
the statement generalizes basically unchanged to multitype voter models.

\begin{theorem}[Process started in product law]
Let\label{T:prodvot} $X$ be a nearest neighbor or range $R$ voter model on
$\Z^d$. Assume that the $(X_0(i))_{i\in\Z^d}$ are i.i.d.\ with intensity
$\P[X_0(i)=1]=p\in[0,1]$. Then
\begin{equation}\label{nutet}
\dis\P[X_t\in\cdot\,]\Asto{t}\nu_p,
\end{equation}
where $\nu_p$ is an invariant law of the process. If $d=1,2$, then
\begin{equation}
\nu_p=(1-p)\de_{\un 0}+p\de_{\un 1}.
\end{equation}
On the other hand, if $d\geq 3$ and $0<p<1$, then the measures
$\nu_p$ are concentrated on configurations that are not constant.
\end{theorem}

\begin{Proof}
As in the proof of Proposition~\ref{P:votclust}, let $(\xi^{(i,t)}_s)_{s\geq
  0}$ be the backward random walk in the graphical representation starting at
$(i,t)$. Define a random equivalence relation $\sim$ on $\Z^d$ by
\[
i\sim j\quad\mbox{iff}\quad\xi^{(i,0)}_s=\xi^{(j,0)}_s\mbox{ for some }s\geq 0.
\]
We claim that if we color the equivalence classes of $\sim$ in an
i.i.d.\ fashion such that each class gets the color 1 with probability
$p$ and the color 0 with probability $1-p$, then this defines an
invariant law $\nu_p$ such that (\ref{nutet}) holds. Since random walk
in dimensions $d=1,2$ is recurrent, there is a.s.\ only one
equivalence class, and $\nu_p=(1-p)\de_{\un 0}+p\de_{\un 1}$. On the
other hand, since random walk in dimensions $d\geq 3$ is transient,
there are a.s.\ infinitely many\footnote{Although this is intuitively
plausible, it requires a bit of work to prove this. A quick proof,
that however requires a bit of ergodic theory, is as follows: since
Poisson point processes are spatially ergodic, and the number $N$ of
equivalence classes is a translation-invariant random variable, this
random number $N$ must in fact be a.s.\ constant. Since the
probability that two paths coalesce tends to zero as the distance
between their starting points tends to infinity, for each finite $n$
we can find $n$ starting points sufficiently far from each other so
that with positive probability, none of the paths started at these
points coalesce. This implies that $\P[N\geq n]>0$ for each finite $n$
and hence by the fact that $N$ is a.s.\ constant $\P[N=\infty]=1$.}
equivalence classes and hence for $p\neq 0,1$ the measure $\nu_p$ is
concentrated on configurations that are not constant.

To prove (\ref{nutet}), we use coupling. Let $(\chi(i))_{i\in\Z^d}$ be
i.i.d.\ $\{0,1\}$-valued with $\P[\chi(i)=1]=p$. For each $t\geq 0$, we define
a random equivalence relation $\sim_t$ on $\Z^d$ by
\[
i\sim_t j\quad\mbox{iff}\quad\xi^{(i,0)}_s=\xi^{(j,0)}_s\mbox{ for some }0\leq
s\leq t.
\]
We enumerate the elements of $\Z^d$ in some arbitrary way and define
\begin{equation}
\ti X_t(i):=\chi(j)\quad\mbox{where $j$ is the smallest element of }
\{k\in\Z^d:i\sim_t k\}.
\end{equation}
Then $\ti X_t$ is equally distributed with $X_t$ and converges a.s.\ as
$t\to\infty$ to a random variable with law $\nu_p$.
\end{Proof}

\noi
\textbf{Remark} In dimensions $d\geq 3$, it is in fact known that the
measures $\nu_p$ are extremal, and each extremal invariant law of the
voter model is of this form. See \cite[Thm~V.1.8]{Lig85}.

\begin{Exercise}
Let $(Y_t)_{t\geq 0}$ be coalescing random walks with generator as in
(\ref{Grw}). Show that the upper invariant law is $\de_{\un 0}$, the
delta-measure on the all-zero configuration. \emph{Hint:} Use
Lemma~\ref{L:difP} to derive a differential equation for
$\E^{\ov\nu}[Y_t(0)]$. To complete the argument, you will need to
argue that if $\ov\nu\neq\de_{\un 0}$, then the event
$\{y:y(0)=y(1)=1\}$ has positive probability under $\ov\nu$.
\end{Exercise}

\begin{Exercise}
Prove that the voter model started in a finite initial state dies out:
\[
\P^x\big[\exists t\geq 0\mbox{ s.t.\ }X_t=\un 0\big]=1
\qquad\big(x\in\Si_{\rm fin}(\Z^d)\big).
\]
\emph{Hint:} You can use the previous exercise and
duality. Alternatively, you can use martingale convergence.
\end{Exercise}

Note that the statements of both previous exercises are not true
on finite lattices.

\section{Homogeneous invariant laws}\label{S:hom}

In the present section, we show how the self-duality\index{self-duality}
of the contact process can be used to prove that for contact processes
with some sort of translation invariant structure, the upper invariant
law is the limit law started from any nontrivial translation invariant
initial law, and we will show that this in turn implies that the
function $\tet(\la)$ from (\ref{tetdef}) is continuous everywhere,
except possibly at the critical point. The methods of the present
section are not restricted to additive particle systems. Applications
of the technique to cancellative systems can be found in
\cite{BDD91,SS08,CP14} while \cite{LS23a} treats a coupling of an
additive and a cancellative system.

We work in the set-up of Section~\ref{S:contmon}, so we consider
contact processes with generator of the form (\ref{Gcont3}) where
$\La$ is a countable set, $p$ is a symmetric probability kernel on
$\La$ such that $p(i,i)=0$ $(i\in\La)$, and the pair $(\La,p)$ is
vertex transitive as defined in (\ref{ptrans}). We will also assume
that $p$ is irreducible in the sense that for each $i,j\in\La$, there
exists an $n\geq 0$ such that $p^n(i,j)>0$. We start with a simple
observation, that has been anticipated before, and that says that the
functions $\tet(\la)$ from (\ref{surprob}) and (\ref{tetdef}) are the
same. We continue to use the notation $|x|:=\sum_ix(i)$ and let
$\Si_{\rm fin}(\La):=\{x\in\{0,1\}^\La:|x|<\infty\}$ denote the space
of finite configurations.

\begin{lemma}[The function theta]\index{0zzhthetalambda@$\tet(\la)$}
Let\label{L:theta} $X$ denote the contact process with infection rate
$\la$ on a graph $\La$ and let $\ov\nu$ denote its upper invariant
law. Then
\[
\int\!\ov\nu(\di x)\,x(i)=\P^{e_i}[X_t\neq\un 0\ \forall t\geq 0]
\qquad(i\in\La).
\]
More generally, for any $y\in\Si_{\rm fin}(\La)$,
\[
\int\!\ov\nu(\di x)\,1_{\txt\{x\wedge y\neq\un 0\}}
=\P^y[X_t\neq\un 0\ \forall t\geq 0].
\]
\end{lemma}

\begin{Proof}
By Theorem~\ref{T:addu}, the contact process $X$ is self-dual with
respect to the additive systems duality function, that is,
\[
\P^x[X_t\wedge y=\un 0]=\P^y[x\wedge X_t=\un 0]
\qquad\big(t\geq 0,\ x,y\in\{0,1\}^\La\big).
\]
In particular, setting $x=\un 1$, we see that for any $y\in\Si_{\rm fin}(\La)$,
\[
\int\!\ov\nu(\di x)\,1_{\txt\{x\wedge y\neq\un 0\}}
=\lim_{t\to\infty}\P^{\un 1}[X_t\wedge y\neq\un 0]
=\lim_{t\to\infty}\P^y[\un 1\wedge X_t\neq\un 0]
=\P^y[X_t\neq\un 0\ \forall t\geq 0].
\]
Note that since $|y|<\infty$, the function $x\mapsto 1_{\txt\{x\wedge y\neq\un 0\}}$
is continuous, which together with the weak convergence
$\de_{\un 1}P_t\Rightarrow\ov\nu$ implies the first equality above. The
condition $|y|<\infty$ can be removed by doing this step more
carefully, using monotone convergence instead of weak convergence, but
since contact processes started in infinite initial states a.s.\ do
not die out in finite time, this case is less interesting.
\end{Proof}

We will be interested in processes that are started in a translation
invariant initial law. Since we are working in a rather general
set-up, we have to say more precisely what we mean by translation
invariance. Let ${\rm Aut}(\La,p)$ be the group of all automorphisms
of $(\La,p)$, in the sense defined in Section~\ref{S:contmon}. We say
that a subgroup $\Ti\sub{\rm Aut}(\La,p)$ is \emph{vertex transitive}
\index{transitivity!of group} \index{vertex transitive!group} if for
each $i,j\in\La$, there exists a $\psi\in\Ti$ such that
$\psi(i)=j$. To see an example, consider the case that $\La=\Z^d$ and
$p$ is the transition kernel of a symmetric nearest-neighbor random
walk on $\Z^d$. For each $j\in\Z^d$, let $\psi_j$ be the translation
defined as $\psi_j(i):=i+j$ $(i\in\Z^d)$. Then the group of
translations $\Ti:=\{\psi_j:j\in\Z^d\}$ is vertex transitive. In this
example, $\Ti$ is smaller than the group ${\rm Aut}(\Z^d,p)$ of all
automorphisms of $(\Z^d,p)$ which also contains rotations, inversions,
and more.

From now on, we fix a vertex transitive subgroup $\Ti\sub{\rm Aut}(\La,p)$.
We say that a probability law $\mu$ on $\{0,1\}^\La$
is \emph{homogeneous}\index{homogeneous} or \emph{translation
invariant}\index{translation invariance} if $\mu\circ\psi^{-1}=\mu$
for all $\psi\in\Ti$. The main aim of the present section is to prove
the following result, which is originally due to Harris \cite{Har76},
with a similar result for a one-dimensional discrete time process
already having been proved by Vasil'ev \cite{Vas69}. We can think of
this result as a sort of spatial analogue of the observation in
Section~\ref{S:contmean} that for the mean-field contact process,
solutions of the differential equation (\ref{meancont}) started in any
nonzero initial state converge to the upper fixed point. Recall from
Section~\ref{S:contmon} that a probability law $\mu$ on $\{0,1\}^\La$
is \emph{nontrivial}\index{nontrivial law} if $\mu(\{\un 0\})=0$, that
is, if $\mu$ gives zero probability to the all-zero configuration.

\begin{theorem}[Convergence to upper invariant law]
Let\label{T:homcon} $(X_t)_{t\geq 0}$ be a contact process started in a
homogeneous nontrivial initial law $\P[X_0\in\cdot\,]$. Then
\[
\P[X_t\in\cdot\,]\Asto{t}\ov\nu,
\]
where $\ov\nu$ is the upper invariant law.
\end{theorem}

We start with two preparatory lemmas. We will use the graphical
representation of the contact process as an additive particle system
(see Section~\ref{S:add}) and use the shorthand
\[\index{0Xxt@$X^x_t$}
X^x_t:=\Xb_{0,t}(x)\qquad\big(t\geq 0,\ x\in\{0,1\}^\La\big),
\]
where $(\Xb_{s,t})_{s\leq t}$ is the stochastic flow constructed from
the graphical representation.

\begin{lemma}[Extinction versus unbounded growth]\label{L:exgro}
For each $x\in\Si_{\rm fin}(\La)$, one has
\begin{equation}\label{exgro}
X^x_t=\un 0\ \mbox{ for some }t\geq 0\quad\mbox{or}\quad
|X^x_t|\asto{t}\infty\quad{\rm a.s.}
\end{equation}
\end{lemma}

\begin{Proof}
Define
\[\index{0zzrrho@$\rho$}
\rho(x):=\P\big[X^x_t\neq\un 0\ \forall t\geq 0\big]
\qquad\big(x\in\Si_{\rm fin}(\La)\big).
\]
Since there is a positive probability that each infected site dies
before it manages to reproduce, it is not hard to see that for each
$N\geq 0$ there exists an $\eps>0$ such that
\begin{equation}\label{Next}
|x|\leq N\quad\mbox{implies}\quad\rho(x)\leq1-\eps.
\end{equation}
We first argue why it is plausible that this implies (\ref{exgro}) and then
give a rigorous proof. Imagine that $|X^x_t|\not\to\infty$. Then, in view
of (\ref{Next}), the process infinitely often gets a chance of at least $\eps$
to die out, hence eventually it should die out.

To make this rigorous, let
\[
\Ai_x:=\{X^x_t\neq\un 0\ \forall t\geq 0\}
\qquad\big(x\in\Si_{\rm fin}(\La)\big)
\]
denote the event that the process $(X^x_t)_{t\geq 0}$ survives
and let $\Fi_t$ be the \si-field generated by the Poisson point processes used
in our graphical representation till time~$t$. Then
\begin{equation}\label{marco}
\rho(X^x_t)=\P\big[\Ai_x\,\big|\,\Fi_t\big]\asto{t}1_{\Ai_x}\quad{\rm a.s.},
\end{equation}
where we have used an elementary result from probability theory that
says that if $\Fi_n$ is an increasing sequence of \si-fields and
$\Fi_\infty=\sig(\bigcup_n\Fi_n)$, then
$\lim_{n\to\infty}\P[\Ai|\Fi_n]=\P[\Ai|\Fi_\infty]$ a.s.\ for each
measurable event $\Ai$. (See \cite[\S~29, Complement~10~(b)]{Loe63}.)
In view of (\ref{Next}), formula (\ref{marco}) implies (\ref{exgro}).
\end{Proof}

\begin{lemma}[Nonzero intersection]
Let\label{L:nonempt} $(X_t)_{t\geq 0}$ be a contact process with a
homogeneous nontrivial initial law $\P[X_0\in\cdot\,]$. Then for each
$s,\eps>0$ there exists an $N\geq 1$ such that for any
$y\in\Si_{\rm fin}(\La)$
\[
|y|\geq N\quad\mbox{implies}\quad
\P\big[X_s\wedge y=\un 0\big]\leq\eps.
\]
\end{lemma}

\begin{Proof}
We construct $(X_t)_{t\geq 0}$ as $X_t:=\Xb_{0,t}(X_0)$ $(t\geq 0)$,
where $X_0$ is independent of the graphical representation $\om$. We
let $(\Yb_{t,s})_{t\geq s}$ denote the backward stochastic flow of the
dual process as in Theorem~\ref{T:addu}. For each $\de>0$, we let
$p_\de$ denote the matrix defined by $p_\de(i,j):=1_{\{p(i,j)\geq\de\}}p(i,j)$
$(i,j\in\La)$, and for $m\geq 1$ we set
\[
\La_m(i):=\big\{j\in\La:p^n_{1/m}(i,j)>0\mbox{ for some }0\leq n\leq m\big\}.
\]
We fix an arbitrary reference point $0\in\La$. By vertex transitivity,
$|\La_{m}(i)|=|\La_{m}(0)|$ does not depend on $i\in\La$. It
is not hard to see that for each $y\in\{0,1\}^\La$ with $|y|\geq N$ we
can find a $y'\leq y$ with $|y'|\geq N/|\La_{m}(0)|$ such that
the sets $\La_{m}(i)$ where $i$ ranges through $\{i:y'(i)=1\}$
are disjoint. We let $(\Yb^{m,i}_{t,s})_{t\geq s}$ denote the
backward stochastic flow of the dual process restricted to
$\La_{m}(i)$. More precisely, this is the stochastic flow
associated with the modified graphical representation obtained by
removing all branchings from inside $\La_{m}(i)$ to its
complement. Then, using
H\"older's inequality\footnote{Recall that H\"older's inequality
\index{H\"older's inequality} says that $1/p+1/q=1$ implies
$\|fg\|_1\leq\|f\|_p\|g\|_q$, where
$\|f\|_p:=(\int|f|^p\di\mu)^{1/p}$. By induction, this gives
$\|\prod_{i=1}^nf_i\|_1\leq\prod_{i=1}^n\|f_i\|_n$.} in the inequality
marked with an exclamation mark, we have
\[\begin{array}{l}
\dis\P\big[X_s\wedge y=\un 0\big]
=\P\big[X_0\wedge\Yb_{s,0}(y)=\un 0\big]
=\int\P[X_0\in\di x]\,\P\big[x\wedge\Yb_{s,0}(y)=\un 0\big]\\[5pt]
\dis\quad\leq\int\P[X_0\in\di x]
\,\P\big[x\wedge\bigvee_{i:\,y'(i)=1}\Yb^{m,i}_{s,0}(e_i)=\un 0\big]\\[5pt]
\dis\quad=\int\P[X_0\in\di x]
\prod_{i:\,y'(i)=1}\P\big[x\wedge\Yb^{m,i}_{s,0}(e_i)=\un 0\big]\\[5pt]
\dis\quad\stackrel{!}\leq\prod_{i:\,y'(i)=1}\Big(\int\P[X_0\in\di x]
\,\P\big[x\wedge\Yb^{m,i}_{s,0}(e_i)=\un 0\big]^{|y'|}\Big)^{1/|y'|}\\[5pt]
\dis\quad=\prod_{i:\,y'(i)=1}\Big(\int\P[X_0\in\di x]
\,\P\big[x\wedge\Yb^{m,0}_{s,0}(e_0)=\un 0\big]^{|y'|}\Big)^{1/|y'|}\\[5pt]
\dis\quad=\int\P[X_0\in\di x]
\,\P\big[x\wedge\Yb^{m,0}_{s,0}(e_0)=\un 0\big]^{|y'|},
\end{array}\]
where we have used the homogeneity of $\P[X_0\in\cdot\,]$ in the
penultimate equality. Our arguments so far show that $|y|\geq N$
implies that
\[
\P\big[X_s\wedge y=\un 0\big]
\leq\int\P[X_0\in\di x]
\P\big[x\wedge\Yb^{m,0}_{s,0}(e_0)=\un 0\big]^{N/|\La_{m}(0)|}=:f(N,m).
\]
Here, using the fact that
\[
\P\big[x\wedge\Yb^{m,0}_{s,0}(e_0)=\un 0\big]<1\quad\mbox{if }
x(i)=1\mbox{ for some }i\in\La_{m}(0),
\]
we see that
\[
\lim_{N\up\infty}f(N,m)=\int\P[X_0\in\di x]1_{\{x(i)=0\ \forall i\in\La_{m}(0)\}}
=\P[X_0(i)=0\ \forall i\in\La_{m}(0)].
\]
Since the kernel $p$ is irreducible we have that $\La_m(0)\up\La$, so by the nontriviality of $\P[X_0\in\cdot\,]$ it follows that
\[
\lim_{m\up\infty}\P[X_0(i)=0\ \forall i\in\La_m(0)]=\P[X_0=\un 0]=0.
\]
Together with our previous equation, this shows that
\[
\lim_{m\to\infty}\lim_{N\to\infty}f(N,m)=0.
\]
For each $\eps>0$, we can first choose $m$ large enough such that
$\lim_{N\to\infty}f(N,m)\leq\eps/2$ and then $N$ large enough such
that $f(N,m)\leq\eps$, proving our claim.
\end{Proof}

\begin{Exercise}
Show by counterexample that the statement of Lemma~\ref{L:nonempt} is false
for $s=0$.
\end{Exercise}

\begin{Proof}[of Theorem~\ref{T:homcon}]
As in the proof of Lemma~\ref{L:exgro}, we set
\[\index{0zzrrho@$\rho$}
\rho(x):=\P\big[X^x_t\neq\un 0\ \forall t\geq 0\big]
\qquad\big(x\in\Si_{\rm fin}(\La)\big).
\]
We construct $(X_t)_{t\geq 0}$ as $X_t:=\Xb_{0,t}(X_0)$ $(t\geq 0)$,
where $X_0$ is independent of the graphical representation. By
Lemmas~\ref{L:weaksuf}, \ref{L:adddist}, and \ref{L:theta}, it
suffices to show that
\[
\lim_{t\to\infty}\P\big[\Xb_{0,t}(X_0)\wedge y\neq\un 0\big]=\rho(y)
\qquad\big(y\in\Si_{\rm fin}(\La)\big).
\]
By duality, for any $s>0$, this is equivalent to
\[
\lim_{t\to\infty}\P\big[\Xb_{0,s}(X_0)\wedge\Yb_{t,s}(y)\neq\un 0\big]=\rho(y)
\qquad\big(y\in\Si_{\rm fin}(\La)\big).
\]
Setting $X_s:=\Xb_{0,s}(X_0)$ and $X^y_t:=\Yb_{s+t,s}(y)$, we may
equivalently show that
\[
\lim_{t\to\infty}\P\big[X_s\wedge X^y_t\neq\un 0\big]=\rho(y)
\qquad\big(y\in\Si_{\rm fin}(\La)\big),
\]
where $X_s$ and $X^y_t$ are independent and $s>0$ is some fixed
constant. For each $\eps>0$, we can choose $N$ as in
Lemma~\ref{L:nonempt}, and write
\[\begin{array}{r@{\,}c@{\,}l}
\dis\P\big[X_s\wedge X^y_t\neq\un 0\big]
&=&\dis\P\big[X_s\wedge X^y_t\neq\un 0\,\big|\,|X^y_t|=0\big]
\,\P\big[|X^y_t|=0\big]\\[5pt]
&&\dis+\P\big[X_s\wedge X^y_t\neq\un 0\,\big|\,0<|X^y_t|<N\big]
\,\P\big[0<|X^y_t|<N\big]\\[5pt]
&&\dis+\P\big[X_s\wedge X^y_t\neq\un 0\,\big|\,|X^y_t|\geq N\big]
\,\P\big[|X^y_t|\geq N\big].
\end{array}\]
Here, by Lemma~\ref{L:exgro} and our choice of $N$,
\[\begin{array}{rl}
{\rm(i)}&\dis
\P\big[X_s\wedge X^y_t\neq\un 0\,\big|\,|X^y_t|=0\big]=0,\\[5pt]
{\rm(ii)}&\dis\lim_{t\to\infty}\P\big[0<|X^y_t|<N\big]=0,\\[5pt]
{\rm(iii)}&\dis\liminf_{t\to\infty}
\P\big[X_s\wedge X^y_t\neq\un 0\,\big|\,|X^y_t|\geq N\big]
\geq 1-\eps,\\[5pt]
{\rm(iv)}&\dis\lim_{t\to\infty}\P\big[|X^y_t|\geq N\big]=\rho(x),
\end{array}\]
from which we conclude that
\[
(1-\eps)\rho(x)
\leq\liminf_{t\to\infty}\P\big[X_s\wedge X^y_t\neq\un 0\big]
\leq\limsup_{t\to\infty}\P\big[X_s\wedge X^y_t\neq\un 0\big]
\leq\rho(x).
\]
Since $\eps>0$ is arbitrary, our proof is complete.
\end{Proof}

Theorem~\ref{T:homcon} has a simple corollary.

\begin{corollary}[Homogeneous invariant laws]
All\label{C:hom} homogeneous invariant laws of a contact process are
convex combinations of $\de_{\un 0}$ and $\ov\nu$.
\end{corollary}

\begin{Proof}
Let $\nu$ be any homogeneous invariant law. We will show that $\nu$ is a
convex combination of $\de_{\un 0}$ and $\ov\nu$. If $\nu=\de_{\un 0}$ we are
done. Otherwise, as in the proof of Lemma~\ref{L:nontriv}, we can write
$\nu=(1-p)\de_{\un 0}+p\mu$ where $p\in(0,1]$ and $\mu$ is a
nontrivial homogeneous invariant law. But now Theorem~\ref{T:homcon}
implies that
\[
\mu=\mu P_t\Asto{t}\ov\nu,
\]
so we conclude that $\mu=\ov\nu$.
\end{Proof}

Recall from Exercise~\ref{E:tetright} that the function $\la\mapsto\tet(\la)$
from (\ref{tetdef}) is right-continuous everywhere. We let
\begin{equation}\label{lac}\index{0zzllambdac@$\la_{\rm c}$}
\la_{\rm c}:=\inf\{\la\in\R:\tet(\la)>0\}
\end{equation}
denote the \emph{critical point}\index{critical point}
\index{critical!infection rate} of the contact process.
As an application of Theorem~\ref{T:homcon}, we prove the following result.

\begin{proposition}[Continuity above the critical point]\label{P:leftco}
The function $\la\mapsto\tet(\la)$ is left-continuous on $(\la_{\rm c},\infty)$.
\end{proposition}

\begin{Proof}
Let $\ov\nu_\la$ denote the upper invariant law of the contact process with
infection rate $\la$. Fix $\la>\la_{\rm c}$ and choose $\la_{\rm c}<\la_n\up\la$.
Since the space $\Mi_1(\{0,1\}^\La)$ of probability measures on $\{0,1\}^\La$,
equipped with the topology of weak convergence, is compact, it suffices to
show that each subsequential limit $\nu_\ast$ of the measures $\ov\nu_{\la_n}$
equals $\ov\nu_\la$. By Proposition~\ref{P:invlim}, each subsequential
limit $\nu_\ast$ is an invariant law. It is clearly also homogeneous. Since
$\la>\la_{\rm c}$, by Lemma~\ref{L:nontriv}, the measures $\ov\nu_{\la_n}$ are
nontrivial for all $n$, and hence, using also
Proposition~\ref{P:moninla}, the same is true for $\nu_\ast$. By
Corollary~\ref{C:hom}, we conclude that $\nu_\ast=\ov\nu$. This argument shows
that the map
\[
(\la_{\rm c},\infty)\ni\la\mapsto\ov\nu_\la
\]
is left-continuous w.r.t.\ the topology of weak convergence. Since $x\mapsto
x(i)$ is a continuous function and $\tet(\la)$ is its expectation under
$\ov\nu_\la$, the claim follows.
\end{Proof}

\begin{Exercise}
Let $(X_t)_{t\geq 0}$ be a additive interacting particle system and
let $(Y_t)_{t\geq 0}$ be its additive dual. Show that the upper
invariant law $\ov\nu$ of $(X_t)_{t\geq 0}$ is uniquely characterized
by
\begin{equation}\label{uppchar}
\int\!\ov\nu(\di x)\,1_{\txt\{x\wedge y\neq\un 0\}}
=\P^y[Y_t\neq\un 0\ \forall t\geq 0]\qquad\big(y\in\Si_{\rm fin}(\La)\big).
\end{equation}
\end{Exercise}

\begin{Exercise}
Let $(X_t)_{t\geq 0}$ be a cancellative interacting particle system
and let $(Y_t)_{t\geq 0}$ be it cancellative dual. Let $\pi_{1/2}$
denote product measure with intensity $1/2$. Show that
\[
\P^{\pi_{1/2}}\big[X_t\in\,\cdot\,\big]\Asto{t}\nu_{1/2},
\]
where $\nu_{1/2}$ is an invariant law that is uniquely characterized
by the relation
\[
\int\nu_{1/2}(\di x)1_{\txt\{|x\wedge y|\mbox{ is odd}\}}=\ha\P^y\big[Y_t\neq\un 0\ \forall t\geq 0\big]\qquad\big(y\in\Si_{\rm fin}(\La)\big).
\]
Because of the similarity of this formula to the characterization of
the upper invariant law of an additive interacting particle system in
(\ref{uppchar}), the measure $\nu_{1/2}$ is sometimes called the
\emph{odd upper invariant law}.\index{odd upper invariant law}
\index{upper invariant law!odd}
\end{Exercise}

\begin{Exercise}
If we drop the assumption that the probability kernel in
(\ref{Gcont3}) is symmetric, then the contact process is no longer
self-dual. Show that in such a setting, vertex transitivity implies
that the constant
\begin{equation}\label{Kback}
K:=\sum_{i\in\La}p(i,j)
\end{equation}
does not depend on $j\in\La$. Assuming that $K<\infty$, show that the
dual process is a contact process with kernel $p'(i,j):=K^{-1}p(j,i)$
and infection rate $\la':=K\la$. Give an example of a transitive pair
$(\La,p)$ for which $K\neq 1$. \emph{Hint:} Consider an infinite tree
in which each vertex has three neighbors. Give each edge an
orientation so that at each vertex, there are two incoming edges and
one outgoing edge, and the oriented paths starting at any two vertices
eventually meet. (Compare Exercise~\ref{E:Kupdo}.)
\end{Exercise}

\begin{Exercise}
Show that Theorem~\ref{T:homcon} remains true if we drop the
assumption that $p$ is symmetric but assume that the constant in
(\ref{Kback}) is finite.
\end{Exercise}

\section{Equality of critical points}\label{S:eqcrit}

The contact-voter model $X$ on $\Z^d$, that has a mixture of contact
process and voter model dynamics, has been introduced in
Section~\ref{S:convot}. It has two parameters: the infection rate
$\la$ and the voter rate $\al$. We say that $X$ \emph{survives} if
\[
\P^{e_0}[X_t\neq\un 0\ \forall t\geq 0]>0.
\]
For each $\al\geq 0$, we define critical\index{critical point}
\index{critical!infection rate} infection rates $\la_{\rm c}(\al)$
and $\la'_{\rm c}(\al)$ by
\[\index{0zzllambdac@$\la_{\rm c}$}
\begin{array}{r@{\,}c@{\,}l}
\dis\la_{\rm c}(\al)&:=&\inf\big\{\la\in\R:
\mbox{the upper invariant law is nontrivial}\big\},\\[5pt]
\dis\la'_{\rm c}(\al)&:=&\inf\big\{\la\in\R:
\mbox{the process survives}\big\}.
\end{array}\]
The paper \cite{DLZ14} studies the asymptotics of $\la_{\rm c}(\al)$ as
$\al\to\infty$. Here, we will use duality to prove a more simple statement,
namely, that $\la_{\rm c}(\al)=\la'_{\rm c}(\al)$ for all $\al\geq 0$.

For $\al=0$ (that is, the pure contact process), we already know this,
as it is a direct consequence of Lemma~\ref{L:theta}, which follows
from the self-duality of the contact process. We will use a similar
argument here using Proposition~\ref{P:covo}, which says that the
contact-voter model is self-dual with respect to the duality function
$\psi_r$ from (\ref{LSpsi}) with $r:=\la/(\al+\la)$. Note that if
$\al=0$ (the pure contact process), then $r=1$ which corresponds to
additive systems duality.

\begin{proposition}[Characterization of the upper invariant law]
Let $r:=\la/(\al+\la)$. The upper invariant law $\ov\nu$ of the contact-voter
model satisfies
\begin{equation}\label{covoup}
\int\!\ov\nu(\di x)\psi_r(x,y)
=\P^y\big[X_t=\un 0\mbox{ for some }t\geq 0\big]
\end{equation}
for all finite $y\in\{0,1\}^{\Z^d}$. Moreover, the upper invariant law
is nontrivial if and only if the process survives. As a consequence,
$\la_{\rm c}(\al)=\la'_{\rm c}(\al)$ for all $\al\geq 0$.
\end{proposition}

\begin{Proof}
Letting $X^{\un 1}$ and $X^y$ denote the processes started in $X^{\un 1}_0=\un 1$
and $X^y_0=y$, we observe that by Proposition~\ref{P:covo},
\[
\int\!\ov\nu(\di x)\,\psi_r(x,y)
=\lim_{t\to\infty}\E\big[\psi_r(X^{\un 1}_t,y)\big]
=\lim_{t\to\infty}\E\big[\psi_r(\un 1,X^y_t)\big]
=\lim_{t\to\infty}\E\big[(1-r)^{\txt|X^y_t|}\big].
\]
The proof of Lemma~\ref{L:exgro} carries over without a change to the
contact-voter model, so
\[
X^y_t=\un 0\mbox{ for some }t\geq 0\quad\mbox{or}\quad
|X^y_t|\asto{t}\infty\quad{\rm a.s.}
\]
Using this, we see that
\[
\lim_{t\to\infty}\E\big[(1-r)^{\txt|X^y_t|}\big]
=\P\big[X^y_t=\un 0\mbox{ for some }t\geq 0\big],
\]
completing the proof of (\ref{covoup}).

Inserting $y=e_0$ into (\ref{covoup}), we see that
\[
\int\!\ov\nu(\di x)\,(1-r)^{\txt x(0)}
=\P^{e_0}\big[X_t=\un 0\mbox{ for some }t\geq 0\big],
\]
or equivalently, using the fact that $1-(1-r)^{x(0)}=rx(0)$ with
$r=\la/(\al+\la)$,
\[
\frac{\la}{\al+\la}\int\!\ov\nu(\di x)\,x(0)
=\P^{e_0}\big[X_t\neq\un 0\ \forall t\geq 0\big].
\]
This shows that $\ov\nu\neq\de_{\un 0}$ if and only if the process survives.
\end{Proof}

\chapter{Oriented percolation}\label{C:percol}

\section{Introduction}\label{S:percint}

Although we have seen phase transitions in our simulations of interacting
particle systems in Chapter~\ref{C:intro}, and we have seen how phase
transitions are defined and can be calculated in the mean-field limit in
Chapter~\ref{C:meanfield}, we have not yet proved the existence of a phase
transition for any of the spatial models that we have seen so far.

In the present chapter, we fill this gap by proving that the contact process
on $\Z^d$ with generator as in (\ref{Gcontact}) and death rate $\de=1$
undergoes a phase transition. We will show that the critical point
$\la_{\rm c}$ defined in (\ref{lac}) is nontrivial in the sense that
$0<\la_{\rm c}<\infty$. Note that by Lemma~\ref{L:theta},
\[\begin{array}{r@{\,}c@{\,}l}
\la_{\rm c}&=&\inf\{\la\in\R:\mbox{ the contact process survives}\}\\[5pt]
&=&\inf\{\la\in\R:\mbox{ the upper invariant law is nontrivial}\}.
\end{array}\]

In Exercise~\ref{E:context}, which is based on Theorem~\ref{T:ergo}, we have
already proved that\footnote{The apparent difference between this formula and the formula in Exercise~\ref{E:context} is due to the different conventions in (\ref{Gcontact}) and (\ref{Gcont3}).}
\[
\frac{1}{|\Ni_0|}\leq\la_{\rm c},
\]
where $|\Ni_0|=2d$ or $=(2R+1)^d-1$ is the size of the neighborhood of the
origin for the nearest-neighbor process and for the range $R$ process,
respectively. In view of this, it suffices to prove that $\la_{\rm c}<\infty$.
A simple comparison argument (Exercise~\ref{E:dimcomp}) shows that if the
nearest-neighbor one-dimensional contact process survives for some value of
$\la$, then the same is true for the nearest-neighbor and range $R$ processes
in dimensions $d\geq 2$. Thus, it suffices to show that $\la_{\rm c}<\infty$
for the nearest-neighbor process in dimension one.

The method we will use is comparison with oriented percolation. This neither
leads to a particularly short proof nor does it yield a very good upper bound
on $\la_{\rm c}$, but it has the advantage that it is a very robust method
that can be applied to many other interacting particle systems. For example,
in \cite{SS08} and \cite{SS15}, the method is applied to rebellious voter
models and systems with cooperative branching and coalescing random walk
dynamics, respectively. An important paper for propagating the technique was
\cite{Dur91}, where this was for the first time applied to non-monotone
systems and it was shown that ``basically, all one needs'' to prove survival
is that a particle system spreads into empty areas at a positive speed.

\section{Oriented percolation}

\begin{figure}[htb]
\begin{center}
\inputtikz{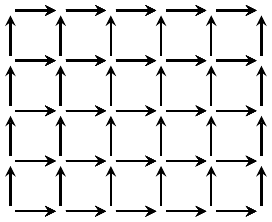}
\caption{$\Z^2$ as an oriented graph.}
\label{fig:orperc}
\commentAlt{Figure~\ref{fig:orperc}}{Square lattice with arrows
  between lattice points, pointing to the right (east) or upwards
  (north).}
\end{center}
\end{figure}

In order to prepare for the proof that the critical infection rate of the
contact process is finite, in the present section, we will study 
\emph{oriented} (or \emph{directed}) \emph{bond percolation} on
$\Z^d$.\index{oriented percolation}\index{percolation!oriented}
\index{directed percolation}\index{percolation!directed} For $i,j\in\Z^d$, we
write $i\leq j$ if $i=(i_1,\ldots,i_d)$ and $j=(j_1,\ldots,j_d)$ satisfy
$i_k\leq j_k$ for all $k=1,\ldots,d$. Let
\begin{equation}\label{Aidef}
\Ai:=\big\{(i,j):i,j\in\Z^d,\ i\leq j,\ |i-j|=1\big\}.
\end{equation}
We view $\Z^d$ as an infinite directed graph, where elements $(i,j)\in\Ai$
represent arrows (or \emph{directed bonds}) between neighboring sites. Note
that all arrows point ``upwards'' in the sense of the natural order on
$\Z^d$. See Figure~\ref{fig:orperc}.

\begin{figure}[htb]
\begin{center}
\inputtikz{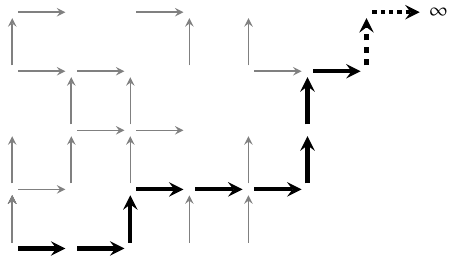}
\caption{An infinite path in oriented percolation.}
\label{fig:orperc2}
\commentAlt{Figure~\ref{fig:orperc2}}{A random subset of the arrows
  from the previous figure is drawn and an infinite path leading to
  infinity is highlighted in green.}
\end{center}
\end{figure}

Now fix some \emph{percolation parameter}\index{percolation!parameter}
$p\in[0,1]$ and let $(\om_{(i,j)})_{(i,j)\in\Ai}$ be a collection of
i.i.d.\ Bernoulli random variables with $\P[\om_{(i,j)}=1]=p$. We say that
there is an \emph{open path} \index{open path!oriented percolation} from
a site $i\in\Z^d$ to $j\in\Z^d$ if there exist $n\geq 0$ and a function
$\ga\cn\{0,\ldots,n\}\to\Z^d$ such that $\ga(0)=i$, $\ga(n)=j$, and
\[
(\ga(k-1),\ga(k))\in\Ai\quad\mbox{and}\quad\om_{(\ga(k-1),\ga(k))}=1
\qquad(k=1,\ldots,n).
\]
We denote the presence of an open path by $\leadsto$. Note that open paths
must walk upwards in the sense of the order on $\Z^d$. We write
$0\leadsto\infty$ to indicate the existence of an infinite open path starting
at the origin $0\in\Z^d$. See Figure~\ref{fig:orperc2}.

\begin{Exercise}
Show that the number of vertices that can be reached by an open path from the
origin is infinite if and only if there starts an infinite open path in the
origin.
\end{Exercise}

\begin{theorem}[Critical percolation parameter]
For\label{T:pcrit} oriented percolation in dimensions $d\geq 2$, there
exists a critical parameter $p_{\rm c}=p_{\rm c}(d)$ such that
$\P[0\leadsto\infty]=0$ for $p<p_{\rm c}$ and $\P[0\leadsto\infty]>0$
for $p>p_{\rm c}$. One has
\[
\frac{1}{d}\leq p_{\rm c}(d)\leq\frac{8}{9}.
\]
\end{theorem}

\begin{Proof}
Set
\[
p_{\rm c}:=\inf\big\{p\in[0,1]:\P[0\leadsto\infty]>0\big\}.
\]
A simple monotone coupling argument shows that $\P[0\leadsto\infty]=0$ for
$p<p_{\rm c}$ and $\P[0\leadsto\infty]>0$ for $p>p_{\rm c}$.

To prove that $0<p_{\rm c}$, let $N_n$ denote the number of open paths
of length $n$ starting in $0$. Since there are $d^n$ different upward
paths of length $n$ starting at the origin, and each path has
probability $p^n$ to be open, we see that
\[
\P[N_n\neq 0]\leq\E[N_n]=d^np^n.
\]
Since the events $\{N_n\neq 0\}$ decrease as $n\to\infty$ to the event
$\{0\leadsto\infty\}$, taking the limit, we see that
$\P[0\leadsto\infty]=0$ for all $p<1/d$, and therefore
$1/d\leq p_{\rm c}(d)$.

To prove that $p_{\rm c}(d)\leq 8/9$ for $d\geq 2$ it suffices to consider the
case $d=2$, for we may view $\Z^2$ as a subset of $\Z^d$ $(d\geq 3)$ and then,
if there is an open path that stays in $\Z^2$, then certainly there is an open
path in $\Z^d$. (Note, by the way, that in $d=1$ one has
$\P[0\leadsto\infty]=0$ for all $p<1$ and hence $p_{\rm c}(1)=1$.)

We will use a Peierls argument,\index{Peierls argument} named after
R.~Peierls who used a similar argument in 1936 for the Ising model
\cite{Pei36}. In Figure~\ref{fig:Peierls}, we have drawn a piece of
$\Z^2$ with a random collection of open arrows.
Sites $i\in\Z^2$ such that $0\leadsto i$ are drawn
green. These sites are called \emph{wet}.\index{wet sites} Consider the 
\emph{dual lattice}\index{dual lattice}
\[
\hat\Z^2:=\{(n+\ffrac{1}{2},m+\ffrac{1}{2}):(n,m)\in\Z^2\}.
\]
If there are only finitely many wet sites, then the set of all non-wet sites
in $\N^2$ contains one infinite connected component. (Here ``connected'' is to
be interpreted in terms of the undirected graph $\N^2$ with nearest-neighbor
edges.) Let $\ga$ be the boundary of this infinite component. Then $\ga$ is a
nearest-neighbor path in $\hat\Z^2$, starting in some point
$(k+\frac{1}{2},-\frac{1}{2})$ and ending in some point
$(-\frac{1}{2},m+\frac{1}{2})$ with $k,m\geq 0$, such that all sites
immediately to the left of $\ga$ are wet, and no open arrows starting at these
sites cross $\ga$. In Figure~\ref{fig:Peierls}, we have indicated $\ga$ with
red arrows.

\begin{figure}
\begin{center}
\inputtikz{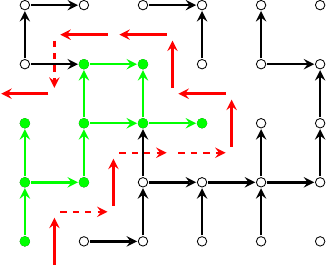}
\caption{Peierls argument for oriented percolation. The green cluster of
  points reachable from the origin is surrounded by a red contour. The
  \emph{north} and \emph{west} steps of this contour cannot cross open arrows.}
\label{fig:Peierls}
\commentAlt{Figure~\ref{fig:Peierls}}{Drawn are the random arrows of
  oriented percolation. Points reachable from the origin are
  green. The outer boundary of this set of points is a path drawn with
  red arrows. Some east and south red arrows cross percolation
  arrows.}
\end{center}
\end{figure}

From these considerations, we see that the following statement is true: one
has $0\not\leadsto\infty$ if and only if there exists a path in $\hat\Z^2$,
starting in some point $(k+\frac{1}{2},-\frac{1}{2})$ $(k\geq 0)$, ending in
some point $(-\frac{1}{2},m+\frac{1}{2})$ $(m\geq 0)$, and passing to the
northeast of the origin, such that all arrows of $\ga$ in the north and west
directions (solid red arrows in the figure) are not crossed by an open
arrow. Let $M_n$ be the number of paths of length $n$ with these
properties. Since there are $n-1$ dual sites from where such a path of length
$n$ can start, and since in each step, there are three directions where it can
go, there are less than $n3^n$ paths of length $n$ with these properties. Since
each path must make at least half of its steps in the north and west
directions, the expected number of these paths satisfies
\[
\E\big[\sum_{n=2}^\infty M_n\big]\leq\sum_{n=2}^\infty n3^n(1-p)^{n/2}<\infty
\qquad(p>\ffrac{8}{9})
\]
and therefore
\[
\P[0\not\leadsto\infty]\leq\P\big[\sum_{n=2}^\infty M_n\geq 1\big]
\leq\E\big[\sum_{n=2}^\infty M_n\big]<\infty.
\]
This does not quite prove what we want yet, since we need the
right-hand side of this equation to be less than one. To fix this, we
use a trick. (This part of the argument comes from
\cite[Section~5a]{Dur88}.) Set $D_m:=\{0,\ldots,m\}^2$ and let us say
that a set $i$ is ``wet'' if $j\leadsto i$ for some $j\in D_m$. If
$D_m\not\leadsto\infty$, then the set of wet sites must be finite,
and, just as before, there must be a dual path surrounding this set of
wet sites. Then, by the same arguments as before
\[
\P[D_m\not\leadsto\infty]\leq\P\big[\sum_{n=2m}^\infty M_n\geq 1\big]
\leq\E\big[\sum_{n=2m}^\infty M_n\big]
\leq\sum_{n=2m}^\infty n3^n(1-p)^{n/2},
\]
where now the sum starts at $2m$ since the dual paths must surround
$D_m$ and hence have length $2m$ at least. If $p>\ffrac{8}{9}$, then
the sum is finite so it can be made arbitrarily small by choosing $m$
sufficiently large. It follows that $\P[D_m\leadsto\infty]>0$ for some
$m$, hence $\P[i\leadsto\infty]>0$ for at least one $i\in D_m$, and
therefore, by translation invariance, also $\P[0\leadsto\infty]>0$.
\end{Proof}

\section{Survival}

The main result of the present chapter is the following theorem, which
rigorously establishes the existence of a phase transition for the contact
process on $\Z^d$.

\begin{theorem}[Nontrivial critical point]
For\label{T:lacont} the nearest-neighbor or range $R$ contact process
on $\Z^d$ $(d\geq 1)$, the critical infection rate satisfies
$0<\la_{\rm c}<\infty$.
\end{theorem}

\begin{Proof}
As already mentioned in Section~\ref{S:percint}, the fact that $0<\la_{\rm c}$
has already been proved in Exercise~\ref{E:context}. By
Exercise~\ref{E:dimcomp}, to prove that $\la_{\rm c}<\infty$, it suffices to
consider the one-dimensional nearest-neighbor case.

\begin{figure}
\begin{center}
\inputtikz{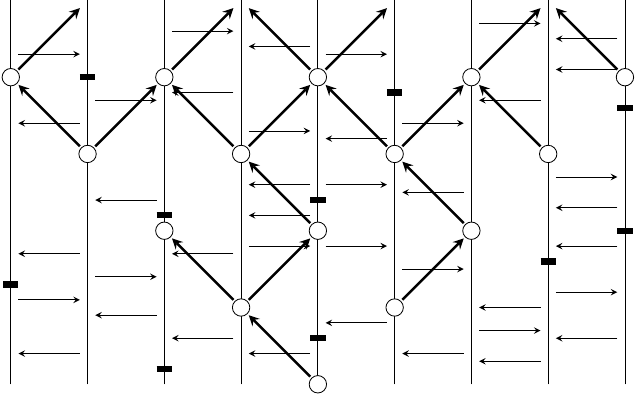}
\caption{Comparison with oriented percolation. Good events in the graphical
  representation of the contact process (blue) correspond to open percolation
  arrows (black). An infinite open path along percolation arrows implies an
  infinite open path in the graphical representation of the contact process.}
\label{fig:percomp}
\commentAlt{Figure~\ref{fig:percomp}}{An oriented percolation process
  is drawn on top of the graphical representation of a contact
  process, with percolation arows indicating the presence of an open
  path in the graphical representation.}
\end{center}
\end{figure}

We will set up a comparison between the graphical representation of the
one-dimensional nearest-neighbor contact process and oriented bond percolation
on $\Z^2$; see Figure~\ref{fig:percomp}.

We fix $T>0$ and define a map $\psi\cn\Z^2\to\Z\times\R$ by
\[
\psi(i)=\big(\kappa_i,\sig_i\big):=\big(i_1-i_2,T(i_1+i_2)\big)
\qquad\big(i=(i_1,i_2)\in\Z^2\big).
\]
The points $(\kappa_i,\sig_i)$ with $i\in\N^2$ are indicated by open circles in
Figure~\ref{fig:percomp}. As before, we make $\Z^2$ into an directed graph by
defining a collection of arrows $\Ai$ as in (\ref{Aidef}).
We wish to define a collection $(\om_{(i,j)})_{(i,j)\in\Ai}$ of
Bernoulli random variables such that
\[
\om_{(i,j)}=1\quad\mbox{implies}\quad
(\kappa_i,\sig_i)\leadsto(\kappa_j,\sig_j)\qquad\big((i,j)\in\Ai\big).
\]
For each $i\in\Z^2$ we let
\[
\tau^\pm_i:=\inf\{t\geq\sig_i:
\mbox{at time $t$ there is an infection arrow from
$\kappa_i$ to $\kappa_i\pm1$}\}
\]
denote the first time after $\sig_i$ that an arrow points out of $\kappa_i$ to
the left or right, respectively, and we define ``good events''
\[\begin{array}{r@{\,}l}
\dis\Gi^\pm_i:=\big\{
&\dis\tau^\pm_i<\sig_i+T\mbox{ and there are no blocking symbols on}\\[3pt]
&\dis\{\kappa_i\}\times(\sig_i,\tau^\pm_i]\mbox{ and }
\{\kappa_i\pm 1\}\times(\tau^\pm_i,\sig_i+T]\big\}.
\end{array}\]
Clearly,
\[\begin{array}{rl}
&\dis\Gi^-_i\quad\mbox{implies}\quad\psi(i_1,i_2)\leadsto\psi(i_1,i_2+1),\\[5pt]
\dis\mbox{and}\quad
&\dis\Gi^+_i\quad\mbox{implies}\quad\psi(i_1,i_2)\leadsto\psi(i_1+1,i_2).
\end{array}\]
In view of this, we set
\[
\om_{\txt((i_1,i_2),(i_1,i_2+1))}:=1_{\txt\Gi^-_i}
\quand
\om_{\txt((i_1,i_2),(i_1+1,i_2))}:=1_{\txt\Gi^+_i}.
\]
Then the existence of an infinite open path in the oriented percolation model
defined by the $(\om_{(i,j)})_{(i,j)\in\Ai}$ implies the existence of an
infinite open path in the graphical representation of the contact process, and
hence survival of the latter.

We observe that
\begin{equation}\label{plaT}
p:=\P[\om_{(i,j)}=1]=\P(\Gi^\pm_i)=(1-e^{-\la T})e^{-T}
\qquad\big((i,j)\in\Ai\big),
\end{equation}
which tends to one as $\la\to\infty$ while $T\to 0$ in such a way that $\la
T\to\infty$. It follows that for $\la$ sufficiently large, by a suitable
choice of $T$, we can make $p$ as close to one as we wish. We would like to
conclude from this that $\P[(0,0)\leadsto\infty]>0$ for the oriented
percolation defined by the $\om_{(i,j)}$, and therefore also
$\P[(0,0)\leadsto\infty]>0$ for the contact process.

\begin{figure}
\begin{center}
\inputtikz{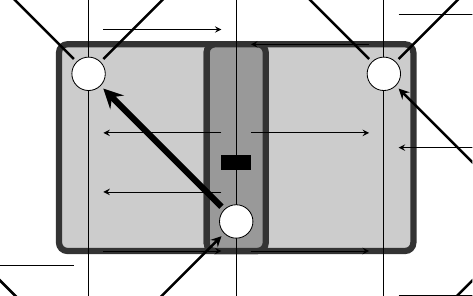}
\caption{Good events use information from partially overlapping regions of
  space-time.}
\label{fig:overlap}
\commentAlt{Figure~\ref{fig:overlap}}{Zoomed-in part of the previous
  picture indicating the overlapping space-time areas to which the
  good events relate.}
\end{center}
\end{figure}

Unfortunately, life is not quite so simple, since as shown in
Figure~\ref{fig:overlap}, the good events $\Gi^\pm_i$
have been defined using information from partially overlapping space-time
regions of the graphical representation of the contact process,
and in view of this are not independent. They are, however, 3-dependent in the
sense of Theorem~\ref{T:kdep} below, so by applying that result we can
estimate the Bernoulli random variables $(\om_{(i,j)})_{(i,j)\in\Ai}$
from below by i.i.d.\ Bernoulli random variables
$(\ti\om_{(i,j)})_{(i,j)\in\Ai}$ whose success probability $\ti p$ can be made
arbitrarily close to one, so we are done.
\end{Proof}

\section{K-dependence}\label{S:gener}

To finish the proof of Theorem~\ref{T:lacont} we need to provide the proof of
Theorem~\ref{T:kdep} below, which states that $K$-dependent random variables
with success probability $p$ can be estimated from below by i.i.d.\ random
variables with a success probability $\ti p$ that tends to one as $p\to 1$.

Traditionally, for $k\geq 0$, one says that a collection $(X_i)_{i\in\Z^d}$ of
random variables, indexed by the integer square lattice, is $k$-dependent if
for any $A,B\sub\Z^d$ with
\[
\inf\big\{|i-j|:i\in A,\ j\in B\big\}>k,
\]
the collections of random variables $(X_i)_{i\in A}$ and $(X_j)_{j\in B}$ are
independent of each other. Note that in particular, $0$-dependence means
independence.

It is a bit unfortunate that the traditional definition of $k$-dependence is
strictly tied to the integer lattice $\Z^d$, while the structure of $\Z^d$ has
little to do with the essential idea. Therefore, in this book, we
will deviate from tradition and replace(!) the definition above by the
following definition.

Let $\La$ be countable and let $(X_i)_{i\in\La}$ be a countable collection of
random variables. Then we will say that the $(X_i)_{i\in\La}$ are
\emph{$K$-dependent}\index{k-dependence@$k$-dependence}
if for each $i\in\La$ there exists a $\De_i\sub\La$
with $i\in\De_i$ and $|\De_i|\leq K$, such that
\[
\chi_i\mbox{ is independent of }(\chi_j)_{j\in\La\beh\De_i}.
\]
Note that according to our new definition, $1$-dependence means independence. The next theorem is a reformulation of \cite[Thm~B26]{Lig99}, who in turn cites \cite{LSS97}.

\begin{theorem}[$K$-dependence]\label{T:kdep}
Let $\La$ be a countable set and let $p\in(0,1)$, $K<\infty$. Assume that
$(\chi_i)_{i\in\La}$ are $K$-dependent Bernoulli random variables with
$\P[\chi_i=1]\geq p$ $(i\in\La)$, and that
\[
\ti p:=\big(1-(1-p)^{1/K}\big)^2\geq\ffrac{1}{4}.
\]
Then it is possible to couple $(\chi_i)_{i\in\La}$ to a collection of
independent Bernoulli random variables $(\ti\chi_i)_{i\in\La}$ with
\begin{equation}\label{tichi}
\P[\ti\chi_i=1]=\ti p\qquad(i\in\La),
\end{equation}
in such a way that $\ti\chi_i\leq\chi_i$ for all $i\in\La$.
\end{theorem}

\begin{Proof}
In the language of Theorem~\ref{T:stochord}, we must show that the law
of $(\ti\chi_i)_{i\in\La}$ lies below the law of $(\chi_i)_{i\in\La}$
in the stochastic order. Since we can always choose some arbitrary
denumeration of $\La$, we may assume that $\La=\N$. Our strategy will
be as follows. We will choose $\{0,1\}$-valued random variables
$(\psi_i)_{i\in\La}$ with $\P[\psi_i=1]=r$, independent of each other
and of the $(\chi_i)_{i\in\N}$, and put
\[
\chi'_i:=\psi_i\chi_i\qquad(i\in\N).
\]
Note that the $(\chi'_i)_{i\in\N}$ are a thinning\index{thinning}
of the $(\chi_i)_{i\in\N}$. In particular, $\chi'_i\leq\chi_i$
$(i\in\N)$, so the law of $(\chi'_i)_{i\in\La}$ lies below the law of
$(\chi_i)_{i\in\La}$ in the stochastic order. We will show that for an
appropriate choice of $r$,
\begin{equation}\label{conbd}
\P[\chi'_n=1\,|\,\chi'_0,\ldots,\chi'_{n-1}]\geq\ti p
\end{equation}
for all $n\geq 0$, and we will show that this implies that the law of
$(\ti\chi_i)_{i\in\La}$ lies below the law of $(\chi'_i)_{i\in\La}$ in
the stochastic order. Thus
\[
\P\big[(\ti\chi_i)_{i\in\La}\in\,\cdot\,\big]
\leq\P\big[(\chi'_i)_{i\in\La}\in\,\cdot\,\big]
\leq\P\big[(\chi_i)_{i\in\La}\in\,\cdot\,\big],
\]
which implies by Theorem~\ref{T:stochord} that the
$(\chi_i)_{i\in\La}$ can be coupled to $(\ti\chi_i)_{i\in\La}$ such
that $\ti\chi_i\leq\chi_i$ for all $i\in\La$.

We start by showing that (\ref{conbd}) implies that the
$(\ti\chi_i)_{i\in\La}$ and $(\chi'_i)_{i\in\La}$ can be coupled such
that $\ti\chi_i\leq\chi'_i$ for all $i\in\La$. Set
$p'_0:=\P[\chi'_0=1]$ and
\[
p'_n(\eps_0,\ldots,\eps_{n-1})
:=\P[\chi'_n=1\,|\,\chi'_0=\eps_0,\ldots,\chi'_{n-1}=\eps_{n-1}]
\]
whenever $\P[\chi'_0=\eps_0,\ldots,\chi'_{n-1}=\eps_{n-1}]>0$. Let
$(U_n)_{n\in\N}$ be independent, uniformly distributed $[0,1]$-valued
random variables. Set
\[
\ti\chi_n:=1_{\txt\{U_n<\ti p\}}\qquad(n\in\N)
\]
and define inductively
\[
\chi'_n:=1_{\txt\{U_n<p'_n(\chi'_0,\ldots,\chi'_{n-1})\}}\qquad(n\in\N).
\]
Then
\[
\P[\chi'_n=\eps_n,\ldots,\chi'_0=\eps_0]
=p'_n(\eps_0,\ldots,\eps_{n-1})\cdots p'_1(\eps_0)\cdot p'_0.
\]
This shows that these new $\chi'_n$ have the same distribution as the old
ones, and they are coupled to $\ti\chi_i$ as in (\ref{tichi}) in such a way
that $\ti\chi_i\leq\chi'_i$.

What makes life complicated is that (\ref{conbd}) does not always hold
for the original $(\chi_i)_{i\in\N}$, which is why we have to work
with the thinned variables $(\chi'_i)_{i\in\N}$.\footnote{Indeed, let
$(\phi_n)_{n\geq 0}$ be independent $\{0,1\}$-valued random variables
with $\P[\phi_n=1]=\sqrt p$ for some $p<1$, and put
$\chi_n:=\phi_n\phi_{n+1}$. Then the $(\chi_n)_{n\geq 0}$ are
$3$-dependent with $\P[\chi_n=1]=p$, but
$\P[\chi_n=1|\chi_{n-1}=0,\chi_{n-2}=1]=0$.} We observe that
\begin{equation}\label{pI}
\P[\chi'_n=1\,|\,\chi'_0=\eps_0,\ldots,\chi'_{n-1}=\eps_{n-1}]
=r\P[\chi_n=1\,|\,\chi'_0=\eps_0,\ldots,\chi'_{n-1}=\eps_{n-1}].
\end{equation}
We will prove by induction that for an appropriate choice of $r$,
\begin{equation}\label{IH}
\P[\chi_n=0\,|\,\chi'_0=\eps_0,\ldots,\chi'_{n-1}=\eps_{n-1}]\leq 1-r.
\end{equation}
Note that this is true for $n=0$ provided that $r\leq p$. Let us put
\[\begin{array}{c}
\dis E_0:=\big\{i\in\De_n:0\leq i\leq n-1,\ \eps_i=0\big\},\quad
E_1:=\big\{i\in\De_n:0\leq i\leq n-1,\ \eps_i=1\big\},\\[5pt]
\dis F:=\big\{i\not\in\De_n:0\leq i\leq n-1\big\}.
\end{array}\]
Then
\begin{equation}\begin{array}{l}\label{LSS}
\dis \P[\chi_n=0\,|\,\chi'_0=\eps_0,\ldots,\chi'_{n-1}=\eps_{n-1}]\\[7pt]
\quad\dis=\P\big[\chi_n=0\,\big|\,\chi'_i=0\ \forall i\in E_0,
\ \chi_i=1=\psi_i\ \forall i\in E_1,\ \chi'_i=\eps_i\ \forall i\in F\big]\\[7pt]
\quad\dis=\P\big[\chi_n=0\,\big|\,\chi'_i=0\ \forall i\in E_0,
\ \chi_i=1\ \forall i\in E_1,\ \chi'_i=\eps_i\ \forall i\in F\big]\\[8pt]
\quad\dis=\frac{\P\big[\chi_n=0,\ \chi'_i=0\ \forall i\in E_0,
\ \chi_i=1\ \forall i\in E_1,\ \chi'_i=\eps_i
\ \forall i\in F\big]}{\P\big[\chi'_i=0\ \forall i\in E_0,\ \chi_i=1
\ \forall i\in E_1,\ \chi'_i=\eps_i\ \forall i\in F\big]}\\[13pt]
\quad\dis\leq\frac{\P\big[\chi_n=0,\ \chi'_i=\eps_i
\ \forall i\in F\big]}{\P\big[\psi_i=0\ \forall i\in E_0,
\ \chi_i=1\ \forall i\in E_1,\ \chi'_i=\eps_i\ \forall i\in F\big]}\\[13pt]
\quad\dis=\frac{\P\big[\chi_n=0\,\big|\,\chi'_i=\eps_i
\ \forall i\in F\big]}{\P\big[\psi_i=0\ \forall i\in E_0,\ \chi_i=1
\ \forall i\in E_1\,\big|\,\chi'_i=\eps_i\ \forall i\in F\big]}\\[13pt]
\quad\dis\leq\frac{1-p}{(1-r)^{|E_0|}\P\big[\chi_i=1
\ \forall i\in E_1\,\big|\,\chi'_i=\eps_i\ \forall i\in F\big]}
\leq\frac{1-p}{(1-r)^{|E_0|}\,r^{|E_1|}},
\end{array}\end{equation}
where in the last step we have used $K$-dependence and the (nontrivial) fact
that
\begin{equation}\label{IHc}
\P\big[\chi_i=1\ \forall i\in E_1\,\big|\,\chi'_i=\eps_i
\ \forall i\in F\big]\geq r^{|E_1|}.
\end{equation}
We claim that (\ref{IHc}) is a consequence of the induction hypothesis
(\ref{IH}). Indeed, we may assume that the induction hypothesis (\ref{IH})
holds regardless of the ordering of the first $n$ elements, so without loss of
generality we may assume that $E_1=\{n-1,\ldots,m\}$ and $F=\{m-1,\ldots,0\}$,
for some $m$. Then the left-hand side of (\ref{IHc}) may be written as
\[\begin{array}{l}
\dis\prod_{k=m}^{n-1}\P\big[\chi_k=1\,\big|\,\chi_i=1\ \forall m\leq i<k,
\ \chi'_i=\eps_i\ \forall 0\leq i<m\big]\\[5pt]
\dis\quad=\prod_{k=m}^{n-1}\P\big[\chi_k=1\,\big|\,\chi'_i=1
\ \forall m\leq i<k,\ \chi'_i=\eps_i\ \forall 0\leq i<m\big]\geq r^{n-m}.
\end{array}\]
If we assume moreover that $r\geq\frac{1}{2}$, then
$r^{|E_1|}\geq(1-r)^{|E_1|}$ and therefore the right-hand side of
(\ref{LSS}) can be further estimated as
\[
\frac{1-p}{(1-r)^{|E_0|}\,r^{|E_1|}}\leq\frac{1-p}{(1-r)^{|\De_n\cap\{0,\ldots,n-1\}|}}
\leq\frac{1-p}{(1-r)^{K-1}}.
\]
We see that in order for our proof to work, we need $\frac{1}{2}\leq
r\leq p$ and
\begin{equation}\label{pr}
\frac{1-p}{(1-r)^{K-1}}\leq 1-r.
\end{equation}
In particular, choosing $r=1-(1-p)^{1/K}$ yields equality in
(\ref{pr}). Having proved (\ref{IH}), we see by (\ref{pI}) that
(\ref{conbd}) holds provided that we put $\ti p:=r^2$.
\end{Proof}

\begin{Exercise}
Combine Theorem~\ref{T:pcrit} and formulas (\ref{plaT}) and (\ref{tichi}) to
derive an explicit upper bound on the critical infection rate $\la_{\rm c}$ of
the one-dimensional contact process.
\end{Exercise}


\begin{Exercise}
The one-dimensional contact process with double deaths has been introduced
just before Exercise~\ref{E:doubdeath}. Use comparison with oriented
percolation to prove that the one-dimensional contact process with double
deaths survives with positive probability if its branching rate $\la$ is large
enough. When you apply Theorem~\ref{T:kdep}, what value of $K$ do you (at
least) need to use?
\end{Exercise}

\begin{Exercise}
Use the previous exercise and Exercise~\ref{E:doubdeath} to conclude
that for the cooperative branching process considered there, if $\la$
is large enough, then: I.\ If the process is started with at least two
particles on neighboring sites, then there is a positive probability
that there will always be pairs of particles on neighboring
sites. II.\ The upper invariant law is nontrivial.
\end{Exercise}

\begin{Exercise}
Assume that there exists some $t>0$ such that the contact
process satisfies
\[
r:=\E^{e_0}\big[|X_t|\big]<1.
\]
Show that this then implies that
\[
\E^{e_0}\big[|X_{nt}|\big]\leq r^n\qquad(n\geq 0)
\]
and the process started in any finite initial state dies out a.s.
Can you use this to improve the lower bound $1/|\Ni_i|\leq\la_{\rm c}$
from Exercise~\ref{E:context}, for example for the one-dimensional
nearest-neighbor process?
\end{Exercise}


\addtocontents{toc}{\vspace{\baselineskip}}

\backmatter

\renewcommand{\refname}{Bibliography}

\cleardoublepage
\phantomsection

\printindex

\end{document}